\documentclass[11pt,openany]{book}
\usepackage{t1enc}
\usepackage[T2A,T1]{fontenc}
\usepackage[utf8]{inputenc}
\usepackage[russian,english]{babel}
\usepackage{amsmath,amssymb,amsfonts,amsthm,mathrsfs,textcomp,url,bbm,cancel,comment}
\usepackage{graphicx}
\usepackage{float}
\usepackage[all,ps]{xy}
\usepackage[colorlinks,urlcolor=cyan,citecolor=blue,linkcolor=blue]{hyperref}
\usepackage[dvipsnames]{xcolor}
\usepackage[margin=3.5cm]{geometry}
\usepackage{longtable}
\usepackage{tikz}
\usetikzlibrary{matrix,calc,decorations,positioning}
\pgfkeys{/tikz/.cd,
    alt double distance/.initial=5pt,
    alt double step/.initial=1pt,
}

\pgfdeclaredecoration{double deco}{initial}
{% initial arrow butt
\state{initial}[width=\pgfkeysvalueof{/tikz/alt double step},next state=cont] {
    \pgfmoveto{\pgfpoint{\pgfkeysvalueof{/tikz/alt double step}}{\pgfkeysvalueof{/tikz/alt double distance}/2}}
    \pgfpathlineto{\pgfpoint{0.3\pgflinewidth}{\pgfkeysvalueof{/tikz/alt double distance}/2}}
    \pgfpathmoveto{\pgfpoint{0.3\pgflinewidth}{-\pgfkeysvalueof{/tikz/alt double distance}/2}}
    \pgfpathlineto{\pgfpoint{1pt}{-\pgfkeysvalueof{/tikz/alt double distance}/2}}
    \pgfcoordinate{lastup}{\pgfpoint{1pt}{\pgfkeysvalueof{/tikz/alt double distance}/2}}
    \pgfcoordinate{lastdown}{\pgfpoint{1pt}{-\pgfkeysvalueof{/tikz/alt double distance}/2}}
  }
  \state{cont}[width=\pgfkeysvalueof{/tikz/alt double step}]{
     \pgfmoveto{\pgfpointanchor{lastup}{center}}
     \pgfpathlineto{\pgfpoint{\pgfkeysvalueof{/tikz/alt double step}}{\pgfkeysvalueof{/tikz/alt double distance}/2}}
     \pgfcoordinate{lastup}{\pgfpoint{\pgfkeysvalueof{/tikz/alt double step}}{\pgfkeysvalueof{/tikz/alt double distance}/2}}
     \pgfmoveto{\pgfpointanchor{lastdown}{center}}
     \pgfpathlineto{\pgfpoint{\pgfkeysvalueof{/tikz/alt double step}}{-\pgfkeysvalueof{/tikz/alt double distance}/2}}
     \pgfcoordinate{lastdown}{\pgfpoint{\pgfkeysvalueof{/tikz/alt double step}}{-\pgfkeysvalueof{/tikz/alt double distance}/2}}
  }
  \state{final}[width=0pt]
  { % perhaps unnecessary but doesn't hurt either
    \pgfmoveto{\pgfpointdecoratedpathlast}
  }
}
\tikzset{alt double/.style={decorate,decoration=double deco}}
\allowdisplaybreaks

\renewcommand{\labelenumi}{\theenumi}
\renewcommand{\theenumi}{\rm{(\arabic{enumi})}}
\newcommand{\toverset}[2]{%
\mathop{#2}\limits^{\vbox to -.1ex{\kern-0.4ex\hbox{$\scriptstyle #1$}\vss}}}
\newcommand{\tightoverset}[2]{%
\mathop{#2}\limits^{\vbox to -.5ex{\kern-0.4ex\hbox{$\scriptstyle #1$}\vss}}}
\renewcommand{\underset}[2]{%
\mathop{#2}\limits_{\vbox to -.5ex{\kern-1.6ex\hbox{$\scriptstyle #1$}\vss}}}
\newcommand{\tightunderset}[2]{%
\mathop{#2}\limits_{\vbox to -.5ex{\kern-1.8ex\hbox{$\scriptstyle #1$}\vss}}}
\newcommand{\xra}[1]{%
\mathop{\xrightarrow{~#1~}}}
\newcommand{\xla}[1]{%
\mathop{\xleftarrow{~#1~}}}
\newcommand{\congq}[0]{%
\mathop{\tightunderset{\mathbb{Q}}{\cong}}}
\newcommand{\congc}[1]{%
\mathop{\tightunderset{\mathscr{C}_{#1}}{\cong}}}

\newcommand{\liminfty}[1]{%
\mathop{\underset{#1\to\infty}\lim}}
\newcommand{\cpt}[1]{%
\mathop{\toverset{^*}{#1}}}

\newcommand{\usqcup}[1]{%
\mathop{\tightunderset{#1}{\sqcup}}}
\newcommand{\utimes}[1]{%
\mathop{\tightunderset{#1}{\times}}}
\newcommand{\sims}[1]{%
\text{\scalebox{#1}[1]{$\sim$}}}
\newcommand{\widesqcup}{%
\text{\scalebox{1.5}[1]{$\bigsqcup$}}}
\newcommand{\lowtilde}[1]{%
\raisebox{-.7ex}{\~{#1}}}
\def\F{\mathbb{F}}
\def\N{\mathbb{N}}
\def\Z{\mathbb{Z}}
\def\Q{\mathbb{Q}}
\def\R{\mathbb{R}}
\def\C{\mathbb{C}}
\def\HH{\mathbb{H}}
\def\FP{\mathbb{F}P}
\def\RP{\mathbb{R}P}
\def\CP{\mathbb{C}P}
\def\HP{\mathbb{H}P}
\def\AA{\mathscr{A}}
\def\CC{\mathscr{C}}
\def\FF{\mathscr{F}}
\def\GG{\mathscr{G}}
\def\NN{\mathfrak{N}}
\def\into{\hookrightarrow}
\def\imto{\looparrowright}
\def\cs{\smallsmile}
\def\ol{\overline}
\def\ul{\underline}
\def\lv{\lVert}
\def\rv{\rVert}
\def\la{\langle}
\def\ra{\rangle}
\def\ua{\mathop{\uparrow}}
\def\da{\mathop{\downarrow}}
\def\ext{\mathrm{ext}}
\def\im{\mathop{\rm im}}
\def\coker{\mathop{\rm coker}}

\def\codim{\mathop{\rm codim}}

\def\id{\mathop{\rm id}}
\def\pr{\mathop{\rm pr}}

\def\mod{\mathop{\rm mod}}

\def\Tp{\mathrm{Tp}}
\def\Sp{\mathrm{Sp}}
\def\Sq{\mathrm{Sq}}
\def\SP{\mathop{\rm SP}}
\def\Spin{\mathrm{Spin}}
\def\O{\mathrm{O}}
\def\SO{\mathrm{SO}}
\def\U{\mathrm{U}}

\def\rk{\mathop{\rm rk}}

\def\Emb{\textstyle{\mathop{\rm Emb}}}
\def\Stab{\textstyle{\mathop{\rm Stab}}}
\def\Imm{\textstyle{\mathop{\rm Imm}}}
\def\Cob{\textstyle{\mathop{\rm Cob}}}
\def\Prim{\textstyle{\mathop{\rm Prim}}}
\def\Bord{\textstyle{\mathop{\rm Bord}}}
\def\Stong{\textstyle{\mathop{\rm Stong}}}

\def\Spin{\mathop{\rm Spin}}
\def\Diff{\textstyle{\mathop{\rm Diff}}}
\def\Aut{\textstyle{\mathop{\rm Aut}}}
\def\Hom{\textstyle{\mathop{\rm Hom}}}

\newenvironment{prf}%
{\par\noindent\textbf{Proof.\enspace\ignorespaces}}%
{~$\square$\par\medskip}%
\newenvironment{sclaim}%
{\par\noindent\textit{Claim.\enspace\ignorespaces}\begin{em}}%
{\end{em}\par}%
\newenvironment{sprf}%
{\medskip\par\noindent\textit{Proof.\enspace\ignorespaces}}%
{~$\diamond$\par\medskip}%
{\medskip\par\noindent\textit{Remark.\enspace\ignorespaces}}%
{\par\medskip}%
{\medskip\par\noindent\textbf{Acknowledgement.\enspace\ignorespaces}}%
{\par\medskip}%
\newtheorem{thm}{Theorem}[section]%
\newtheorem{lemma}[thm]{Lemma}%
\newtheorem{prop}[thm]{Proposition}%
\newtheorem{crly}[thm]{Corollary}%
\newtheorem*{thm*}{Theorem}%
\newtheorem*{add*}{Addendum}%
\newtheorem*{claim*}{Claim}%
\newtheorem*{case}{The Eucledian case}%
\theoremstyle{definition}
\newtheorem{defi}[thm]{Definition}%
\newtheorem{rmk}[thm]{Remark}%
\newtheorem{rmks}[section]{Remark}%
\newtheorem{ex}[thm]{Example}%
\title{}
\author{}
\date{}
\begin{document}
%\maketitle

\newpage
\thispagestyle{empty}
\begin{center}
\textsc{\Large Eötvös Loránd University}\\[.4cm]
\textsc{\large Faculty of Science}\\[.4cm]
\textsc{\large Institute of Mathematics}\\[.4cm]
\textsc{\large Department of Analysis}\\[.4cm]
\hrulefill\\[1.5cm]
		
{\Huge \textbf{Cobordisms of singular maps}}\\[.6cm]
\textsc{\large --- master's thesis ---}\\
\end{center}

\vfill

\hspace{-.5cm}
\begin{tabular}[t]{l} 
\textit{Author:}\\
András Csépai
\end{tabular}
\hfill
\begin{tabular}[t]{l}
\textit{Supervisor:}\\
András Szűcs
\end{tabular}

\begin{center}\begin{figure}[h!]
\centering
\includegraphics[width=0.35\textwidth]{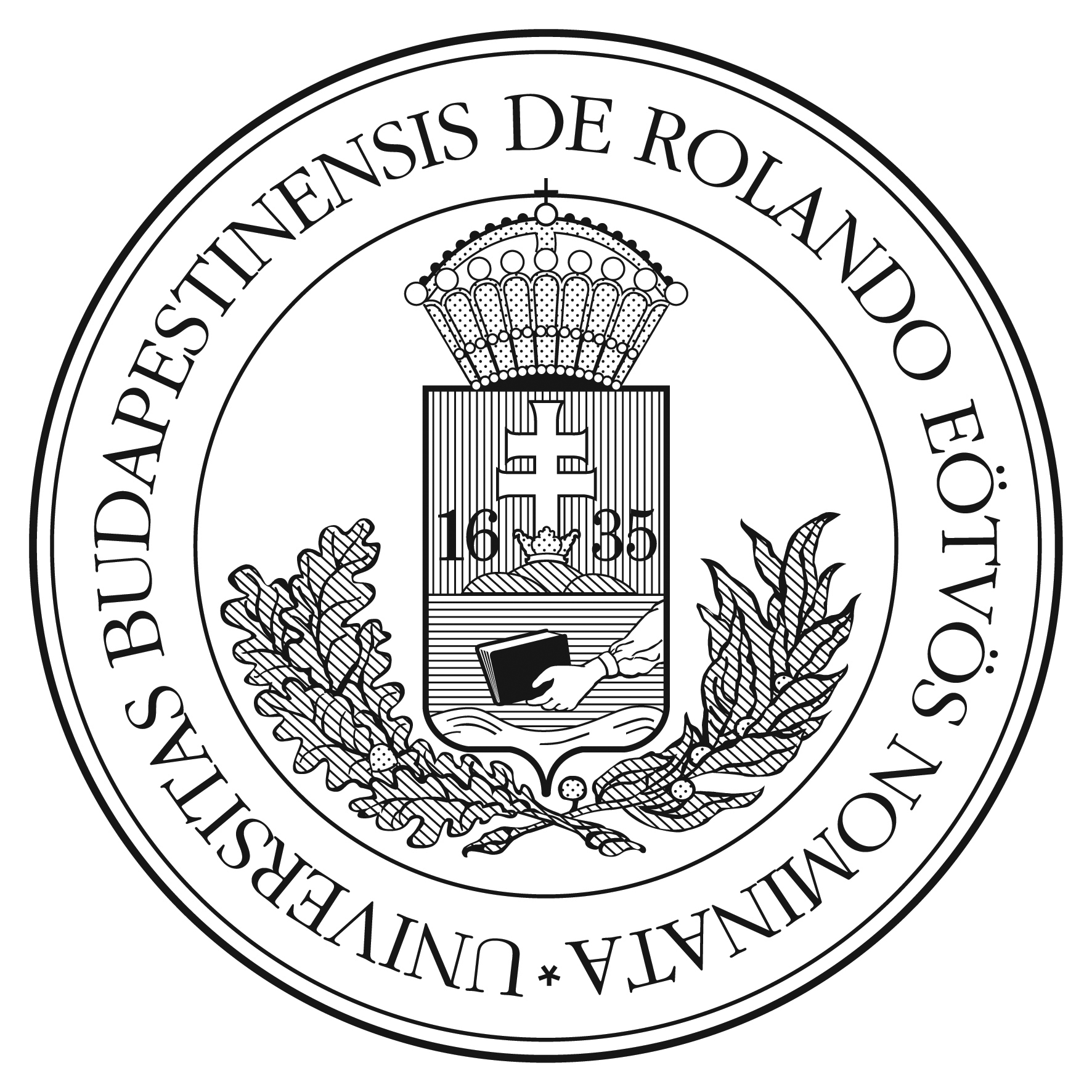}
\end{figure}~\\
2020.
\end{center}

\newpage
\pagenumbering{roman}

\tableofcontents
%\addcontentsline{toc}{chapter}{Contents}

%-------------------------------------------------------------------------------------------------------------------------------
%-------------------------------------------------------------------------------------------------------------------------------

\chapter*{Introduction}
\addcontentsline{toc}{chapter}{Introduction}

The main motivating question of the present thesis is the following:
\medskip\begin{center}
\emph{What do smooth maps between manifolds look like?}
\end{center}\medskip

Of course, this question is not at all precise, however, the description of smooth maps is a central problem in differential topology. There are many variants of what we mean by ``look like'' and precisely which sorts of smooth maps we are considering; these variants give rise to quite a few large theories of topology (for example Morse theory or the immersion theory of Hirsch, Smale, Gromov and others). Similarly to them, we will also narrow down the problem of investigating smooth maps to a specific theory in differential topology.

When thinking about mapping a manifold to another, probably the first kind of maps that comes to mind is the embeddings. These are very easy to visualise, since a manifold embedded into another one has its manifold structure as the subspace topology of the containing manifold. Another relatively simple kind of smooth maps is the immersions, which are just locally embeddings.

\begin{figure}[H]
\begin{center}
\centering\includegraphics[scale=0.1]{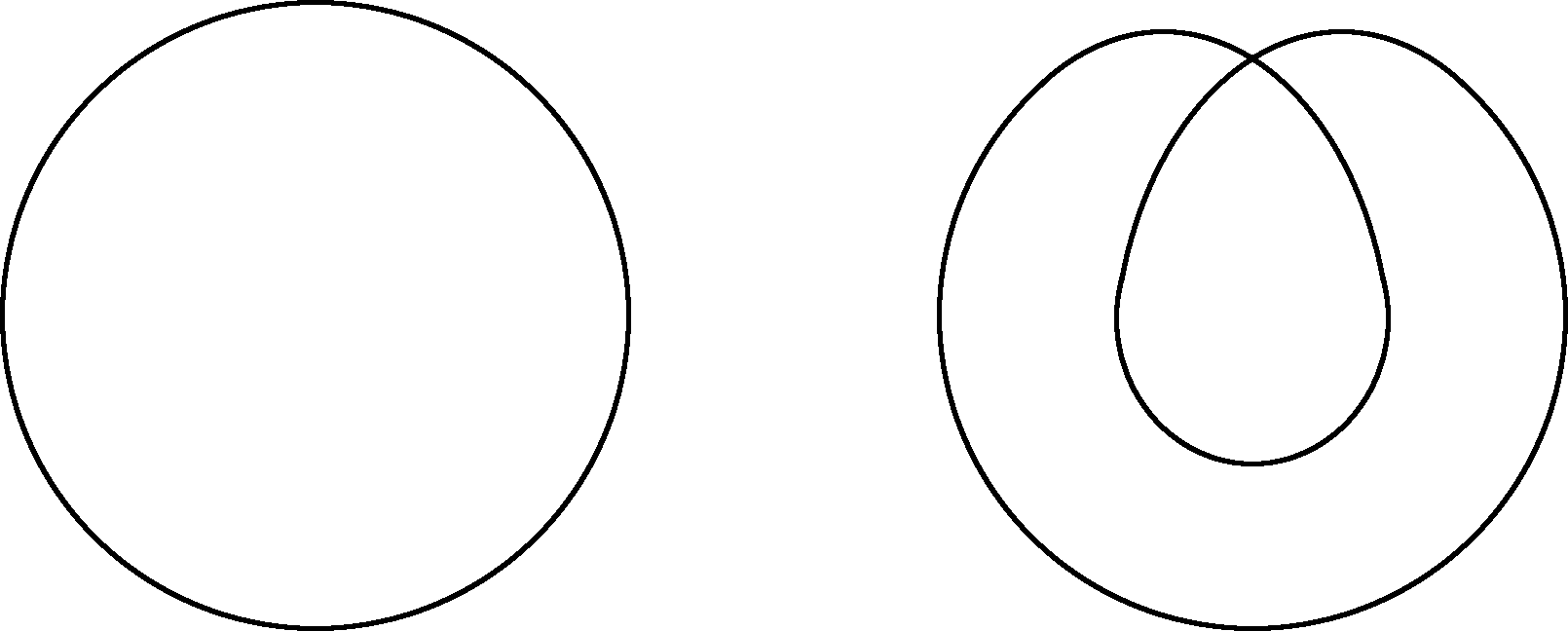}
\begin{changemargin}{2cm}{2cm} 
\caption{\hangindent=1.4cm\small Here we show the images of an embedding and an immersion of the circle into the plane.}
\end{changemargin} 
\vspace{-1.5cm}
\end{center}
\end{figure}

Embeddings and immersions share the same property that their differentials have maximal ranks in all points, however, generic smooth maps can be more complicated. There are many ways that the rank of the differential of a map can decrease in a point, which is then called a singular point of that map; maps with singular points are called singular maps. One way to think about describing general (singular) smooth maps is that it has the following three levels:
\begin{enumerate}
\item Local description: This concerns the local forms of map germs, which describe the map in small neighbourhoods of points. Being a singular point of a map is a local property and this talks about the singularities of maps individually. For a long time this was the main subject of singularity theory; see for example the results of Whitney \cite{wh}, Morin \cite{mor}, Mather \cite{math}, Arnold \cite{arnold} and others.
\item Semiglobal description: This investigates the automorphism groups of local forms, which provide an uderstanding of the map in neighbourhoods of its singularity strata (a stratum is a set of points with equal local forms). For results on this see Jänich \cite{jan}, Wall \cite{wall} and Rimányi \cite{rrthesis}.
\item\label{lazasag} Global description: This means describing how the different singularity strata fit together, how the more complicated strata are attached to the simpler strata to form the map. A very important result on this is a universal singular map constructed by Szűcs and Rimányi \cite{rsz} (and preceeded by a number of similar results of Szűcs, e.g. \cite{analog}) that spawns a theory of studying singular maps globally.
\end{enumerate}

We remark here that global singularity theory does not restrict to the line of investigation described above; for example one can also study the Thom polynomials which stand for the homology classes represented by singularity strata (see Thom \cite{tp}, Fehér, Rimányi \cite{fr}, Kazarian \cite{kazthom} and others). However, the above point \ref{lazasag} is what we will particularly be interested in now.

\begin{figure}[H]
\begin{center}
\centering\includegraphics[scale=0.1]{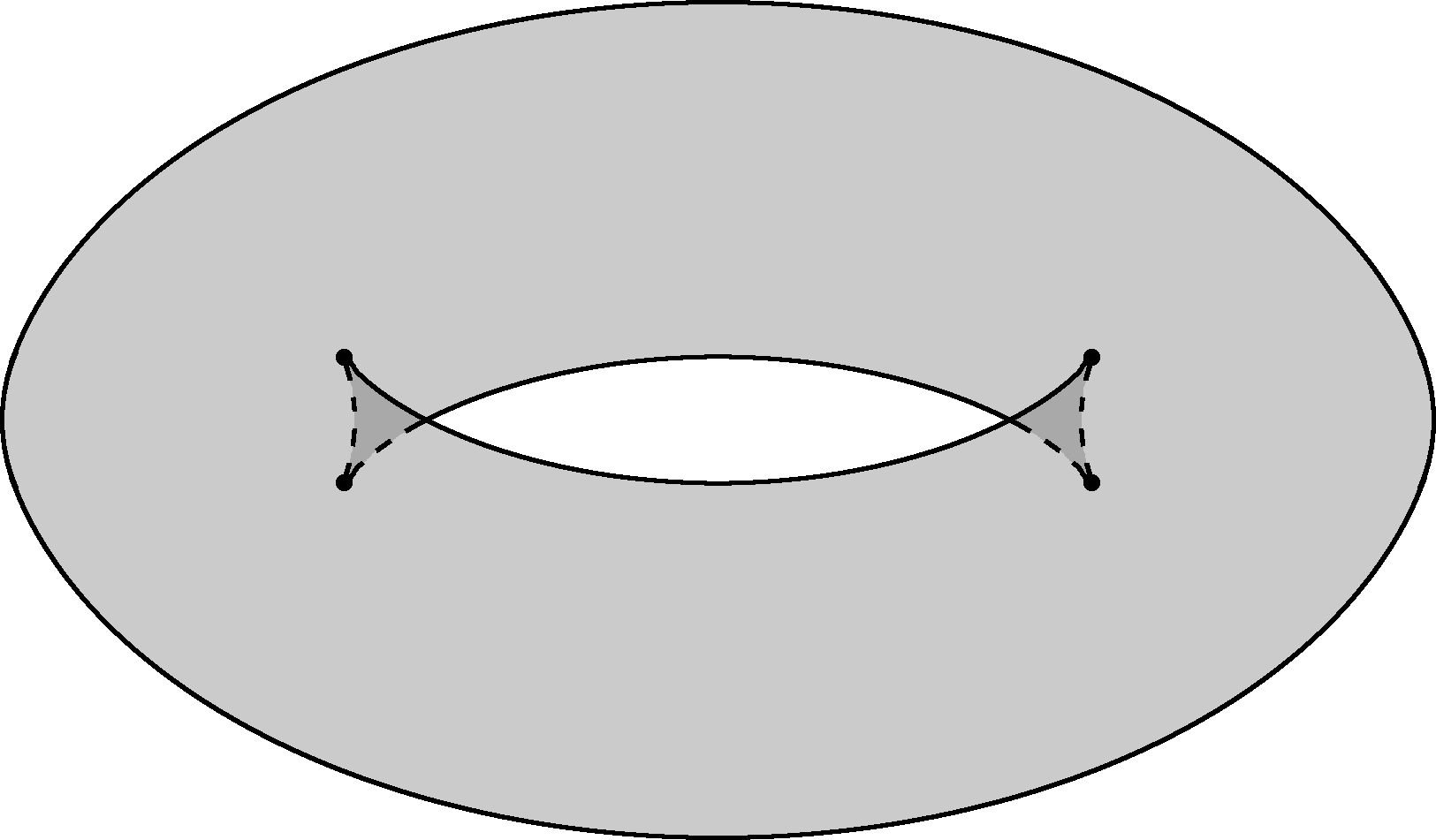}
\begin{changemargin}{2cm}{2cm} 
\caption{\hangindent=1.4cm\small The grey area represents the image of a torus mapped to the plane and all the curves bounding it represent where the surface was ``folded'' by the mapping; these are the singular points of the map. Note that the inside curve is not smooth at the four ``corner'' points; we will see that the map at these points has a more complicated local form than it has at the rest of the curves.}
\end{changemargin} 
\vspace{-1.5cm}
\end{center}
\end{figure}

Now that we have some understanding about the phrase ``singular maps'' let us get to the other half of the title of this thesis and explain the word ``cobordism''.

As usual in mathematics and in particular in topology, we do not want to differentiate between objects that are very similar in some sense, which in topology means that they can be deformed into each other in some way. The theory we are about to investigate considers maps which have singularities from a fixed list $\tau$ and the deformation of one such map to another is given by deforming the source manifold of the first map to that of the second by a higher dimensional manifold (called a cobordism of manifolds) and connect the two maps by a map of this higher dimensional manifold that also has singularities from $\tau$. This yields an equivalence relation called the cobordism of maps with singularities from $\tau$.

\begin{figure}[H]
\begin{center}
\centering\includegraphics[scale=0.1]{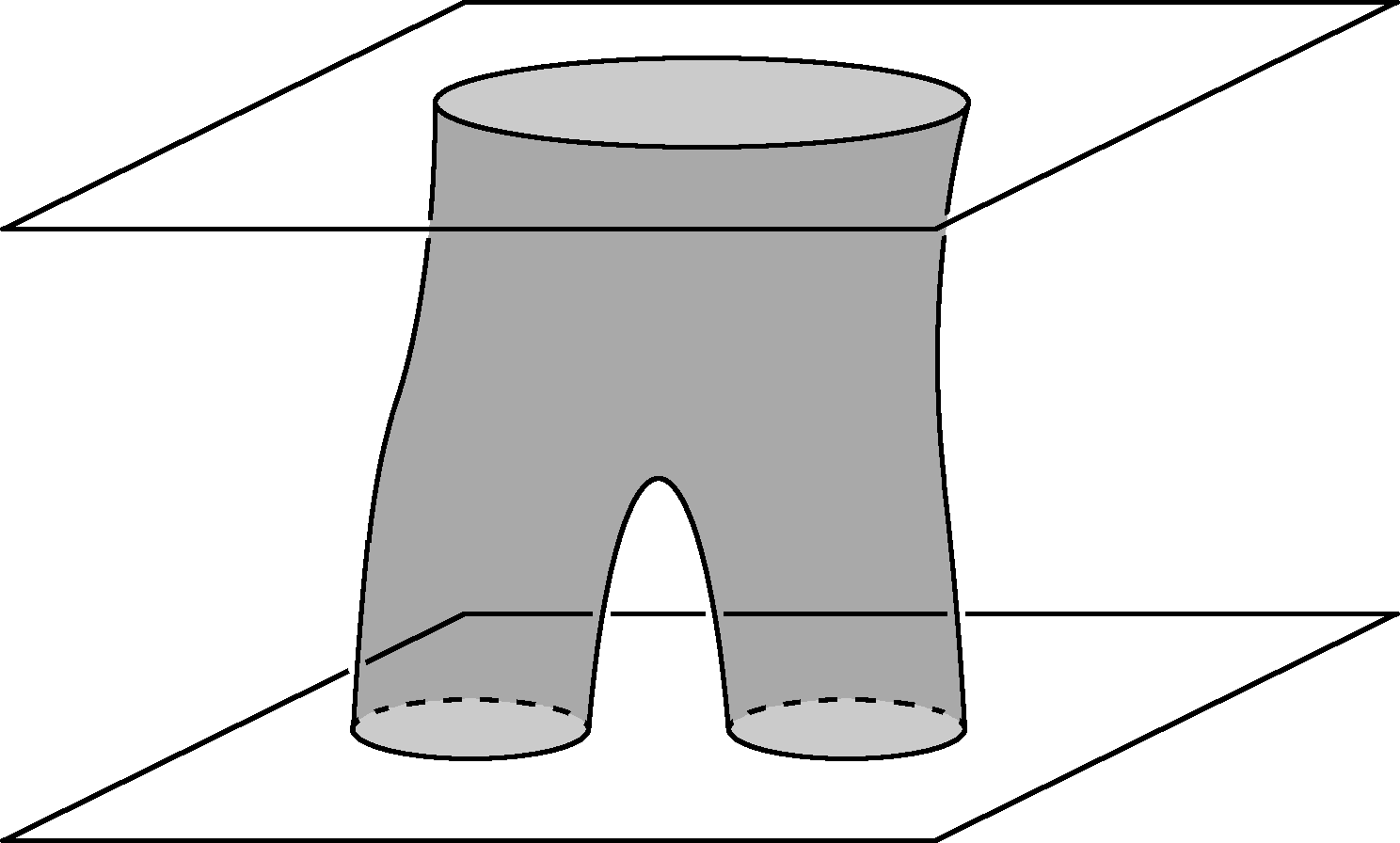}
\begin{changemargin}{2cm}{2cm} 
\caption{\hangindent=1.4cm\small This ``pair of pants'' is a cobordism between an embedding of one circle into the plane (on the top) and an embedding of two circles into the plane (on the bottom).}
\end{changemargin} 
\vspace{-1.5cm}
\end{center}
\end{figure}

The equivalence classes given by the above relation on singular maps are called cobordism classes. The set of cobordism classes of maps to a fixed manifold and with singularities from the list $\tau$ admits a natural group operation; the cobordism theory of singular maps concerns the computation of these groups for various target manifolds and lists of singularities.

These groups were and are investigated among others by Koschorke \cite{kosvf}, Ando \cite{ando}, Saeki \cite{saeki}, Ikegami \cite{ike}, Kalmár \cite{kal} and Sadykov \cite{sad}. However, one of the first people to consider them is Szűcs \cite{analog}; the aim of the present thesis is to collect and organise the results towards and on the computation of cobordism groups by his method published in papers of Szűcs, Terpai and other coauthors \cite{rsz}, \cite{cobmor}, \cite{elimcob}, \cite{eszt}, \cite{hosszu}, \cite{2k+2}, \cite{key}, \cite{szszt}, \cite{nulladik}, \cite{nszt}, \cite{nehez}, \cite{hominv} and \cite{ctrl}.

This thesis contains no individual results, however, I introduce many notions in a more general setting than the original papers did, which yields definitions and theorems that already existed implicitly in literature but were never written down. It is my hope that this will help in expanding the theory and obtaining results for cobordism groups with more general structures.

%-------------------------------------------------------------------------------------------------------------------------------
%-------------------------------------------------------------------------------------------------------------------------------

\chapter*{Conventions and notations}
\addcontentsline{toc}{chapter}{Conventions and notations}

Throughout this thesis, when we talk about manifolds, we are in the smooth category, that is, all manifolds and maps are assumed to be $C^\infty$. If we want to indicate the dimension of a manifold, we put it in a superindex (i.e. $M^n$ means that $M$ is a manifold of dimension $n$), but in most cases we omit this index. There will be many times when we use certain transversality conditions without mentioning it; the validity of these always follows from standard approximation and compactness reasons. For homologies and cohomologies, if we do not indicate the coefficient group, then it is assumed to be $\Z$. By fibration we usually mean Serre fibration.

\subsection*{Notations}\vspace{-.2cm}

\setlength\LTleft{-.2cm}
\renewcommand{\arraystretch}{1.3}
\begin{longtable}{p{2.9cm}p{11.05cm}}
$\R_+$ & non-negative real numbers\\
$\RP^n,\CP^n,\HP^n$ & $n$-dimensional real, complex and quaternionic projective spaces\\
$S^n$ & the sphere $\{x\in\R^{n+1}\mid\lv x\rv=1\}$\\
$D^n$ & the disk $\{x\in\R^n\mid\lv x\rv\le1\}$\\
$\O,\SO,\Spin,\U,\Sp$ & direct limits of the orthogonal groups $\O(n)$, special orthogonal groups $\SO(n)$, spin groups $\Spin(n)$, unitary groups $\U(n)$ and symplectic groups $\Sp(n)$ as $n\to\infty$\\
$*$ & one-point space\\
$X^+$ & the space $X\sqcup*$\\
$\cpt X$ & one-point compactification of the space $X$ (if $X$ is compact, then this is $X^+$)\\
$X\utimes GY$ & factor space of $X\times Y$ by the diagonal $G$-action, i.e. the space $X\times Y/(x,y)\sim(gx,gy)~\forall\,x\in X,y\in Y, g\in G$ (for some fixed actions of the group $G$ on the spaces $X$ and $Y$)\\
$X\usqcup\rho Y$ & the space $Y$ glued to $X$ by $\rho$, i.e. the factor space $X\sqcup Y/y\sim\rho(y)$ (where $\rho$ is a map from a subspace of $Y$ to $X$)\\
$CX,SX$ & cone and suspension of the space $X$\\
$PX,\Omega X$ & path space and loop space of the pointed space $X$\\
$\Gamma X$ & the space $\Omega^\infty S^\infty X:=\liminfty{n}\Omega^n S^n X$\\
$\SP X$ & the infinite symmetric product $\liminfty n\underbrace{X\times\ldots\times X}_{n\text{ times}}\big/S_n$\\
$[X]_p$ & $p$-localisation of the space $X$ for a prime $p$\\
$\partial X$ & boundary of $X$\\
$\overline X$ & closure of $X$\\
$\toverset{_\circ} X$ & interior of $X$\\
$\id_X$ & the identity map $X\to X$\\
$\pr_{X_j}$ & the projection $\underset{i\in I}\prod X_i\to X_j$ (for any $j\in I$)\\
$EG\xra GBG$ & universal principal $G$-bundle\\
$\gamma_n^G$ & universal $n$-dimensional vector bundle with structure group $G$; if $G=\O(n)$, then we omit the superindex\\
$\varepsilon^n$ & $n$-dimensional trivial vector bundle (over any base space)\\
$\nu_f$ & virtual normal bundle of the map $f$\\
$D\zeta,S\zeta,T\zeta$ & the unit disk bundle and sphere bundle of the vector bundle $\zeta$ (with respect to any metric) and the Thom space of $\zeta$\\
$\tilde G_\zeta$ & Local coefficient system from the group $G$ twisted by the determinant bundle of the vector bundle $\zeta$\\
$e(\zeta),u(\zeta)$ & Euler and Thom class of the vector bundle $\zeta$\\
$w_i(\zeta),c_i(\zeta),p_i(\zeta)$ & Stiefel--Whitney, Chern and Pontryagin classes of $\zeta$\\
$[X,Y]$ & homotopy classes of maps $X\to Y$\\
$X\nrightarrow Y$ & stable map, i.e. a map $S^nX\to S^nY$ for some $n$\\
$\{X,Y\}$ & stable homotopy classes of maps $X\nrightarrow Y$, i.e. $\liminfty n[S^nX,S^nY]$\\
$\pi^s(n)$ & $n$-th stable homotopy group of spheres, i.e. $\liminfty k\pi_{n+k}(S^k)$\\
$[G]_p$ & $p$-primary part of the Abelian group $G$, i.e. the quotient of $G$ by the subgroup of torsion elements of order coprime to $p$\\
$[R]^n$ & degree-$n$ part of the graded ring $R$\\
$S_n$ & symmetric group on $n$ elements\\
$\Z_n$ & cyclic group of $n$ elements\\
$E_{i,j}^n\implies G_{i+j}$ & a spectral sequence with $n$-th page $(E_{i,j}^n)_{i,j\in\Z}$ converges to $G_*$, i.e. $\underset{i+j=k}\bigoplus E^\infty_{i,j}$ is associated to the group $G_k$ for all indices $k$\\
$\CC_\varphi$ & Serre class of groups of order some combination of the primes that satisfy the condition $\varphi$\\
$\congc{}$ & isomorphism modulo the Serre class $\CC$\\
$\cong,\congq$ & homotopy equivalence or isomorphism and rational homotopy equivalence or rational isomorphism\\
$\approx$ & diffeomorphism\\
$df_p$ & differential of the map $f$ in the point $p$\\
$\partial f$ & the restriction $f|_{\partial X}$ of the map $f\colon(X,\partial X)\to(Y,\partial Y)$\\
$f_\#,f_*,f^*$ & homomorphisms induced by the map $f$ in the homotopy groups, homologies and cohomologies\\
$f_!$ & Gysin pushforward in cohomologies induced by the map $f$\\
$\Tp_\eta$ & Thom polynomial of the singularity $\eta$\\
$p_k(n)$ & number of partitions of $n$ to sums of positive integers not greater than $k$\\
$\alpha_k(n)$ & sum of digits of $n$ in $k$-adic system\\
$\nu_p(n)$ & exponent of the highest power of the prime $p$ that divides $n$
\end{longtable}

%-------------------------------------------------------------------------------------------------------------------------------
%-------------------------------------------------------------------------------------------------------------------------------

\newpage
\chapter{Preliminaries}
\pagenumbering{arabic}

%-------------------------------------------------------------------------------------------------------------------------------

\section{Singular maps}

We will consider cobordisms of maps with fixed stable singularity classes. As the word ``singularity'' is defined in slightly different ways depending on where it is used, we first give the definition we will use.

\begin{defi}\label{sing}
By singularity class (or simply singularity) we mean the equivalence class of a germ
$$\eta\colon(\R^c,0)\to(\R^{c+k},0)$$
where two germs are defined to be equivalent if some combination of the following two conditions holds:
\begin{itemize}
\item[(i)] The germs $\eta$ and $\eta'$ are $\mathscr{A}$-equivalent (also called left-right equivalent). That is, there is a commutative diagram
$$\xymatrix{
(\mathbb{R}^c,0)\ar[r]^\eta\ar[d]_\varphi & (\mathbb{R}^{c+k},0)\ar[d]^\psi \\
(\mathbb{R}^c,0)\ar[r]^{\eta'} & (\mathbb{R}^{c+k},0)
}$$
where $\varphi$ is the germ of a diffeomorphism of $(\mathbb{R}^c,0)$ and $\psi$ is the germ of a diffeomorphism of $(\mathbb{R}^{c+k},0)$.
\item[(ii)] The germ $\eta'$ is the suspension of $\eta$, that is,
$$\textstyle\eta'=\eta\times\id_{\mathbb{R}^1}\colon(\mathbb{R}^c\times\mathbb{R}^1,0)\to(\mathbb{R}^{c+k}\times\mathbb{R}^1,0);~(x,t)\mapsto(\eta(x),t).$$
\end{itemize}
The singularity class of $\eta$ will be denoted by $[\eta]$.
\end{defi}

We observe that in any singularity class the dimension $c$ in the germs $(\mathbb{R}^c,0)\to(\mathbb{R}^{c+k},0)$ is not fixed, but the codimension $k$ is fixed. Throughout this thesis we consider singularities of positive codimension (i.e. $k>0$ in the above formula), therefore, if not specified otherwise, the codimensions will always assumed to be positive.

Moreover, we will always assume that the germs are stable in the sense that a small perturbation in the space of $C^\infty$-maps does not change their type (an exact definition of stable germs can be found for example in the beginning of \cite{rrthesis}, see also \cite{math}).

\begin{defi}
We say that the germ $\eta\colon(\mathbb{R}^c,0)\to(\mathbb{R}^{c+k},0)$ has an isolated singularity at the origin if there is a neighbourhood $U$ of the origin such that at no point in $U\setminus\{0\}$ the germ of $\eta$ is equivalent to that at the origin. Such an $\eta$ is called the root of the singularity class $[\eta]$.
\end{defi}

\begin{rmk}
The root of a singularity is also characterised by having the smallest possible dimension $c$ for any germ $(\mathbb{R}^c,0)\to(\mathbb{R}^{c+k},0)$ in its class.
\end{rmk}

Now we define the singularities of a map between manifolds by looking at the germs of this map at each point of the source (where the map is locally identified with a map $(\mathbb{R}^n,0)\to(\mathbb{R}^{n+k},0)$).

\begin{defi}
For a smooth map $f\colon M^n\to P^{n+k}$ a point $p\in M$ is called an $[\eta]$-point (or simply $\eta$-point) if the germ of $f$ at $p$ is equivalent to the germ $\eta\colon(\mathbb{R}^c,0)\to(\mathbb{R}^{c+k},0)$. The set of $\eta$-points in $M$ is denoted by $\eta(f)$.
\end{defi}

\begin{prop}
If $\eta\colon(\mathbb{R}^c,0)\to(\mathbb{R}^{c+k},0)$ is the root of the singularity class $[\eta]$, then for any map $f\colon M^n\to P^{n+k}$ the set $\eta(f)\subset M$ is a submanifold of dimension $n-c$.
\end{prop}

\begin{prf}
For any point $p\in\eta(f)$, the map $f$ restricted to a small neighbourhood of $p$ is such that there is a commutative diagram
$$\xymatrix{
(U,p)\ar[r]^{f|_U}\ar[d]_\varphi & (V,f(p))\ar[d]^\psi \\
(\mathbb{R}^n,0)\ar[r]^{\eta'} & (\mathbb{R}^{n+k},0)
}$$
where $U\subset M^n$, $V\subset P^{n+k}$, $\varphi$ and $\psi$ are diffeomorphisms and $\eta'$ is the $(n-c)$-times suspension of $\eta$. This means $\eta'=\eta\times\id_{\mathbb{R}^{n-c}}\colon(\mathbb{R}^c\times\mathbb{R}^{n-c},0)\to(\mathbb{R}^{c+k}\times\mathbb{R}^{n-c},0)$, hence all points of the factor $\mathbb{R}^{n-c}$ in the domain are $\eta$-points and there are no other $\eta$-points in a neighbourhood of $\mathbb{R}^{n-c}$.

Therefore the neighbourhood $U\subset M$ of $p$ can be chosen so that $\eta(f)\cap U=\varphi^{-1}(\mathbb{R}^{n-c})$. We can of course apply this construction for all points $p\in\eta(f)$, which gives an $(n-c)$-manifold structure on $\eta(f)$.
\end{prf}

\begin{defi}\label{taumap}
Given a set $\tau$ of singularities of a fixed codimension $k$, we say that a map $f\colon M^n\to P^{n+k}$ is a $\tau$-map if the germ of $f$ at any point of $M$ belongs to a singularity class in $\tau$. If $M$ is a manifold with boundary, then we assume the following extra property for any $\tau$-map: The target manifold $P$ also has a boundary and $f$ maps $\partial M$ to $\partial P$, moreover, the boundaries have collar neighbourhoods $\partial M\times[0,\varepsilon)$ and $\partial P\times[0,\varepsilon)$ such that
$$\textstyle f|_{\partial M\times[0,\varepsilon)}=g\times\id_{[0,\varepsilon)}\colon\partial M\times[0,\varepsilon)\to\partial P\times[0,\varepsilon)$$
for a $\tau$-map $g\colon\partial M\to\partial P$.
\end{defi}

\subsection{Singularity strata}\label{strata}

An important categorisation of singularities is the so-called Thom--Boardman type defined as follows.

\begin{defi}\label{thomboar}
A singularity class $[\eta]$ is of type $\Sigma^i$ (for $i=0,1,\ldots$) if the rank of the differential of the germ $\eta\colon(\mathbb{R}^c,0)\to(\mathbb{R}^{c+k},0)$ at the origin drops by $i$, that is $\rk d\eta_0=c-i$. Singularities of type $\Sigma^1$ are also called Morin singularities.
\end{defi}

\begin{rmk}
Clearly the Thom--Boardman type does not depend on the choice of representative of $[\eta]$ and the $i$ above can be at most the dimension of the source manifold of the root of $[\eta]$.
\end{rmk}

Now we introduce the Thom--Boardman stratification; for the proofs see \cite{boar}.

\begin{defi}
For a map $f\colon M^n\to P^{n+k}$ we define for $i=0,\ldots,n$
$$\Sigma^i(f):=\{p\in M^n\mid\rk df_p=n-i\}$$
the set of points where the singularity of $f$ is of type $\Sigma^i$. We call the points of $\Sigma^0(f)$ regular and the points of $\Sigma(f):=\Sigma^1(f)\cup\ldots\cup\Sigma^n(f)$ singular.
\end{defi}

This way we get a stratification of the source manifold
$$M=\Sigma^0(f)\cup\ldots\cup\Sigma^n(f).$$
By resticting $f$ to the strata we get further stratifications with $\Sigma^{i,j}(f):=\Sigma^j(f|_{\Sigma^i(f)})$ and we can iterate this process to obtain $\Sigma^{i,j,k}(f)$, $\Sigma^{i,j,k,l}(f)$, etc.

On the set of stable singularities with a fixed codimension $k$ a natural partial order arises in the following way.

\begin{defi}\label{pos}
For two singularity classes $[\eta]\ne[\vartheta]$ of codimension $k$, we call $[\eta]$ more complicated than $[\vartheta]$ if in any neighbourhood of any $\eta$-point there is a $\vartheta$-point. We denote this relation by $[\eta]>[\vartheta]$ and we put $\partial[\eta]:=\{[\vartheta]<[\eta]\}$.
\end{defi}

Throughout this thesis we will use the technical assumption that for any singularity $[\eta]$ any decreasing chain in $\partial[\eta]$ is finite (i.e. the partial order is well-founded).

\begin{rmk}
If $f\colon M^n\to P^{n+k}$ is a $\tau$-map (for some set of singularities $\tau$) and $[\eta]\in\tau$, then we also have $[\vartheta]\in\tau$ for all singularities $[\vartheta]<[\eta]$. Hence from now on we will assume all singularity sets to be decreasing.
\end{rmk}

\begin{ex}
\begin{itemize}
\item[]
\item[(1)] The lowest ``singularity'' is $\Sigma^0$, the class of germs of regular (non-singular) maps.
\item[(2)] Fold singularity (i.e. singularity of type $\Sigma^{1,0}$) is lower than cusp singularity (i.e. that of type $\Sigma^{1,1,0}$).
\end{itemize}
\end{ex}

\begin{figure}[H]
\begin{center}
\centering\includegraphics[scale=0.1]{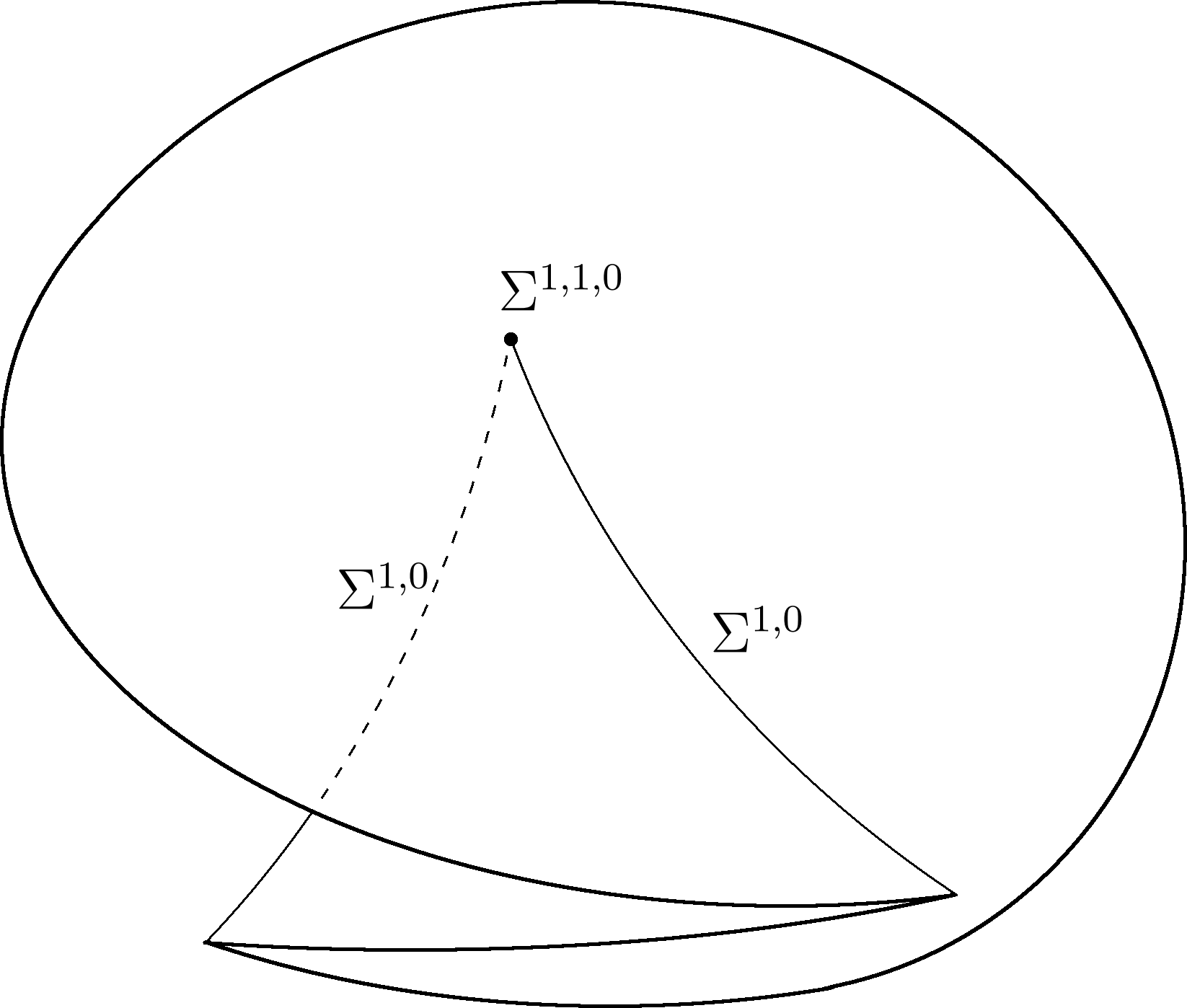}
\begin{changemargin}{2cm}{2cm} 
\caption{\hangindent=1.4cm\small Here we show a cusp germ $(\mathbb{R}^2,0)\to(\mathbb{R}^2,0)$ (with one cusp point). The two indicated curves are the fold points and everything else is regular.}\label{kep1}
\end{changemargin} 
\vspace{-1.5cm}
\end{center}
\end{figure}

Now if we have again a map $f\colon M^n\to P^{n+k}$, then $\eta(f)\subset\overline{\vartheta(f)}$ whenever $[\eta]>[\vartheta]$ and $M$ is the union of submanifolds of the form $\eta(f)$ (where the $\eta$ are the singularities of $f$). This way we get another stratification of $M$ with strata $\eta(f)$ that refines the Thom--Boardman stratification.

\begin{ex}\label{morsing}
If the $f$ above is a Morin map (i.e. it only has singularities of types $\Sigma^0$ and $\Sigma^1$), then the singularity strata are $\Sigma^0(f)$ (regular points), $\Sigma^{1,0}(f)$ (fold points), $\Sigma^{1,1,0}(f)$ (cusp points), etc. The symbol $\Sigma^{1,\ldots,1,0}$ where the number of $1$-s is $r$ will be shortened to $\Sigma^{1_r}$. The above stratification follows from the fact that each type $\Sigma^{1_r}$ only contains one singularity class and if $r\ne s$, then the singularity classes $\Sigma^{1_r}$ and $\Sigma^{1_s}$ are different (see \cite{mor}) and Morin singularities form a single increasing sequence $\Sigma^{1,0}<\Sigma^{1,1,0}<\ldots$
\end{ex}

\subsection{Multisingularities}

Until this point we only considered local restrictions on maps (the singularities of $f$ only describe how $f$ behaves locally) but later on we will sometimes also need global restrictions. This gives rise to the notion of multisingularities which describe how many points of each singularity type a map can have.

\begin{defi}
A multisingularity is a formal sum of singularities (of a fixed codimension) $m_1[\eta_1]+\ldots+m_r[\eta_r]$ with coefficients $m_1,\ldots,m_r\in\mathbb{N}$. A multisingularity of the form $1[\eta]$ will be called monosingularity.
\end{defi}

In order to distinguish the notation of multisingularities and the notation of singularities (with no prescribed multiplicites), we will typically use the convention to underline letters that denote multisingularities or sets of multisingularities.

\begin{defi}
Let $\underline\eta=m_1[\eta_1]+\ldots+m_r[\eta_r]$ be a multisingularity. For a smooth map $f\colon M^n\to P^{n+k}$ a point $p\in M$ is called an $\underline\eta$-point if $f^{-1}(f(p))$ consists of $m_1+\ldots+m_r$ points and out of all germs of $f$ at these points $m_i$ are equivalent to the germ $\eta_i$ (for $i=1,\ldots,r$). The set of $\underline\eta$-points in $M$ is denoted by $\underline\eta(f)$.
\end{defi}

The following are direct analogues of the notions we introduced earlier.

\begin{defi}
Given a set $\underline\tau$ of multisingularities of a fixed codimension $k$, we say that a map $f\colon M^n\to P^{n+k}$ is a $\underline\tau$-map if each point of $M$ belongs to $\underline\eta(f)$ for some $\underline\eta\in\underline\tau$. For manifolds with boundaries we have the analogous extra condition as in definition \ref{taumap}.
\end{defi}

\begin{ex}
The $\{1\Sigma^0\}$-maps are the embeddings, $\{1\Sigma^0,2\Sigma^0\}$-maps are the immersions with no triple points, etc. The $\{1\Sigma^0,2\Sigma^0,\ldots\}$-maps are all immersions, which is the same as $\{\Sigma^0\}$-maps (where no multiplicity is given).
\end{ex}

\begin{rmk}
If $\tau$ is a set of singularities and $\underline\tau$ is the set of all possible formal sums of the elements of $\tau$ (with coefficients in $\mathbb{N}$), then $\underline\tau$-maps are the same as $\tau$-maps.
\end{rmk}

\begin{defi}\label{poms}
For two multisingularities $\underline\eta\ne\underline\vartheta$ of codimension $k$, we call $\underline\eta$ more complicated than $\underline\vartheta$ if in any neighbourhood of any $\underline\eta$-point there is a $\underline\vartheta$-point. We denote this relation by $\underline\eta>\underline\vartheta$ and we put $\partial\underline\eta:=\{\underline\vartheta<\underline\eta\}$.
\end{defi}

\begin{rmk}
\begin{itemize}
\item[]
\item[(1)] We have $\underline\eta\ge\underline\vartheta$ whenever $\underline\eta=m_1[\eta_1]+\ldots+m_r[\eta_r],\underline\vartheta=n_1[\vartheta_1]+\ldots+n_r[\vartheta_r]$ and $m_i\ge n_i,[\eta_i]\ge[\vartheta_i]$ (for $i=1,\ldots,r$).
\item[(2)] If $f\colon M^n\to P^{n+k}$ is a $\underline\tau$-map (for some set of multisingularities $\underline\tau$) and $\underline\eta\in\underline\tau$, then we also have $\underline\vartheta\in\underline\tau$ for all multisingularities $\underline\vartheta<\underline\eta$. Hence from now on we will assume all multisingularity sets to be decreasing.
\item[(3)] For a map $f\colon M^n\to P^{n+k}$ the sets of points of the same multisingularities stratify the manifold $M$ in the same way as the singularities (with no prescribed multiplicites) do. This stratification refines the singularity stratification.
\end{itemize}
\end{rmk}

\begin{figure}[H]
\begin{center}
\centering\includegraphics[scale=0.1]{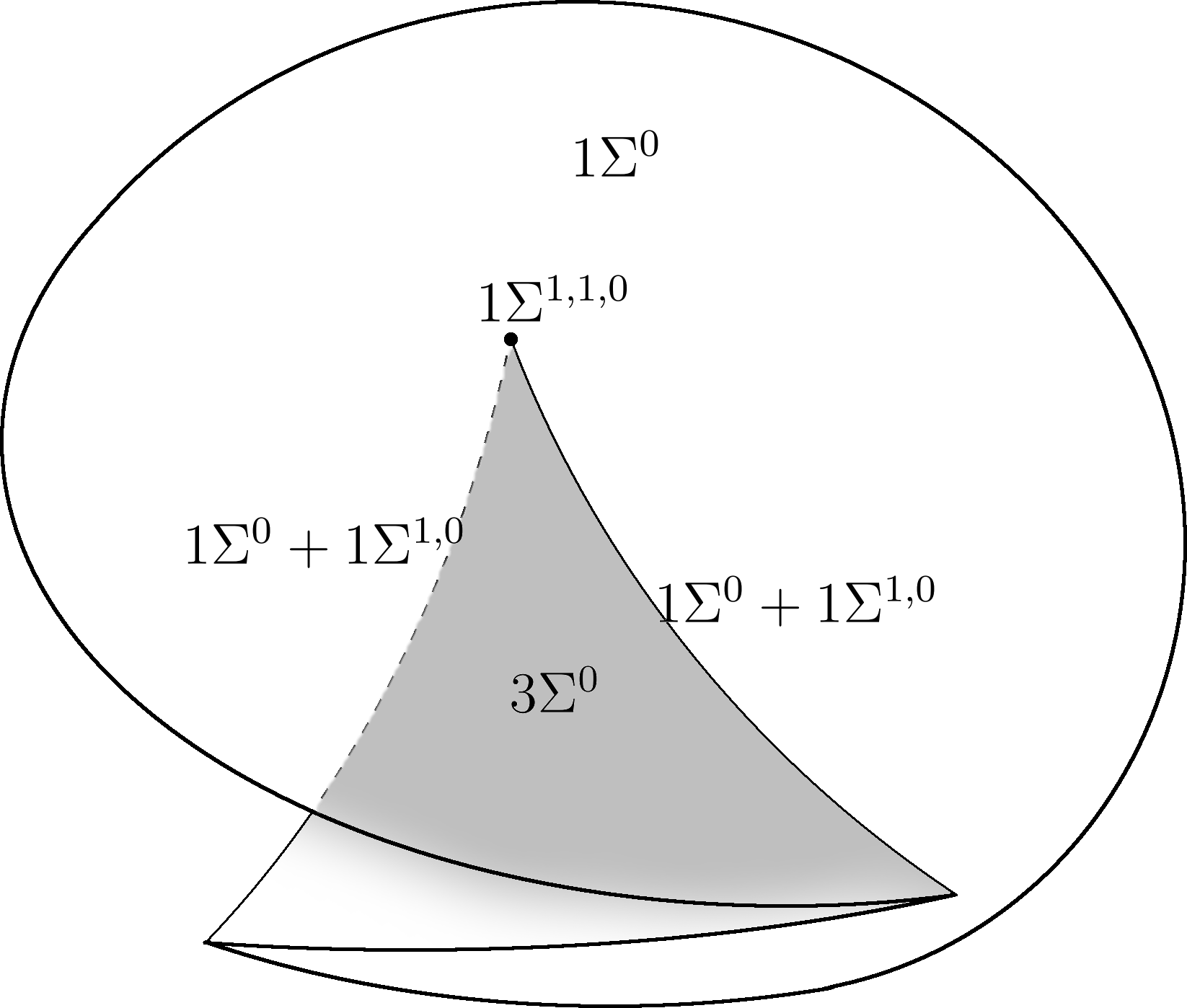}\label{kep2}
\begin{changemargin}{2cm}{2cm} 
\caption{\hangindent=1.4cm\small Here we show the same germ $(\mathbb{R}^2,0)\to(\mathbb{R}^2,0)$ as earlier, but with all multisingularities indicated.}
\end{changemargin} 
\vspace{-1.5cm}
\end{center}
\end{figure}

%-------------------------------------------------------------------------------------------------------------------------------

\section{Semiglobal description of singular maps}\label{semiglob}

In this section we will investigate how a singular map behaves in a small neighbourhood of a (multi)singularity stratum. This is not a central question in the present thesis, therefore the proofs will be omitted here.

\begin{defi}\label{A}
Put $\mathscr{A}:=\Diff_0(\mathbb{R}^c)\times\Diff_0(\mathbb{R}^{c+k})$ where $\Diff_0(\mathbb{R}^n)$ denotes the group of diffeomorphism germs of $(\mathbb{R}^n,0)$. We define the (left) action of the group $\mathscr{A}$ on the set of stable germs $\eta\colon(\mathbb{R}^c,0)\to(\mathbb{R}^{c+k},0)$ by the formula
$$(\varphi,\psi)\colon\eta\mapsto\psi\circ\eta\circ\varphi^{-1}~~~~(\varphi,\psi)\in\mathscr{A}.$$
\end{defi}

\begin{defi}\label{orbstab}
The orbits and stabilisers of the above $\mathscr{A}$-action also have the following names:
\begin{itemize}
\item[(1)] Two germs $(\mathbb{R}^c,0)\to(\mathbb{R}^{c+k},0)$ are called $\mathscr{A}$-equivalent if they are in the same $\mathscr{A}$-orbit (see definition \ref{sing}).
\item[(2)] The stabiliser of the germ $\eta\colon\mathbb{R}^c\to\mathbb{R}^{c+k}$ is called the automorphism group of $\eta$, that is
$$\Aut_\mathscr{A}\eta:=\{(\varphi,\psi)\in\mathscr{A}\mid\psi\circ\eta=\eta\circ\varphi\}.$$
\end{itemize}
\end{defi}

Let us fix a singularity $[\eta]$ where $\eta\colon\mathbb{R}^c\to\mathbb{R}^{c+k}$ is the root of its class. In order to investigate how a map behaves around an $[\eta]$-stratum we need some knowledge on $\Aut_\mathscr{A}\eta$. This group does not admit any convenient topology, however it still shares some properties with Lie groups; we will now introduce these.

\begin{defi}\label{cpct}
A subgroup of $\Aut_\mathscr{A}\eta$ is called compact if it is conjugate in $\mathscr{A}$ with a compact linear group.
\end{defi}

The following important theorem is the result of Jänich \cite{jan} and Wall \cite{wall} and we will not prove it here.

\begin{thm}\label{janwall}
Any compact subgroup of $\Aut_\mathscr{A}\eta$ is contained in a maximal compact subgroup and any two maximal compact subgroups are conjugate in $\Aut_\mathscr{A}\eta$.
\end{thm}

We will denote a maximal compact subgroup of $\Aut_\mathscr{A}\eta$ by $G_\eta$. Then $G_\eta$ has natural representations $\lambda$ and $\tilde\lambda$ on the source and target spaces respectively, defined by
$$\begin{array}{l}
\lambda(\varphi,\psi):=\varphi\colon(\mathbb{R}^c,0)\to(\mathbb{R}^c,0)\\
\tilde\lambda(\varphi,\psi):=\psi\colon(\mathbb{R}^{c+k},0)\to(\mathbb{R}^{c+k},0)
\end{array}
~~~~(\varphi,\psi)\in G_\eta.$$
By possibly choosing another representative in the $\mathscr{A}$-equivalence class of $\eta$, we can assume that the images of $\lambda$ and $\tilde\lambda$ are subgroups of $\O(c)$ and $\O(c+k)$ respectively.

\begin{defi}\label{phieta}
Denote the universal principal $G_\eta$-bundle by $EG_\eta\xrightarrow{~G_\eta~}BG_\eta$. Using the above representations we associate to them the vector bundles
$$\xi_\eta:=EG_\eta\tightunderset{\lambda}\times\mathbb{R}^c~~~~\text{and}~~~~\tilde\xi_\eta:=EG_\eta\tightunderset{\tilde\lambda}\times\mathbb{R}^{c+k}.$$
We define the fibrewise (but generally not linear) map $\Phi_\eta\colon\xi_\eta\to\tilde\xi_\eta$ using the following diagram (where we obtain the inner square by factoring the outer square by the $G_\eta$-action):
$$\xymatrix@R=.333pc{
EG_\eta\times\mathbb{R}^c\ar[rrr]^{\id_{EG_\eta}\times\eta}\ar[dr]^(.65)\lambda\ar[ddddd]_{\pr_{EG_\eta}}&&&EG_\eta\times\mathbb{R}^{c+k}\ar[dl]_(.65){\tilde\lambda}\ar[ddddd]^{\pr_{EG_\eta}}\\
&\xi_\eta\ar@{-->}[r]^{\Phi_\eta}\ar[ddd]&\tilde\xi_\eta\ar[ddd]&\\ \\ \\
&BG_\eta\ar[r]^{\id_{BG_\eta}}&BG_\eta&\\
EG_\eta\ar[ur]^{G_\eta}\ar[rrr]_{\id_{EG_\eta}}&&&EG_\eta\ar[ul]_{G_\eta}
}$$
Here $\eta$ is invariant under the $G_\eta$-action, therefore $\Phi_\eta$ is well-defined and its restriction to any fibre is $\mathscr{A}$-equivalent to $\eta$. This is called the global normal form of $[\eta]$.
\end{defi}

\begin{rmk}
Actually $\id_{EG_\eta}\times\eta$ and therefore also $\Phi_\eta$ is only defined on a neighbourhood of the zero-section, but we can assume (for simplicity) that this neighbourhood is the whole total space.
\end{rmk}

The map $\Phi_\eta$ is the key to the semiglobal description of a map around the $1[\eta]$-stratum because of the following property. This is the main theorem of this section, but again we will not give a proof here; for the proof see \cite{rrthesis} or \cite{rsz}.

\begin{thm}\label{univ}
If the map $f\colon M^n\to P^{n+k}$ restricted to the $\eta$-stratum is an embedding (i.e. $\eta(f)=1[\eta](f)$), then there are tubular neighbourhoods $\eta(f)\subset T\subset M$ and $f(\eta(f))\subset\tilde T\subset P$ of the $\eta$-stratum and its image such that there is a commutative diagram of fibrewise isomorphisms (the vertical arrows are the projections)
$$\xymatrix@R=.333pc{
\xi_\eta\ar[rrr]^{\Phi_\eta}\ar[ddddd]&&&\tilde\xi_\eta\ar[ddddd]\\
&T\ar[r]^{f|_T}\ar[ul]\ar[ddd]&\tilde T\ar[ur]\ar[ddd]&\\ \\ \\
&\eta(f)\ar[r]^{f|_{\eta(f)}}\ar[dl]&f(\eta(f))\ar[dr]&\\
BG_\eta\ar[rrr]_{\id_{BG_\eta}}&&&BG_\eta
}$$
\end{thm}

In other words this theorem states that $\Phi_\eta\colon\xi_\eta\to\tilde\xi_\eta$ is the universal object for the maps of tubular neighbourhoods of the $1[\eta]$-stratum into those of its image. This theorem has an analogue for more complicated multisingularities as we will see in the next subsection.

Although the proof is omitted here, we remark that it is based on the contractibility of the quotient $\Aut_\mathscr{A}\eta/G_\eta$ (in a generalised sense) which is the result of Rimányi \cite{rrthesis}. This is again a property of $\Aut_\mathscr{A}\eta$ common with Lie groups.

%Before we get to the proof we will first need some preparations.

%\begin{defi}
%Let $G$ be any subgroup of $\Aut_\mathscr{A}\eta$.
%\begin{itemize}
%\item[(1)] If $M$ is a manifold (with boundary), we call a map $a\colon M\to\Aut_\mathscr{A}\eta/G$ smooth if $M$ can be covered by open sets with the following property: For any $U\subset M$ in the cover, $a|_U$ can be represented by pairs $(\alpha_c,\alpha_{c+k})$ such that $\alpha_i$ is a diffeomorphism germ of $(U\times\mathbb{R}^i,U\times\{0\})$ that maps all fibres $\{p\}\times\mathbb{R}^i$ to themselves ($p\in U, i=c,c+k$).
%\item[(2)] We call $\Aut_\mathscr{A}\eta/G$ contractible if for any manifold $M$ with boundary any smooth map $\partial M\to\Aut_\mathscr{A}\eta/G$ can be extended to a smooth map $M\to\Aut_\mathscr{A}\eta/G$.
%\end{itemize}
%\end{defi}

%The next theorem also shows that $\Aut_\mathscr{A}\eta$ has a lot of properties common with Lie groups. This was proved by Rimányi in \cite{rrthesis} and we will not give a proof here.

%\begin{thm}
%\begin{itemize}
%\item[]
%\item[\rm{(1)}] If $G$ is a maximal compact subgroup of $\Aut_\mathscr{A}\eta$, then $\Aut_\mathscr{A}\eta/G$ is contractible.
%\item[\rm{(2)}] There is a section $\sigma\colon\Aut_\mathscr{A}\eta/G\to\Aut_\mathscr{A}\eta$ such that any smooth map $M\to\Aut_\mathscr{A}\eta/G$ composed with $\sigma$ is also smooth.
%\end{itemize}
%\end{thm}

\subsection{Automorphism groups of multisingularities}

The previous constructions show how we can give a semiglobal description of a map around a monosingularity stratum. Now we will see that everything works similarly if we take a multisingularity stratum instead.

Fix a multisingularity $\underline\eta=m_1[\eta_1]+\ldots+m_r[\eta_r]$ where $m_i\in\N$ and $\eta_i\colon(\mathbb{R}^{c_i},0)\to(\mathbb{R}^{c_i+k},0)$ is the root of its class (for $i=1,\ldots,r$). Observe that $\underline\eta$ can be identified with the equivalence class of a germ
$$\eta_I\colon(\mathbb{R}^{\sum_{i=1}^rm_ic_i+\left(\sum_{i=1}^rm_i-1\right)k},S)\to(\mathbb{R}^{\sum_{i=1}^rm_ic_i+\sum_{i=1}^rm_ik},0)$$
where $S\subset\mathbb{R}^{\sum_{i=1}^rm_ic_i+\left(\sum_{i=1}^rm_i-1\right)k}$ contains $m_1+\ldots+m_r$ points and at $m_i$ of these points the germ $\eta_I$ is equivalent to the germ $\eta_i$ (for $i=1,\ldots,r$).

Now we can define the analogues of the notions we had for monosingularities. To make notations simpler we put $c:=\overset{r}{\underset{i=1}\sum}m_ic_i+\left(\overset{r}{\underset{i=1}\sum}m_i-1\right)k$. If we define $\mathscr{A}:=\Diff_S(\mathbb{R}^c)\times\Diff_0(\mathbb{R}^{c+k})$ where $\Diff_S(\mathbb{R}^c)$ denotes the group of germs $(\mathbb{R}^c,S)\to(\mathbb{R}^c,S)$ that are diffeomorphisms near the points of $S$, then definitions \ref{A}, \ref{orbstab} and \ref{cpct} are the same for $\eta_I$. Then $\underline\eta$ is the equivalence class of $\eta_I$ in the equivalence relation generated by $\mathscr{A}$-equivalence and suspension.

Theorem \ref{janwall} is also true in this setting and the analogous maximal compact subgroup $G_{\underline\eta}$ of the automorphism group also has representations on the source and target spaces. This way we can define the vector bundles $\xi_{\underline\eta}$ and $\tilde\xi_{\underline\eta}$ and the fibrewise map $\Phi_{\underline\eta}\colon\xi_{\underline\eta}\to\tilde\xi_{\underline\eta}$ as earlier, and the analogue of theorem \ref{univ} is also true. A more detailed description of this can be found in \cite{rrthesis}.

%-------------------------------------------------------------------------------------------------------------------------------

\section{Cobordism groups}\label{cobsec}

Now we have arrived to the point where we can begin the global description of singular maps, which is the main subject of this thesis. We will define everything for singularities with no prescribed multiplicities, but we refer to remark \ref{singmulti}: Almost everything in the following two sections has a direct analogue for multisingularities.

\begin{defi}\label{cobtau}
Let $\tau$ be a set of singularities of a fixed codimension $k$. We call two $\tau$-maps $f_0\colon M_0^n\to P^{n+k}$ and $f_1\colon M_1^n\to P^{n+k}$ (with closed source manifolds $M_0$ and $M_1$) $\tau$-cobordant if there is
\begin{itemize}
\item[(i)] a compact manifold with boundary $W^{n+1}$ such that $\partial W=M_0\sqcup M_1$,
\item[(ii)] a $\tau$-map $F\colon W^{n+1}\to P\times[0,1]$ such that for $i=0,1$ we have $F^{-1}(P\times\{i\})=M_i$ and $F|_{M_i}=f_i$.
\end{itemize}
The set of $\tau$-cobordism classes of $\tau$-maps to the manifold $P^{n+k}$ will be denoted by $\Cob_\tau(n,P^{n+k})$, and in the case $P^{n+k}=\mathbb{R}^{n+k}$ we abbreviate this notation to $\Cob_\tau(n,k)$.
\end{defi}

We will use the following conventions throughout this thesis:
\begin{enumerate}
\item\label{conv1} Observe that putting the $n$ in the notation of $\Cob_\tau(n,P^{n+k})$ is a bit redundant, as the codimension $k$ is fixed in $\tau$, so only $n$-manifolds can have $\tau$-maps to $(n+k)$-manifolds. Therefore in most cases we will omit it from the notation and put $\Cob_\tau(P^{n+k}):=\Cob_\tau(n,P^{n+k})$.
\item\label{conv2} A large part of this thesis consists of the investigation of Morin maps (see definition \ref{thomboar} and example \ref{morsing}), therefore we will use the following abbreviation:
$$\Cob_{\{\Sigma^0,\Sigma^{1_1},\ldots,\Sigma^{1_r}\}}(n,P)=:\Cob_r(n,P).$$
Here $r=\infty$ is also allowed, the symbol $\Cob_\infty(n,P)$ denotes the cobordism group of all Morin maps.
\item\label{conv3} We will always denote the cobordism relations by $\sim$ and the cobordism class of a map $f$ by $[f]$; the context will make it clear which cobordism relation they stand for.
\end{enumerate}

The aim of this thesis is to collect some of the more recent computations concerning the cobordism groups $\Cob_\tau(P)$ with various restrictions and extra conditions (the group structure will be defined in the next subsection).

\subsection{The group operation}\label{grp}

Fix a set of singularities $\tau$ of a fixed codimension $k$ and a manifold $P^{n+k}$. The cobordism set $\Cob_\tau(P)$ trivially admits a commutative semigroup operation by the disjoint union: If $f_0\colon M_0\to P$ and $f_1\colon M_1\to P$ are $\tau$-maps, then
$$f_0\sqcup f_1\colon M_0\sqcup M_1\to P$$
is also a $\tau$-map and $[f_0]+[f_1]:=[f_0\sqcup f_1]$ is well-defined. This operation also has a null element represented by the empty map.

In order to make this an Abelian group structure, we have to prove that all elements of $\Cob_\tau(P)$ have inverse elements. We fix a $\tau$-map $f\colon M^n\to P^{n+k}$ and look for the inverse of $[f]$.

\begin{prop}\label{inv}
If $P=N\times\mathbb{R}^1$ for a manifold $N^{n+k-1}$, then $[f]$ has an inverse in $\Cob_\tau(P)$.
\end{prop}

\begin{prf}
Let $\rho\colon N\times\mathbb{R}^1\to N\times\mathbb{R}^1$ denote the reflection to a hypersurface $N\times\{t\}$ for a $t\in\mathbb{R}^1$. Then the composition $\rho\circ f$ is such that $f\sqcup(\rho\circ f)\sim\varnothing$, since the manifold $M\times[0,1]$ can be mapped by a $\tau$-map to the path of the rotation of the image of $f$ to the image of $\rho\circ f$ around $N\times\{t\}$ in $P\times\mathbb{R}_+\approx P\times[0,1)$. Hence $[\rho\circ f]$ is the inverse of $[f]$.
\end{prf}

\begin{figure}[H]
\begin{center}
\centering\includegraphics[scale=0.1]{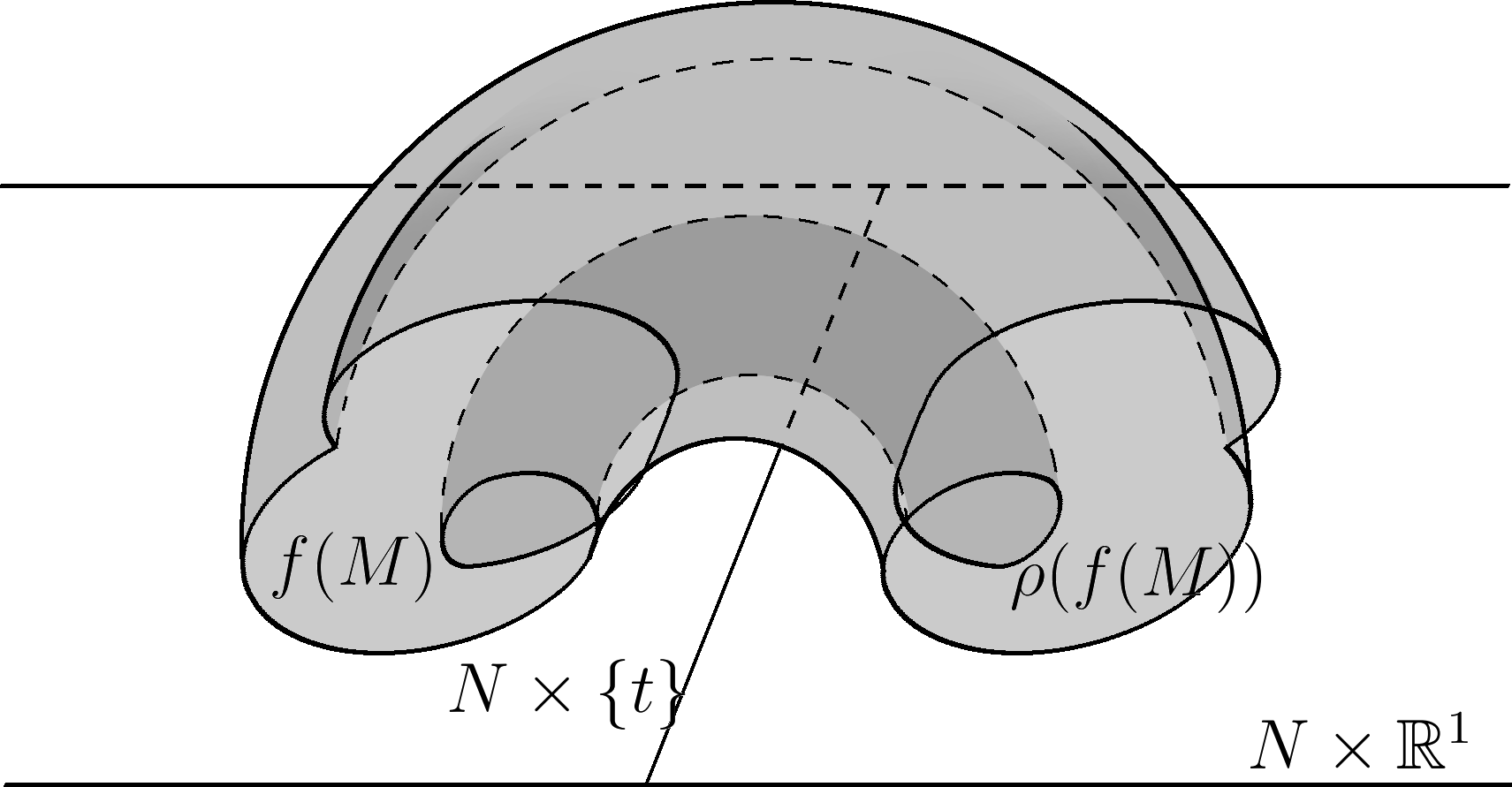}\label{kep3}
\begin{changemargin}{2cm}{2cm} 
\caption{\hangindent=1.4cm\small We indicated the rotation of the image of $f$ to the image of $\rho\circ f$ around $N\times\{t\}$ in $P\times\mathbb{R}^1_+=N\times\mathbb{R}^1\times\mathbb{R}^1_+$.}
\end{changemargin} 
\vspace{-1.3cm}
\end{center}
\end{figure}

If $P$ is an arbitrary manifold, then we need the following notion.

\begin{defi}\label{lframed}
An $l$-framed $\tau$-map from a manifold $M^n$ to a manifold $Q^{n+k+l}$ is the germ along $M\approx M\times\{0\}$ of a map
$$\tilde f\colon M^n\times\mathbb{R}^l\to Q^{n+k+l}$$
such that for all points $p\in M$ there are coordinate neighbourhoods $(p,0)\in U\times\tightoverset{_\circ}{D^l_\varepsilon}\subset M\times\mathbb{R}^l$ and $\tilde f(p,0)\in V\subset Q$ of the point $(p,0)\in M\times\mathbb{R}^l$ and its image, with the following property: $U\approx\mathbb{R}^n,V\approx\mathbb{R}^{n+k}\times\tightoverset{_\circ}{D^l_\varepsilon}$ and with these identifications we have
$$\textstyle\tilde f|_{U_p\times\toverset{_\circ}{D^l_\varepsilon}}=g\times\id_{\toverset{_\circ}{D^l_\varepsilon}}\colon\mathbb{R}^n\times\tightoverset{_\circ}{D^l_\varepsilon}\to\mathbb{R}^{n+k}\times\tightoverset{_\circ}{D^l_\varepsilon}$$
for a $\tau$-map $g\colon\mathbb{R}^n\to\mathbb{R}^{n+k}$.
\end{defi}

\begin{ex}
An $l$-framed $\{\Sigma^0\}$-map is an immersion with $l$ pointwise independent normal vector fields.
\end{ex}

The cobordism set of $l$-framed $\tau$-maps to the manifold $Q^{n+k+l}$ can be defined by an obvious modification of the definition of that of $\tau$-maps (\ref{cobtau}). This will be denoted by $\Cob_{\tau\oplus l}(n,Q^{n+k+l})$ or by $\Cob_{\tau\oplus l}(Q^{n+k+l})$.

\begin{thm}\label{susp}
For any $l\in\mathbb{N}$, if we assign to a $\tau$-map $f\colon M^n\to P^{n+k}$ the germ of the map $f\times\id_{\mathbb{R}^l}$ along $M$, we get a bijective correspondence
$$\Cob_\tau(n,P^{n+k})\to\Cob_{\tau\oplus l}(n,P^{n+k}\times\mathbb{R}^l).$$
\end{thm}

The proof of this theorem relies on the notion of $\tau$-embeddings (see section \ref{szpt}) which does not fit with the topic of the present section. Therefore we will only prove it later in remark \ref{tauembrmk}.

Now the cobordism class $[f]$ corresponds to its suspension in the cobordism set $\Cob_{\tau\oplus 1}(P\times\mathbb{R}^1)$, where the target manifold is of the form we needed in proposition \ref{inv}. Therefore in $\Cob_{\tau\oplus 1}(P\times\mathbb{R}^1)$ the inverse element exists, so it exists in $\Cob_\tau(P)$ too. Moreover, the construction of the inverse of $[f]$ will also be clear from the proof in subsection \ref{tauemb}. Hence we obtained the following.

\begin{crly}
$\Cob_\tau(P)$ is an Abelian group with the disjoint union operation.
\end{crly}

\subsection{Cobordisms with additional structures}

\begin{defi}
Let $L=\underset{k\to\infty}\lim L(k)$ be a stable linear group (see \cite[8.2]{stabgr}). We say that a map $f\colon M^n\to P^{n+k}$ is equipped with a normal $L$-structure, if the virtual normal bundle of $f$ has structure group $L$.
\end{defi}

\begin{rmk}
We can define the global normal form of a multisingularity $\underline\eta$ for maps with normal $L$-structures similarly to section \ref{semiglob}. Recall that the maximal compact subgroup $G_{\underline\eta}$ was a inside $\O(c)\times\O(c+k)$, hence we can take
$$G_{\underline\eta}^L:=\big\{(A,B)\in G_{\underline\eta}\mid\big(\begin{smallmatrix}
A&0\\
0&B
\end{smallmatrix}\big)\in L(2c+k)\big\}.$$
The corresponding fibrewise map between the universal $G_{\underline\eta}^L$-bundles will be denoted by $\Phi_{\underline\eta}^L\colon\xi_{\underline\eta}^L\to\tilde\xi_{\underline\eta}^L$. This way the virtual bundle $\tilde\xi_{\underline\eta}^L-\xi_{\underline\eta}^L$ has structure group $L$, so any virtual bundle that we get as a pullback of $\tilde\xi_{\underline\eta}^L-\xi_{\underline\eta}^L$ has it too. Now the map $\Phi_{\underline\eta}^L$ has the same universal property for the maps with normal $L$-structures around an $\underline\eta$-stratum as described in theorem \ref{univ}.
\end{rmk}

\begin{defi}\label{taug}
Let $\tau$ be a set of singularities of a fixed codimension $k$ and $L$ a stable linear group. The $(\tau,L)$-cobordism of two $\tau$-maps $f_0$ and $f_1$ equipped with normal $L$-structures is the same as in definition \ref{cobtau} with the extra condition that the map $F$, which joins $f_0$ and $f_1$, also has a normal $L$-structure which matches that of $f_0$ and $f_1$. The set of $(\tau,L)$-cobordism classes to the manifold $P^{n+k}$ will be denoted by $\Cob_\tau^L(n,P^{n+k})$.
\end{defi}

\begin{rmk}
The normal $L$-structure can be added to every notion in the previous subsection, hence $\Cob_\tau^L(n,P)$ is also a group.
\end{rmk}

In the remaining part of the present section we will describe those additional structures to cobordism groups which we will use throughout this thesis.

\begin{ex}
\begin{itemize}
\item[]
\item[(1)] If $L(k)=\O(k)$, then $\Cob_\tau^L(P)=\Cob_\tau(P)$ is the unoriented $\tau$-cobordism group.
\item[(2)] If $L(k)=\SO(k)$, then $\Cob_\tau^L(P)=\Cob_\tau^{\SO}(P)$ is the cooriented $\tau$-cobordism group.
\item[(3)] If $L(k)=\Spin(k)$, then $\Cob_\tau^L(P)=\Cob_\tau^{\Spin}(P)$ is the cobordism group of $\tau$-maps equipped with spin normal structures.
\end{itemize}
\end{ex}

\begin{rmk}
In the latter two examples above, if the manifold $P$ is oriented (resp. spin), then an orientation (resp. spin structure) of the virtual normal bundle of a map $f\colon M\to P$ is equivalent to an orientation (resp. spin structure) of the tangent bundle $TM$ and the same is true for maps to $P\times[0,1]$. Hence for an oriented (resp. spin) $P$ the group $\Cob_\tau^{\SO}(P)$ (resp. $\Cob_\tau^{\Spin}(P)$) is the cobordism group of $\tau$-maps of oriented (resp. spin) manifolds to $P$.
\end{rmk}

\begin{defi}
A map between manifolds $f\colon M\to P$ is called prim (projected immersion), if it is the cobination of an immersion $i_f\colon M\looparrowright P\times\mathbb{R}^1$ and the projection $\pr_P\colon P\times\mathbb{R}^1\to P$.
\end{defi}

\begin{rmk}\label{kerbundle}
Observe that a prim map $f$ is always a Morin map (i.e. $\dim\ker df\le1$) and the line bundle $\ker df|_{\Sigma(f)}$ over the set of singular points is trivialised. It is not hard to see that the converse is also true: A Morin map equipped with a trivialisation of its kernel line bundle is a prim map.
\end{rmk}

\begin{defi}\label{primcob}
Given a set $\tau$ of Morin singularities of a fixed codimension $k$ and a stable linear group $L$, the cobordism of prim $\tau$-maps to a manifold $P^{n+k}$ equipped with normal $L$-structures is the analogue of definition \ref{taug} for prim maps. The set (and group) of these cobordism classes will be denoted by $\Prim_\tau^L(n,P^{n+k})$.
\end{defi}

\begin{rmk}\label{primrmk}
We can give analogous definitions for the cobordisms of those prim $\tau$-maps $f\colon M\to P$ for which the normal bundle of the immersion $i_f\colon M\looparrowright P\times\mathbb{R}^1$ is equipped with a complex or quaternionic structure; we will denote these cobordism groups by $\Prim_\tau^{\U}(P)$ and $\Prim_\tau^{\Sp}(P)$ respectively.
\end{rmk}

We will use the same conventions for the notations of cobordisms with additional structures as described after definition \ref{cobtau}.

\begin{rmks}\label{singmulti}
Observe that almost everything in the section above works in exactly the same way when we replace the set $\tau$ of singularities by a set $\underline\tau$ of multisingularities. The only exception is theorem \ref{susp} which is not necessarily true for multisingularities, hence the cobordism set $\Cob_{\underline\tau}^L(P)$ is only a group if $P$ is a manifold of the form $N\times\mathbb{R}^1$.

We also remark that for a set $\underline\tau$ of multisingularities there is a natural forgetful map $\Cob_{\underline\tau}^L(P)\to\Cob_\tau^L(P)$ where $\tau$ is the set of all singularities in the elements of $\underline\tau$.
\end{rmks}

%-------------------------------------------------------------------------------------------------------------------------------
%-------------------------------------------------------------------------------------------------------------------------------

\chapter{Classifying spaces}\label{classp}

One of the most important elements of any type of cobordism theory is an analogue of the Pontryagin--Thom construction, which means constructing a bijection between the cobordism set and a set of homotopy classes to a specific space. This space is then called the classifying space for that cobordism theory.

In our setting the cobordism theory is $\Cob_{\underline\tau}^L(n,P^{n+k})$ where $\underline\tau$ is a set of (multi)singu-larities of codimension $k$ and $L$ is a stable linear group. We will obtain a space $X_{\underline\tau}^L$ which has the property
$$\Cob_{\underline\tau}^L(P)=[\toverset{^*}{P},X_{\underline\tau}^L]_*$$
where $\toverset{^*}{P}$ denotes the one-point compactification of $P$ and $[\cdot,\cdot]_*$ denotes the set of based homotopy classes of maps that fix the image of $\infty\in\toverset{^*}{P}$ (if $P$ was initially compact, then this homotopy set is just $[P,X_{\underline\tau}^L]$). Hence $X_{\underline\tau}^L$ is the classifying space for the cobordisms of singular maps.

Although this classifying space can also be obtained using Brown's representability theorem (see \cite{hosszu}), it is useful to see its construction. We only remark that this classifying space is homotopically unique by the Brown representability.

In section \ref{rszpt} we construct the classifying space $X_{\underline\tau}^L$ if $\underline\tau$ is a set of multisingularities (this is taken from \cite{rrthesis} and \cite{rsz}); in section \ref{szpt} the classifying space $X_\tau^L$ is constructed in a different way for a set $\tau$ of singularities without fixed multiplicities; then in section \ref{kazspace} we construct another type of classifying space, which will be called Kazarian's space (although it was already considered by Thom), then we show a connection between $X_\tau^L$ and Kazarian's space (the latter two sections are based on \cite{hosszu}).

%It is not hard to see (by a modified version of the Pontryagin--Thom construction) that the above classifying space also has the properties
%$$\Bord_{\underline\tau}^L(n,\mathfrak{N}_{n+k})=\mathfrak{N}_{n+k}(X_{\underline\tau}^L)~~~~\text{and}~~~~\Bord_{\underline\tau}^L(n,\Omega_{n+k})=\Omega_{n+k}(X_{\underline\tau}^L).$$

Analogous classifying spaces can be constructed for cobordisms of prim maps (see definition \ref{primcob}) by the same methods. The classifying space for $\Prim_{\underline\tau}^L(P)$ will be denoted by $\overline X_{\underline\tau}^L$.

%\section{Brown representability}

%-------------------------------------------------------------------------------------------------------------------------------

\section{The construction of Szűcs and Rimányi}\label{rszpt}

Throughout this section we fix a set $\underline\tau$ of multisingularities of codimension $k$ and a stable linear group $L$. We also fix a complete ordering of the elements of $\underline\tau$ which extends the natural partial order (see definition \ref{poms}); it will be denoted by $\underline\eta_0<\underline\eta_1<\ldots$ (where the $\underline\eta_i$ are the elements of $\underline\tau$). We can also assume this to be a well-ordering, as the natural partial order is well-founded (as a consequence of the assumption following definition \ref{pos}).

We will show a generalisation of the Pontryagin--Thom construction to $(\underline\tau,L)$-maps due to Szűcs and Rimányi. This means constructing a ``universal $(\underline\tau,L)$-map'' in the same sense as the inclusion $B\O(k)\hookrightarrow M\O(k)$ of the Grassmannian manifold into the Thom space of the universal rank-$k$ bundle is the universal embedding. In other words the following is true.

\begin{thm}\label{ptpullback}
There are topological spaces $X_{\underline\tau}^L$, $Y_{\underline\tau}^L$ and a continuous map $f_{\underline\tau}^L\colon Y_{\underline\tau}^L\to X_{\underline\tau}^L$ with the following properties:
\begin{enumerate}
\item\label{pt1} For any manifold with boundary $P^{n+k}$ and any closed manifold $N^{n-1}$, if there are maps $f$, $\kappa_N$ and $\kappa_P$ such that the outer square in the diagram below is a pullback diagram, then there is a manifold $M^n$ with boundary $\partial M=N$, an extension $\kappa_M$ of $\kappa_N$ and a $(\underline\tau,L)$-map $F\colon M\to P$ such that the upper inner square is a pullback diagram as well.
$$\xymatrix{
Y_{\underline\tau}^L\ar[rr]^{f_{\underline\tau}^L} && X_{\underline\tau}^L\\
& M\ar@{-->}[ul]_(.35){\kappa_M}\ar@{-->}[r]^F & P\ar[u]_{\kappa_P}\\
N\ar[uu]^{\kappa_N}\ar@{^(-->}[ur]\ar[rr]^f && \partial P\ar@{^(->}[u] 
}$$
\item\label{pt2} For any $(\underline\tau,L)$-map $F\colon M^n\to P^{n+k}$ between manifolds with boundaries, if there are maps $\kappa_{\partial M}$ and $\kappa_{\partial P}$ such that the outer square in the diagram below is a pullback diagram, then there are extensions $\kappa_M$ and $\kappa_P$ of $\kappa_{\partial M}$ and $\kappa_{\partial P}$ respectively such that the upper inner square is a pullback diagram as well.
$$\xymatrix{
Y_{\underline\tau}^L\ar[rrr]^{f_{\underline\tau}^L} &&& X_{\underline\tau}^L\\
& M\ar@{-->}[ul]_(.35){\kappa_M}\ar[r]^F & P\ar@{-->}[ur]^(.35){\kappa_P} &\\
\partial M\ar[uu]^{\kappa_{\partial M}}\ar@{^(->}[ur]\ar[rrr]^{F|_{\partial M}} &&& \partial P\ar@{_(->}[ul]\ar[uu]_{\kappa_{\partial P}}
}$$
\end{enumerate}
\end{thm}

\begin{rmk}\label{gentau}
We call this $f_{\underline\tau}^L$ the universal $(\underline\tau,L)$-map although it is not really a $\underline\tau$-map, as its source and target spaces are not finite dimensional manifolds. However, it will be constructed from direct limits of $(\underline\tau,L)$-maps, so it is a ``generalised $(\underline\tau,L)$-map''. It will follow from the proof that it does not only classify $(\underline\tau,L)$-maps (in the sense of the above theorem), but also generalised $(\underline\tau,L)$-maps obtained from direct limits of $(\underline\tau,L)$-maps.
\end{rmk}

Before proving the above theorem, we prove its most important corollary, which states that $X_{\underline\tau}^L$ is indeed the classifying space for the cobordisms of $(\underline\tau,L)$-maps.

\begin{thm}\label{ptclass}
For any manifold $P^{n+k}$ there is a bijective correspondence between $\Cob_{\underline\tau}^L(P)$ and $[\toverset{^*}{P},X_{\underline\tau}^L]_*$.
\end{thm}

We will only prove this if $P$ is a closed manifold, but the next subsection will make it clear how the same proof works for an open $P$ too; see remark \ref{1cptrmk}.

\medskip\begin{prf}
If $M^n$ is a closed manifold and $f\colon M\to P$ is a $(\underline\tau,L)$-map, then by \ref{pt2} we can assign to the map $f$ the map $\alpha(f):=\kappa_P$ (and the map $\beta(f):=\kappa_M$).

In order to make this $\alpha$ a correspondence between cobordism classes and homotopy classes, we have to prove the following: If $f_0\colon M_0\to P$ and $f_1\colon M_1\to P$ are $(\underline\tau,L)$-cobordant, then $\alpha(f_0)$ and $\alpha(f_1)$ are homotopic. This is indeed so, since a $(\underline\tau,L)$-cobordism $F\colon W\to P\times[0,1]$ is a map that satisfies the conditions of \ref{pt2}, hence there is a map $\alpha(F)=\kappa_{P\times[0,1]}$ which is clearly a homotopy between $\alpha(f_0)$ and $\alpha(f_1)$.

Hence we found a map $\tilde\alpha\colon\Cob_{\underline\tau}^L(P)\to[P,X_{\underline\tau}^L]$; we just have to prove that this $\tilde\alpha$ is bijective.

\medskip\noindent I. \emph{$\tilde\alpha$ is surjective.}\medskip

If we are given a map $\kappa\colon P^{n+k}\to X_{\underline\tau}^L$, then by \ref{pt1} we obtain a manifold $M^n$, a map $\kappa_M\colon M\to Y_{\underline\tau}^G$ and a $(\underline\tau,L)$-map $f\colon M\to P$, so the surjectivity is proved if $\alpha(f)$ is homotopic to $\kappa$. To prove this, we apply \ref{pt2} to the map 
$$\textstyle f\times\id_{[0,1]}\colon N\times[0,1]\to P\times[0,1]$$
where the boundaries are mapped by $\kappa\sqcup\alpha(f)\colon P\sqcup P\to X_{\underline\tau}^L$ and $\kappa_M\sqcup\beta(f)\colon N\sqcup N\to Y_{\underline\tau}^L$, which yields the desired homotopy.

\medskip\noindent II. \emph{$\tilde\alpha$ is injective.}\medskip

Suppose we have two maps $f_0\colon M_0\to P$ and $f_1\colon M_1\to P$ such that $\alpha(f_0)$ is homotopic to $\alpha(f_1)$. If $\kappa$ is a homotopy between them, then we can apply \ref{pt1} to the maps
$$f_0\sqcup f_1\colon M_0\sqcup M_1\to P\sqcup P=\partial(P\times[0,1])$$
and $\kappa\colon P\times[0,1]\to X_{\underline\tau}^L$ and $\beta(f_0)\sqcup\beta(f_1)\colon M_0\sqcup M_1\to Y_{\underline\tau}^L$. This gives a $(\underline\tau,L)$-cobordism between $f_0$ and $f_1$.
\end{prf}

\subsection{Proof of theorem \ref{ptpullback}}\label{prf}

The proof proceeds by (transfinite) induction on the elements of $\underline\tau$. We will glue blocks together to obtain the spaces $X_{\underline\tau}^L$ and $Y_{\underline\tau}^L$; one block corresponds to one multisingularity in $\underline\tau$. The block corresponding to $\underline\eta\in\underline\tau$ is glued precisely to the blocks corresponding to multisingularities less complicated than $\underline\eta$.

This, of course, only works if the index of $\underline\eta$ is a successor ordinal, but in the case when the order-type of $\underline\tau$ is a limit, we can define $f_{\underline\tau}^L\colon Y_{\underline\tau}^L\to X_{\underline\tau}^L$ as a direct limit. Hence from now on we will assume that $\underline\tau$ has a maximal element.

The first step in the induction is when we have $\underline\tau=\{1\Sigma^0\}$, that is, $(\underline\tau,L)$-maps are the embeddings with normal $L$-structures. In this case the classical Thom construction yields the inclusion
$$f_{\{1\Sigma^0\}}^L\colon Y_{\{1\Sigma^0\}}^L:=BL(k)\hookrightarrow T\gamma_k^L=:X_{\{1\Sigma^0\}}^L$$
of the zero-section into the Thom space of the universal rank $k$ vector bundle with structure group $L$. This way $T\gamma_k^L$ and $BL(k)$ are the first blocks of $X_{\underline\tau}^L$ and $Y_{\underline\tau}^L$ respectively (for any multisingularity set $\underline\tau$).

Now suppose we know the theorem for $\underline\tau':=\underline\tau\setminus\{\underline\eta\}$ (where $\underline\eta$ is the maximal element of $\underline\tau$) and we want to prove for $\underline\tau$.

Recall definition \ref{phieta} (and its extensions later) where we defined the universal map $\Phi_{\underline\eta}^L\colon\xi_{\underline\eta}^L\to\tilde\xi_{\underline\eta}^L$ for tubular neighbourhoods of the $\underline\eta$-strata. Now if $D\tilde\xi_{\underline\eta}^L$ is the disk bundle of $\tilde\xi_{\underline\eta}^L$ (of a sufficiently small radius) and we set $D\xi_{\underline\eta}^L:=(\Phi_{\underline\eta}^L)^{-1}(D\tilde\xi_{\underline\eta}^L)$, which is homeomorphic to the disk bundle of $\xi_{\underline\eta}^L$, then the restriction
$$\Phi_{\underline\eta}^L|_{S\xi_{\underline\eta}^L}\colon S\xi_{\underline\eta}^L=\partial D\xi_{\underline\eta}^L\to\partial D\tilde\xi_{\underline\eta}^L=S\tilde\xi_{\underline\eta}^L$$
is a $(\underline\tau',L)$-map. Hence, by the induction hypothesis on \ref{pt2}, there are inducing maps $\rho,\tilde\rho$ which make the following diagram commutative:
$$\xymatrix@C=4pc{
Y_{\underline\tau'}^L\ar[r]^{f_{\underline\tau'}^L} & X_{\underline\tau'}^L\\
S\xi_{\underline\eta}^L\ar[u]^{\rho}\ar[r]^{\Phi_{\underline\eta}^L|_{S\xi_{\underline\eta}^L}} & S\tilde\xi_{\underline\eta}^L\ar[u]_{\tilde\rho}
}$$
We define
$$X_{\underline\tau}^L:=X_{\underline\tau'}^L\tightunderset{\tilde\rho}{\sqcup}D\tilde\xi_{\underline\eta}^L,~~~~Y_{\underline\tau}^L:=Y_{\underline\tau'}^L\tightunderset{\rho}{\sqcup}D\xi_{\underline\eta}^L~~~~\text{and}~~~~f_{\underline\tau}^L:=f_{\underline\tau'}^L\cup\Phi_{\underline\eta}^L|_{D\xi_{\underline\eta}^L}.$$
Now the only thing left is to prove that \ref{pt1} and \ref{pt2} hold for this newly constructed map.

\medskip\noindent\textbf{Proof of \ref{pt1}.\enspace\ignorespaces}
Let $P,N,f,\kappa_P,\kappa_N$ be as in \ref{pt1}. The zero-section $BG_{\underline\eta}^L$ is inside $D\tilde\xi_{\underline\eta}^L$ and $D\xi_{\underline\eta}^L$, hence it is also in $X_{\underline\tau}^L$ and $Y_{\underline\tau}^L$. Moreover, $f_{\underline\tau}^L$ maps $BG_{\underline\eta}^L\subset Y_{\underline\tau}^L$ onto $BG_{\underline\eta}^L\subset X_{\underline\eta}^L$ homeomorphically.

The submanifold $\kappa_N^{-1}(BG_{\underline\eta}^L)$ is mapped by $f$ diffeomorphically onto $\kappa_P|_{\partial P}^{-1}(BG_{\underline\eta}^L)$ because of the pullback property. We will denote both preimages by $K_N$. Now the submanifold $K_P:=\kappa_P^{-1}(BG_{\underline\eta}^L)$ has boundary $K_N$. If we set $U:=\kappa_P^{-1}(\toverset{^{~\circ}}D\tilde\xi_{\underline\eta}^L)$, then we may assume that $U$ is a tubular neighbourhood of $K_P$, hence its closure $\overline U$ can be identified with the disk bundle of $\kappa_P|_{K_P}^*\tilde\xi_{\underline\eta}^L$. The boundary $\partial U$ is now the union of a sphere bundle $\partial_SU$ over $K_P$ and a disk bundle $\partial_DU$ over $K_N$. We put $P':=P\setminus U$ and $Q:=\overline{\partial P\setminus\partial_DU}$.

Now let $V$ be the disk bundle of $\kappa_P|_{K_P}^*\xi_{\underline\eta}^L$ and let $\kappa_V\colon V\to D\xi_{\underline\eta}^L$ be the map defined by the inducing map $\kappa_P|_{K_P}$. Then the boundary $\partial V$ is the union of a sphere bundle $\partial_SV$ over $K_P$ and a disk bundle $\partial_DV$ over $K_N$. Observe that $\partial_DV$ can be identified with a (closed) tubular neighbourhood of $K_N\subset N$; we will put $N':=\overline{N\setminus\partial_DV}$.

\begin{figure}[H]
\begin{center}
\centering\includegraphics[scale=0.1]{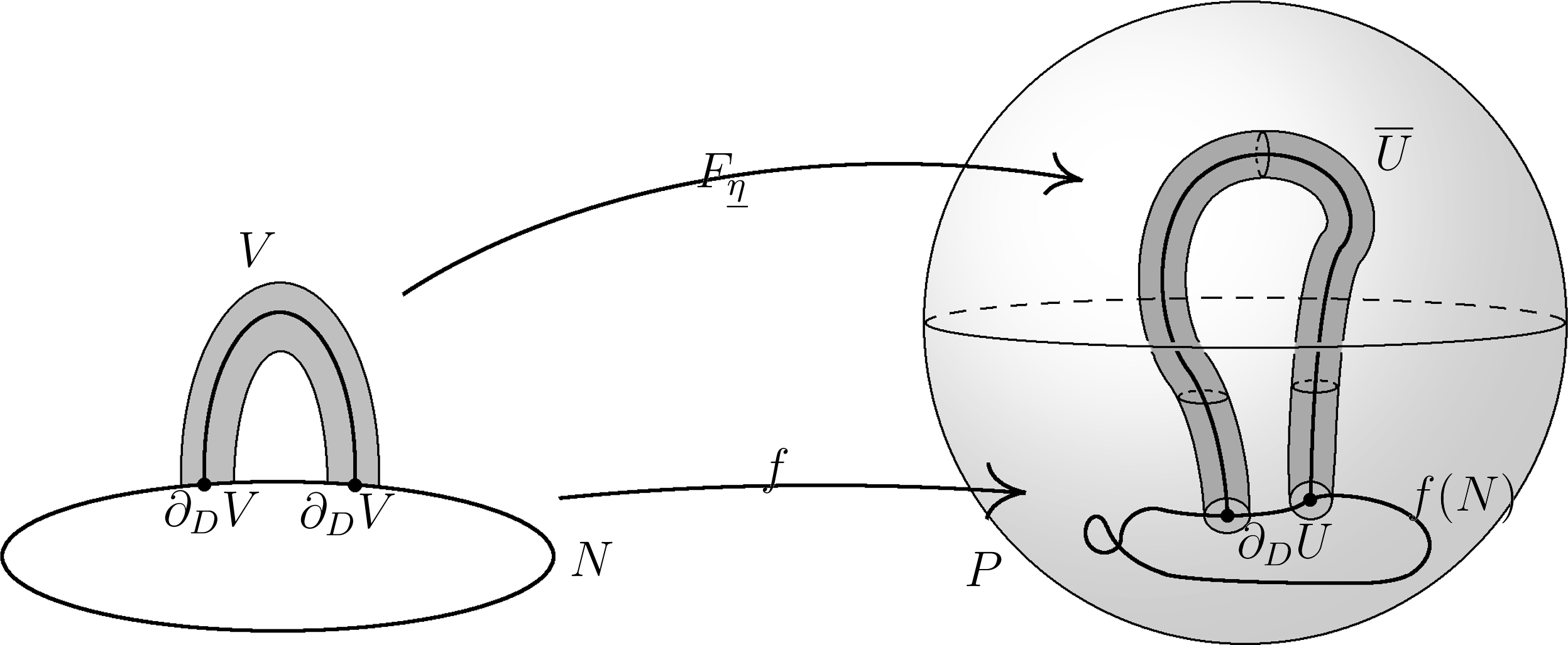}\label{kep4}
\begin{changemargin}{2cm}{2cm} 
\caption{\hangindent=1.4cm\small The ellipse on the left stands for $N$ and the (closed) ball on the right for $P$; the two indicated points on $N$ and the other two on $f(N)$ represent $K_N$ and the arc joining them is $K_P$; the disk bundles $V$ and $\overline U$ are also indicated.}
\end{changemargin} 
\vspace{-1.3cm}
\end{center}
\end{figure}

The map $\Phi_{\underline\eta}^L\colon\xi_{\underline\eta}^L\to\tilde\xi_{\underline\eta}^L$ induces a map $F_{\underline\eta}\colon V\to\overline U$ between the pullbacks, which is a $(\underline\tau,L)$-map. By definition, the restriction $F_{\underline\eta}|_{\partial_DV}$ coincides with $f|_{\partial_DV}$.

Now $P'$ is a manifold with boundary $Q\cup\partial_SU$ and the (closed) manifold $N'\cup\partial_SV$ is mapped to this boundary by the $(\underline\tau',L)$-map $f|_{N'}\cup F_{\underline\eta}|_{\partial_SV}$. Moreover, if we also add the map $\kappa_N|_{N'}\cup\kappa_V|_{\partial_SV}\colon N'\cup\partial_SV\to Y_{\underline\tau'}^L$, then we can apply the induction hypothesis on \ref{pt1} to the lower square in the following diagram to obtain the manifold $M'$ and the maps $F'$ and $\kappa_{M'}$:
$$\xymatrix@C=.7pc{
& V\ar[rrrrrr]^{F_{\underline\eta}}\ar[d]_{\kappa_V} &&&&&& \overline U\ar[d]^{\kappa_P|_{\overline U}} &\\
Y_{\underline\tau}^L\ar@{=}[r] & D\xi_{\underline\eta}^L\ar@{{}{}{}}[r]|\widesqcup_*+\txt{$_\rho$}\ar@/^1.5pc/@<.3pc>[rrrrrr]^(.7){\Phi_{\underline\eta}^L} & Y_{\underline\tau'}^L\ar[rrrr]^{f_{\underline\tau'}^L} &&&& X_{\underline\tau'}^L\ar@{{}{}{}}[r]|\widesqcup_*+\txt{$_{\tilde\rho}$} & D\tilde\xi_{\underline\eta}^L\ar@{=}[r] & X_{\underline\eta}^L\\
&&&& M'\ar@{-->}[ull]_(.35){\kappa_{M'}}\ar@{-->}[rr]^{F'} && P'\ar[u]_{\kappa_P|_{P'}} &&\\
&& N'\cup\partial_SV\ar[uu]^{\kappa_N|_{N'}\cup\kappa_V|_{\partial_SV}}\ar@{^(-->}[urr]\ar[rrrr]^{f|_{N'}\cup F_{\underline\eta}|_{\partial_SV}} &&&& Q\cup\partial_SU\ar@{^(->}[u] &&
}$$
We can glue $M'$ to $V$ along the common boundary part $\partial_SV$. The maps
$$F'\cup F_{\underline\eta}\colon M'\cup V\to P'\cup\overline U=P~~~~\text{and}~~~~\kappa_{M'}\cup\kappa_V\colon M'\cup V\to Y_{\underline\tau'}^L\tightunderset{\rho}{\sqcup}D\xi_{\underline\eta}^L=Y_{\underline\tau}^L$$
are well-defined, as the maps used in both cases coincide on $\partial_SV$. Hence we got a manifold $M:=M'\cup V$ with boundary $\partial M=N$ and maps $F:=F'\cup F_{\underline\eta}$ and $\kappa_M:=\kappa_{M'}\cup\kappa_V$ and it is easy to see that these satisfy \ref{pt1}.
~$\square$\medskip

\noindent\textbf{Proof of \ref{pt2}.\enspace\ignorespaces}
Let $M,P,F,\kappa_{\partial M},\kappa_{\partial P}$ be as in \ref{pt2}. We put $K:=\underline\eta(F)$ and $\tilde K:=f(\underline\eta(f))$; this way $f|_K\colon K\to\tilde K$ is a diffeomorphism. We now use theorem \ref{univ} (and its extensions later) to obtain closed tubular neighbourhoods $K\subset T\subset M$ and $\tilde K\subset\tilde T\subset P$ such that $f$ maps $T$ to $\tilde T$ and $\partial T$ to $\partial\tilde T$ and there is a commutative diagram
$$\xymatrix@R=.333pc{
D\xi_{\underline\eta}^L\ar[rrr]^{\Phi_{\underline\eta}^L}\ar[ddddd]&&&D\tilde\xi_{\underline\eta}^L\ar[ddddd]\\
&T\ar[r]^{f|_T}\ar[ul]_(.35){\kappa_T}\ar[ddd]&\tilde T\ar[ur]^(.35){\kappa_{\tilde T}}\ar[ddd]&\\ \\ \\
&K\ar[r]^{f|_K}\ar[dl]&\tilde K\ar[dr]&\\
BG_{\underline\eta}^L\ar[rrr]_{\id_{BG_{\underline\eta}^L}}&&&BG_{\underline\eta}^L
}$$
where $T$ and $\tilde T$ are identified with disk bundles over $K$ and $\tilde K$ respectively.

We have $\partial T=\partial_ST\cup\partial_DT$ and $\partial\tilde T=\partial_S\tilde T\cup\partial_DT$, where $\partial_ST$ and $\partial_S\tilde T$ are sphere bundles over $K\subset M$ and $\tilde K\subset P$ respectively and $\partial_DT$ and $\partial_D\tilde T$ are disk bundles over $\partial K\subset\partial M$ and $\partial\tilde K\subset\partial P$ respectively. We put $M':=\overline{M\setminus T}$, $N:=\overline{\partial M\setminus\partial_DT}$, $P':=\overline{P\setminus\tilde T}$ and $Q:=\overline{\partial P\setminus\partial_D\tilde T}$.

Now $F|_{M'}\colon M'\to P'$ is a $(\underline\tau',L)$-map and the boundaries $\partial M'=N\cup\partial_ST$ and $\partial P'=Q\cup\partial_S\tilde T$ are mapped to $Y_{\underline\tau'}^L$ and $X_{\underline\tau'}^L$ respectively by the maps $\kappa_{\partial M}|_N\cup(\rho\circ\kappa_T|_{\partial_ST})$ and $\kappa_{\partial P}|_Q\cup(\tilde\rho\circ\kappa_{\tilde T}|_{\partial_S\tilde T})$. Hence by the induction hypothesis on \ref{pt2} there are extensions $\kappa_{M'}$ and $\kappa_{P'}$ of these maps to $M'$ and $P'$ respectively, and a pullback diagram
$$\xymatrix{
Y_{\underline\tau'}^L\ar[r]^{f_{\underline\tau'}^L} & X_{\underline\tau'}^L\\
M'\ar[u]^{\kappa_{M'}}\ar[r]^{F|_{M'}} & P'\ar[u]_{\kappa_{P'}}
}$$

The maps $\kappa_{M'}$ and $\kappa_T$ coincide on $\partial_ST$ and the maps $\kappa_{P'}$ and $\kappa_{\tilde T}$ coincide on $\partial_S\tilde T$. Therefore, the maps
$$\kappa_M:=\kappa_{M'}\cup\kappa_T\colon M'\cup T=M\to Y_{\underline\tau}^L~~~~\text{and}~~~~\kappa_P:=\kappa_{P'}\cup\kappa_{\tilde T}\colon P'\cup\tilde T=P\to X_{\underline\tau}^L$$
are well-defined and it is not hard to see that they satisfy \ref{pt2}.
~$\square$

\begin{rmk}\label{1cptrmk}
The first block in the construction of $X_{\underline\tau}^L$ was the Thom space $T\gamma_k^L$ and it is clear from the proofs that the points of $P$ which are ``far away'' from the image of $M$ are mapped (according to the Thom construction) to the special point of this Thom space. Hence $\infty\in\toverset{^*}P$ is always mapped to this point.
\end{rmk}

%-------------------------------------------------------------------------------------------------------------------------------

\section{Szűcs's construction by $\tau$-embeddings}\label{szpt}

In the rest of this chapter we will have a fixed stable linear group $L$ and a fixed set $\tau$ of singularities of codimension $k$ with a fixed complete ordering $[\eta_0]<[\eta_1]<\ldots$ of its elements which extends the natural partial order (see definition \ref{pos}). As before, we can assume this to be a well-ordering and that $\tau$ has a maximal element with respect to this order.

We will again show a Pontryagin--Thom type construction (due to Szűcs) to obtain the classifying space for cobordisms of $(\tau,L)$-maps, but now the techniques of section \ref{rszpt} do not entirely work. The reason for this is that in the proof above we used that the image of each multisingularity stratum is embedded into the target space, however, this is not generally true for singularity strata. This will be one of the main problems to solve now.

\begin{defi}
By a $(\tau,L)$-embedding we mean a triple $(e,V,\mathscr{F})$ with the following properties:
\begin{itemize}
\item[(i)] $e\colon M^n\hookrightarrow Q^q$ is an embedding of a manifold $M$ to a manifold $Q$. If $M$ is a manifold with boundary, then we also assume (as in definition \ref{taumap}) that $Q$ has a boundary too, $e^{-1}(\partial Q)=\partial M$ and $e$ is transverse to $\partial Q$.
\item[(ii)] $V=(v_1,\ldots,v_m)$ where $m=q-n-k$ and the $v_i$-s are pointwise independent vector fields in $Q$ along $e(M)$ (i.e. sections of the bundle $TQ|_{e(M)}$). For manifolds with boundaries we require that the vector fields on $e(\partial M)$ are tangent to $\partial Q$. We will identify $V$ with the trivialised subbundle generated by the $v_i$-s.
\item[(iii)] $\mathscr{F}$ is a foliation of dimension $m$ on a neighbourhood of $e(M)$ and it is tangent to $V$ along $e(M)$.
\item[(iv)] Any point $p\in M$ has a neighbourhood on which the composition of $e$ with the projection along the leaves of $\mathscr{F}$ to a small $(n+k)$-dimensional transverse slice is a map that has at $p$ a singularity which belongs to $\tau$.
\item[(v)] The normal bundle $\nu_e$ of the embedding $e$ has structure group $L$.
\end{itemize}
\end{defi}

\begin{rmk}
If $(e,V,\mathscr{F})$ is a $(\tau,L)$-embedding, a stratification of $M$ arises by the submanifolds
$$\eta(e):=\eta(e,V,\mathscr{F}):=\{p\in M\mid p\in\eta(\pi\circ e)\}$$
where $\pi$ denotes the local projection around $e(p)$ along the leaves of $\mathscr{F}$.
\end{rmk}

\begin{ex}
If $f\colon M^n\to P^{n+k}$ is a $(\tau,L)$-map and $i\colon M^n\hookrightarrow\R^m$ is an embedding, then we can define a $(\tau,L)$-embedding $(e,V,\mathscr{F})$ of $M$ into $P\times\R^m$: We put $e:=f\times i$, the vector fields $v_i$ arise from a basis in $\R^m$ and $\mathscr{F}$ is composed of the leaves $\{p\}\times\R^m$ (for all $p\in P$).
\end{ex}

\begin{defi}
\begin{itemize}
\item[]
\item[(1)] The vector fields $V=(v_1,\ldots,v_m)$ and the foliation $\mathscr{F}$ in the above example are called vertical in $P\times\R^m$.
\item[(2)] The subsets $P\times\{x\}\subset P\times\R^m$ (for any $x\in\R^m$) are called horizontal sections.
\end{itemize}
\end{defi}

\begin{rmk}
If $(e,V,\mathscr{F})$ is a $(\tau,L)$-embedding of $M^n$ into $P^{n+k}\times\R^m$ such that $V$ and $\mathscr{F}$ are vertical, then $\pr_P\circ e\colon M\to P$ is a $(\tau,L)$-map.
\end{rmk}

\subsection{Cobordism of $\tau$-embeddings}\label{tauemb}

\begin{defi}
We call two $(\tau,L)$-embeddings $(e_0,V_0,\mathscr{F}_0)$ and $(e_1,V_1,\mathscr{F}_1)$ of closed source manifolds $M_0^n$ and $M_1^n$ into the manifold $Q^q$ cobordant if there is
\begin{itemize}
\item[(i)] a compact manifold with boundary $W^{n+1}$ such that $\partial W=M_0\sqcup M_1$,
\item[(ii)] a $(\tau,L)$-embedding $(E,U,\mathscr{G})$ of $W$ into $Q\times[0,1]$ such that for $i=0,1$ the restriction of $(E,U,\mathscr{G})$ to the boundary part $M_i$ is the $(\tau,L)$-embedding $(e_i,V_i,\mathscr{F}_i)$ of $M_i$ into $Q\times\{i\}$.
\end{itemize}
The set of cobordism classes of $(\tau,L)$-embeddings of $n$-manifolds into the manifold $Q$ will be denoted by $\Emb_\tau^L(n,Q)$.
\end{defi}

\begin{rmk}
If $Q^q$ is a manifold of the form $N^{q-1}\times\R^1$, then $\Emb_\tau^L(n,Q)$ is also an Abelian group with almost the same operation as the one in section \ref{grp}. The difference is that the group operation is now induced by the ``far away'' disjoint union: When forming the sum of two cobordism classes we take such representatives $(e_0,V_0,\mathscr{F}_0)$ and $(e_1,V_1,\mathscr{F}_1)$ of them that the images of $e_0$ and $e_1$ can be separated by a hypersurface $N\times\{t\}$ for a $t\in\R^1$. Such representatives always exist, because by translating any $(\tau,L)$-embedding along the lines $\R^1$ in $N\times\R^1$ we get another representative of the same cobordism class.
\end{rmk}

\begin{lemma}\label{vert}
Let $(e,V,\mathscr{F})$ be a $(\tau,L)$-embedding of $M^n$ into $P^{n+k}\times\R^m$ where $M$ is a compact manifold and $P$ is any manifold. Then there is a diffeotopy $\varphi_t~(t\in[0,1])$ of $P\times\R^m$ such that $\varphi_0$ is the identity and (the differential of) $\varphi_1$ takes $V$ to the vertical vector fields $V'$ and $\mathscr{F}$ to the vertical foliation $\mathscr{F}'$ around the image of $M$. 
\end{lemma}

\begin{prf}
The manifold $M$ is finitely stratified by the submanifolds
$$S_i:=\bigcup_{\eta\in\tau\atop\dim\eta(e)=i}\eta(e)~~~~i=0,\ldots,n.$$
By the stratified compression theorem \ref{sct} there is a diffeotopy of $P\times\R^m$ which turns the vector fields $V$ into vertical vector fields. Therefore we may assume that $V$ is already vertical, and so we only need to find a diffeotopy that takes the foliation $\FF$ into the vertical foliation $\FF'$ and its differential keeps the vector fields $V$ vertical.

We will recursively deform $\FF$ into $\FF'$ around the images of the strata $S_i~(i=0,\ldots,n)$. First we list some trivial general observations.
\begin{enumerate}
\item\label{tr1} If $R\subset\R^m$ and $K,K'\subset P\times R$ are such that each of $K$ and $K'$ intersects each horizontal section $P\times\{x\}~(x\in R)$ exactly once, then a bijective correspondence $K\to K'$ arises by associating the points on the same horizontal section to each other.
\item\label{tr2} If $A\subset P\times\R^m$ is such that for each $a\in A$ subsets $R_a,K_a,K'_a$ are given as in \ref{tr1}, then a family of bijective maps $\{K_a\to K'_a\mid a\in A\}$ arises. If we have $K_{a_1}\cap K_{a_2}=\varnothing=K'_{a_1}\cap K'_{a_2}$ for any two different points $a_1,a_2\in A$, then the union of these bijections gives a continuous bijective map
$$\alpha\colon U:=\bigcup_{a\in A}K_a\to\bigcup_{a\in A}K'_a=:U'$$
\item\label{tr3} If the subsets $A,K_a,K'_a$ in \ref{tr2} are submanifolds of $P\times\R^m$ such that $U$ and $U'$ are also submanifolds, then the map $\alpha$ is smooth.
\end{enumerate}

Denote by $V^\perp$ the orthogonal complement of the bundle $V|_{e(S_0)}\oplus Te(S_0)$ in $T(P\times\R^m)|_{e(S_0)}$ (with respect to a Riemannian metric). Choose a small neighbourhood $A$ of $e(S_0)$ in $\exp(V^\perp)$ (where $\exp$ denotes the exponential map of $P\times\R^m$) and for all $a\in A$ let $K_a$ and $K'_a$ be the intersections of a small neighbourhood of $a$ and the leaves of $\FF$ and $\FF'$ respectively.

If the neighbourhoods were chosen sufficiently small, then we are in the setting of \ref{tr3}, hence a diffeomorphism $\alpha\colon U\to U'$ arises (with the same notations as above). Note that $U$ and $U'$ are both neighbourhoods of $e(S_0)$, the map $\alpha$ fixes $e(S_0)$ and for all $a\in e(S_0)$ we have $d\alpha_a=\id_{T_a(P\times\R^m)}$.

For all points $(p,x)\in U$ we can join $(p,x)$ and $\alpha(p,x)$ by a minimal geodesic in the horizontal section $P\times\{x\}$, and using these we can extend $\alpha$ to an isotopy $\alpha_t~(t\in[0,1])$ of $U$ (for which $\alpha_0=\id_U$ and $\alpha_1=\alpha$). Observe that where the foliations $\FF$ and $\FF'$ initially coincide, this method just gives the identity for all $t\in[0,1]$. Of course this isotopy can be extended to a diffeotopy of $P\times\R^m$ (by the isotopy extension theorem) and it takes the leaves of $\FF$ to the leaves of $\FF'$ around the image of $S_0$.

Next we repeat the same procedure around $e(S_1)$, the image of the next stratum, to get a new diffeotopy (that leaves a neighbourhood of $e(S_0)$ unchanged), and so on. In the end we obtain a diffeotopy of $P\times\R^m$ which turns $\FF$ into the vertical foliation $\FF'$ around the image of $M$ and does not change the vertical vector fields $V$.
\end{prf}

\begin{rmk}
It is clear from the proof above, that the relative version of lemma \ref{vert} is also true, that is, if the vector fields $V$ and the foliation $\FF$ are already vertical on a neighbourhood of a compact subset $C\subset e(M)$, then the diffeotopy $\varphi_t~(t\in[0,1])$ is fixed on a neighbourhood of $C$.
\end{rmk}

Now we can prove the key observation on $\tau$-embeddings, namely, that the computation of cobordisms of $\tau$-maps can be reduced to that of $\tau$-embeddings.

\begin{thm}\label{cob=emb}
For any manifold $P^{n+k}$, if the number $m$ is sufficiently large (compared to $n$), then
$$\Cob_\tau^L(n,P^{n+k})\cong\Emb_\tau^L(n,P^{n+k}\times\R^m).$$
\end{thm}

\begin{prf}
Take any number $m\ge2n+4$, so any manifold of dimension at most $n+1$ can be embedded into $\R^m$ uniquely up to isotopy. We will define two homomorphisms $\varphi\colon\Cob_\tau^L(n,P)\to\Emb_\tau^L(n,P\times\R^m)$ and $\psi\colon\Emb_\tau^L(n,P\times\R^m)\to\Cob_\tau^L(n,P)$ which will turn out to be each other's inverses.

\medskip\noindent I. \emph{Construction of $\varphi$.}\medskip

For a $(\tau,L)$-map $f\colon M^n\to P^{n+k}$ we can choose any embedding $i\colon M\hookrightarrow\R^m$ and define a $(\tau,L)$-embedding $(e,V,\FF)$ by putting $e:=f\times i$ and defining $V$ and $\FF$ as the vertical vector fields and foliation. Define the map $\varphi$ to assign to the cobordism class of $f$ the cobordism class of $(e,V,\FF)$. In order to prove that $\varphi$ is well-defined, we have to show that the cobordism class of $(e,V,\FF)$ does not depend on the choice of the embedding $i$ and the representative of the cobordism class $[f]$.

\medskip\begin{sclaim}
If $i_0\colon M\hookrightarrow\R^m$ and $i_1\colon M\hookrightarrow\R^m$ are two embeddings and the above method assigns to them the $(\tau,L)$-embeddings $(e_0,V_0,\mathscr{F}_0)$ and $(e_1,V_1,\mathscr{F}_1)$ respectively, then
$$(e_0,V_0,\mathscr{F}_0)\sim(e_1,V_1,\mathscr{F}_1).$$
\end{sclaim}

\begin{sprf}
Because of the dimension condition, $i_0$ and $i_1$ can be connected by an isotopy $i_t~(t\in[0,1])$. We define a $(\tau,L)$-embedding
$$E\colon M\times[0,1]\hookrightarrow P\times\R^m\times[0,1];~(p,t)\mapsto (f(p),i_t(p),t)$$
(again with the vertical vector fields and foliation), which is clearly a cobordism between $(e_0,V_0,\mathscr{F}_0)$ and $(e_1,V_1,\mathscr{F}_1)$.
\end{sprf}

\begin{sclaim}
If $f_0\colon M_0\to P$ and $f_1\colon M_1\to P$ are cobordant $(\tau,L)$-maps and the above method assigns to them the $(\tau,L)$-embeddings $(e_0,V_0,\mathscr{F}_0)$ and $(e_1,V_1,\mathscr{F}_1)$ respectively, then
$$(e_0,V_0,\mathscr{F}_0)\sim(e_1,V_1,\mathscr{F}_1).$$
\end{sclaim}

\begin{sprf}
Let $F\colon W^{n+1}\to P\times[0,1]$ be a cobordism between $f_0$ and $f_1$. Again by the dimension condition, the embedding $i_0\sqcup i_1\colon M_0\sqcup M_1=\partial W\hookrightarrow\R^m$ extends to an embedding $I\colon W\hookrightarrow\R^m$. Hence the map $E:=F\times I$ is a $(\tau,L)$-embedding of $W$ into $P\times\R^m\times[0,1]$ (the vector fields and foliation are again vertical) and it is easy to see that this is a cobordism between $(e_0,V_0,\mathscr{F}_0)$ and $(e_1,V_1,\mathscr{F}_1)$.
\end{sprf}

\noindent II. \emph{Construction of $\psi$.}\medskip

If $(e,V,\FF)$ is a $(\tau,L)$-embedding of a manifold $M^n$ into $P^{n+k}\times\R^m$, then by lemma \ref{vert} we obtain a diffeotopy of $P\times\R^m$ that makes $V$ and $\FF$ vertical. A diffeotopy of $P\times\R^m$ trivially yields a cobordism of $(\tau,L)$-embeddings, hence we can assume that $V$ and $\FF$ were initially vertical. Now we can define $\psi$ to assign to the cobordism class of $(e,V,\FF)$ the cobordism class of the $(\tau,L)$-map $f:=\pr_P\circ e$. In order to prove that $\psi$ is well-defined, we have to show that the cobordism class of $f$ does not depend on the choice of the representative of the cobordism class $[(e,V,\FF)]$.

\medskip\begin{sclaim}
If $(e_0,V_0,\mathscr{F}_0)$ and $(e_1,V_1,\mathscr{F}_1)$ are cobordant $(\tau,L)$-embeddings of the manifolds $M_0$ and $M_1$ respectively into $P\times\R^m$ and the above method assigns to them the $(\tau,L)$-maps $f_0\colon M_0\to P$ and $f_1\colon M_1\to P$ respectively, then $f_0\sim f_1$.
\end{sclaim}

\begin{sprf}
We apply a diffeotopy $\varphi_t^i~(t\in[0,1])$ of $P\times\R^m\times\{i\}$ to turn the vector fields $V_i$ and foliation $\FF_i$ vertical (for $i=0,1$), this way we obtain the $(\tau,L)$-map $f_i\colon M_i\to P\times\{i\}$. If $(e_0,V_0,\mathscr{F}_0)$ and $(e_1,V_1,\mathscr{F}_1)$ are connected by a cobordism $(E,U,\mathscr{G})$, which is a $(\tau,L)$-embedding of a manifold $W^{n+1}$ into $P\times\R^m\times[0,1]$, then we can apply (the relative version of) lemma \ref{vert} to obtain a diffeotopy $\Phi_t~(t\in[0,1])$ of $P\times\R^m\times[0,1]$ that extends the given diffeotopies $\varphi_t^0,\varphi_t^1$ on the boundary and turns the vector fields $U$ and the foliation $\mathscr{G}$ vertical. Now combining $E$ with the final diffeomorphism $\Phi_1$ and the projection to $P\times[0,1]$, we obtain a $(\tau,L)$-cobordism $F:=\pr_{P\times[0,1]}\circ\Phi_1\circ E$ between $f_0$ and $f_1$.
\end{sprf}

It is trivial from the constructions of $\varphi$ and $\psi$ that they are homomorphisms and clearly $\psi$ is the inverse of $\varphi$, hence they are both isomophisms between $\Cob_\tau^L(n,P)$ and $\Emb_\tau^L(n,P\times\R^m)$.
\end{prf}

\begin{rmk}\label{tauembrmk}
The same proof also shows the following:
\begin{itemize}
\item[(1)] If $m$ is sufficiently large (compared to $n$ and $l$), then
$$\Cob_{\tau\oplus l}^L(n,Q^{n+k+l})\cong\Emb_\tau^L(n,Q^{n+k+l}\times\R^m).$$
\item[(2)] The cobordism group of $(\tau,L)$-embeddings stabilises, that is, for large numbers $m$ the inclusion $\R^m\subset\R^{m+1}$ induces an isomorphism
$$\Emb_\tau^L(n,Q^q\times\R^m)\cong\Emb_\tau^L(n,Q^q\times\R^{m+1}).$$
\end{itemize}
\end{rmk}

\begin{defi}
For any manifold $Q^q$ we put $\Stab_\tau^L(n,Q):=\underset{m\to\infty}\lim\Emb_\tau^L(n,Q\times\R^m)$ and call this the stable cobordism group of $(\tau,L)$-embeddings.
\end{defi}

\subsection{Semiglobal description}

In order to describe stable cobordisms of $(\tau,L)$-embeddings, we will need an analogue of section \ref{semiglob} in this setting. We fix a natural number $n$ and consider $[\eta]$-strata of $(\tau,L)$-embeddings of $n$-manifolds to $q$-manifolds where $q$ is sufficiently large and $[\eta]\in\tau$ is the maximal element.

\begin{defi}\label{approx}
Choose a finite (say $r$-)dimensional approximation $BG_\eta^L(n)$ of $BG_\eta^L$ such that the pair $(BG_\eta^L,BG_\eta^L(n))$ is $(n+1)$-connected (the existence of such a space is clear from the Milnor construction of $BG_\eta^L$); we can assume that $BG_\eta^L(n)$ is an $r$-dimensional manifold. The restriction of the global normal form of $[\eta]$ will be denoted by
$$\Phi_\eta^L(n):=\Phi_\eta^L|_{\xi_\eta^L(n)}\colon\xi_\eta^L(n):=\xi_\eta^L|_{BG_\eta^L(n)}\to\tilde\xi_\eta^L|_{BG_\eta^L(n)}=:\tilde\xi_\eta^L(n).$$
Fix a number $m\ge2r+2n+2$ and let $i\colon\xi_\eta^L(n)\hookrightarrow\R^m$ be any embedding (the dimensions were chosen such that $i$ is unique up to isotopy) and put
$$\tilde\Phi_\eta^L(n):=\Phi_\eta^L(n)\times i\colon\xi_\eta^L(n)\hookrightarrow\tilde\xi_\eta^L(n)\times\R^m,$$
which is a $(\tau,L)$-embedding (with the vertical vector fields and foliation). This is called the global normal form of $(\tau,L)$-embeddings of $[\eta]$ in dimensions at most $n$.
\end{defi}

\begin{thm}\label{univemb}
Let $(e\colon M^n\hookrightarrow Q^q,V,\mathscr{F})$ be a $(\tau,L)$-embedding where $m=q-n-k\ge2r+2n+2$. Then there are tubular neighbourhoods $\eta(e)\subset T\subset M$ and $e(\eta(e))\subset\tilde T\subset Q$ of the $\eta$-stratum and its image such that there is a commutative diagram of fibrewise isomorphisms (the vertical arrows are the projections and $BG_\eta^L(n)\subset\tilde\xi_\eta^L(n)\times\R^m$ is identified with the graph of $i|_{BG_\eta^L(n)}$)
$$\xymatrix@R=.333pc{
\xi_\eta^L(n)\ar[rrr]^{\tilde\Phi_\eta^L(n)}\ar[ddddd]&&&\tilde\xi_\eta^L(n)\times\R^m\ar[ddddd]\\
&T\ar[r]^{e|_T}\ar[ul]\ar[ddd]&\tilde T\ar[ur]\ar[ddd]&\\ \\ \\
&\eta(e)\ar[r]^{e|_{\eta(e)}}\ar[dl]&e(\eta(e))\ar[dr]&\\
BG_\eta^L(n)\ar[rrr]_{\id_{BG_\eta^L(n)}}&&&BG_\eta^L(n)
}$$
\end{thm}

\begin{prf}
Trivial from theorem \ref{univ}.
\end{prf}

\subsection{The classifying space for $\tau$-embeddings}

Our aim in the present subsection will be to construct a virtual complex (see appendix \ref{virtc}) $V_\tau^L$ with the property
$$\Stab_\tau^L(n,Q^{n+k+m})\cong\{\cpt{Q},S^mV_\tau^L\}_*$$
for any manifold $Q$ (where the notation is $\{X,Y\}_*:=\liminfty{r}[S^rX,S^rY]_*$). Recall section \ref{rszpt} where we constructed the space $X_{\underline\tau}^L$ by attaching the disk bundles $D\tilde\xi_{\underline\eta}^L$ to each other by appropriate gluing maps. Here we will follow a similar method; namely, we glue the disk bundles $D\tilde\xi_\eta^L$ together to form $V_\tau^L$, but now the gluing maps are only defined in stable sense, hence the resulting space will only be a virtual complex.

\begin{thm}\label{embclass}
There is a virtual complex $V_\tau^L$ such that for all $n,m\in\N$ where $m$ is sufficiently large (compared to $n$ and $k$), for any manifold $Q^{n+k+m}$ there is an isomorphism
$$\Emb_\tau^L(n,Q^{n+k+m})\cong\{\cpt{Q},S^mV_\tau^L\}_*.$$
\end{thm}

\begin{add*}
For any $l\in\N$ there is an approximation $V_\tau^L(l)$ of $V_\tau^L$ for which the space $S^mV_\tau^L(l)$ exists and its $(m+l)$-homotopy type is that of $S^mV_\tau^L$ whenever $m$ is large enough, and there is a subspace $K_\tau^L(l)\subset S^mV_\tau^L(l)$ with the following properties for all $n\le l$:
\begin{enumerate}
\item\label{ad1} For any manifold with boundary $Q^{n+k+m}$, if there is a $(\tau,L)$-embedding $(e,V,\FF)$ of a closed manifold $N^{n-1}$ into $\partial Q$, a map $\kappa_N$ and a stable map $\kappa_Q$ (i.e. a map between the $j$-th suspensions of the spaces involved for some $j$) such that the outer square in the diagram below is a pullback diagram, then there is a manifold $M^n$ with boundary $\partial M=N$, an extension $\kappa_M$ of $\kappa_N$ and a $(\tau,L)$-embedding $(E,U,\mathscr{G})$ that extends $(e,V,\FF)$ and the upper inner square is a pullback diagram as well.
$$\xymatrix{
K_\tau^L(l)\ar@{^(->}[rr] && S^mV_\tau^L(l)\\
& M\ar@{-->}[ul]_(.35){\kappa_M}\ar@{^(-->}[r]^E & Q\ar[u]|{\object@{/}}_{\kappa_Q}\\
N\ar[uu]^{\kappa_N}\ar@{^(-->}[ur]\ar@{^(->}[rr]^e && \partial Q\ar@{^(->}[u] 
}$$
\item\label{ad2} For any $(\tau,L)$-embedding $(E\colon M^n\hookrightarrow Q^{n+k+m},U,\mathscr{G})$ of a compact manifold with boundary, if there is a map $\kappa_{\partial M}$ and a stable map $\kappa_{\partial Q}$ such that the outer square in the diagram below is a pullback diagram, then there are extensions $\kappa_M$ and $\kappa_Q$ of $\kappa_{\partial M}$ and $\kappa_{\partial Q}$ respectively such that the upper inner square is a pullback diagram as well.
$$\xymatrix{
K_\tau^L(l)\ar@{^(->}[rrr] &&& S^mV_\tau^L(l)\\
& M\ar@{-->}[ul]_(.35){\kappa_M}\ar@{^(->}[r]^E & Q\ar@{-->}[ur]|{\object@{/}}^(.35){\kappa_Q} &\\
\partial M\ar[uu]^{\kappa_{\partial M}}\ar@{^(->}[ur]\ar@{^(->}[rrr]^{E|_{\partial M}} &&& \partial Q\ar@{_(->}[ul]\ar[uu]|{\object@{/}}_{\kappa_{\partial Q}}
}$$
\end{enumerate}
\end{add*}

\begin{prf}
We prove the theorem and the addendum together by (transfinite) induction on the elements of $\tau$ analogously to the proof of theorem \ref{ptpullback}.

The starting step is when we have $\tau=\{\Sigma^0\}$, hence $(\tau,L)$-embeddings $M^n\hookrightarrow Q^{n+k+m}$ are embeddings equipped with normal $L$-structures and $m$ pointwise independent normal vector fields. If we choose $V_{\{\Sigma^0\}}^L$ as the Thom space $T\gamma_k^L$, then the $m$-th suspension is $S^mV_{\{\Sigma^0\}}^L=S^mT\gamma_k^L=T(\gamma_k^L\oplus\varepsilon^m)$. Thus the usual Pontryagin--Thom construction yields the statement of the theorem and the addendum by putting
$$V_{\{\Sigma^0\}}^L(l):=V_{\{\Sigma^0\}}^L:=T\gamma_k^L~~~~\text{and}~~~~K_{\{\Sigma^0\}}^L(l):=K_{\{\Sigma^0\}}^L:=BL(k).$$

Now suppose we know the theorem and the addendum for $\tau':=\tau\setminus\{[\eta]\}$ (where $[\eta]$ is the maximal element of $\tau$) and we want to prove for $\tau$.

Fix a natural number $n$ and consider (as in definition \ref{approx}) the $r$-dimensional approximation $BG_\eta^L(n)$ of $BG_\eta^L$ such that $(BG_\eta^L,BG_\eta^L(n))$ is an $(n+1)$-connected pair. Let $m\ge2r+2n+2$ be such a number that the $(m+r+n)$-homotopy type of $S^mV_{\tau'}^L$ is well-defined and $V_{\tau'}^L(r+n)$ is an approximation of $V_{\tau'}^L$ for which the space $S^mV_{\tau'}^L(r+n)$ exists and its $(m+r+n)$-type is that of $S^mV_{\tau'}^L$. Then $S^mV_{\tau'}^L(r+n)$ classifies the $(\tau',L)$-embeddings with source dimensions at most $r+n$.

Take the global normal form 
$$\tilde\Phi_\eta^L(n)=\Phi_\eta^L(n)\times i\colon\xi_\eta^L(n)\hookrightarrow\tilde\xi_\eta^L(n)\times\R^m,$$
(see definition \ref{approx}) with the disk bundle $D\tilde\xi_\eta^L(n)$ (of a sufficiently small radius) and its preimage $D\xi_\eta^L(n):=(\Phi_\eta^L(n))^{-1}(D\tilde\xi_\eta^L(n))$, which is diffeomorphic to the disk bundle of $\xi_\eta^L(n)$. Now we use the restriction of $\Phi_\eta^L(n)$ to the boundary sphere bundle to define the $(\tau',L)$-embedding
$$\Phi_S:=\Phi_\eta^L(n)|_{S\xi_\eta^L(n)}\times i|_{S\xi_\eta^L(n)}\colon S\xi_\eta^L(n)\hookrightarrow S\tilde\xi_\eta^L(n)\times\R^m.$$

The image of $i|_{D\xi_\eta^L(n)}$ is contained in a ball $D^m\subset\R^m$ of a sufficiently large radius, hence the image of $\Phi_S$ is contained in $S\tilde\xi_\eta^L(n)\times D^m$ which is part of the boundary $S:=\partial(D\tilde\xi_\eta^L(n)\times D^m)$. Therefore $\Phi_S$ is a $(\tau',L)$-embedding of a manifold of dimension not greater than $r+n$ into $S$, hence by the induction hypothesis on \ref{ad2}, there is a map $\rho(n)$ and a stable map $\tilde\rho(n)$ which make the following diagram commutative:
$$\xymatrix{
K_{\tau'}^L(r+n)\ar@{^(->}[r] & S^mV_{\tau'}^L(r+n)\\
S\xi_\eta^L(n)\ar[u]^{\rho(n)}\ar@{^(->}[r]^(.55){\Phi_S} & S\ar[u]|{\object@{/}}_{\tilde\rho(n)}
}$$

Observe that by construction, the points of $S$ that are ``far away'' from the image of $\Phi_S$ are mapped to the special point of $T(\gamma_k^L\oplus\varepsilon^m)=S^mT\gamma_k^L\subset S^mV_{\tau'}^L(r+n)$. Thus the stable map $\tilde\rho(n)$ factors through the quotient $q\colon S\to S^mS\tilde\xi_\eta^L(n)$ by the subspace $D\tilde\xi_\eta^L(n)\times S^m$, that is, $\tilde\rho(n)=\hat\rho\circ q$ for a stable map $\hat\rho\colon S^mS\tilde\xi_\eta^L(n)\nrightarrow S^mV_{\tau'}^L(r+n)$. But now $m$ can be chosen large enough so that $\hat\rho$ is an existing map (not only stable), hence $\tilde\rho(n)$ is too.

We define
\begin{gather*}
S^{m}V_\tau^L(n):=S^{m}V_{\tau'}^L(r+n)\usqcup{\tilde\rho(n)}(D\tilde\xi_\eta^L(n)\times D^m)\\
\text{and}~~~~K_\tau^L(n):=K_{\tau'}^L(r+n)\usqcup{\rho(n)}D\xi_\eta^L(n).
\end{gather*}

Now we can see that these newly constructed spaces satisfy conditions \ref{ad1} and \ref{ad2} of the addendum completely analogously to the proofs of \ref{pt1} and \ref{pt2} of theorem \ref{ptpullback}, therefore we omit the proofs here.

We obtain a virtual complex by putting
$$V_\tau^L:=\liminfty{n}V_\tau^L(n).$$
Now by \ref{ad1}, \ref{ad2} and a complete analogue of the proof of theorem \ref{ptclass} we see that the claim of the theorem is also true.
\end{prf}

By this theorem we have obtained the following description of the classifying space for cobordisms of $(\tau,L)$-maps.

\begin{crly}\label{gammav}
$X_\tau^L\cong\Gamma V_\tau^L$.
\end{crly}

\begin{prf}
The Brown representability theorem implies that $X_\tau^L$ is unique up to homotopy equivalence (see \cite{hosszu} for a detailed description of this), hence we only need to prove that $\Gamma V_\tau^L$ also classifies cobordisms of $(\tau,L)$-maps. This is indeed so, since for any manifold $P^{n+k}$ we have
\begin{alignat*}2
\Cob_\tau^L(n,P)&\toverset{^*}{\cong}\Stab_\tau^L(n,P)\toverset{^{**}}{\cong}\{\cpt{P},V_\tau^L\}_*=\liminfty{r}[S^r\cpt{P},S^rV_\tau^L]_*\cong\\
&\cong\liminfty{r}[\cpt{P},\Omega^rS^rV_\tau^L]_*=[\cpt{P},\Gamma V_\tau^L]_*
\end{alignat*}
where $^*$ follows from theorem \ref{cob=emb} and $^{**}$ is an easy consequence of theorem \ref{embclass}.
\end{prf}

\begin{rmk}\label{h}
Any loop space (i.e. space of the form $\Omega Y$) is an H-space, hence we have also proved that $X_\tau^L$ is an H-space. This also means that the functors $[\cdot,X_\tau^L]$ and $[\cdot,X_\tau^L]_*$ coincide.
\end{rmk}

\begin{rmk}\label{lfr}
Similarly to corollary \ref{gammav}, it follows from remark \ref{tauembrmk} and theorem \ref{embclass} that the classifying space for cobordisms of $l$-framed $(\tau,L)$-maps (see definition \ref{lframed}) is $X_{\tau\oplus l}^L\cong\Gamma S^lV_\tau^L$.
\end{rmk}

%-------------------------------------------------------------------------------------------------------------------------------

\section{Kazarian's space}\label{kazspace}

In the proof of theorem \ref{embclass} we defined the spaces $K_\tau^L(n)$, but we did not emphasise any of their important properties. We will do so now.

\begin{rmk}
The recursive definition $K_\tau^L(n):=K_{\tau'}^L(r+n)\usqcup{\rho(n)}D\xi_\eta^L(n)$ shows that $K_\tau^L(n-1)\subset K_\tau^L(n)$, since we can assume the inclusion $K_{\tau'}^L(r+n-1)\subset K_{\tau'}^L(r+n)$ by induction, the inclusion $D\xi_\eta^L(n-1)\subset D\xi_\eta^L(n)$ is trivially true and the restriction of the gluing map $\rho(n)$ can be assumed to be $\rho(n-1)$.
\end{rmk}

\begin{defi}
The space $K_\tau^L:=\liminfty nK_\tau^L(n)$ will be called (following Szűcs \cite{hosszu}) the Kazarian space of $\tau$-maps with normal $L$-structures.
\end{defi}

For later use, we define an important tool (a spectral sequence) now. Recall that we fixed a complete order of $\tau$ that extends the natural partial order.

\begin{defi}\label{kazspec}
Suppose that the singularity set $\tau$ has order-type $\omega$ or less, let its elements be $\eta_0<\eta_1<\ldots$ and use the notations $K_i:=K_{\{\eta_j\mid j\le i\}}^L$, $\xi_i:=\xi_{\eta_i}^L$, $G_i:=G_{\eta_i}^L$ and $c_i:=\dim\xi_i$. We define the Kazarian spectral sequence corresponding to this (ordered) set of singularities as follows. The filtration $K_0\subset K_1\subset\ldots\subset K_\tau^L$ defines a homological spectral sequence with first page
$$E^1_{i,j}=H_{i+j}(T\xi_i;G)\cong H_{i+j-c_i}(BG_i;\tilde G_{\xi_i}),$$
where $G$ is any coefficient group, the isomorphism is the homological Thom isomorphism and $\tilde G_{\xi_i}$ means $G$ twisted by the orientation of $\xi_i$. This spectral sequence converges to $H_*(K_\tau^L;G)$, that is, $\underset{i+j=n}\bigoplus E^\infty_{i,j}$ is associated to $H_n(K_\tau^L;G)$.
\end{defi}

\subsection{Generalities on the Borel construction}

We will give a different way to obtain Kazarian's space, which was used by Kazarian \cite{kaz} but even considered by Thom \cite{thomthm}. More precisely, we will construct a space $\hat K_\tau^L$ and then prove that it is homotopy equivalent to $K_\tau^L$. First we recall the Borel construction which will be used in the definition of $\hat K_\tau^L$.

\begin{defi}
If $G$ is a topological group and $J$ is a space with a given $G$-action on it, then the Borel construction on $J$ is the space $BJ:=EG\tightunderset{G}{\times}J$.
\end{defi}

\begin{rmk}
The Borel construction $BJ$ has the following properties:
\begin{enumerate}
\item If $\Sigma\subset J$ is a $G$-invariant subspace, then $B\Sigma\subset BJ$.
\item If $J=\underset{i}{\bigcup}\Sigma_i$ is the decomposition of $J$ to $G$-orbits, then $BJ=\underset{i}{\bigcup}B\Sigma_i$.
\item If $J$ is contractible, $\Sigma\subset J$ is a $G$-orbit and $G_x$ is the stabiliser of $x\in\Sigma\subset J$, then $BG_x=EG_x/G_x=EG\tightunderset{G}{\times}(G/G_x)=EG\tightunderset{G}{\times}\Sigma=B\Sigma$.
\end{enumerate}
\end{rmk}

\begin{lemma}\label{borel}
Let $G$ be a compact Lie group, $J$ a contractible manifold with a smooth (left) $G$-action on it, $\Sigma\subset J$ a $G$-orbit and $G_x$ the stabiliser of $x\in\Sigma\subset J$. Fix a Riemannian metric on $J$ such that the $G$-action is isometric (such a metric exists) and choose a small orthogonal slice $S_x$ of $\Sigma$ at $x$. This way $G_x$ acts orthogonally on the tangent space $T_xS_x$; denote by $\rho_x\colon G_x\to\O(T_xS_x)$ the arising representation. Then a neighbourhood of $B\Sigma$ in $BJ$ can be identified with the universal bundle $EG_x\utimes{\rho_x}T_xS_x$ over $BG_x$.
\end{lemma}

\begin{prf}
Let $S\subset J$ be a $G$-invariant tubular neighbourhood of $\Sigma$ composed of the fibres $S_{gx}~(g\in G)$. Then for all $g\in G$ diffeomorphisms $S_x\to S_{gx}$ arise (of course, if $g_1x=g_2x$, then $S_{g_1x}=S_{g_2x}$ and we obtain different diffeomorphisms).

Lift the bundle $S$ over $\Sigma$ to a bundle $\tilde S$ over $G$ by the following pullback diagram:
$$\xymatrix@C=.75pc{
\tilde S\ar@{-->}[rrr]\ar@{-->}[d]_{G_x} &&& G\ar[d]^{G_x} \\
S\ar[rr] && \Sigma\ar@{=}[r] & G/G_x
}$$
Here the fibre $\tilde S_g$ of $\tilde S$ over $g\in G$ is identified with the fibre $S_{[g]}=S_{gx}$ of $S$ by the pullback definition (where $[g]\in G/G_x$ is the coset of $g$ and it is identified with the point $gx\in\Sigma$). This way the diffeomorphisms $S_x\to S_{gx}$ give canonical identifications $\tilde S_e\to\tilde S_g$ for all $g\in G$ (where $e$ is the neutral element of $G$), hence $\tilde S$ is a trivial bundle.

The fibration $\tilde S\to S$ is $G$-equivariant, hence applying the Borel construction results in a map
$$B\tilde S=EG\utimes{G}\tilde S\to EG\utimes{G}S=BS.$$
We have a decomposition $\tilde S=G\times\tilde S_e$ which is $G$-equivariant (using the trivial $G$-action on $\tilde S_e$), hence the domain is $EG\utimes{G}\tilde S=EG\utimes{G}G\times\tilde S_e=EG\times\tilde S_e$. The map is the quotient by the diagonal $G_x$-action, therefore the target (i.e. $BS$) can be identified with $EG_x\utimes{G_x}\tilde S_e$. The fibre $\tilde S_e$ is identified with $S_x$ which can be identified with $T_xS_x$ and the $G_x$-action on it is the representation $\rho_x$, thus $BS$ (which is a neighbourhood of $B\Sigma$ in $BJ$) is identified with $EG_x\utimes{\rho_x}T_xS_x$.
\end{prf}

\begin{crly}\label{transv}
If $G,J,\Sigma$ are as above, $R\subset J$ is transverse to $\Sigma$ and we identify $J$ with the fibre of $BJ\xra{J}BG$, then $R$ is transverse to $B\Sigma\subset BJ$.
\end{crly}

\subsection{Constructing Kazarian's space}

In this subsection we will construct the space $\hat K_\tau^L$ and prove its homotopy equivalence with $K_\tau^L$. First we construct a preliminary ``unstable'' version; to do this, we fix natural numbers $n$ and $r$ and suppose that (the roots of) all elements of $\tau$ have source dimensions at most $n$ and are $r$-determined (that is, for all $[\vartheta]\in\tau$, all germs that have the same $r$-jet as $\vartheta$ are $\AA$-equivalent to $\vartheta$).

\begin{defi}
We put $J(n):=J^r_0(\R^n,\R^{n+k})$ the $r$-jet space of germs $(\R^n,0)\to(\R^{n+k},0)$ at the origin and
$$G^L(n):=\big\{(J^r\varphi,J^r\psi)\in J^r_0(\Diff_0(\R^n))\times J^r_0(\Diff_0(\R^{n+k}))\mid\big(\begin{smallmatrix}
d\varphi_0&0\\
0&d\psi_0
\end{smallmatrix}\big)\in L(2n+k)\big\}$$
the group of $r$-jets of those diffeomorphism germs of $(\R^n,0)$ and $(\R^{n+k},0)$ for which the effect of their derivatives on the virtual normal bundle fixes the group $L$.
\end{defi}

\begin{defi}\label{taujet}
In $J(n)$ the $r$-jets of the germs that belong to the same singularity class form a $G^L(n)$-orbit. If $\Sigma_\vartheta\subset J(n)$ is the orbit corresponding to the singularity class $[\vartheta]$, then we put
$$B\vartheta(n):=B\Sigma_\vartheta=EG^L(n)\tightunderset{G^L(n)}{\times}\Sigma_\vartheta.$$
The set of $r$-jets of all singularity classes in $\tau$ is also $G^L(n)$-invariant. This set will be denoted by $J_\tau^L(n)\subset J(n)$ and we put 
$$\tilde K_\tau^L(n):=BJ_\tau^L(n)=EG^L(n)\tightunderset{G^L(n)}{\times}J_\tau^L(n).$$
We will denote the fibration $\tilde K_\tau^L(n)\xra{J_\tau^L(n)}BG^L(n)$ by $p(n)$.
\end{defi}

\begin{rmk}\label{miafasz}
$~$
\begin{itemize}
\item[(1)] By the decomposition of $J_\tau^L(n)$ to $G^L(n)$-orbits we have $\tilde K_\tau^L(n)=\underset{[\vartheta]\in\tau}\bigcup B\vartheta(n)$.
\item[(2)] Note that $BG^L(n)$ is homotopy equivalent to $BH^L(n)$ where
$$H^L(n):=\big\{(A,B)\in\O(n)\times\O(n+k)\mid\big(\begin{smallmatrix}
A&0\\
0&B
\end{smallmatrix}\big)\in L(2n+k)\big\},$$
so we can lift the universal bundles $\gamma_n$ and $\gamma_{n+k}$ to bundles over $BG^L(n)$.
\end{itemize}
\end{rmk}

\begin{defi}\label{nuniv}
Pull back the lift of the bundles $\gamma_n$ and $\gamma_{n+k}$ by $p(n)$ to bundles over $\tilde K_\tau^L(n)$ and denote their formal difference (which is a $k$-dimensional virtual vector bundle with structure group $L$) by $\nu_\tau^L(n)$. We call this the $n$-universal virtual normal bundle for $(\tau,L)$-maps.
\end{defi}

The following property was proved by Thom, see \cite{thomthm} or \cite{kaz}. We will not prove it here, only sketch the construction that proves it.

\begin{thm}\label{kazthm}
For any $(\tau,L)$-map $f\colon M^n\to P^{n+k}$ there is a map $\tilde\kappa_f\colon M\to\tilde K_\tau^L(n)$ such that $\tilde\kappa_f^{-1}(B\vartheta(n))=\vartheta(f)$ for all singularities $[\vartheta]\in\tau$. Moreover, this map also has the property that the pullback $\tilde\kappa_f^*\nu_\tau^L(n)$ is the virtual normal bundle $\nu_f$.
\end{thm}

\noindent\textbf{Sketch of the proof.\enspace\ignorespaces}
Fix Riemannian metrics on $M$ and $P$, denote by $\exp_M$ and $\exp_P$ their exponential maps respectively and let $\toverset{_\circ}{T}M\subset TM$ be a neighbourhood of the zero-section such that $\exp_M|_{{\overset{_\circ}{T}}_pM}$ is a diffeomorphism onto its image for all $p\in M$. Then there is a fibrewise map $F$ that makes the following diagram commutative:
$$\xymatrix{
\toverset{_\circ}{T}M\ar@{-->}[r]^{F}\ar[d]_{\exp_M} & TP\ar[d]^{\exp_P} \\
M\ar[r]^f & P
}$$

If $J^r_0(TM,f^*TP)$ denotes the $r$-jet bundle of germs of fibrewise maps $TM\to f^*TP$ along the zero-section and $J_\tau^L(M)\subset J^r_0(TM,f^*TP)$ is the subspace corresponding to the singularities in $\tau$, then $J_\tau^L(M)\xra{J_\tau^L(n)}M$ is a fibre bundle with structure group $G^L(n)$. Hence it can be induced from the universal bundle by maps $\kappa,\tilde\kappa$ as shown in the diagram
$$\xymatrix{
J_\tau^L(M)\ar[r]^{\tilde\kappa}\ar[d]_{J_\tau^L(n)} & \tilde K_\tau^L(n)\ar[d]^{J_\tau^L(n)} \\
M\ar@{-->}@/^4pc/[u]^{\tilde f}\ar[r]^(.45)\kappa & BG^L(n)
}$$
Moreover, the map $F$ defines a section $\tilde f$ of that bundle (also indicated in the diagram) and we define $\tilde\kappa_f:=\tilde\kappa\circ\tilde f$.
~$\square$\medskip

Until this point we did not use that when forming the singularity classes we identified each germ with its suspension (see definition \ref{sing}). Now we will use it in the following way: Observe that for any $(\tau,L)$-map $f\colon M^n\to P^{n+k}$ we can add an arbitrary vector bundle $\zeta$ over $M$ to both $TM$ and $f^*TP$ and replace the map $F$ in the above proof by
$$F\oplus{\id}_\zeta\colon\toverset{_\circ}TM\oplus\zeta\to TP\oplus\zeta,$$
this way obtaining the same singularity stratification of $M$. 

Now we put $\zeta:=\nu_M$, the stable normal bundle of $M$ (i.e. $\zeta:=\varepsilon^m-TM$ for a sufficiently large number $m$ such that $M$ can be embedded into $\R^m$ uniquely up to isotopy) and consider the $r$-jet bundle $J^r_0(\varepsilon^m,f^*TP\oplus\nu_M)\to M$ of germs of fibrewise maps $\varepsilon^m\to f^*TP\oplus\nu_M$ along the zero-section. This jet bundle can be induced from the bundle $J^r_0(\varepsilon^m,\gamma_{m+k}^L)\to BL(m+k)$.

\begin{defi}\label{kazj}
In each fibre $J^r_0(\R^m,\R^{m+k})$ of the bundle $J^r_0(\varepsilon^m,\gamma_{m+k}^L)\to BL(m+k)$ we can take the subspace $J_\tau^L(m)$ (see definition \ref{taujet}). The union of these in all fibres will be denoted by $\hat K_\tau^L(m)$ and the partition of $\hat K_\tau^L(m)$ corresponding to singularities in $\tau$ (by the decomposition of $J_\tau^L(m)$ to the $G^L(m)$-orbits) will be denoted by a slight abuse of notation by $\hat K_\tau^L(m)=\underset{[\vartheta]\in\tau}\bigcup B\vartheta(m)$.
\end{defi}

\begin{rmk}
$~$
\begin{enumerate}
\item Another description of the space $\hat K_\tau^L(m)$ is the following: The space $\tilde K_\tau^L(m)$ was the total space of a bundle over $BG^L(m)$, which is homotopy equivalent to $BH^L(m)$ (see remark \ref{miafasz}); the inclusion $L(m)\times L(m+k)\hookrightarrow H^L(m)$ induces a map $BL(m)\times BL(m+k)\to BH^L(m)$, that we compose with the projection $BL(m+k)\cong EL(m)\times BL(m+k)\to BL(m)\times BL(m+k)$. The pullback of $\tilde K_\tau^L(m)$ by this map is $\hat K_\tau^L(m)$, that is, we have the following pullback diagram:
$$\xymatrix@C=.75pc{
\hat K_\tau^L(m)\ar[rrrr]\ar[d]_{J_\tau^L(m)} &&&& \tilde K_\tau^L(m)\ar[dd]^{J_\tau^L(m)} \\
BL(m+k)\ar@{=}[r]^(.37){\sims{1.6}} & EL(m)\times BL(m+k)\ar[d] \\
& BL(m)\times BL(m+k)\ar[rr] && BH^L(m)\ar@{=}[r]^{\sims{1.6}} & BG^L(m)
}$$
\item Recall the definition of the $m$-universal virtual normal bundle (\ref{nuniv} and change $n$ to $m$) that we get by pulling back the virtual bundle $\gamma_{m+k}-\gamma_m$ over $BH^L(m)$ to $\tilde K_\tau^L(m)$. Using the above diagram we can pull this back to a virtual bundle over $\hat K_\tau^L(m)$ too; from now on this will be called $\nu_\tau^L(m)$. The reason for this redefinition is part \ref{harom} of this remark.
\item\label{harom} The considerations above the definition of $\hat K_\tau^L(m)$ imply that the analogue of theorem \ref{kazthm} is true for the space $\hat K_\tau^L(m)$ in the following sense: Given a $(\tau,L)$-map $f\colon M^n\to P^{n+k}$ we can construct a map $\hat\kappa_f\colon M\to\hat K_\tau^L(m)$ with the properties $\hat\kappa_f^{-1}(B\vartheta(m))=\vartheta(f)$ (for all $[\vartheta]\in\tau$) and $\hat\kappa_f^*\nu_\tau^L(m)=\nu_f$.
\end{enumerate}
\end{rmk}

Observe that the number $r$ that we fixed in the beginning of this subsection can be replaced by any larger number, hence we can even use the infinite jet space $J^\infty_0(\R^m,\R^{m+k})$ (i.e. the space of all polynomial maps $\R^m\to\R^{m+k}$ with $0$ constant term) instead of $J^r_0(\R^m,\R^{m+k})$ (i.e. the polynomial maps $\R^m\to\R^{m+k}$ with $0$ constant term and degree at most $r$). Of course for any $(\tau,L)$-map $f\colon M^n\to P^{n+k}$ we can choose the map $\hat\kappa_f$ such that its image is in a finite dimensional approximation of this space.

This way the condition we assumed above, that all elements of $\tau$ are $r$-determined, can be omitted. We may also get rid of the other condition we assumed, that all elements of $\tau$ have source dimensions at most $n$, by converging with $m$ to $\infty$. Hence we get the following.

\begin{defi}\label{kaznu}
The (second version of the) Kazarian space, the stratum of the Kazarian space corresponding to $[\vartheta]\in\tau$ and the universal virtual normal bundle for $(\tau,L)$-maps are respectively defined as
$$\hat K_\tau^L:=\liminfty{m}\hat K_\tau^L(m),~~~~B\vartheta:=\liminfty{m}B\vartheta(m)~~~~\text{and}~~~~\nu_\tau^L:=\liminfty{m}\nu_\tau^L(m)$$
\end{defi}

\begin{rmk}\label{kappaf}
To formulate the following theorem, we need to recall that a $(\tau,L)$-map $f\colon M^n\to P^{n+k}$ corresponds to a $(\tau,L)$-embedding $M\into P\times\R^m$ for any sufficiently large number $m$ (by theorem \ref{cob=emb}), so part \ref{ad2} of the addendum of theorem \ref{embclass} assigns to $f$ a map $\kappa_M\colon M\to K_\tau^L$. This map will be denoted by $\kappa_f$ from now on.
\end{rmk}

\begin{thm}
We have $\hat K_\tau^L\cong K_\tau^L$ and for any $(\tau,L)$-map $f\colon M^n\to P^{n+k}$ the maps $\hat\kappa_f$ and $\kappa_f$ coincide up to this homotopy equivalence.
\end{thm}

\begin{prf}
The proof is an induction on the elements of $\tau$. We will use manifold-like properties of the spaces involved, which all work because they hold for finite dimensional approximations and for their direct limits as well.

The starting step is $\tau=\{\Sigma^0\}$. Observe that in this case we have $K_{\{\Sigma^0\}}^L=BG_{\Sigma^0}^L=BL(k)$, the bundle $\xi_{\Sigma^0}^L$ is $0$-dimensional and $\tilde\xi_{\Sigma^0}^L$ is the universal bundle $\gamma_k^L$. The inclusion $i\colon BL(k)\into T\gamma_k^L$ defines a map $\hat\kappa_i\colon BL(k)\to\hat K_{\{\Sigma^0\}}^L$ (by defining $\hat\kappa_i$ for finite dimensional approximations and taking the direct limit) for which $\hat\kappa_i^*\nu_{\{\Sigma^0\}}^L=\gamma_k^L$. But $\nu_{\{\Sigma^0\}}^L$ admits $L$ as a structure group, hence it can be induced from the universal bundle $\gamma_k^L$ by a map $g\colon\hat K_{\{\Sigma^0\}}^L\to BL(k)$. Now $g\circ\hat\kappa_i$ induces the bundle $\gamma_k^L$ from $\gamma_k^L$, so the uniqueness of the inducing map implies that $g\circ\hat\kappa_i$ is homotopic to the identity.

Observe that the homology groups of $K_{\{\Sigma^0\}}^L=BL(k)$ and $\hat K_{\{\Sigma^0\}}^L$ are finitely generated, since their finite dimensional approximations are homotopy equivalent to compact spaces. Hence $g$ and $\hat\kappa_i$ induce isomorphisms in the homologies, which implies (using the homological Whitehead theorem) that $\hat\kappa_i$ is a homotopy equivalence. This proves the theorem for $\tau=\{\Sigma^0\}$.

Now suppose we know the theorem for $\tau':=\tau\setminus\{[\eta]\}$ (where $[\eta]$ is the maximal element of $\tau$) and we want to prove for $\tau$. The induction step will be similar to the first step, only in a bit different setting.

Recall that $K_\tau^L$ was constructed recursively with the last step being the attachment of $D\xi_\eta^L$ to $K_{\tau'}^L$ by the map $\kappa_{\Phi_\eta^L|_{S\xi_\eta^L}}$ (more precisely, the construction was using finite dimensional approximations and then taking a direct limit). Let $h'\colon K_{\tau'}^L\to\hat K_{\tau'}^L$ be the homotopy equivalence constructed by the earlier steps in the induction and define
$$h\colon K_\tau^L=K_{\tau'}^L\usqcup{\kappa_{\Phi_\eta^L|_{S\xi_\eta^L}}}D\xi_\eta^L\to\hat K_\tau^L$$
by $h|_{K_{\tau'}^L}:=h'$ and $h|_{D\xi_\eta^L}:=\hat\kappa_{\Phi_\eta^L|_{D\xi_\eta^L}}$ (again using finite dimensional approximations and direct limits). This is well-defined, as by the induction hypothesis we have $h'\circ\kappa_{\Phi_\eta^L|_{S\xi_\eta^L}}=\hat\kappa_{\Phi_\eta^L|_{S\xi_\eta^L}}$.

By the homological Whitehead theorem it is enough to show that the map $h$ induces isomorphisms in the homologies and to do this, using the 5-lemma on the long exact sequence, it is enough to show that the map
$$\hat\kappa\colon(D\xi_\eta^L,S\xi_\eta^L)\to(\hat K_\tau^L,\hat K_{\tau'}^L)$$
defined by $\hat\kappa_{\Phi_\eta^L|_{D\xi_\eta^L}}$ induces isomorphisms in the homologies.

Observe that $\hat\kappa$ maps the $\eta$-stratum $BG_\eta^L\subset D\xi_\eta^L$ to $B\eta\subset\hat K_\tau^L$.

\medskip\begin{sclaim}
The map $\hat\kappa$ maps each fibre of $D\xi_\eta^L$ transverse to $B\eta$ in $\hat K_\tau^L$.
\end{sclaim}

\begin{sprf}
Let $r$ be such a number that $[\eta]$ is $r$-determined and let $m$ be an arbitrary large number. We will use corollary \ref{transv} with the substitutions $G:=G_\eta^L$, $J:=J^r_0(\R^m,\R^{m+k})$ (the space of polynomial maps $\R^m\to\R^{m+k}$ with $0$ constant term and degree at most $r$) and $\Sigma:=[\eta]\cap J$ (those polynomial maps in $J$ whose germs at $0$ are equivalent to $\eta$). We may of course assume that $\eta\colon\R^c\to\R^{c+k}$ is itself a polynomial map and that it is the root of its type.

Recall that $\Phi_\eta^L$ restricted to each fibre $\R^c$ of $\tightoverset{\,_\circ~~~}{D\xi_\eta^L}$ is the map (germ) $\eta\colon\R^c\to\R^{c+k}$, so $\hat\kappa$ restricted to this $\R^c$ maps each point $a\in\R^c$ to the polynomial map $p_a\in J$ (in some fibre of $\hat K_\tau^L$) for which
$$p_a(x,y)=(\eta(x+a)-\eta(a),y)\in\R^{c+k}\times\R^{m-c}~~~~(x,y)\in\R^c\times\R^{m-c}.$$
Now corollary \ref{transv} implies that it is enough to show that the image of $\R^c$ in $J$ is transverse to $[\eta]\cap J$.

%If $U\subset J$ is a sufficiently small neighbourhood of $p_0$, then for any element $q\in U$ there is a linear subspace $L_q\subset\R^m$ (recall that $q$ is a map $\R^m\to\R^{m+k}$) with the properties $\ker dq\subset L_q$ and $\dim L_q=\dim\ker dp_0=:l$. If $U$ is small enough, then all subspaces $L_q~(q\in U)$ can be projected isomorphically onto $L_{p_0}=\ker dp_0$, so all of these subspaces can be identified with $\R^l$. This defines a restriction map
%$$\rho\colon U\to J^r_0(\R^l,\R^{m+k});~q\mapsto q|_{L_q},$$
%where $\rho(q)$ is either the genotype of $q$ (i.e. restriction to the kernel of its differential) or a multiple suspension of its genotype.

If we put $\tilde\eta:=\eta\times\id_{\R^{m-c}}$, then for any $a\in\R^c$ the polynomial $p_a$ can be identified with $J^r\tilde\eta(a,b)\in J^r(\R^c\times\R^{m-c},\R^{m+k})$ for any $b\in\R^{m-c}$. Hence the map $a\mapsto p_a$ can be identified with any section $s\colon\R^c\to\R^c\times\R^{m-c}$ composed with $J^r\tilde\eta$ and the projection $\pr_J$ to the fibre of the trivial bundle $J^r(\R^m,\R^{m+k})$. By a theorem of Mather \cite{math5}, the section $J^r\tilde\eta$ is transverse to the fibrewise $\AA$-equivalence classes in $J^r(\R^m,\R^{m+k})$ (since $\tilde\eta$ is stable), which implies that $\pr_J\circ J^r\tilde\eta\circ s$ is transverse to $[\eta]\cap J$ and this is what we wanted to prove.
\end{sprf}

We may assume that $\hat\kappa$ maps each fibre of $\tightoverset{\,_\circ~~~}{D\xi_\eta^L}\approx\xi_\eta^L$ into a fibre of the normal bundle $\nu_\eta$ of the stratum $B\eta$ in $\hat K_\tau^L$. Now putting $f:=\hat\kappa|_{BG_\eta^L}\colon BG_\eta^L\to B\eta$ we have $f^*\nu_\eta=\xi_\eta^L$. The bundle $\nu_\eta$ admits a $G_\eta^L$-structure (a consequence of lemma \ref{borel}), hence it can be induced from the universal such bundle, which is $\xi_\eta^L$. In other words there is a map $g\colon B\eta\to BG_\eta^L$ with the property $g^*\xi_\eta^L=\nu_\eta$. Hence $g\circ f$ induces the bundle $\xi_\eta^L$ from $\xi_\eta^L$, so by the uniqueness of the inducing map, $g\circ f$ is homotopic to the identity.

The homology groups of $BG_\eta^L$ and $B\eta$ are finitely generated in each dimension, hence $f$ (as well as $g$) induces isomorphisms in the homologies. This implies that the map $\hat f\colon T\xi_\eta^L\to T\nu_\eta$ defined by $f$ between the Thom spaces also induces isomorphisms in the homologies. Hence we have
$$H_*(D\xi_\eta^L,S\xi_\eta^L)\cong H_*(T\xi_\eta^L)\xra{\hat f_*}H_*(T\nu_\eta)\cong H_*(\hat K_\tau^L,\hat K_{\tau'}^L)$$
and the composition of $\hat f_*$ with these isomorphisms is $\hat\kappa_*$, thus $\hat\kappa$ also induces isomorphisms in the homologies.
\end{prf}

From now on we will identify the spaces $\hat K_\tau^L$ and $K_\tau^L$.

\subsection{Kazarian's conjecture}

In this subsection we show a theorem of Szűcs which was formulated as a conjecture by Kazarian and is as follows.

\begin{thm}\label{kazc}
If $L$ is positive (i.e. all of its elements have positive determinants), then $X_\tau^L\congq\Gamma S^k(K_\tau^L)^+$.
\end{thm}

In other words, the space $X_\tau^L$ gives the stable rational homotopy type of the $k$-th suspension of $(K_\tau^L)^+:=K_\tau^L\sqcup*$. The proof is the composition of several individual statements and it will actually show more. Namely, if we modify the right-hand side properly, then we get a true homotopy equivalence (not only a rational one) relating $X_\tau^L$ to $K_\tau^L$, which remains true even if $L$ is not positive.

Observe that $S^k(K_\tau^L)^+$ is the Thom space of the trivial $k$-dimensional vector bundle over $K_\tau^L$. The proper modification we mentioned above is to replace this trivial bundle with the universal virtual normal bundle $\nu_\tau^L$ (see definition \ref{kaznu}), which makes sense by corollary \ref{gammatnu}.

\begin{thm}\label{kula}
$X_\tau^L\cong\Gamma T\nu_\tau^L$.
\end{thm}

\begin{prf}
By corollary \ref{gammav} it is enough to prove that the stable homotopy type of $T\nu_\tau^L$ is the same as that of $V_\tau^L$.

Note that it is actually enough to prove this homotopy equivalence between sufficiently large suspensions of finite dimensional approximations of the (virtual) spaces involved. Indeed, if we get sequences of finite dimensional approximations and numbers $m_i~(i\in\N)$ converging to $\infty$ such that the $m_i$-th suspensions of the $i$-th approximations exist and are homotopy equivalent, then we also have $\Gamma V_\tau^L\cong\Gamma T\nu_\tau^L$.

Now fix natural numbers $n,m$ such that the $(m+n)$-homotopy type of $S^mV_\tau^L(n)$ exists and is that of $S^mV_\tau^L$ (see theorem \ref{embclass}), the virtual bundle $\nu:=\nu_\tau^L|_{K_\tau^L(n)}$ can be represented by $\alpha-\varepsilon^m$ for a vector bundle $\alpha$ over $K_\tau^L(n)$, and so $S^mT\nu$ exists (see corollary \ref{gammatnu}). We want to prove that $S^mV_\tau^L(n)$ is homotopy equivalent to $S^mT\nu$; by the cohomological Whitehead theorem it is enough to show a map between them which induces isomorphisms in the cohomologies.

Recall that $K_\tau^L$ is the union of the spaces $D\xi_\vartheta^L~([\vartheta]\in\tau)$ and the universal virtual normal bundle is such that
$$\nu_\tau^L|_{D\xi_\vartheta^L}=\pi^*\tilde\xi_\vartheta^L-\pi^*\xi_\vartheta^L$$
where $\pi$ is the projection of the disk bundle $D\xi_\vartheta^L\to BG_\vartheta^L$. Recall also that the inclusion $K_\tau^L(n)\subset S^mV_\tau^L(n)$ is by construction such that for all $[\vartheta]\in\tau$ and any $l$ for which $D\xi_\vartheta^L(l)$ took part in the construction, $D\xi_\vartheta^L(l)\subset K_\tau^L(n)$ is embedded into $D\tilde\xi_\vartheta^L(l)\times D^m$ by the restriction of the global normal form $\tilde\Phi_\vartheta^L(l)$ to $D\xi_\vartheta^L(l)$ (see the proof of theorem \ref{embclass}). Hence the normal bundle of this embedding is the restriction $(\nu\oplus\varepsilon^m)|_{D\xi_\vartheta^L(l)}$ (which is a non-virtual bundle).

Now if $[\eta]$ is the maximal element of $\tau$ and $\tau':=\tau\setminus\{[\eta]\}$, then the inclusion of $K_\tau^L(n)$ into $S^mV_\tau^L(n)$ can be written as
$$(K_\tau^L(n),K_{\tau'}^L(n+r))\hookrightarrow(S^mV_\tau^L(n),S^mV_{\tau'}^L(n+r))$$
where $r\in\N$ and $K_{\tau'}^L(n+r)$ and $S^mV_{\tau'}^L(n+r)$ are finite dimensional approximations of $K_{\tau'}^L$ and $S^mV_{\tau'}^L$ respectively. Hence this inclusion factors to an embedding between the quotient spaces
$$T\xi_\eta^L(n)\hookrightarrow S^mT\tilde\xi_\eta^L(n),$$
which is the identity on the base space $BG_\eta^L(n)$.

Let $U$ be a tubular neighbourhood of $K_\tau^L(n)\subset S^mV_\tau^L(n)$ and factor $S^mV_\tau^L(n)$ by the complement of $U$. If we take the quotient of this factorisation by $S^mV_{\tau'}^L(n+r)$ on the left-hand side and by the image of $U|_{K_{\tau'}^L(n+r)}$ on the right-hand side, then we get the diagram
$$\xymatrix@C=.75pc{
S^mV_\tau^L(n)\ar[rrrrr]^{/(S^mV_\tau^L(n)\setminus U)}\ar[d]_{/S^mV_{\tau'}^L(n+r)} &&&&& S^mT\nu\ar[d]^{/S^mT(\nu|_{K_{\tau'}^L(n+r)})}\\
S^mT\tilde\xi_\eta^L(n)\ar[rrrrr]^(.36){/(S^mT\tilde\xi_\eta^L(n)\setminus U)} &&&&& S^mT(\nu|_{D\xi_\eta^L(n)})/S^mT(\nu|_{S\xi_\eta^L(n)})
}$$

The lower horizontal arrow in this diagram maps the space $T(\tilde\xi_\eta^L(n)\oplus\varepsilon^m)$ to the space $S^mT(\nu|_{BG_\eta^L(n)}\oplus\xi_\eta^L(n))=T(\nu|_{BG_\eta^L(n)}\oplus\xi_\eta^L(n)\oplus\varepsilon^m)$ and it is again the identity on the base space $BG_\eta^L(n)$, hence it induces isomorphisms in the cohomologies by the Thom isomorphisms. Now an induction on the elements of $\tau$ shows that the upper horizontal arrow also induces isomorphisms in the cohomologies and this is what we wanted to prove.
\end{prf}

In the following we will use the infinite symmetric product operation, which will be denoted by $\SP$.

\begin{lemma}\label{doldth}
If $A$ and $B$ are connected CW-complexes such that $H_*(A)\cong H_*(B)$, then $\SP A$ is homotopy equivalent to $\SP B$.
\end{lemma}

\begin{prf}
We use a few known properties of the infinite symmetric product, which can be read in \cite[4.K]{hatcher}. The space $\SP A$ is homotopy equivalent to a product of Eilenberg--MacLane spaces and by the Dold--Thom theorem we have $\pi_*(\SP A)\cong\tilde H_*(A)$. Of course the same is true for $B$, hence the equivalence
$$\SP A\cong\prod_{i=1}^\infty K(\tilde H_i(A),i)=\prod_{i=1}^\infty K(\tilde H_i(B),i)\cong\SP B$$
holds.
\end{prf}

\begin{crly}
For any virtual complex $V$ the space $\SP V$ is well-defined.
\end{crly}

\begin{prf}
For each natural number $i$ fix a number $n(i)$ and a subcomplex $V_i$ of $V$ such that for all $n\ge n(i)$ the space $S^nV_i$ exists and its $(n+i)$-homotopy type is that of $S^nV$. Then the same is true for the space $\SP S^nV_i$ and the virtual complex $\SP S^nV$, because the $\SP$ functor turns cofibrations $B\hookrightarrow A\to A/B$ into quasi-fibrations $\SP B\to\SP A\to\SP(A/B)$ for which the homotopy exact sequence holds (see \cite[4.K]{hatcher}), hence if the pair $(A,B)$ was $(n+i)$-connected, then so is the pair $(\SP A,\SP B)$. By the previous lemma we have $\Omega\SP SA\cong\SP A$ for all connected complexes $A$, thus putting
$$\SP V:=\liminfty i\Omega^{n(i)}\SP S^{n(i)}V_i$$
works for all virtual complexes and extends the usual definition of $\SP$.
\end{prf}

\begin{lemma}\label{hur}
For any virtual complex $V$ there is a map $h\colon\Gamma V\to\SP V$ that induces a rational homotopy equivalence.
\end{lemma}

\begin{prf}
First suppose that $V$ is an existing (non-virtual) space. Then we have
$$\Gamma V=\left(\coprod_{i=1}^\infty WS_i\utimes{S_i}(V\times\ldots\times V)\right)\bigg/\sim$$
where $WS_i$ is a contractible space with a free $S_i$-action on it and ``$/\sim~$'' means gluing by some natural equivalences; for the precise definition see \cite{barec}. Now the projections
$$WS_i\utimes{S_i}(V\times\ldots\times V)\to(V\times\ldots\times V)/S_i$$
are consistent with the gluings of these spaces for different $i$-s forming $\Gamma V$ on the left-hand side and $\SP V$ on the right-hand side, hence the union of these projections forms a map $\Gamma V\to\SP V$ that we denote by $h$. This map induces in the homotopy groups
$$h_\#\colon\pi_*(\Gamma V)=\pi_*^s(V)\to H_*(V)=\pi_*(\SP V),$$
which is the stable Hurewicz homomorphism, hence by Serre's theorem (see for example \cite[theorem 18.3]{charclass}) it is a rational isomorphism.

Now the case when $V$ is just a virtual complex can be solved again by approximating $V$ with subcomplexes $V_i$ and taking numbers $n(i)$ such that for all $i$ the space $S^{n(i)}V_i$ exists and has the same $(n(i)+i)$-type as $S^{n(i)}V$. Then we replace $\Gamma V$ with $\Omega^{n(i)}\Gamma S^{n(i)}V_i$ and $\SP V$ with $\Omega^{n(i)}\SP S^{n(i)}V_i$ and repeat the same proof.
\end{prf}

Now we only have to combine the above observations to prove the initial statement (i.e. Kazarian's conjecture).

\medskip\par\noindent\textbf{Proof of theorem \ref{kazc}.\enspace\ignorespaces}
We have a sequence of equivalences
$$X_\tau^L\toverset{^*}\cong\Gamma T\nu_\tau^L\toverset{^{**}}\congq\SP T\nu_\tau^L\toverset{^{***}}\cong\SP T\varepsilon^k=\SP S^k(K_\tau^L)^+\toverset{^{**}}\congq\Gamma S^k(K_\tau^L)^+$$
where $^*$ is theorem \ref{kula}, $^{**}$ follows from lemma \ref{hur} and $^{***}$ is a consequence of lemma \ref{doldth}, since using again sufficiently large suspensions of finite dimensional approximations of Kazarian's space we obtain $H_*(T\nu_\tau^L)\cong H_{*-k}(K_\tau^L)\cong H_*(T\varepsilon^k)$ by the homological Thom isomorphism (note that this is the only part we use the positivity of $L$, which is the structure group of $\nu_\tau^L$).
~$\square$\par\medskip

\subsection{Corollaries}

Here we gather a few important consequences of theorems \ref{kazc} and \ref{kula} that will turn out to be very useful in the study of cobordism groups. We will consider $(\tau,L)$-maps to an arbitrary manifold $P^{n+k}$, which is a bit complicated, so it is useful to separately see the simpler case of maps to a Eucledian space (i.e. when $P=\R^{n+k}$).

\begin{prop}\label{cobhomo}
$\Cob_\tau^L(n,P^{n+k})\otimes\Q\cong\overset\infty{\underset{i=1}{\bigoplus}}H^i(\cpt P;H_{i-k}(K_\tau^L;\tilde\Q_{\nu_\tau^L}))$.
\end{prop}

\begin{case}
$\Cob_\tau^L(n,k)\otimes\Q\cong H_n(K_\tau^L;\tilde\Q_{\nu_\tau^L})$.
\end{case}

\begin{prf}
We have $\Cob_\tau^L(n,P^{n+k})=[\cpt P,X_\tau^L]=[\cpt P,\Gamma V_\tau^L]=[\cpt P,\Gamma T\nu_\tau^L]$, for which we can use the homomorphism induced by the map $h$ of lemma \ref{hur}. This homomorphism maps to $[\cpt P,\SP T\nu_\tau^L]$, which is the same as
$$\left[\cpt P,\prod_{i=1}^\infty K(H_i(T\nu_\tau^L),i)\right]=\left[\cpt P,\prod_{i=1}^\infty K(H_{i-k}(K_\tau^L;\tilde\Z_{\nu_\tau^L}),i)\right]=\bigoplus_{i=1}^\infty H^i(\cpt P;H_{i-k}(K_\tau^L;\tilde\Z_{\nu_\tau^L})).$$
If we compose the map induced by $h$ with these identifications and take its tensor product with $\id_\Q$, then we get an isomorphism because $h$ induces a rational isomophism.
\end{prf}

%The considerations above give a way to compute the ranks of the cobordism groups $\Cob_\tau^L(P)$. In the following we show an approach to the computation of the torsion parts as well. %We will assume here that $\tau$ has finitely many elements.
The above proposition can be extended to show details on the torsion parts of the cobordism groups $\Cob_\tau^L(P)$ as well: We will see that the map $h$ not only induces rational isomorphism, but also isomorphism in the $p$-components for sufficiently large primes $p$. To prove this, we need the following theorem of Arlettaz \cite{arl}, which will not be proved here.

\begin{thm}\label{arlthm}
For any $(l-1)$-connected spectrum $V$, the stable Hurewicz homomorphism $h_m\colon\pi^s_m(V)\to H_m(V)$ has the properties
\begin{itemize}
\item[\rm{(i)}] $\rho_1\ldots\rho_{m-l}\cdot\ker h_m=0$ for all $m\ge l+1$,
\item[\rm{(ii)}] $\rho_1\ldots\rho_{m-l-1}\cdot\coker h_m=0$ for all $m\ge l+2$,
\end{itemize}
where $\rho_i$ is the exponent of the $i$-th stable homotopy group of spheres $\pi^s(i)$.
\end{thm}

Now we can get to the torsion parts of the cobordism groups $\Cob_\tau^L(P)$.

\begin{prop}\label{cobtors}
For any prime $p>\frac n2+1$ (where $n\ge2$) there is an isomorphism of $p$-primary parts
$$\big[\Cob_\tau^L(n,P^{n+k})\big]_p\cong\displaystyle\bigoplus_{i=1}^\infty\big[H^i(\cpt P;H_{i-k}(K_\tau^L;\tilde\Z_{\nu_\tau^L}))\big]_p.$$
\end{prop}

\begin{case}
For any prime $p>\frac n2+1$ (where $n\ge2$) there is an isomorphism of $p$-primary parts
$$\big[\Cob_\tau^L(n,k)\big]_p\cong\big[H_n(K_\tau^L;\tilde\Z_{\nu_\tau^L})\big]_p.$$
\end{case}

\begin{prf}
Recall that the map $h\colon\Gamma T\nu_\tau^L\to\SP T\nu_\tau^L$ of lemma \ref{hur} was construced such that in the homotopy groups it induces the stable Hurewicz homomorphism for the virtual complex $T\nu_\tau^L$. Now any spectrum defined by $T\nu_\tau^L$ (according to remark \ref{spec}) is $(k-1)$-connected, so we can apply the theorem of Arlettaz with the substitutions $l:=k$ and $m:=n+k$.

Serre proved in \cite{serre} that the number $\rho_i$ is not divisible by the prime $p$, if $p>\frac i2+1$, hence we obtain that $h$ induces an isomorphism between the $p$-components of the homotopy groups for $p>\frac n2+1$ until dimension $n+k$. Now (as in the proof of proposition \ref{cobhomo}) we use the isomorphisms
$$\Cob_\tau^L(n,P^{n+k})\displaystyle\cong[\cpt P,\Gamma T\nu_\tau^L]~~~~\text{and}~~~~[\cpt P,\SP T\nu_\tau^L]\cong\bigoplus_{i=1}^\infty H^i(\cpt P;H_{i-k}(K_\tau^L;\tilde\Z_{\nu_\tau^L}))$$
and we get the statement of the present proposition.
\end{prf}

In the remaining part of this section we assume that $L$ is positive (i.e. the bundle $\nu_\tau^L$ is orientable).

We will give a (rational) analogue for $(\tau,L)$-maps of the Pontryagin--Thom theorem claiming that two manifolds are cobordant iff their characteristic numbers coincide. This will also be analogous to the Conner--Floyd theorem in \cite{conflo} about the characteristic numbers of rational bordism classes. The first task is to define the characteristic numbers of a $(\tau,L)$-map.

\begin{defi}
For a $(\tau,L)$-map $f\colon M^n\to P^{n+k}$ consider the map $\kappa_f\colon M\to K_\tau^L$ (see remark \ref{kappaf}). For any cohomology class $x\in H^*(K_\tau^L;\Q)$,
\begin{enumerate}
\item the class $\kappa_f^*(x)\in H^*(M;\Q)$ is called the $x$-characteristic class of $f$,
\item for any $y\in H^{n-*}(P;\Q)$, the number $\la\kappa_f^*(x)\smallsmile f^*(y),[M]\ra\in\Q$ is called the $(x,y)$-characteristic number of $f$.
\end{enumerate}
\end{defi}

\begin{thm}\label{charcob}
There is an isomorphism
$$\Cob_\tau^L(n,P)\otimes\Q\cong\Hom\left(\displaystyle\bigoplus_{i=0}^nH^i(K_\tau^L;\Q)\otimes H^{n-i}(P;\Q),\Q\right)$$
given by assigning to a $(\tau,L)$-map $f\colon M^n\to P^{n+k}$ the homomorphism $x\otimes y\mapsto\la\kappa_f^*(x)\smallsmile f^*(y),[M]\ra$.
\end{thm}

\begin{case}
There is an isomorphism
$$\Cob_\tau^L(n,k)\otimes\Q\cong\Hom(H^n(K_\tau^L;\Q),\Q)$$
given by assigning to a $(\tau,L)$-map $f\colon M^n\to\R^{n+k}$ the homomorphism $x\mapsto\la\kappa_f^*(x),[M]\ra$.
\end{case}

In other words this theorem states that the characteristic numbers are well-defined invariants of rational cobordism classes of $(\tau,L)$-maps and they completely determine the class they stand for. Moreover, given any possible collection of rational numbers assigned to the corresponding cohomology classes, there is a unique rational $(\tau,L)$-cobordism class with those characteristic numbers.

\medskip\begin{prf}
%We have $\la\kappa_f^*(x),[M]\ra=\la x,\kappa_{f*}[M]\ra$ and the set of these numbers for all $x\in H^n(K_\tau^L;\Q)$ determines the homology class $\kappa_{f*}[M]$. We will prove the identity $\kappa_{f*}[M]=\varphi([f])$ in a more general setting in lemma \ref{kappafi}, hence (using that $\varphi$ is a rational isomorphism) we have $\kappa_{f_0*}[M]=\kappa_{f_1*}[M]$ iff $f_0\sim f_1$. Thus we get a monomorphism $\Cob_\tau^L(n,k)\otimes\Q\to\Hom(H^n(K_\tau^L);\Q)$, but the source and target have equal dimensions over $\Q$, so it is an isomorphism.
Let $f\colon M\to P$ be a $(\tau,L)$-map. Recall the Pontryagin--Thom type construction in theorem \ref{ptpullback} or in the addendum of theorem \ref{embclass} (also using theorem \ref{cob=emb}), which assigns to $f$ a map
$$\kappa'_f:=\kappa_P\colon\cpt P\to X_\tau^L.$$

Now $\kappa'_f$ maps to a finite dimensional approximation $\Omega^mS^m T\nu_\tau^L$ of $X_\tau^L=\Gamma T\nu_\tau^L$ (here $\nu_\tau^L$ is restricted to a finite dimensional approximation of Kazarian's space, but to make notations simpler, we omit this) such that the bundle $\nu_\tau^L\oplus\varepsilon^m$ is non-virtual. Recall that the maps $A\to\Omega^mB$ bijectively correspond to maps $S^mA\to B$ for any (pointed) spaces $A,B$; we call the maps that correspond to each other adjoint. Hence this $\kappa'_f$ can be identified with its adjoint map $S^m\cpt P\to S^mT\nu_\tau^L=T(\nu_\tau^L\oplus\varepsilon^m)$, which we will also denote by $\kappa'_f$.

We put $e:=f\times i\colon M\hookrightarrow P\times\R^m\subset S^m\cpt P$ the $(\tau,L)$-embedding corresponding to $f$ (with the vertical vector fields and foliation; we may assume that $m$ is large enough so that this exists uniquely up to isotopy). Observe that $\kappa'_f$ can be written (by its construction) as a composition $\kappa'_f=t\circ q$, where $q\colon S^m\cpt P\to T\nu_e$ is the quotient by the complement of a tubular neighborhood of $e(M)\subset P\times\R^m$ (hence a degree-$1$ map) and $t\colon T\nu_e\to T(\nu_\tau^L\oplus\varepsilon^m)$ is a fibrewise isomorphism. Moreover, the restriction of $t$ to the zero-section is the map $\kappa_f\colon M\to K_\tau^L$, hence we have a commutative diagram
$$\xymatrix{
S^m\cpt P\ar[r]^(.55){q}\ar@/^1.3pc/[rr]^{\kappa'_f} & T\nu_e\ar[r]^(.35){t} & T(\nu_\tau^L\oplus\varepsilon^m)\\
& M\ar@{_(->}[ul]^{f\times i=e}\ar@{_(->}[u]\ar[r]^{\kappa_f} & K_\tau^L\ar@{^(->}[u]
}$$

The image of the cobordism class $[f]\in\Cob_\tau^L(n,P)$ in
$$\bigoplus_{i=1}^{n+k} H^i(\cpt P;H_{i-k}(K_\tau^L;\Q))\congq\bigoplus_{i=1}^{n+k}\Hom(H_i(\cpt P;\Q),H_{i-k}(K_\tau^L;\Q))$$
(according to proposition \ref{cobhomo}) is the homomorphism that we get by the composition
$$H_i(\cpt P;\Q)\cong H_{i+m}(S^m\cpt P;\Q)\xra{\kappa'_{f*}}H_{i+m}(T(\nu_\tau^L\oplus\varepsilon^m);\Q)\xra{\varphi^{-1}}H_{i-k}(K_\tau^L;\Q)$$
for all $i$, where $\varphi$ denotes the Thom isomorphism. The $\Hom$-dual of this map is a map between the compactly supported cohomology groups given by the formula
\pagebreak
\begin{alignat*}2
H_c^{i-k}(K_\tau^L;\Q)\to&~H_c^i(P;\Q)\\
x\mapsto&~S^{-m}q^*t^*\varphi(x)=S^{-m}q^*\varphi\kappa_f^*(x)=S^{-m}(f\times i)_!\kappa_f^*(x)=\\
&=S^{-m}(S^mf)_!\kappa_f^*(x)=f_!\kappa_f^*(x).
\end{alignat*}

Now the Poincaré duality implies that the cohomology classes $f_!\kappa_f^*(x)$ (for any $x\in H_c^{i-k}(K_\tau^L;\Q)$) are completely determined by the numbers $\la f_!\kappa_f^*(x)\smallsmile y,[P]\ra~(y\in H_c^{n-i+k}(P;\Q))$. We can rewrite these numbers in the following way:
\begin{alignat*}2
\la f_!\kappa_f^*(x)\smallsmile y,[P]\ra&=\la f_!(\kappa_f^*(x)\smallsmile f^*(y)),[P]\ra=\la S^m(f_!(\kappa_f^*(x)\smallsmile f^*(y))),[S^m\cpt P]\ra=\\
&=\la q^*\varphi(\kappa_f^*(x)\smallsmile f^*(y)),[S^m\cpt P]\ra=\la \varphi(\kappa_f^*(x)\smallsmile f^*(y)),[T\nu_e]\ra=\\
&=\la \kappa_f^*(x)\smallsmile f^*(y),[M]\ra.
\end{alignat*}

We can set this to be an arbitrary rational number depending on $x$ and $y$ and the compact support in the cohomologies does not play a role here, since the image of $M$ is compact in both $P$ and $K_\tau^L$. This proves the statement of the theorem.
\end{prf}

\begin{ex}
If $\tau$ is the set of all possible singularities of $k$-codimensional germs, then we have $K_\tau^L\cong BL$, because the fibre of the bundle $K_\tau^L\to BL$ (see definition \ref{kazj}) is the space of all polynomial maps, which is contractible. Now setting $L(k):=\SO(k)$, we get $\Cob_\tau^{\SO}(n,k)=\Omega_n(\R^{n+k})=\Omega_n$ and the above propositions mean
$$\Omega_n\otimes\Q\cong H_n(B\SO;\Q)\cong\Hom(H^n(B\SO;\Q),\Q),$$
which is a well-known theorem of Thom.
\end{ex}

%-------------------------------------------------------------------------------------------------------------------------------
%-------------------------------------------------------------------------------------------------------------------------------

\chapter{The key fibration}

In the previous chapter we saw that the investigation of cobordisms of singular maps is equivalent to the homotopical investigation of classifying spaces. One of the key tools of this investigation is a fibre bundle
$$X_\tau^L\xra{X_{\tau'}^L}\Gamma T\tilde\xi_\eta^L,$$
where $[\eta]$ is a top singularity in $\tau$ and $\tau'$ denotes $\tau\setminus\{[\eta]\}$ as before (as usual, $\tau$ denotes a set of singularities and $L$ is a stable linear group). By ``top singularity'' we mean such a singularity that is maximal in $\tau$ with the natural partial order. We call this the key fibration (or key bundle) to indicate its importance.

The existence of such a fibre bundle was conjectured by Szabó and proved by Szűcs \cite{hosszu} as a corollary of the description of classifying spaces with virtual complexes (corollary \ref{gammav}), later a constructive proof was given by Terpai \cite{key} in a bit more general setting. In this chapter we recall these proofs, then see some important properties and corollaries of this fibration from \cite{hominv}, \cite{hosszu} and \cite{ctrl}.

%-------------------------------------------------------------------------------------------------------------------------------

\section{Existence of the key fibration}

\subsection{Existence derived from virtual complexes}

It was proved in \cite{barec} that the $\Gamma$ functor turns cofibrations into fibrations. We need an analogue of this for virtual complexes.

\begin{lemma}\label{cofibfib}
Given virtual complexes $V,V'$ and a cofibration $V'\into V\to V/V'$, we have a fibration $\Gamma V'\into\Gamma V\to\Gamma(V/V')$.
\end{lemma}

\begin{prf}
For each natural number $i$ fix a number $n(i)$ and subcomplexes $V_i$ and $V'_i$ of $V$ and $V'$ respectively such that for all $n\ge n(i)$ the spaces $S^nV_i$ and $S^nV'_i$ exist and their $(n+i)$-homotopy types are that of $S^nV$ and $S^nV'$ respectively. Now applying the $\Gamma$ functor, we get a fibration
$$\Gamma S^{n(i)}V'_i\hookrightarrow\Gamma S^{n(i)}V_i\to\Gamma S^{n(i)}(V_i/V'_i).$$

Recall that whenever we have a fibration $E\xra{F}B$, we can (homotopically) continue it to the left infinitely, obtaining a sequence of fibrations (called the resolvent of $E\xra{F}B$)
$$\ldots\to\Omega^2F\to\Omega^2E\to\Omega^2B\to\Omega F\to\Omega E\to\Omega B\to F\to E\to B.$$
Now applying the resolvent of the fibration above, we can consider the fibration
$$\Omega^{n(i)}\Gamma S^{n(i)}V'_i\hookrightarrow\Omega^{n(i)}\Gamma S^{n(i)}V_i\to\Omega^{n(i)}\Gamma S^{n(i)}(V_i/V'_i).$$
Of course, $\Omega^{n(i)}\Gamma S^{n(i)}=\Gamma$, hence this is a fibration $\Gamma V_i\xra{\Gamma V'_i}\Gamma(V_i/V'_i)$. Thus by setting $i\to\infty$, we obtain the fibration as claimed.
\end{prf}

\begin{thm}\label{keythm}
If $L$ is a stable linear group, $\tau$ is a set of singularities (without fixed multiplicites), $[\eta]$ is a top singularity in $\tau$ and we put $\tau':=\tau\setminus\{[\eta]\}$, then there is a fibration
$$X_\tau^L\xra{X_{\tau'}^L}\Gamma T\tilde\xi_\eta^L.$$
\end{thm}

\begin{prf}
Observe that by the construction of the virtual complex $V_\tau^L$ we have a cofibration of virtual complexes
$$V_{\tau'}^L\into V_\tau^L\to V_\tau^L/V_{\tau'}^L=T\tilde\xi_\eta^L.$$
Now we get the statement of the theorem by the previous lemma.
\end{prf}

\begin{rmk}\label{keyrmk}
We can define the cobordism group of immersions with normal bundles induced from a given bundle in a standard way (analogously to definition \ref{cobtau}). The classifying space of cobordisms of immersions with normal bundles induced from the bundle $\zeta$ is the space $\Gamma T\zeta$ (see \cite{wells}), that is, for any manifold $P$, the cobordism group $\Imm^\zeta(P)$ of immersions to $P$ with normal bundles induced from $\zeta$ is isomorphic to $[P,\Gamma T\zeta]$.

Now if we apply the functor $[P,\cdot]$ to the key fibration $X_\tau^L\to\Gamma T\tilde\xi_\eta^L$, then we get a map $\Cob_\tau^L(P)\to\Imm^{\tilde\xi_\eta^L}(P)$. It is not hard to see from the proofs above (and it also follows from the proof in the next subsection) that this map is the same as the one we get by assigning to the cobordism class of a $(\tau,L)$-map $f\colon M\to P$ the cobordism class of the immersion $f|_{\eta(f)}\colon\eta(f)\imto P$ (clearly $f|_{\eta(f)}$ is an immersion with normal bundle induced from $\tilde\xi_\eta^L$).

This observation has a very nice consequence about the problem of eliminating certain singularities of a given map. Namely that the only obstruction for a $(\tau,L)$-map $f$ to be $(\tau,L)$-cobordant to a $(\tau',L)$-map is the cobordism class of $f|_{\eta(f)}$ in $\Imm^{\tilde\xi_\eta^L}(P)$.
\end{rmk}

\subsection{Terpai's construction}

In this subsection we will be in a slightly more general setting. Namely, let $L,\tau,\eta$ and $\tau'$ be as before, let $\underline\tau'$ denote the set of all possible multisingularities composed of the elements of $\tau'$ and put
$$\underline\tau_r:=\{\underline\vartheta+i[\eta]\mid\underline\vartheta\in\underline\tau',i=0,\ldots,r\}$$
for all $r\in\N$. We will use the notation $X_r:=X_{\underline\tau_r}^L$ and denote by $\Gamma_r$ the classifying space for cobordisms of immersions with normal bundles induced from $\tilde\xi_\eta^L$ and with at most $r$-tuple points (see \cite{madar}, \cite{limm} or \cite{s2}).

\begin{thm}\label{terpthm}
For all $r\in\N$ there is a fibration
$$X_r\xra{X_{\tau'}^L}\Gamma_r.$$
\end{thm}

This implies the statement of theorem \ref{keythm}, since by converging with $r$ to $\infty$, we get the same fibration with $X_\infty$ and $\Gamma_\infty$ in the place of $X_r$ and $\Gamma_r$ respectively, and of course $X_\infty=X_\tau^L$ and $\Gamma_\infty=\Gamma T\tilde\xi_\eta^L$.

\medskip\begin{prf}
Recall the universal $(\tau,L)$-map $f_\tau^L\colon Y_\tau^L\to X_\tau^L$ of theorem \ref{ptpullback}. This is not a real $(\tau,L)$-map itself, as it is a map between infinite dimensional spaces, but it was constructed using direct limits; in particular its restriction to the $\eta$-stratum $\overset\infty{\underset{r=1}\bigcup}BG_{r[\eta]}^L$ is a direct limit of immersions.

Thus the direct limit of the classifying maps of the cobordism classes of these immersions is a well-defined map $\chi_\tau^L\colon X_\tau^L\to\Gamma T\tilde\xi_\eta^L$ and its restriction to $X_r$ is a map $X_r\to\Gamma_r$. We will prove the theorem for a fixed $r$ and we proceed by directly checking the homotopy lifting property for the map $\chi_\tau^L|_{X_r}$.

Fix a finite CW-complex $P$ and maps $g,h$ such that the square below commutes; we want to find a map $H$ for which the diagram is still commutative.
$$\xymatrix{
P\ar@{_(->}[d]_{\id_P\times0}\ar[r]^{g} & X_r\ar[d]^{\chi_\tau^L|_{X_r}}\\
P\times[0,1]\ar@{-->}[ur]^H\ar[r]^(.6)h & \Gamma_r
}$$

We may assume that $P$ is a manifold with boundary, since we can embed $P$ into a Eucledian space and then replace it with a small closed neighbourhood that $P$ is a deformation retract of. By part \ref{pt1} of theorem \ref{ptpullback}, the map $g\colon P\to X_r$ corresponds to a $(\underline\tau_r,L)$-map $f\colon M^n\to P^{n+k}=P\times\{0\}$. By the analogous pullback property of the classifying space $\Gamma_r$, the map $h\colon P\times[0,1]\to\Gamma_r$ corresponds to an immersion $f_\eta\colon N^{n-c+1}\imto P\times[0,1]$ (if the dimension of $\tilde\xi_\eta^L$ is $c+k$) with at most $r$-tuple points and normal bundle induced from $\tilde\xi_\eta^L$ with the properties $f_\eta^{-1}(P\times\{0\})=\eta(f)$ and $f_\eta|_{\eta(f)}=f|_{\eta(f)}$. The proposed lift $H$ of the map $h$ would correspond to a $(\underline\tau_r,L)$-map $F\colon W^{n+1}\to P\times[0,1]$ that extends $f\colon M\to P\times\{0\}$ and for which $\eta(F)=N$ and $F|_{\eta(F)}=f_\eta$.

We remark here that $P\times[0,1]$ (and so $N$ and the proposed $W$ too) is not a manifold with boundary but a manifold with corners (that is, it can be covered by Eucledian neighbourhoods diffeomorphic to $[0,\varepsilon)^m\times\R^{n+k+1-m}~(m\in\N)$). Fortunately the notions of $(\underline\tau_r,L)$-maps and immersions with normal bundles induced from $\tilde\xi_\eta^L$ extend to manifolds with corners and the pullback properties of the classifying spaces are also true in this extended sense, hence there is no problem here.

We will assume that the images of $M$ and $N$ in $P\times[0,1]$ are disjoint from $\partial P\times[0,1]$, so they do not have corners. We can assume this, as otherwise we can apply each of the following steps with $\partial P$ in the place of $P$, obtaining an extension over the boundary, then apply the same steps again to extend this over the whole manifold.

The normal structure of the map $f_\eta$ and theorem \ref{univ} (together with its analogues later) allow us to extend the map $f_\eta$ to a disk bundle $U$ over $N$, which is induced from the union of the disk bundles $D\xi_{i[\eta]}^L$ and its mapping to $P\times[0,1]$ is induced from the global normal forms $\Phi_{i[\eta]}^L$ (for $i=1,\ldots,r$), so it maps to a tubular neighbourhood of $f_\eta(N)$. This extended map will be denoted by $F_\eta$.

Now the boundary of $U$ is the union of a sphere bundle over $N$ denoted by $\partial_SU$, a disk bundle over $f_\eta^{-1}(P\times\{0\})$ denoted by $\partial_0U$ and a disk bundle over $f^{-1}_\eta(P\times\{1\})$ denoted by $\partial_1U$. Here $\partial_0U$ is a (closed) tubular neighbourhood of $\eta(f)$ and we can assume that $F_\eta|_{\partial_0U}$ is the same as $f|_{\partial_0U}$.

Of course, nothing changes if we replace the map $f\colon M\to P\times[0,1]$ with $f\times\id_{[-\varepsilon,0]}\colon M\times[-\varepsilon,0]\to P\times[-\varepsilon,0]$ for a number $\varepsilon>0$ (i.e. if we add a $[-\varepsilon,0]$-collar to the map $f$). Now we are searching for a map $F$ that extends this $f$ with the added $[-\varepsilon,0]$-collar and $F_\eta$ such that no new $\eta$-points of $F$ arise besides the $\eta$-strata in $M\times[-\varepsilon,0]$ and in $U$.

\begin{figure}[H]
\begin{center}
\centering\includegraphics[scale=0.1]{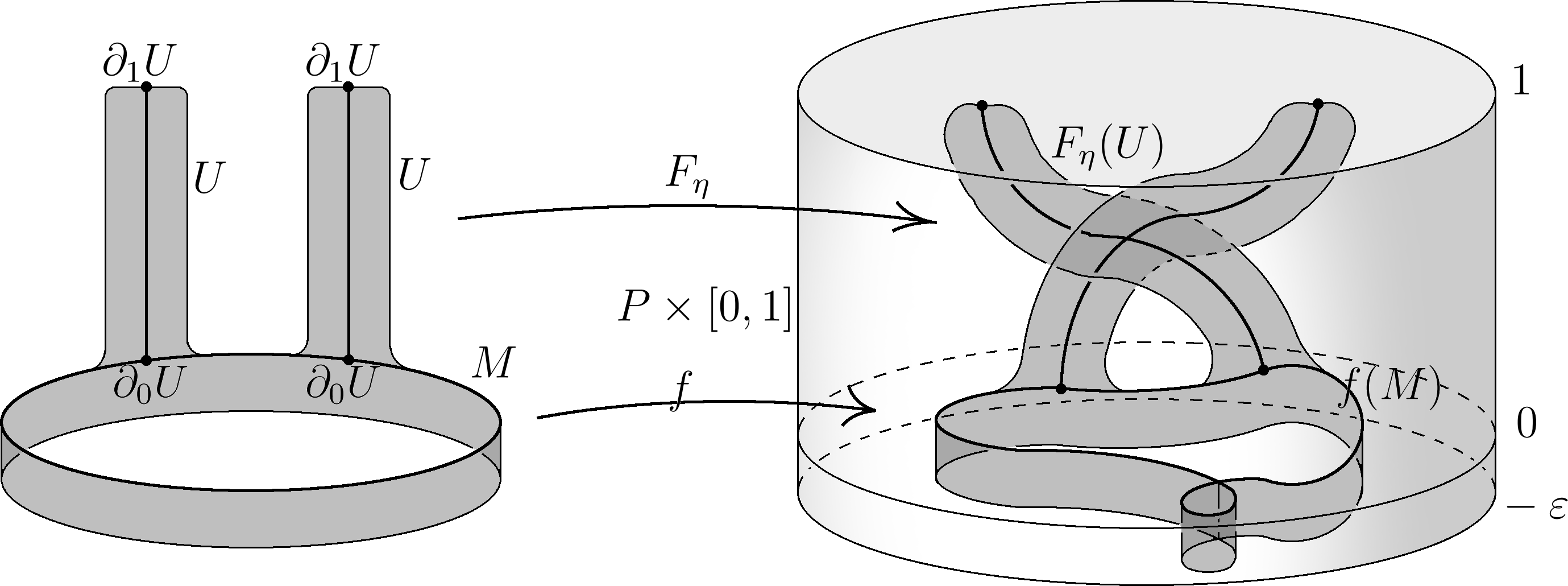}\label{kep5}
\begin{changemargin}{2cm}{2cm} 
\caption{\hangindent=1.4cm\small We represent $M$ and $P\times[0,1]$ by the ellipse on the left and the cylinder on the right with the $[-\varepsilon,0]$-collars added; the two indicated points on $M$ represent $\eta(f)$ and the vertical segments starting from them are $N$; the disk bundle $U$ and its image are also indicated.}
\end{changemargin} 
\vspace{-1.3cm}
\end{center}
\end{figure}

Now $V:=(M\times[-\varepsilon,0])\cup U$ is a manifold with boundary $\partial V=(M\times\{-\varepsilon\})\cup(M\setminus\partial_0U)\cup\partial_SU\cup\partial_1U$. The map $F_V:=(f\times\id_{[-\varepsilon,0]})\cup F_\eta\colon V\to P\times[-\varepsilon,1]$ is almost a $(\underline\tau_r,L)$-map that extends $f$ and $F_\eta$, the only differece is that it does not map the boundary of the source to the boundary of the target. In order to change $F_V$ to a map that maps boundary to boundary (hence proving the theorem), we first prove the following.

\medskip\begin{sclaim}
There is a diffeotopy $\varphi_t~(t\in[0,1])$ of $P\times[0,1]$ for which $\varphi_0=\id_{P\times[0,1]}$ and for any $p\in\partial_SU$, the subspace $d\varphi_1(dF_V(T_p\partial_SU))\subset T_{\varphi_1(F_V(p))}(P\times[0,1])$ does not contain the tangent line of $[0,1]$.
\end{sclaim}

\begin{sprf}
The restriction $F_\eta|_{\partial_SU}$ maps to the ``boundary'' of a tubular neighbourhood of $f_\eta(N)$; actually this ``boundary'' is just a $1$-codimensional subcomplex of $P\times[0,1]$ which is the image of an immersion. We can choose a vector field on this subcomplex (i.e. a section of $T(P\times[0,1])$ restricted to the subcomplex) that is nowhere tangent to it.

If we denote the restriction of such a vector bundle to $F_\eta(\partial_SU)$ by $u$, then $u$ is pointwise independent of the subspaces $dF_\eta(T_p\partial_SU)$ for all $p\in\partial_SU$. Now the stratified compression theorem \ref{sct} applied to the multisingularity stratification of $\partial_SU$ yields a diffeotopy of $P\times[0,1]$ that turns $u$ vertical (i.e. everywhere parallel to $T[0,1]$). Using that $u$ was everywhere independent of the subspaces $dF_\eta(T_p\partial_SU)=dF_V(T_p\partial_SU)$, the images of these subspaces cannot contain the vertical tangent lines and this is what we wanted to achieve.
\end{sprf}

We can trivially extend $\varphi_1$ to a diffeomorphism of $P\times[-\varepsilon,1]$. We will use the notations $\tilde f:=\varphi_1\circ f,\tilde f_\eta:=\varphi_1\circ f_\eta,\tilde F_\eta:=\varphi_1\circ F_\eta$ and $\tilde F_V:=\varphi_1\circ F_V$. Define
$$V':=\{(p,t)\in\partial V\times[0,1]\mid\tilde F_V(p)=(q,s)\in P\times[0,1],t\in[s,1]\}$$
and $\tilde F_{V'}(p,t):=(\tilde F_V(p),t)$ for all $(p,t)\in V'$. Now we can glue $V'$ to $V$ by the map $(p,s)\mapsto p$ (with the above notation) to form the manifold $W$, which is mapped to $P\times[-\varepsilon,1]$ by the map $\tilde F:=\tilde F_V\cup\tilde F_{V'}$. This map is not necessarily smooth where the gluing happened, but it can easily be smoothed, so we may assume that it is.

Now $\tilde F$ is a $(\underline\tau_r,L)$-map, it extends $\tilde f\times\id_{[-\varepsilon,0]}$ and $\tilde f_\eta$ and its $\eta$-stratum is $(\eta(f)\times[-\varepsilon,0])\cup N$. Thus putting $F:=\varphi_1^{-1}\circ\tilde F$ satisfies every condition we need.
\end{prf}

In the rest of this thesis (as in the above proof) we will use the symbol $\chi_\tau^L$ to denote the key fibration $X_\tau^L\to\Gamma T\tilde\xi_\eta^L$.

%-------------------------------------------------------------------------------------------------------------------------------

\section{Inducing map of the key fibration}

In the present section we will give a purely homotopy theoretical description of the key fibration $\chi_\tau^L$ (where $L,\tau,\eta$ and $\tau'$ are fixed as before) based on \cite{hominv}. This fibration can be induced from the universal $X_{\tau'}^L$-bundle, that is, there is a pullback diagram
$$\xymatrix{
X_\tau^L\ar[d]_{X_{\tau'}^L}\ar[r] & EX_{\tau'}^L\ar[d]^{X_{\tau'}^L}\\
\Gamma T\tilde\xi_\eta^L\ar[r]^{b_\eta^L} & BX_{\tau'}^L
}$$
where $EX_{\tau'}^L$ is a contractible space, $BX_{\tau'}^L$ is a space with the property $\Omega BX_{\tau'}^L=X_{\tau'}^L$ and the map $b_\eta^L$ is called the inducing map of the key fibration.

We will describe $b_\eta^L$ using the attaching map $\tilde\rho_\eta^L$ of $D\tilde\xi_\eta^L$ to $X_{\tau'}^L$ in the construction of the classifying space for cobordisms of $\tau$-maps with no multiple $\eta$-points (see subsection \ref{prf}; the attaching map was only denoted by $\tilde\rho$ there).

\subsection{Formulation of the theorem}

The description of the map $b_\eta^L$ will be summarised by a theorem, but to state this theorem we need some introduction. In this subsection we give this introduction using a series of lemmas and then we formulate the theorem. These lemmas will not be proved here in order to get to the statement of the theorem faster; the proofs will be given in the next subsection.

First we consider the key fibration $\chi_\tau^L\colon X_\tau^L\to\Gamma T\tilde\xi_\eta^L$ over the open disk bundle $\tightoverset{\,_\circ~~~}{D\tilde\xi_\eta^L}\subset T\tilde\xi_\eta^L\subset\Gamma T\tilde\xi_\eta^L$ (i.e. the Thom space without its special point).

\begin{lemma}\label{l1}
$(\chi_\tau^L)^{-1}\big(\tightoverset{\,_\circ~~~}{D\tilde\xi_\eta^L}\big)=\tightoverset{\,_\circ~~~}{D\tilde\xi_\eta^L}\times X_{\tau'}^L$.
\end{lemma}

In other words, the key bundle is trivial over this open disk bundle. Now we add the extra point to $\tightoverset{\,_\circ~~~}{D\tilde\xi_\eta^L}$ to form the Thom space $T\tilde\xi_\eta^L$ and consider $X_T:=(\chi_\tau^L)^{-1}(T\tilde\xi_\eta^L)$. Recall that $X_{\tau'}^L$ is an H-space (see remark \ref{h}); we will denote the multiplication in this H-space by $x\cdot y~(x,y\in X_{\tau'}^L)$ and we define
$$\tilde\rho_T\colon S\tilde\xi_\eta^L\times X_{\tau'}^L\to X_{\tau'}^L;~(s,x)\mapsto\tilde\rho_\eta^L(s)\cdot x.$$

\begin{lemma}\label{l3}
$X_T=(D\tilde\xi_\eta^L\times X_{\tau'}^L)\usqcup{\tilde\rho_T}X_{\tau'}^L$.
\end{lemma}

Next we describe the fibration
$$\chi_\tau^L|_{X_T}\colon X_T\xra{X_{\tau'}^L}T\tilde\xi_\eta^L.$$
We denote its inducing map by $b_T\colon T\tilde\xi_\eta^L\to BX_{\tau'}^L$ and note that $b_T$ is just the restriction of $b_\eta^L$ to the Thom space $T\tilde\xi_\eta^L$.

Since the key fibration was trivial over the open disk bundle, it is also trivial over the zero-section $BG_\eta^L$, hence we can assume that the inducing map $b_T$ maps $BG_\eta^L$ to a single point. This means that $b_T$ factors through the quotient map $q\colon T\tilde\xi_\eta^L\to T\tilde\xi_\eta^L/BG_\eta^L$ that contracts the zero-section, that is, there is a map $\sigma\colon T\tilde\xi_\eta^L/BG_\eta^L\to BX_{\tau'}^L$ with the property $b_T=\sigma\circ q$.

Observe that $T\tilde\xi_\eta^L/BG_\eta^L$ coincides with the suspension $S(S\tilde\xi_\eta^L)$ of the sphere bundle $S\tilde\xi_\eta^L$, hence $\sigma$ is a map $S(S\tilde\xi_\eta^L)\to BX_{\tau'}^L$. Recall again that there is a bijective correspondence, called the adjoint correspondence, between the maps $A\to\Omega B$ and the maps $SA\to B$ for any (pointed) spaces $A,B$. Now the adjoint of the map $\sigma$ is a map $S\tilde\xi_\eta^L\to\Omega BX_{\tau'}^L=X_{\tau'}^L$.

\begin{lemma}\label{l2}
The adjoint of $\sigma$ is the attaching map $\tilde\rho_\eta^L\colon S\tilde\xi_\eta^L\to X_{\tau'}^L$.
\end{lemma}

Recall that $X_{\tau'}^L=\Gamma V_{\tau'}^L$ (by theorem \ref{gammav}), hence we have $\Omega BX_{\tau'}^L=X_{\tau'}^L=\Omega\Gamma SV_{\tau'}^L$ and we can define $BX_{\tau'}^L$ as $\Gamma SV_{\tau'}^L$, since its path fibration has contractible total space and fibre $X_{\tau'}^L$. Now the inducing map $b_\eta^L$ is such that both its source space $\Gamma T\tilde\xi_\eta^L$ and its target space $BX_{\tau'}^L=\Gamma SV_{\tau'}^L$ are infinite loop spaces.

\begin{lemma}\label{l4}
The map $b_\eta^L\colon\Gamma T\tilde\xi_\eta^L\to\Gamma SV_{\tau'}^L$ is an infinite loop map. 
\end{lemma}

This implies that $b_\eta^L$ is completely determined by its restriction $b_T$ to $T\tilde\xi_\eta^L$, since any map $g\colon A\to B$ from a (pointed) space to an infinite loop space extends to an infinite loop map $g_\ext\colon\Gamma A\to B$ in a homotopically unique way by \cite[pp. 42--43]{may}.

\begin{lemma}\label{l5}
If $f\colon A\to A'$ is a map between any (pointed) spaces and $g\colon A\to B$ and $g'\colon A'\to B$ are maps to an infinite loop space $B$ for which $g=g'\circ f$, then $g_\ext=g'_\ext\circ\Gamma f$.
\end{lemma}

A direct consequence of these lemmas is a description of the inducing map $b_\eta^L$ as follows.

\begin{thm}\label{induc}
Let $\sigma\colon S(S\tilde\xi_\eta^L)\to BX_{\tau'}^L$ be the adjoint of the gluing map $\tilde\rho_\eta^L$ and let $q\colon T\tilde\xi_\eta^L\to T\tilde\xi_\eta^L/BG_\eta^L=S(S\tilde\xi_\eta^L)$ be the quotient map. Then we have
$$b_\eta^L=(\sigma\circ q)_\ext.$$
\end{thm}

\begin{rmk}
This theorem can also be visualised on a commutative diagram:
$$\xymatrix@C=.75pc{
X_T\ar@{^(->}[rr]\ar[d]_{X_{\tau'}^L} && X_\tau^L\ar[rrrr]\ar[d]^{X_{\tau'}^L} &&&& EX_{\tau'}^L\ar@{=}[r]^(.45){\sims{1.8}}\ar[dd]^{X_{\tau'}^L} & \ast~~~~~~~~ \\
T\tilde\xi_\eta^L\ar@{^(->}[rr]\ar[drrrrrr]^{b_T}\ar[d]_q && \Gamma T\tilde\xi_\eta^L\ar[drrrr]^{(b_T)_\ext=b_\eta^L} &&&&&\\
T\tilde\xi_\eta^L/BG_\eta^L\ar@{=}[rr] && S(S\tilde\xi_\eta^L)\ar[rrrr]^(.53)\sigma && \ar@{<->}[d]^{\text{adjoint}} && BX_{\tau'}^L\ar@{=}[r] & \Gamma SV_{\tau'}^L~ \\
&& S\tilde\xi_\eta^L\ar[rrrr]_(.53){\tilde\rho_\eta^L} &&&& X_{\tau'}^L\ar@{=}[r] & \Omega BX_{\tau'}^L~
}$$
\end{rmk}

\subsection{Proofs of the lemmas}

\medskip\par\noindent\textbf{Proof of lemma \ref{l1}.\enspace\ignorespaces}
Recall the construction of the space $X_{\tau'}^L$ (section \ref{prf}): We glue together the disk bundles
$$D_{\underline\vartheta}:=D\tilde\xi_{\underline\vartheta}^L=\prod_{i=1}^r(D\tilde\xi_{\vartheta_i}^L)^{m_i}$$
along their boundaries for all multisingularities $\underline\vartheta=m_1[\vartheta_1]+\ldots+m_r[\vartheta_r]$ composed of the elements of $\tau'$ according to an order that is compatible with the natural partial order. The attaching map $\tilde\rho_{\underline\vartheta}^L$ of $D_{\underline\vartheta}$ maps to the subspace corresponding to multisingularities less complicated than $\underline\vartheta$. 

To form the space $X_\tau^L$, we follow the same procedure using the disk bundles $D_{\underline\vartheta}\times(D\tilde\xi_\eta^L)^r$ for all natural numbers $r$. Now if we project each such disk bundle onto the factor $(D\tilde\xi_\eta^L)^r$ and attach these correspondingly, then we get the space $\Gamma T\tilde\xi_\eta^L$. Clearly the projections are compatible with the gluings that form the space $X_\tau^L$, hence we get a map $X_\tau^L\to\Gamma T\tilde\xi_\eta^L$, which is $\chi_\tau^L$.

This way it is easy to see that we obtain $(\chi_\tau^L)^{-1}\big(\tightoverset{\,_\circ~~~}{D\tilde\xi_\eta^L}\big)$ by gluing the spaces $D_{\underline\vartheta}\times\tightoverset{\,_\circ~~~}{D\tilde\xi_\eta^L}$ together along their boundaries by the maps $\tilde\rho_{\underline\vartheta}^L\times\id_{\toverset{_\circ~~}{D\tilde\xi_\eta^L}}$. This procedure results in the space $X_{\tau'}^L\times\tightoverset{\,_\circ~~~}{D\tilde\xi_\eta^L}$ as claimed.
~$\square$\par

\medskip\par\noindent\textbf{Proof of lemma \ref{l3}.\enspace\ignorespaces}
The space $X_T$ is a direct product over the open disk bundle $\tightoverset{\,_\circ~~~}{D\tilde\xi_\eta^L}$, hence over the Thom space it can be obtained as
$$X_T=(D\tilde\xi_\eta^L\times X_{\tau'}^L)\usqcup{\tilde\rho}X_{\tau'}^L,$$
where the attaching map is some map $\tilde\rho\colon S\tilde\xi_\eta^L\times X_{\tau'}^L\to X_{\tau'}^L$ for which the restriction $\{s\}\times X_{\tau'}^L\to X_{\tau'}^L$ is a homeomorphism for all $s\in S\tilde\xi_\eta^L$. We can describe this map using the construction of the spaces $X_{\tau'}^L$ and $X_T$; note that $X_T$ is the classifying space for those $(\tau,L)$-maps where the $\eta$-points are not multiple.

Recall the universal $(\tau',L)$-map $f_{\tau'}^L\colon Y_{\tau'}^L\to X_{\tau'}^L$ in theorem \ref{ptpullback}. The pullback property of this map (i.e. that all $(\tau',L)$-maps can be induced from it by a map of the target space to $X_{\tau'}^L$) remains true for generalised $(\tau',L)$-maps obtained as direct limits of $(\tau',L)$-maps as well (see remark \ref{gentau}). Note that the operation in the classifying space $X_{\tau'}^L$ is the same as the operation in the cobordism group classified by it (by a version of the Brown representability theorem), that is, the disjoint union of (generalised) $(\tau',L)$-maps corresponds to the pointwise multiplication of their classifying maps.

Now the attaching map $\tilde\rho$ is the same as the inducing map of a generalised $(\tau',L)$-map as shown on the following diagram:
$$\xymatrix@C=11pc{
Y_{\tau'}^L\ar[r]^{f_{\tau'}^L} & X_{\tau'}^L \\
(S\tilde\xi_\eta^L\times Y_{\tau'}^L)\sqcup(S\xi_\eta^L\times X_{\tau'}^L)\ar[u]^\rho\ar[r]^(.575){\big(\id_{S\tilde\xi_\eta^L}\times f_{\tau'}^L\big)\sqcup\big(\Phi_\eta^L|_{S\xi_\eta^L}\times\id_{X_{\tau'}^L}\big)} & S\tilde\xi_\eta^L\times X_{\tau'}^L\ar[u]_{\tilde\rho}
}$$
Here the restriction $\id_{S\tilde\xi_\eta^L}\times f_{\tau'}^L\colon S\tilde\xi_\eta^L\times Y_{\tau'}^L\to S\tilde\xi_\eta^L\times X_{\tau'}^L$ is induced by the projection $\pr_{X_{\tau'}^L}\colon S\tilde\xi_\eta^L\times X_{\tau'}^L\to X_{\tau'}^L$ and the other restriction $\Phi_\eta^L|_{S\xi_\eta^L}\times\id_{X_{\tau'}^L}\colon S\xi_\eta^L\times X_{\tau'}^L\to S\tilde\xi_\eta^L\times X_{\tau'}^L$ is induced by the map $\tilde\rho_\eta^L\circ\pr_{S\tilde\xi_\eta^L}\colon S\tilde\xi_\eta^L\times X_{\tau'}^L\to X_{\tau'}^L$. Hence the disjoint union is induced by the pointwise product that maps $(s,x)\in S\tilde\xi_\eta^L\times X_{\tau'}^L$ to the point $\tilde\rho_\eta^L(s)\cdot x$.
~$\square$\par

\medskip\par\noindent\textbf{Proof of lemma \ref{l2}.\enspace\ignorespaces}
To prove that $\sigma\colon S(S\tilde\xi_\eta^L)\to BX_{\tau'}^L$ and $\tilde\rho_\eta^L\colon S\tilde\xi_\eta^L\to X_{\tau'}^L$ are adjoint maps, we will use a general statement.

Let $\mathscr{H}$ be an H-space and for any (pointed) space $A$ denote by $\mathscr{H}(A)$ the set of principal $\mathscr{H}$-bundles over $A$. There are two well-known bijective correspondences
$$\varphi\colon\mathscr{H}(SA)\to[SA,B\mathscr{H}]~~~~\text{and}~~~~\psi\colon\mathscr{H}(SA)\to[A,\mathscr{H}],$$
where $\varphi$ sends each bundle to (the homotopy class of) its inducing map and $\psi$ sends each bundle $E\xra{\mathscr{H}}SA$ to (the homotopy class of) the map $\alpha\colon A\to\mathscr{H}$ for which
$$E=(CA\times\mathscr{H})\usqcup{\beta}\mathscr{H}$$
and the attaching map is $\beta(a,h)=\alpha(a)\cdot h~(a\in A,h\in\mathscr{H})$. Such an $\alpha$ exists and is unique up to homotopy.

\medskip\begin{sclaim}
The bijection $\psi\circ\varphi^{-1}$ coincides with the adjoint correspondence.
\end{sclaim}

\begin{sprf}
The symbol $*$ will denote a fixed point in any space throughout the proof. Let $f\colon(SA,*)\to(B\mathscr{H},*)$ be a map that induces the bundle $f^*\pi$ from the universal bundle $\pi\colon E\mathscr{H}\to B\mathscr{H}$. Lift the map $f$ to a map $\tilde f$ as indicated in the diagram
$$\xymatrix@C=.75pc{
(CA,A)\ar[rr]^(.4){\tilde f}\ar[d]_{A\times\{1\}} && (PB\mathscr{H},\Omega B\mathscr{H})\ar@{=}[r]^(.58){\sims{1.5}}\ar[d]^{\Omega B\mathscr{H}} & (E\mathscr{H},\mathscr{H})\ar[d]^{\mathscr{H}}\\
(SA,*)\ar[rr]^f && (B\mathscr{H},*)\ar@{=}[r] & (B\mathscr{H},*)
}$$
where $CA=(A\times[0,1])/(A\times\{0\})$, its map to $SA$ is the quotient by $A\times\{1\}$ and the lifted map $\tilde f$ to the path space is defined by $\tilde f(a,t):=f|_{\{a\}\times[0,t]}$. It is not hard to see that $\tilde f$ maps $A=A\times\{1\}$ to $\Omega B\mathscr{H}$ by the adjoint of $f$ that we denote by $f'$.

Observe that lifting the universal bundle $\pi\colon E\mathscr{H}\xra{\mathscr{H}}B\mathscr{H}$ to its contractible total space results in the trivial bundle $E\mathscr{H}\times\mathscr{H}$; conversely, identifying each point $(e,h)\in E\mathscr{H}\times\mathscr{H}$ with the points $(e\cdot h',h\cdot h')$ for all $h'\in\mathscr{H}$ (i.e. factoring by the diagonal $\mathscr{H}$-action) recovers the universal bundle. Hence the pullback bundle $f^*\pi$ can also be obtained by pulling back the trivial bundle $E\mathscr{H}\times\mathscr{H}$ by the map $\tilde f$ to a bundle over $CA$, then identifying the points $(a,1,h)$ and $(a,1,f'(a)\cdot h)$ for all $(a,1)\in A\times\{1\}\subset CA$ and $h\in\mathscr{H}$. So we have $\varphi(f^*\pi)=f$ and $\psi(f^*\pi)=f'$.
\end{sprf}

Let $\Sigma\xra{X_{\tau'}^L}S(S\tilde\xi_\eta^L)$ denote the bundle induced by $\sigma$. We will use the claim above with the substitutions $\mathscr{H}:=X_{\tau'}^L$ and $A:=S\tilde\xi_\eta^L$ on the element $\Sigma\in\mathscr{H}(SA)$. Clearly we have $\varphi(\Sigma)=\sigma$, so we only have to see the identity $\psi(\Sigma)=\tilde\rho_\eta^L$. This follows from the observaion
$$\Sigma=(C(S\tilde\xi_\eta^L)\times X_{\tau'}^L)\usqcup{\tilde\rho_T}X_{\tau'}^L,$$
where lemma \ref{l3} implies that the gluing map is $\tilde\rho_T\colon S\tilde\xi_\eta^L\times X_{\tau'}^L\to X_{\tau'}^L$, which maps each point $(s,x)$ to the point $\tilde\rho_\eta^L(s)\cdot x$.
~$\square$\par

\medskip\par\noindent\textbf{Proof of lemma \ref{l4}.\enspace\ignorespaces}
Now we will use the notion of $l$-framed $(\tau,L)$-maps (see definition \ref{lframed}). We saw in remark \ref{lfr} that the classifying space for cobordisms of $l$-framed $(\tau,L)$-maps is of the form $X_{\tau\oplus l}^L\cong\Gamma S^lV_\tau^L$. Hence we have $BX_\tau^L=X_{\tau\oplus1}^L$ and we can set $B^lX_\tau^L:=X_{\tau\oplus l}^L$ to denote the space with the property $\Omega^lB^lX_\tau^L=X_\tau^L$.

We can now consider the key fibration for $l$-framed $(\tau,L)$-maps, as the proof of theorem \ref{keythm} works for the spaces $B^lX_\tau^L$ without change. This is a fibration
$$\chi_l\colon B^lX_\tau^L\xra{B^lX_{\tau'}^L}\Gamma S^lT\tilde\xi_\eta^L$$
and it can be induced from the universal $B^lX_{\tau'}^L$-bundle $*\xra{B^lX_{\tau'}^L}B^{l+1}X_{\tau'}^L$ by a map $b_l\colon\Gamma S^lT\tilde\xi_\eta^L\to B^{l+1}X_{\tau'}^L$. The following diagram shows this inducing map on the rightmost square and the resolvents of both horizontal arrows to the left:
$$\xymatrix{
\ldots\ar[r] & \Omega B^lX_\tau^L\ar[r]^{\chi_{l-1}}\ar[d] & \Omega\Gamma S^lT\tilde\xi_\eta^L\ar[r]\ar[d]^{b_{l-1}} & B^lX_{\tau'}^L\ar[r]\ar[d]^{\id_{B^lX_{\tau'}^L}} & B^lX_\tau^L\ar[r]^{\chi_l}\ar[d] & \Gamma S^lT\tilde\xi_\eta^L\ar[d]^{b_l}\\
\ldots\ar[r] & \ast\ar[r] & \Omega B^{l+1}X_{\tau'}^L\ar[r] & B^lX_{\tau'}^L\ar[r] & \ast\ar[r] & B^{l+1}X_{\tau'}^L
}$$

The leftmost square in the diagram above is precisely the pullback diagram of the fibration $\chi_{l-1}$ with the map $b_{l-1}$, hence we obtained that $\chi_{l-1}=\Omega\chi_l$ and $b_{l-1}=\Omega b_l$. We can iterate this process as $l\to\infty$ and we obtain that $b_0=b_\eta^L$ is indeed an infinite loop map.
~$\square$\par

\medskip\par\noindent\textbf{Proof of lemma \ref{l5}.\enspace\ignorespaces}
We want to show that the diagram below remains commutative (up to homotopy) after adding the dashed arrow.
$$\xymatrix{
A\ar[rr]^f\ar[dr]^g\ar@{_(->}[dd]_i && A'\ar[dl]_{g'}\ar@{^(->}[dd]^{i'}\\
& B &\\
\Gamma A\ar@{-->}[rr]^{\Gamma f}\ar[ur]^{g_\ext} && \Gamma A'\ar[ul]_{g'_\ext}
}$$
The square commutes by the construction of the map $\Gamma f$, hence the composition $g'_\ext\circ\Gamma f\circ i$ coincides with the map $g=g'\circ f=g'_\ext\circ i'\circ f$. Thus $g'_\ext\circ\Gamma f$ is an extension of $g$ to the space $\Gamma A$ and it is also an infinite loop map, since it is the composition of such maps. Now the uniqueness of the extension implies that $g_\ext$ and $g'_\ext\circ\Gamma f$ coincide up to homotopy.
~$\square$\par

%-------------------------------------------------------------------------------------------------------------------------------

\section{Corollaries}

Here we will see some applications of the key fibration from \cite{hosszu} and \cite{hominv} that give nice descriptions of the classifying spaces $X_\tau^L$ in some cases. Again $\tau$ will be a set of singularities, $L$ a stable linear group, $[\eta]\in\tau$ a top singularity and $\tau':=\tau\setminus\{[\eta]\}$.

\subsection{Rational decomposition}

In this subsection we show rational decomposition theorems on the key bundle. First we will see that given some restrictions on the allowed singularity classes in $\tau$, the key fibration is rationally trivial, i.e. $X_\tau^L\congq X_{\tau'}^L\times\Gamma T\tilde\xi_\eta^L$.

\begin{rmk}
The rational triviality of the key bundle means that for any manifold $P$ we have
$$\Cob_\tau^L(P)\otimes\Q\cong(\Cob_{\tau'}^L(P)\oplus\Imm^{\tilde\xi_\eta^L}(P))\otimes\Q.$$
\end{rmk}

\begin{thm}\label{ratrivi1}
If $\nu_\tau^L,\xi_\eta^L$ and $\tilde\xi_\eta^L$ are orientable and either $\xi_\eta^L$ or $\tilde\xi_\eta^L$ has non-trivial rational Euler class $e(\xi_\eta^L)$ or $e(\tilde\xi_\eta^L)$ in $H^*(BG_\eta^L;\Q)$, then we have 
$$X_\tau^L\congq X_{\tau'}^L\times\Gamma T\tilde\xi_\eta^L.$$
\end{thm}

\begin{prf}
Consider the Kazarian space $K_\tau^L$ and the virtual complex $V_\tau^L$ and take the cohomological exact sequences obtained from the Puppe sequences
$$K_{\tau'}^L\into K_\tau^L\to T\xi_\eta^L\into SK_{\tau'}^L~~~~\text{and}~~~~V_{\tau'}^L\into V_\tau^L\to T\tilde\xi_\eta^L\into SV_{\tau'}^L.$$
These are dual to the homological long exact sequences of the pairs $(K_\tau^L,K_{\tau'}^L)$ and $(V_\tau^L,V_{\tau'}^L)$. These dual sequences are shown on the upper two rows of the diagram below and the isomorphisms between these rows are the Thom isomorphisms corresponding to the virtual bundle $\nu_\tau^L$.
$$\xymatrix{
H_*(K_{\tau'}^L;\Q)\ar[r]\ar[d]^\cong & H_*(K_{\tau}^L;\Q)\ar[r]\ar[d]^\cong & H_*(T\xi_\eta^L;\Q)\ar[r]^{\partial_K}\ar[d]^\cong & H_*(SK_{\tau'}^L;\Q)\ar[d]^\cong \\
H_{*+k}(V_{\tau'}^L;\Q)\ar[r]\ar[d]^\cong & H_{*+k}(V_{\tau}^L;\Q)\ar[r]\ar[d]^\cong & H_{*+k}(T\tilde\xi_\eta^L;\Q)\ar[r]^{\partial_V}\ar[d]^\cong & H_{*+k}(SV_{\tau'}^L;\Q)\ar[d]^\cong \\
\pi^s_{*+k}(V_{\tau'}^L)\otimes\Q\ar[r]\ar[d]^\cong & \pi^s_{*+k}(V_{\tau}^L)\otimes\Q\ar[r]\ar[d]^\cong & \pi^s_{*+k}(T\tilde\xi_\eta^L)\otimes\Q\ar[r]\ar[d]^\cong & \pi^s_{*+k}(SV_{\tau'}^L)\otimes\Q\ar[d]^\cong \\
\pi_{*+k}(X_{\tau'}^L)\otimes\Q\ar[r] & \pi_{*+k}(X_{\tau}^L)\otimes\Q\ar[r] & \pi_{*+k}(\Gamma T\tilde\xi_\eta^L)\otimes\Q\ar[r]^\partial & \pi_{*+k-1}(X_{\tau'}^L)\otimes\Q
}$$
We obtain the third row by applying the stable Hurewicz homomorphism on the second row (recall that these are also isomorphisms according to Serre's theorem \cite[theorem 18.3]{charclass}) and we get the fourth row by corollary \ref{gammav}.

Now suppose that the vector bundle $\xi_\eta^L$ is orientable and the Euler class $e(\xi_\eta^L)\in H^*(BG_\eta^L;\Q)$ is non-trivial. The dual of the map $\partial_K$ is the homomorphism
$$i^*\colon H^*(SK_{\tau'}^L;\Q)\to H^*(T\xi_\eta^L;\Q)$$
induced by the inclusion $i\colon T\xi_\eta^L\into SK_{\tau'}^L$. We claim that this homomorphism is trivial. This is indeed so, because the cup product in $H^*(SK_{\tau'}^L;\Q)$ is trivial, but in $H^*(T\xi_\eta^L;\Q)$ any element is of the form $u(\xi_\eta^L)\smallsmile x$ for some $x\in H^*(BG_\eta^L;\Q)$ (by the Thom isomophism), hence the cup product here is $(u(\xi_\eta^L)\smallsmile x_1)\smallsmile(u(\xi_\eta^L)\smallsmile x_2)=u(\xi_\eta^L)\smallsmile e(\xi_\eta^L)\smallsmile x_1\smallsmile x_2$. This cannot vanish unless $x_1$ or $x_2$ vanishes, since $H^*(BG_\eta^L;\Q)$ is a subring of the polynomial ring $H^*(BT_\eta^L;\Q)$ where $T_\eta^L$ is the maximal torus in $G_\eta^L$, so it has no zero divisors.

If we suppose that $\tilde\xi_\eta^L$ is an orientable bundle with non-vanishing rational Euler class, then the same reasoning yields that $\partial_V$ is trivial. Either way, the homomorphism $\partial$ in the bottom row of the diagram above vanishes, hence
$$\pi_*(X_\tau^L)\otimes\Q\cong(\pi_*(X_{\tau'}^L)\otimes\Q)\oplus(\pi_*(\Gamma T\tilde\xi_\eta^L)\otimes\Q).$$
It remains to note that $X_\tau^L$, $X_{\tau'}^L$ and $\Gamma T\tilde\xi_\eta^L$ are all H-spaces, so they are rationally homotopy equivalent to products of Eilenberg--MacLane spaces $K(\Q,i)$ (for some $i\in\N$). Thus we obtain the statement of the theorem.
\end{prf}

\begin{thm}\label{ratrivi2}
Suppose that $\tilde\xi_\eta^L=\zeta\oplus\varepsilon^r$ for an orientable vector bundle $\zeta$ of odd dimension and a number $r\in\N$ and suppose also that $\pi_{i}(X_{\tau'}^L)$ is a torsion group whenever $i+r$ is an even number. Then we have 
$$X_\tau^L\congq X_{\tau'}^L\times\Gamma T\tilde\xi_\eta^L.$$
\end{thm}

\begin{prf}
It is enough to show that the inducing map $b_\eta^L\colon\Gamma T\tilde\xi_\eta^L\to BX_{\tau'}^L$ of the key bundle is rationally null-homotopic.

Consider its $r$-th loop map
$$\Omega^rb_\eta^L\colon\Omega^r\Gamma T\tilde\xi_\eta^L=\Gamma T\zeta\to\Omega^{r-1}X_{\tau'}^L=\Omega^rBX_{\tau'}^L$$
using the assumption $T\tilde\xi_\eta^L=T(\zeta\oplus\varepsilon^r)=S^rT\zeta$. This is an infinite loop map by lemma \ref{l4}, hence it is the unique extension of its restriction to $T\zeta$. Recall that $X_{\tau'}^L$ is an infinite loop space, so this also makes sense for $r=0$.

The base space $BG_\eta^L$ has non-trivial rational cohomology only in even degrees (since $G_\eta^L$ is a compact Lie group), hence the Thom isomorphism yields that $T\zeta$ has non-trivial rational cohomology only in odd degrees (because the rank of $\zeta$ was assumed to be odd). The H-space $\Omega^{r-1}X_{\tau'}^L$ is rationally homotopy equivalent to a product of Eilenberg--MacLane spaces $K(\Q,i)$ and none of these numbers $i$ can be odd, as we assumed $\pi_i(\Omega^{r-1}X_{\tau'}^L)\otimes\Q\cong\pi_{i+r-1}(X_{\tau'}^L)\otimes\Q$ to be trivial for odd numbers $i$. This implies that any map $T\zeta\to\Omega^{r-1}X_{\tau'}^L$ induces the trivial homomorphism in rational cohomologies, hence is rationally null-homotopic. Therefore the infinite loop map extension $\Omega^rb_\eta^L$ is also rationally null-homotopic.

The resolvent of the inducing map $b_\eta^L$ has the form
\begin{alignat*}2
\ldots\to\Omega^rX_{\tau'}^L\to\Omega^rX_\tau^L\xra{\Omega^r\chi_\tau^L}\Gamma T\zeta\xra{\Omega^rb_\eta^L}\Omega^{r-1}X_{\tau'}^L\to\ldots&\\
\ldots\to X_{\tau'}^L\to X_\tau^L\xra{\chi_\tau^L}&\Gamma T\tilde\xi_\eta^L\xra{b_\eta^L}BX_{\tau'}^L.
\end{alignat*}
Recall that in a resolvent each map is the inducing map of the bundle whose projection is the previous map, in particular $\Omega^rb_\eta^L$ induces the bundle $\Omega^r\chi_\tau^L\colon\Omega^rX_\tau^L\xra{\Omega^rX_{\tau'}^L}\Omega^r\Gamma T\tilde\xi_\eta^L$ from the universal bundle $*\xra{\Omega^rX_{\tau'}^L}\Omega^{r-1}X_{\tau'}^L$. Now the rational null-homotopy of $\Omega^rb_\eta^L$ implies that this is a rationally trivial bundle, i.e. we have
$$\Omega^rX_\tau^L\congq\Omega^rX_{\tau'}^L\times\Omega^r\Gamma T\tilde\xi_\eta^L.$$
Since the spaces $X_\tau^L$, $X_{\tau'}^L$ and $\Gamma T\tilde\xi_\eta^L$ are H-spaces, their rational homotopy types are products of Eilenberg--MacLane spaces and then $\Omega^rK(\Q,i)=K(\Q,i-r)$ implies that the above splitting of $\Omega^rX_\tau^L$ also holds for $X_\tau^L$, that is, $X_\tau^L\congq X_{\tau'}^L\times\Gamma T\tilde\xi_\eta^L$.
\end{prf}

In general, when the conditions in the above theorems do not apply, there is also a weaker decomposition theorem for the rational homotopy type of $X_\tau^L$. We will denote the rational homotopy types of the spaces $X_\tau^L$, $X_{\tau'}^L$ and $\Gamma T\tilde\xi_\eta^L$ respectively by $X$, $X'$ and $\Gamma$. These are all rational H-spaces, hence they uniquely decompose into products of Eilenberg--MacLane spaces $K(\Q,i)$ (for some $i\in\N$).

\begin{thm}
There is a rational H-space $B$ with the property
$$X\congq\frac{X'\times\Omega B}{\Omega\Gamma}\times B,$$
where dividing by $\Omega\Gamma$ means that all Eilenberg--MacLane factors $K(\Q,i)$ of $\Omega\Gamma$ occur as factors of the numerator and we cancel these factors.
\end{thm}

\begin{prf}
We use the $\Q$-classification of H-bundles \cite[lemma 113]{hosszuv1} on the key bundle to obtain that $X\xra{X'}\Gamma$ is the product of three H-bundles of the forms $*\xra{\Omega A}A$, $B\xra{*}B$ and $C\xra{C}*$. Hence $X=B\times C$, $X'=A\times C$ and $\Gamma=A\times B$ and so $\Omega\Gamma=\Omega A\times\Omega B$, and the statement of the theorem follows.
\end{prf}

\begin{rmk}
If the singularity set is finite, say we have $\tau=\{\eta_0,\ldots,\eta_r\}$ with $\eta_0<\ldots<\eta_r=\eta$, then the space $B$ in the above theorem can be obtained from a spectral sequence as follows.

By the above proof we have $\pi_*(X)=\pi_*(B)\oplus\pi_*(C)$ and $\pi_*(\Gamma)=\pi_*(A)\oplus\pi_*(B)$, the induced map $\chi_\#\colon\pi_*(X)\to\pi_*(\Gamma)$ of the rational key fibration is an isomorphism on $\pi_*(B)$ and this is the maximal such subgroup. Now we have a commutative diagram
$$\xymatrix{
\pi_m(X')\ar[d]^\cong\ar[r] & \pi_m(X)\ar[d]^\cong\ar[r]^{\chi_\#} & \pi_m(\Gamma)\ar[d]^\cong \\
H_{m-k}(K_{\tau'}^L;\tilde\Q)\ar[r] & H_{m-k}(K_\tau^L;\tilde\Q)\ar[r]^\alpha & H_{m-k}(T\xi_\eta^L;\tilde\Q)
}$$
where the vertical arrows are the compositions of the stable Hurewicz homomorphisms (which are rational isomorphisms) and the Thom isomorphisms corresponding to the virtual bundle $\nu_\tau^L$ and $\tilde\Q$ denotes rational coefficients twisted by the orientation of (the restriction of) $\nu_\tau^L$. Now if $E^*_{*,*}$ denotes the Kazarian spectral sequence (see definition \ref{kazspec}) with twisted rational coefficients, then we have
$$H_{m-k}(K_\tau^L;\tilde\Q)=\bigoplus_{i=0}^rE^\infty_{i,m-k-i}~~~~\text{and}~~~~E^\infty_{r,m-k-r}=\im\alpha\subset H_{m-k}(T\xi_\eta^L;\tilde\Q).$$
Hence the rational H-space (i.e. product of Eilenberg--MacLane spaces) $B$ is completely determined by the formula
$$\pi_m(B)=\im\chi_\#\cong\im\alpha=E^\infty_{r,m-k-r}.$$
\end{rmk}

\subsection{A Postnikov-type tower}

In this subsection we assume that $\tau$ is finite and its elements are $\eta_0<\ldots<\eta_r$ with a complete order extending the natural partial order. We will use the simplified notations $G_i:=G_{\eta_i}^L$, $\xi_i:=\xi_{\eta_i}^L$, $c_i:=\dim\xi_i$, $\tilde\xi_i:=\tilde\xi_{\eta_i}^L$, $\Gamma_i:=\Gamma T\tilde\xi_i$, $\tau_i:=\{\eta_j\mid j\le i\}$, $X_i:=X_{\tau_i}^L$, $K_i:=K_{\tau_i}^L$ and $V_i:=V_{\tau_i}^L$.

\begin{prop}\label{postprop}
The classifying space $X_\tau^L=X_r$ can be obtained by a sequence of fibrations
$$\Gamma_r=A_1\xla{\Gamma_{r-1}}A_2\xla{\Gamma_{r-2}}\ldots\xla{\Gamma_1}A_r\xla{\Gamma_0}X_\tau^L.$$
\end{prop}

\begin{prf}
We have $T\tilde\xi_i=V_i/V_{i-1}$ for $i=1,\ldots,r$, so there is a cofibration
$$T\tilde\xi_{r-i+1}\into V_r/V_{r-i}\to V_r/V_{r-i+1}.$$
Now putting $A_i:=\Gamma(V_r/V_{r-i})$ for $i=1,\ldots,r$ and applying that $\Gamma$ turns cofibrations into fibrations (see \cite{barec} and lemma \ref{cofibfib}), we get the proposed sequence of fibrations until $A_r$. Similarly the cofibration $V_0=T\tilde\xi_0\into V_r\to V_r/V_0$ is turned by $\Gamma$ to the last fibration in the sequence.
\end{prf}

\begin{crly}
For any manifold $P^{n+k}$ there is a spectral sequence with first page $\Imm^{\tilde\xi_i}(P\times\R^j)$, which converges to $\Cob_\tau^L(P\times\R^*)$, that is,
$$E_1^{i,j}:=\Imm^{\tilde\xi_i}(P\times\R^j)\implies\Cob_\tau^L(P\times\R^{i+j}).$$
\end{crly}

\begin{prf}
Applying the extraordinary homology theory $h_j:=\{S^j\cpt P,\cdot\}$ to the sequence of cofibrations $V_r\to V_r/V_0\to V_r/V_1\to\ldots\to V_r/V_{r-1}$, we get an exact couple
$$\xymatrix@C=.75pc{
\displaystyle\bigoplus_{i,j}\{S^j\cpt P,V_r/V_{r-i}\}\ar[rr] && \displaystyle\bigoplus_{i,j}\{S^j\cpt P,V_r/V_{r-i}\}\ar[dl] \\
& \displaystyle\bigoplus_{i,j}\{S^j\cpt P,T\tilde\xi_i\}\ar[ul] &
}$$
(with natural indexing by $i,j$ and $V_r/V_{-1}$ is considered to be just $V_r$). Observe that we have $\{S^j\cpt P,V_r/V_{r-i}\}=[S^j\cpt P,A_i]$ and $\{S^j\cpt P,T\tilde\xi_i\}=[S^j\cpt P,\Gamma_i]$, hence the spectral sequence we obtain from this exact couple has first page
$$[S^j\cpt P,\Gamma_i]=[S^j\cpt P,\Gamma T\tilde\xi_i]\cong\Imm^{\tilde\xi_i}(P\times\R^j)$$
and converges to
$$[S^{i+j}\cpt P,\Gamma V_r]=[S^{i+j}\cpt P,\Gamma V_\tau^L]=[S^{i+j}\cpt P,X_\tau^L]\cong\Cob_\tau^L(P\times\R^{i+j}),$$
and this is what we wanted to prove.
\end{prf}

Denote by $0=s_0<s_1<\ldots<s_l<s_{l+1}=r$ the indices where the parity of the numbers $c_i$ changes (for $i=0,\ldots,r$), that is,
$$c_0\equiv\ldots\equiv c_{s_1}\not\equiv c_{s_1+1}\equiv\ldots\equiv c_{s_2}\not\equiv c_{s_2+1}\equiv\ldots~~~~\mod 2.$$
We will say that the indices $i,i'$ are in the same block, if there is a number $0\le t\le l$ for which $s_t<i,i'\le s_{t+1}$.

In the following we consider the (homological) Kazarian spectral sequence $E_{i,j}^t$ (see definition \ref{kazspec}) with rational coefficients and assume that the bundles $\xi_i$ are all orientable.

\begin{lemma}
$~$
\begin{enumerate}
\item If some indices $i,i'$ are in the same block, then the differential $d_{i,j}^t\colon E_{i,j}^t\to E_{i',j'}^t$ (if exists) is trivial.
\item The quotient spaces $K_{s_{t+1}}/K_{s_t}$ have the same rational homologies as the wedge products $T\xi_{s_t+1}\vee\ldots\vee T\xi_{s_{t+1}}$ for $t=0,\ldots,l$.
\end{enumerate}
\end{lemma}

\begin{prf}
The groups $G_i$ are compact Lie groups, hence the spaces $BG_i$ have non-trivial rational cohomologies only in even degrees. Thus (by the homological Thom isomorphism) the groups $E^1_{i,j}=H_{i+j}(T\xi_i;\Q)\cong H_{i+j-c_i}(BG_i;\Q)$ can be non-trivial only if $i+j\equiv c_i~\mod2$. Now the differential $d^t_{i,j}\colon E^t_{i,j}\to E^t_{i-t,j+t-1}$ can be non-trivial only if both groups are non-trivial, which cannot happen if $i$ and $i-t$ are in the same block. This finishes the proof of both claims.
\end{prf}

This lemma will be used in the proof of the following theorem that gives a simplification of proposition \ref{postprop} rationally. We will use $X$ and $\hat\Gamma_t$ (for $t=0,\ldots,l$) respectively to denote the rational homotopy types of $X_\tau^L$ and $\underset{s_t<i\le s_{t+1}}\prod\Gamma_i$.

\begin{thm}
There is a simplified sequence of fibrations for $X$ of the form
$$\hat\Gamma_l=\hat A_1\xla{\hat\Gamma_{l-1}}\hat A_2\xla{\hat\Gamma_{l-2}}\ldots\xla{\hat\Gamma_1}\hat A_l\xla{\hat\Gamma_0}X.$$
\end{thm}

\begin{prf}
The above lemma implies that the quotient spaces $V_{s_{t+1}}/V_{s_t}$ are stably rationally homotopy equivalent to $T\tilde\xi_{s_t+1}\vee\ldots\vee T\tilde\xi_{s_{t+1}}$ for $t=0,\ldots,l$. Now applying the $\Gamma$ functor to the rational homotopy type of a cofibration of the form
$$V_{s_{l-t+1}}/V_{s_{l-t}}\into V_{s_{l+1}}/V_{s_{l-t}}\to V_{s_{l+1}}/V_{s_{l-t+1}},$$
we obtain a fibration with fibre $\hat\Gamma_{l-t}$. Thus denoting by $\hat A_t$ the rational homotopy type of $\Gamma(V_{s_{l+1}}/V_{s_{l-t+1}})$, we obtain the proposed sequence of fibrations similarly to the proof of proposition \ref{postprop}.
\end{prf}

\subsection{Classifying spaces for prim maps}\label{clasprim}

In the beginning of chapter \ref{classp} we noted that the classifying spaces $\overline X_\tau^L$ of the cobordism groups of prim maps $\Prim_\tau^L(n,P^{n+k})$ (see definition \ref{primcob}) can be constructed analogously to the construction of the spaces $X_\tau^L$. However, the key fibration helps in finding a simpler way we can obtain these spaces, as we will soon see (based on \cite{ctrl}). Note, that the key fibration also exists for prim maps, i.e. there is a fibration
$$\ol\chi_\tau^L\colon\ol X_\tau^L\xra{\ol X_{\tau'}^L}\Gamma T\tilde\zeta_\eta^L$$
using the usual notations and denoting by $\tilde\zeta_\eta^L$ the target bundle in the global normal form of prim maps around the $\eta$-stratum.

Recall that prim maps can only have Morin singularities (see definition \ref{thomboar}), hence the singularity set $\tau$ can now only be $\tau_r:=\{\Sigma^0,\Sigma^{1_1},\ldots,\Sigma^{1_r}\}$ for some $r\in\N\cup\{\infty\}$, where $r=\infty$ only means the union of all $\tau_r$-s for $r\in\N$ (see example \ref{morsing}). We will use the notation analogous to \ref{conv2} after definition \ref{cobtau} by putting $\Prim_r^L(n,P):=\Prim_{\tau_r}^L(n,P)$ and $\overline X_r^L:=\overline X_{\tau_r}^L$.

\begin{prop}\label{priminf}
$\overline X_\infty^L\cong\Omega\Gamma T\gamma_{k+1}^L$.
\end{prop}

\begin{prf}
We use a theorem of Wells \cite{wells} which states that the classifying space for cobordisms of $(k+1)$-codimensional immersions with normal $L$-structures is $X_{\{\Sigma^0\}}^L=\Gamma T\gamma_{k+1}^L$. For any manifold $P^{n+k}$, there is a natural bijective correspondence between the cobordisms of $k$-codimensional prim maps to $P$ and the cobordisms of $(k+1)$-codimensional immersions to $P\times\R^1$ (both with normal $L$-structures). Hence we have
$$[\cpt P,\overline X_\infty^L]\cong[S\cpt P,\Gamma T\gamma_{k+1}^L]\cong[\cpt P,\Omega\Gamma T\gamma_{k+1}^L],$$
and now the homotopic uniqueness of the classifying space implies the statement.
\end{prf}

If we restrict the set of allowed singularities for the cobordisms of prim maps, then we get a bit more complicated but similar description of the classifying space. In the rest of this section we will have a fixed number $r\in\N$.

\begin{defi}\label{zetasdef}
Denote by $\pi$ the projection of the sphere bundle of $(r+1)\gamma_{k+1}^L$ (the sum $\gamma_{k+1}^L\oplus\ldots\oplus\gamma_{k+1}^L$ of $(r+1)$ summands). Let $\zeta_S^r$ be the pullback of the universal bundle $\gamma_{k+1}^L$ by $\pi$, that is, we put
$$\xymatrix@C=2pc{
*!<.4cm,0cm>{\zeta_S^r:=\pi^*\gamma_{k+1}^L}\ar[r]\ar[d] & \gamma_{k+1}^L\ar[d] \\
S((r+1)\gamma_{k+1}^L)\ar[r]^(.55){\pi} & BL(k+1)
}$$
\end{defi}

\begin{thm}\label{zetasthm}
$\overline X_r^L\cong\Omega\Gamma T\zeta_S^r$.
\end{thm}

\begin{prf}
We fix a (generic) immersion $g\colon M^n\imto P\times\R^1$ such that the composition $f:=\pr_P\circ g\colon M\to P$ is a $(\tau_r,L)$-map. First we will produce sections $s_1,s_2,\ldots$ of the normal bundle $\nu_g$ such that for all $i$, the points where precisely the first $i$ sections vanish form the $\Sigma^{1_i}$-stratum, i.e. $\underset{j=1}{\overset i\bigcap}s_j^{-1}(0)\setminus s_{i+1}^{-1}(0)=\Sigma^{1_i}(f)$.

The positive unit vector in $\R^1$ defines a constant vector field on $P\times\R^1$; this will be called vertically upwards and denoted by $\ua$. The pullback $g^*\ua$ is a section of $g^*T(P\times\R^1)$ over $M$ and it can be projected to the bundle $\nu_g=g^*T(P\times\R^1)/TM$. We define this projected section to be $s_1$ and note that $s_1^{-1}(0)$ is obviously the singular set $\Sigma(f)$.

The set $\Sigma(f)$ is a submanifold of $M$ of codimension $k+1$; denote its normal bundle by $\nu_2$. Observe that the tangent spaces of the section $s_1$ at the points of $\Sigma(f)$ form the graph of a pointwise linear isomorphism $\iota_2$ between $\nu_2$ and $\nu_g|_{\Sigma(f)}$. We define a section $s'_2$ of $\nu_2$ by projecting $g^*\ua$ to $\nu_2\subset TM$. Now $\iota_2\circ s'_2$ is a section of $\nu_g|_{\Sigma(f)}$ and we can arbitrarily extend it to a section of the whole normal bundle $\nu_g$ and denote this extension by $s_2$. Clearly $\Sigma^{1,1}(f)$ (see definition \ref{thomboar}) coincides with ${s'_2}^{-1}(0)$, hence it is also the same as $s_1^{-1}(0)\cap s_2^{-1}(0)$.

By further applying the same method, we obtain the sections $s_3,s_4,\ldots$ with the desired properties; in particular we have $\underset{i=1}{\overset{r+1}\bigcap}s_i^{-1}(0)=\varnothing$. These sections are not unique, but the difference of any two possible choices of $s_i$ is an arbitrary section of $\nu_g$ that vanishes on $\underset{j=1}{\overset{i-1}\bigcap}s_j^{-1}(0)$, thus the uniqueness holds up to a contractible choice. Hence the collection of the sections $(s_1,\ldots,s_{r+1})$ defines a homotopically unique section $\alpha$ of the sphere bundle $S((r+1)\nu_g)$ over $M$.

The normal bundle $\nu_g$ has structure group $L(k+1)$, hence it can be induced from the universal bundle $\gamma_{k+1}^L$ by a map $b\colon M\to BL(k+1)$ together with a fibrewise isomorphism between $\nu_g$ and $\gamma_{k+1}^L$. This also gives a map $\beta\colon S((r+1)\nu_g)\to S((r+1)\gamma_{k+1}^L)$. The following commutative diagram shows the maps $b$ and $\beta$ and the section $\alpha$ of $S((r+1)\nu_g)$ together with the bundles $\nu_g$ and $\gamma_{k+1}^L$ and their pullbacks by the projections to bundles over $S((r+1)\nu_g)$ and $S((r+1)\gamma_{k+1}^L)$ respectively (observe that this pullback over $S((r+1)\nu_g)$ coincides with $\beta^*\zeta_S^r$).
$$\xymatrix@R=.333pc{
\beta^*\zeta_S^r\ar[dr]\ar[rrr]\ar[ddddd] &&& \zeta_S^r\ar[ddddd]\ar[dl]\\
&S((r+1)\nu_g)\ar[r]^\beta\ar[ddd] & S((r+1)\gamma_{k+1}^L)\ar[ddd] & \\ \\ \\
& M\ar@/^1pc/[uuu]^\alpha\ar[r]^(.42)b & BL(k+1) &\\
\nu_g\ar[rrr]\ar@/^1pc/[uuuuu]\ar[ur] &&& \gamma_{k+1}^L\ar[ul]
}$$

This implies the identity $\nu_g=(\beta\circ\alpha)^*\zeta_S^r$. Because of the homotopic uniqueness of $\alpha$ and $\beta$, we obtain that $\nu_g$ can be induced from $\zeta_S^r$ in a homotopically unique way.

Applying the usual Pontryagin--Thom construction to the diagram above, we obtain a map
$$[\cpt P,\ol X_r^L]\cong\Prim_r^L(P)\to[S\cpt P,\Gamma T\zeta_S^r]\cong[\cpt P,\Omega\Gamma T\zeta_S^r]$$
in the following way: To any cobordism class $[f]$ of prim $(\tau_r,L)$-maps we can assign the cobordism class $[g]\in\Imm^{\zeta_S^r}(P\times\R^1)$ (using the notation above); now \cite{wells} gives a bijective correspondence between $\Imm^{\zeta_S^r}(P\times\R^1)$ and $[S\cpt P,\Gamma T\zeta_S^r]$, hence $[g]$ corresponds to the homotopy class of a map $S\cpt P\to\Gamma T\zeta_S^r$. This can be identified with an element of $[\cpt P,\Omega\Gamma T\zeta_S^r]$ and this is defined to be the image of $[f]$.

The map we have just described is actually a natural transformation of functors $\Prim_r^L\to[\cpt\cdot,\Omega\Gamma T\zeta_S^r]$, hence it is induced by a homotopically unique map between the classifying spaces, that is, by a map
$$\varphi_r\colon\ol X_r^L\to\Omega\Gamma T\zeta_S^r.$$

Recall remark \ref{keyrmk}, which has a direct analogue for prim maps, namely that the key fibration $\ol\chi_r^L\colon\ol X_r^L\xra{\ol X_{r-1}^L}\Gamma T\tilde\zeta_{\Sigma^{1_r}}^L$ induces in the cobordism groups these spaces classify the forgetful map that assigns to a cobordism class of a prim $(\tau_r,L)$-map the cobordism class of its $\Sigma^{1_r}$-stratum. The global normal forms of prim maps are quite simple (see \cite{cobmor} or \cite{hosszu}), in particular the bundle $\tilde\zeta_{\Sigma^{1_r}}^L$ coincides with $(r+1)\gamma_k^L\oplus\varepsilon^r$. Thus the key fibration has the form
\begin{alignat*}2
\ol\chi_r^L\colon\ol X_r^L\xra{\ol X_{r-1}^L}\Gamma T((r+1)\gamma_k^L\oplus\varepsilon^r)&=\Omega\Gamma ST((r+1)\gamma_k^L\oplus\varepsilon^r)=\\
&=\Omega\Gamma T((r+1)(\gamma_k^L\oplus\varepsilon^1)).
\end{alignat*}

Our aim is to prove that $\ol X_r^L$ coincides with $\Omega\Gamma T\zeta_S^r$, which would imply that the above key fibration holds for the space $\Omega\Gamma T\zeta_S^r$. We will now work from the opposite direction and first prove the following.

\medskip\begin{sclaim}
There is a fibration $\psi_r\colon\Omega\Gamma T\zeta_S^r\xra{\Omega\Gamma T\zeta_S^{r-1}}\Omega\Gamma T((r+1)(\gamma_k^L\oplus\varepsilon^1))$.
\end{sclaim}

\begin{sprf}
The following is a trivial general observation: If $A$ is a manifold, $B\subset A$ is a submanifold with a tubular neighbourhood $U$ and normal bundle $\nu$ and $\zeta$ is a vector bundle over $A$, then there is a cofibration of Thom spaces
$$T\zeta|_{A\setminus U}\into T\zeta\to T(\zeta|_B\oplus\nu).$$

We apply this with the substitutions $A:=S((r+1)\gamma_{k+1}^L)$, $B:=S\gamma_{k+1}^L$, $U:=S((r+1)\gamma_{k+1}^L)\setminus S(r\gamma_{k+1}^L)$, $\nu:=r(\gamma_k^L\oplus\varepsilon^1)$ and $\zeta:=\zeta_S^r$. This yields a cofibration
$$T\zeta_S^{r-1}\into T\zeta_S^r\to T((r+1)(\gamma_k^L\oplus\varepsilon^1))$$
and by applying the functor $\Omega\Gamma$ (which turns cofibrations into fibrations) we get the desired fibration.
\end{sprf}

Now the maps $\varphi_{r-1}$ and $\varphi_r$ fit into a diagram as shown below. If we prove that this diagram commutes, then an induction on the number $r$ (together with the 5-lemma) implies that $\varphi_r$ is a homotopy equivalence $\ol X_r^L\cong\Omega\Gamma T\zeta_S^r$ (the starting step is $r=0$, where the diagram degenerates to a map between trivial fibrations).
$$\xymatrix{
\ol X_{r-1}^L\ar[d]^{\varphi_{r-1}}\ar[r] & \ol X_r^L\ar[d]^{\varphi_r}\ar[r]^(.29){\ol\chi_r^L} & \Omega\Gamma T((r+1)(\gamma_k^L\oplus\varepsilon^1))\ar[d]^\cong\\
\Omega\Gamma T\zeta_S^{r-1}\ar[r] & \Omega\Gamma T\zeta_S^r\ar[r]^(.33){\psi_r} & \Omega\Gamma T((r+1)(\gamma_k^L\oplus\varepsilon^1))
}$$

The left-hand square trivially commutes. We prove the commutativity of the right-hand square by showing that the diagram of natural transformations of functors induced by these maps commutes.

The space $\Omega\Gamma T((r+1)(\gamma_k^L\oplus\varepsilon^1))$ classifies cobordisms of immersions to $P\times\R^1$ with normal bundles induced from $(r+1)(\gamma_k^L\oplus\varepsilon^1)$ (for any manifold $P$). When one applies (the transformation induced by) the map $\ol\chi_r^L$ to some $[f]\in\Prim_r^L(P)$, the normal structure induced from $(r+1)(\gamma_k^L\oplus\varepsilon^1)$ is obtained by the splitting of the normal bundle of $\Sigma^{1_r}(f)$ to $r+1$ bundles canonically isomorphic to $\nu_g|_{\Sigma^{1_r}(f)}$ (with the same notation as in the beginning of the proof). On the other hand, applying $\psi_r\circ\varphi_r$ gives the normal structure using the sections $s_1,\ldots,s_{r+1}$ as described earlier. Now these two normal structures may be twisted with respect to each other, but they can be obtained uniquely up to homotopy from each other. Hence the above diagram indeed commutes.
\end{prf}

The following is analogous to the key fibration and also very useful for the investigation of the classifying spaces $\ol X_r^L$.

\begin{thm}\label{prclthm}
$~$
\begin{enumerate}
\item\label{egy} There is a fibration $\overline X_r^L\xra{\Omega^2\Gamma T((r+2)\gamma_{k+1}^L)}\Omega\Gamma T\gamma_{k+1}^L$.
\item\label{ketto} The homotopy exact sequence of this fibration can be identified with that of the pair $(\overline X_\infty^L,\overline X_r^L)$, where we use $\overline X_\infty^L=\liminfty r\overline X_r^L$.
\end{enumerate}
\end{thm}

\begin{prf}
\emph{Proof of \ref{egy}.}\enspace\ignorespaces
We pull back the vector bundle $\gamma_{k+1}^L$ to the disk bundle $D((r+1)\gamma_{k+1}^L)$ by its projection. This pullback is such that
\begin{itemize}
\item[(i)] its restriction to $S((r+1)\gamma_{k+1}^L)$ is $\zeta_S^r$,
\item[(ii)] its Thom space is homotopy equivalent to $T\gamma_{k+1}^L$,
\item[(iii)] its disk bundle is $D((r+2)\gamma_{k+1}^L)$.
\end{itemize}
Now we can apply (as in the proof above) the trivial observation that for a manifold $A$, a submanifold $B\subset A$ with a tubular neighbourhood $U$ and normal bundle $\nu$ and a vector bundle $\zeta$ over $A$, there is a cofibration of Thom spaces $T\zeta|_{A\setminus U}\into T\zeta\to T(\zeta|_B\oplus\nu)$. Hence we obtain a cofibration $T\zeta_S^r\into T\gamma_{k+1}^L\to T((r+2)\gamma_{k+1}^L)$ that is turned by the functor $\Omega\Gamma$ to the fibration
$$\ol X_r^L=\Omega\Gamma T\zeta_S^r\to\Omega\Gamma T\gamma_{k+1}^L\to\Omega\Gamma T((r+2)\gamma_{k+1}^L).$$

Recall that the resolvent of such a fibration is an infinte continuation of the above sequence of maps to the left with the next entry being the loop space of the first entry, i.e. $\Omega^2\Gamma T((r+2)\gamma_{k+1}^L)$. This means we obtain the proposed fibration
$$\Omega^2\Gamma T((r+2)\gamma_{k+1}^L)\to\ol X_r^L\to\Omega\Gamma T\gamma_{k+1}^L.$$

\medskip\noindent\emph{Proof of \ref{ketto}.}\enspace\ignorespaces
We use proposition \ref{priminf}: $\Omega\Gamma T\gamma_{k+1}^L\cong\ol X_\infty^L$. We will construct a map $\varphi\colon\pi_{n+k+1}(\ol X_\infty^L,\ol X_r^L)\to\pi_{n+k+2}(\Gamma T((r+2)\gamma_{k+1}^L))=\pi_{n+k}(\Omega^2\Gamma T((r+2)\gamma_{k+1}^L))$ such that the diagram
$$\xymatrix{
\ldots\ar[r] & \pi_{n+k+1}(\ol X_\infty^L)\ar[r]\ar@{=}[d] & \pi_{n+k+1}(\ol X_\infty^L,\ol X_r^L)\ar[r]\ar[d]^\varphi & \pi_{n+k}(\ol X_r^L)\ar[r]\ar@{=}[d] & \ldots \\
\ldots\ar[r] & \pi_{n+k+1}(\Omega\Gamma T\gamma_{k+1}^L)\ar[r] & \pi_{n+k}(\Omega^2\Gamma T((r+2)\gamma_{k+1}^L))\ar[r] & \pi_{n+k}(\ol X_r^L)\ar[r] & \ldots
}$$
connecting the two long exact sequences commutes. By the 5-lemma, the commutativity of this diagram implies the statement of the theorem, hence it only remains to construct this map $\varphi$.

The relative homotopy group $\pi_{n+k+1}(\ol X_\infty^L,\ol X_r^L)$ can clearly be identified with the relative cobordism group of those prim maps $F\colon(M^{n+1},\partial M^{n+1})\to(\R^{n+k}\times\R_+,\R^{n+k})$ (with normal $L$-structures) for which the restriction $F|_{\partial M}$ is a $\tau_r$-map. We define $\varphi$ to assign to the cobordism class of such a map $F$ the cobordism class in $\Imm^{(r+2)\gamma_{k+1}^L}(\R^{n+k+2})$ of the immersion lift of its $\Sigma^{1_{r+1}}$-stratum (the normal bundle of this immersion splits to the direct sum of $r+1$ isomorphic bundles as before). The classifying space of the latter cobordism group is $\Gamma T((r+2)\gamma_{k+1}^L)$ and this $\varphi$ is clearly such that the diagram above is commutative, thus the proof is complete.
\end{prf}

\begin{crly}
Suppose that the universal bundle $\gamma_{k+1}^L$ is orientable (i.e. the linear group $L(k+1)$ is positive).
\begin{enumerate}
\item\label{ka} If the rational Euler class $e(\gamma_{k+1}^L)$ is non-trivial, then we have
$$\Omega\Gamma T\gamma_{k+1}^L\congq\ol X_r^L\times\Omega\Gamma T((r+2)\gamma_{k+1}^L).$$
\item\label{ki} If the rational Euler class $e(\gamma_{k+1}^L)$ is trivial, then we have
$$\ol X_r^L\congq\Omega^2\Gamma T((r+2)\gamma_{k+1}^L)\times\Omega\Gamma T\gamma_{k+1}^L.$$
\end{enumerate}
\end{crly}

\begin{prf}
Note that $H^*(BL(k+1);\Q)$ has no zero divisors, since it is a subring of the polynomial ring $H^*(BT(k+1);\Q)$ where $T(k+1)$ is the maximal torus in $L(k+1)$. Now if the Euler class of $\gamma_{k+1}^L$ is non-trivial (resp. trivial), then so is the Euler class of $(r+1)\gamma_{k+1}^L$, hence the Gysin sequence yields that the projection $\pi\colon S((r+1)\gamma_{k+1}^L)\to BL(k+1)$ induces an epimorphism (resp. monomorphism) in the cohomologies with rational coefficients.

\medskip\noindent\emph{Proof of \ref{ka}.}\enspace\ignorespaces
The above observation implies that the map of Thom spaces $T\zeta_S^r\to T\gamma_{k+1}^L$ induced by $\pi$ also induces an epimorphism in rational cohomology. Thus the cofibration $T\zeta_S^r\into T\gamma_{k+1}^L\to T((r+2)\gamma_{k+1}^L)$ rationally homologically splits, and so the homotopy long exact sequence of the fibration $\Omega\Gamma T\zeta_S^r\to\Omega\Gamma T\gamma_{k+1}^L\to\Omega\Gamma T((r+2)\gamma_{k+1}^L)$ splits rationally as well. Since the spaces involved are H-spaces, hence rationally products of Eilenberg--MacLane spaces, and $\Omega\Gamma T\zeta_S^r\cong\ol X_r^L$, this implies that we have
$$\Omega\Gamma T\gamma_{k+1}^L\congq\ol X_r^L\times\Omega\Gamma T((r+2)\gamma_{k+1}^L).$$

\medskip\noindent\emph{Proof of \ref{ki}.}\enspace\ignorespaces
Extend the cofibration $T\zeta_S^r\into T\gamma_{k+1}^L\to T((r+2)\gamma_{k+1}^L)$ to the right by the Puppe sequence. We obtain
$$T\zeta_S^r\to T\gamma_{k+1}^L\to T((r+2)\gamma_{k+1}^L)\to ST\zeta_S^r\to ST\gamma_{k+1}^L\to\ldots$$
and the map $ST\zeta_S^r\to ST\gamma_{k+1}^L$ induced by $\pi$ induces a monomorphism in rational cohomology, thus the cofibration $T((r+2)\gamma_{k+1}^L)\to ST\zeta_S^r\to ST\gamma_{k+1}^L$ rationally homologically splits, and so the homotopy long exact sequence of the fibration $\Omega^2\Gamma T((r+2)\gamma_{k+1}^L)\to\Omega^2\Gamma ST\zeta_S^r\to\Omega^2\Gamma ST\gamma_{k+1}^L$ splits rationally as well. Now the statement follows again from the spaces being H-spaces and the equalities $\Omega^2\Gamma ST\zeta_S^r=\Omega\Gamma T\zeta_S^r\cong\ol X_r^L$ and $\Omega^2\Gamma ST\gamma_{k+1}^L=\Omega\Gamma T\gamma_{k+1}^L$.
\end{prf}

\begin{ex}
The conditions of \ref{ka} and \ref{ki} above are satisfied in the case of $L=\SO$ respectively for $k$ odd and $k$ even.
\end{ex}

\begin{rmk}
The corollary above also has a nice geometric interpretation (similarly to remark \ref{keyrmk} about the geometric meaning of the key fibration) using part \ref{ketto} of theorem \ref{prclthm}:

In the case of \ref{ka}, we obtain that the sequence
$$0\to\pi_{n+k}(\ol X_r^L)\to\pi_{n+k}(\ol X_\infty^L)\to\pi_{n+k}(\ol X_\infty^L,\ol X_r^L)\to0$$
is exact at the middle term and has finite homology everywhere else. This means that a prim map $f\colon M^n\to\R^{n+k}$ is rationally prim cobordant to a map having at most $\Sigma^{1_r}$ singularities iff $f|_{\ol{\Sigma^{1_{r+1}}(f)}}$ is rationally prim null-cobordant.

Similarly, in the case of \ref{ki} the same is true for the sequence
$$0\to\pi_{n+k+1}(\ol X_\infty^L,\ol X_r^L)\to\pi_{n+k}(\ol X_r^L)\to\pi_{n+k}(\ol X_\infty^L)\to0.$$
This means that any prim map is rationally prim cobordant to a map having at most $\Sigma^{1_r}$ singularities (for any $r$). Moreover, if a prim map $f\colon M^n\to\R^{n+k}$ has at most $\Sigma^{1_r}$ singularities and represents $0$ in $\Prim_\infty^L(n,k)\otimes\Q$, then its rational prim $(\tau_r,L)$-cobordism class is completely determined by the rational prim cobordism class of
$$F|_{\ol{\Sigma^{1_{r+1}}(F)}}\colon\ol{\Sigma^{1_{r+1}}(F)}\to\R^{n+k+1}$$
for any prim map $F\colon(W^{n+1},\partial W^{n+1})\to(\R^{n+k}\times\R_+,\R^{n+k})$ that it bounds (i.e. for which $\partial W=M$ and $F|_M=f$).
\end{rmk}

%-------------------------------------------------------------------------------------------------------------------------------
%-------------------------------------------------------------------------------------------------------------------------------

\chapter{Computation of cobordism groups}\label{compcobgr}

In the present chapter we gather those cobordism groups $\Cob_\tau^L(n,P^{n+k})$ (and also $\Prim_\tau^L(n,P^{n+k})$) for which (complete or partial) computations already exist in literature using the tools of the chapters before. Although we considered quite general sets of singularities earlier, it should be noted that these computations seem to be rather difficult and most of the cobordism groups computed are such that only Morin singularities are allowed (i.e. the singularity set $\tau$ only contains singularities of coranks $0$ and $1$; see definition \ref{thomboar}).

Recall example \ref{morsing} which implies that if only Morin sungularities are allowed, then the singularity set is $\tau_r:=\{\Sigma^0,\Sigma^{1_1},\ldots,\Sigma^{1_r}\}$ for some $r\in\N\cup\{\infty\}$. Using the notation described in \ref{conv2} after definition \ref{cobtau}, throughout this chapter we put (for any stable linear group $L$), $\Cob_r^L(n,P):=\Cob_{\tau_r}^L(n,P)$, $X_r^L:=X_{\tau_r}^L$, $K_r^L:=K_{\tau_r}^L$, $V_r^L:=V_{\tau_r}^L$, $\nu_r^L:=\nu_{\tau_r}^L$, $\Prim_r^L(n,P):=\Prim_{\tau_r}^L(n,P)$, $\overline X_r^L:=\overline X_{\tau_r}^L$, $\ol K_r^L:=\ol K_{\tau_r}^L$, $\ol\nu_r^L:=\ol\nu_{\tau_r}^L$ (overlined letters denote analogous notions to the non-overlined ones for prim maps) and also $G_r^L:=G_{\Sigma^{1_r}}^L$, $\xi_r^L:=\xi_{\Sigma^{1_r}}^L$, $\tilde\xi_r^L:=\tilde\xi_{\Sigma^{1_r}}^L$, $\ol G_r^L:=\ol G_{\Sigma^{1_r}}^L$, $\zeta_r^L:=\zeta_{\Sigma^{1_r}}^L$, $\tilde\zeta_r^L:=\tilde\zeta_{\Sigma^{1_r}}^L$ (the latter three respectively are the analogues of the notions before them for prim maps). If $L=\O$, then we omit it from the superindex.

In section \ref{arbco} we show computations on cobordism groups in general; those on unoriented and cooriented cobordisms are taken from \cite{szszt} and \cite{hosszu} respectively. In section \ref{larco}, based on \cite{eszt} and \cite{2k+2}, we restrict to the case of maps where the dimension of the target is considerably larger than that of the source; the point of this is that generic maps with this property do not have too complicated singularities. Similarly, section \ref{1co} also considers maps of restricted codimension, namely $1$-codimensional maps; this was described in \cite{nulladik}. Then in section \ref{prico} we turn to cobordisms of prim maps and show results from \cite{nszt}, \cite{nehez} and \cite{ctrl}.

%-------------------------------------------------------------------------------------------------------------------------------

\section{Arbitrary cobordisms}\label{arbco}

The main tool that we will use in this section is the cohomological Kazarian spectral sequence with coefficients in a fixed ring $R$ which contains $\frac12$ (i.e. there is division by $2$). We will compute this spectral sequence for unoriented and cooriented prim maps, from which the sequence for unoriented and cooriented usual (non-prim) Morin maps can be obtained by lemma \ref{cover}. Using this, we compute the ranks of the free parts of the groups $\Prim_r(P)$, $\Cob_r(P)$, $\Prim_r^\SO(P)$ and $\Cob_r^\SO(P)$ for any manifold $P$ and also prove that unoriented cobordism groups are finite $2$-primary in many cases.

The following statement is well-known; see for example \cite{charclass}. Hence we will not give a proof here.

\begin{lemma}\label{cover}
Let $\tilde A$ be a CW-complex equipped with a cellular $\Z_2$-action, put $A:=\tilde A/\Z_2$ and let $q\colon\tilde A\to A$ be the quotient map. Then the homomorphism $q^*\colon H^*(A;R)\to H^*(\tilde A;R)$ is injective and gives an isomorphism between $H^*(A;R)$ and the $\Z_2$-invariant part $H^*(\tilde A;R)^{\Z_2}$ of $H^*(\tilde A;R)$.
\end{lemma}

\begin{rmk}\label{coverr}
In addition to the previous lemma, if $\tilde A'$ and is another space with a $\Z_2$ action as above and $f\colon\tilde A'\to\tilde A$ is a $\Z_2$-invariant map, then $f^*\colon H^*(\tilde A;R)\to H^*(\tilde A';R)$ respects the eigenspace decomposition of the $\Z_2$-action. This implies that if $\tilde A_0\subset\tilde A_1\subset\ldots$ is a filtration of $\tilde A$ by $\Z_2$-invariant subspaces and $A_i:=\tilde A_i/\Z_2$, then the cohomological spectral sequence of the filtration $A_0\subset A_1\subset\ldots$ of $A$ can be identified with the $\Z_2$-invariant part of that of $\tilde A$.
\end{rmk}

We will use the above lemma and remark with the substitutions $\tilde A:=\ol K_r^L$ and $A:=K_r^L$ and the $\Z_2$-action corresponding to changing the orientation of the kernel line bundle (recall remark \ref{kerbundle}: Prim maps are precisely the Morin maps with trivialised kernel line bundles).

\begin{rmk}\label{crutial}
The following observations will be crucial in the computations of the various Kazarian spectral sequences in the next two subsections.
\begin{enumerate}
\item The Kazarian space $\ol K_\infty^L$ for all $k$-codimensional prim maps coincides with the Kazarian space for $(k+1)$-codimensional immersions with normal $L$-structures, which is $BL(k+1)$, since prim maps and their immersion lifts can be identified. Moreover, the virtual normal bundle $\ol\nu_\infty^L$ is stably equivalent to the canonical bundle $\gamma_{k+1}^L$.
\item\label{truncate} If we truncate the Kazarian spectral sequence for any type of Morin maps at the $r$-th column, then we obtain the Kazarian spectral sequence for the same type of maps with at most $\Sigma^{1_r}$ singularities.
\end{enumerate}
\end{rmk}

\begin{rmk}\label{rk}
We claimed above that the ranks of the cobordism groups of unoriented and cooriented prim and usual Morin maps to any manifold $P$ can be computed from the respective Kazarian spectral sequences. The reason for this is proposition \ref{cobhomo} (and its analogue for prim cobordisms), which can be rewritten as
\begin{gather*}
\Cob_r^L(P)\otimes\Q\cong\displaystyle\bigoplus_{i=1}^\infty H_i(\cpt P;\Q)\otimes H^{i-k}(K_r^L;\tilde\Q_{\nu_r^L})\\
\text{and}~~~~\Prim_r^L(P)\otimes\Q\cong\displaystyle\bigoplus_{i=1}^\infty H_i(\cpt P;\Q)\otimes H^{i-k}(\ol K_r^L;\tilde\Q_{\ol\nu_r^L}).
\end{gather*}
Hence computing the ranks of the cobordism groups is equivalent to computing the ranks of the cohomologies of the respective Kazarian spaces twisted by the orientations of the respective virtual normal bundles (and of course the homologies of $\cpt P$, but it is no wonder that we also need such an information). In particular the ranks of the cobordism groups $\Cob_r^L(n,k)$ and $\Prim_r^L(n,k)$ are the same as the ranks of the $n$-th twisted cohomology groups of $K_r^L$ and $\ol K_r^L$ respectively.
\end{rmk}

\subsection{Unoriented cobordism groups}

We will now study the unoriented case $L(k):=\O(k)$ and we use $\tilde R:=\tilde R_{\ol\nu_\infty}$ to denote coefficients in the fixed ring $R$ twisted by the orientation of the universal normal bundle.

\begin{lemma}\label{twisto}
The twisted cohomology of the space $B\O(k+1)$ is
\begin{enumerate}
\item $H^*(B\O(k+1);\tilde R)\cong0$, if $k$ is even,
\item $H^*(B\O(k+1);\tilde R)=e\smallsmile R\big[p_1,\ldots,p_{\frac{k+1}2}\big]$, if $k$ is odd, where $e$ is the twisted Euler class and the $p_i$-s are the Pontryagin classes of $\gamma_{k+1}$.
\end{enumerate}
\end{lemma}

\begin{prf}
We will use the double cover $\pi\colon B\SO(k+1)\to B\O(k+1)$. Applying the Leray spectral sequence to this covering, we get
$$H^*(B\SO(k+1);R)\cong H^*(B\O(k+1);\pi_*(R)),$$
where $\pi_*(R)$ means the pushforward of the untwisted local coefficient system. This is locally $R\oplus R$ at each point and the non-trivial loop (homotopy class) acts on it by interchanging the summands, hence it can be decomposed as the sum $\pi_*(R)=R\oplus\tilde R_\pi$ of the invariant and the anti-invariant part. Thus we have
$$H^*(B\SO(k+1);R)\cong H^*(B\O(k+1);R)\oplus H^*(B\O(k+1);\tilde R_\pi).$$
Now using that the local systems $\tilde R_\pi$ and $\tilde R$ coincide (since $w_1(\ol\nu_\infty)=w_1(\pi)$), we get the statement of the lemma from the well-known facts
$$H^*(B\SO(k+1);R)=\begin{cases}
R\big[p_1\ldots,p_{\frac k2}\big],&\text{if }k\text{ is even}\\
R\big[p_1,\ldots,p_{\frac{k+1}2},e\big]\big/\big(p_{\frac{k+1}2}-e^2\big),&\text{if }k\text{ is odd}
\end{cases}$$
and $H^*(B\O(k+1);R)=R\big[p_1,\ldots,p_{\lfloor\frac{k+1}2\rfloor}\big]$.
\end{prf}

This lemma describes where the twisted Kazarian spectral sequence for unoriented prim maps converges. Now we will describe the spectral sequence itself.

\begin{prop}\label{punorsps}
If $\ol E_*^{*,*}$ is the twisted cohomological Kazarian spectral sequence for $k$-codimensional unoriented prim maps, then it has the following form:
\begin{enumerate}
\item If $k=2l$ is even, then for $i=0,1,\ldots$ the columns $\ol E_1^{2i,*}$ and $\ol E_1^{2i+1,*}$ in the first page are the ring $R[p_1,\ldots,p_l]$ shifted by $(2i+1)k$, i.e. $\ol E_1^{2i,j}\cong\ol E_1^{2i+1,j}$ is the degree-$(j-(2i+1)k)$ part of this polynomial ring, the differential $d_1$ is an isomorphism $\ol E_1^{2i,*}\cong\ol E_1^{2i+1,*}$, hence the second page is $\ol E_2^{*,*}=\ol E_\infty^{*,*}\cong0$.
\item If $k=2l+1$ is odd, then all differentials vanish and for $i=0,1,\ldots$ the columns $\ol E_1^{2i,*}=\ol E_\infty^{2i,*}$ are $0$ and the columns $\ol E_1^{2i+1,*}=\ol E_\infty^{2i+1,*}$ are the ring $R[p_1,\ldots,p_l]$ shifted by $(2i+1)k$.
\end{enumerate}
\end{prop}

\begin{prf}
The first page of this Kazarian spectral sequence (similarly to definition \ref{kazspec}) is $\ol E_1^{i,j}=H^{i+j}(T\zeta_i;\tilde R)$. Since the universal bundle $\zeta_i$ coincides with $i\gamma_k\oplus\varepsilon^i$ over $B\O(k)$ (see \cite{hosszu}) and the virtual normal bundle $\ol\nu_i$ is stably equivalent to $\gamma_k$, we have $w_1(\zeta_i)=rw_1(\gamma_k)=rw_1(\ol\nu_i)$. Hence the Thom isomorphism yields
$$\ol E_1^{i,*}=H^{i+*}(T\zeta_i;\tilde R)\cong H^{*-ik}(B\O(k);\tilde R^{i+1}),$$
where $\tilde R^{i+1}$ means the $(i+1)$-times tensor product $\tilde R\otimes\ldots\otimes\tilde R$.

Now for $k=2l$ or $k=2l+1$, the (untwisted) rational cohomology $H^*(B\O(k);R)$ is $R[p_1,\ldots,p_l]$ and the twisted cohomology $H^*(B\O(k);\tilde R)$ is $0$ for $k=2l+1$ and isomorphic to $R[p_1,\ldots,p_l]$ shifted by $k$ for $k=2l$. This implies everything we claimed about the first page of the sequence and it is not hard to see from this that for odd $k$, all differentials vanish.

It remains to prove that for even $k$, the differential $d_1$ is an isomorphism between the columns $\ol E_1^{2i,*}$ and $\ol E_1^{2i+1,*}$ for all $i$. By lemma \ref{twisto}, we know that in this case all terms $\ol E_\infty^{i,j}$ are $0$, hence the entries of lowest degree in the $0$-th column must be mapped onto those in the next column by $d_1$ isomorphically because this is the only chance for these entries to disappear. Thus $d_1(e)$ (where $e\in e\smallsmile R[p_1,\ldots,p_l]=\ol E_1^{0,*}$) is the twisted Thom class $u(\zeta_1)$.

\medskip\begin{sclaim}
For any monomial $p_I$ of Pontryagin classes and $i=0,1,\ldots$, we have $d_1^{2i,*}(e\smallsmile p_I)=d_1^{2i,*}(e)\smallsmile p_I$.
\end{sclaim}

\begin{sprf}
The differential $d_1^{2i,*}\colon\ol E_1^{2i,*}\to\ol E_1^{2i+1,*}$ is by definition the boundary map $\delta$ in the exact sequence of the triple $(\ol K_{2i+1},\ol K_{2i},\ol K_{2i-1})$. The Pontryagin classes $p_1,\ldots,p_l$ used in $\ol E_1^{2i,*}$ can be identified with those of the universal bundle $\gamma_k$ restricted to the base space of $D(2i\gamma_k\oplus\varepsilon^{2i})\subset\ol K_\infty\cong B\O(k+1)$ (recall that the Kazarian space $\ol K_\infty$ is constructed from blocks $D\zeta_r=D(r\gamma_k\oplus\varepsilon^r)$), thus the map $j^*$ induced by the inclusion $j\colon\ol K_{2i}\into\ol K_{2i+1}$ maps them into each other. Now we have
$$d_1^{2i,*}(e\smallsmile p_I)=\delta(e\smallsmile p_I)=\delta(e\smallsmile j^*(p_I))=\delta(e)\smallsmile p_I+e\smallsmile\delta(j^*(p_I))=\delta(e)\smallsmile p_I,$$
as claimed.
\end{sprf}

Hence the differential $d_1$ is an isomorphism between the first two columns and then an induction shows the same for the columns $2i$ and $(2i+1)$ for all $i$.
\end{prf}

\begin{crly}\label{unorsps}
If $E_*^{*,*}$ is the twisted cohomological Kazarian spectral sequence for $k$-codimensional unoriented Morin maps, then it has the following form:
\begin{enumerate}
\item If $k=2l$ is even, then for $i=0,1,\ldots$ the columns $E_1^{4i,*}$ and $E_1^{4i+1,*}$ in the first page are the ring $R[p_1,\ldots,p_l]$ shifted by $(4i+1)k$, the columns $E_1^{4i+2,*}$ and $E_1^{4i+3,*}$ vanish, the differential $d_1$ is an isomorphism $E_1^{4i,*}\cong E_1^{4i+1,*}$, hence the second page is $E_2^{*,*}=E_\infty^{*,*}\cong0$.
\item If $k=2l+1$ is odd, then the whole spectral sequence is constant $0$.
\end{enumerate}
\end{crly}

\begin{prf}
By lemma \ref{cover}, we need to understand the $\Z_2$-action on the spectral sequence $\ol E_*^{*,*}$ corresponding to changing the orientation of the kernel line bundle to obtain the spectral sequence $E_*^{*,*}$ with coefficients twisted by the orientation of $\nu_\infty$. Considering the $i$-th column (where the first page contains the cohomologies of $B\O(k)$), we have the coefficient system $\tilde R_{\nu_\infty}\otimes\tilde R_{\xi_i}=\tilde R_{\tilde\xi_i}$.

It was proved in \cite{rsz} and \cite{rrthesis} that the maximal compact subgroup $G_i$ of $\Aut_\AA\Sigma^{1_i}$ (see section \ref{semiglob}) is isomorphic to $\O(1)\times\O(k)$ with target representation
$$\tilde\lambda(\varepsilon,A)=\varepsilon^{i+1}\oplus A\oplus\left\lceil\frac{i-1}2\right\rceil1\oplus\left\lfloor\frac{i-1}2\right\rfloor\varepsilon\oplus\left\lfloor\frac{i}2\right\rfloor A\oplus\left\lceil\frac{i}2\right\rceil\varepsilon A$$
for $\varepsilon\in\O(1)$ and $A\in\O(k)$. The $\Z_2$-action discussed above is given by the correspondence $\varepsilon\mapsto-\varepsilon$ on $\O(1)$.

For $k=2l$ this action changes the orientation of $\tilde\lambda$ exactly if $i+1+\lfloor\frac{i-1}2\rfloor$ is odd, i.e. for $i\equiv2,3~\mod4$. This means by remark \ref{coverr} that for $i\equiv0,1~\mod4$ the columns $E_1^{i,*}$ are the same as for prim maps with the differential $d_1$ also being the same; and for $i\equiv2,3~\mod4$ the columns $E_1^{i,*}$ vanish.

For $k=2l+1$ the above action changes the orientation of $\tilde\lambda$ exactly if $i+1+\lfloor\frac{i-1}2\rfloor+\lceil\frac{i}2\rceil$ is odd, i.e. if $i$ is odd. Hence by remark \ref{coverr} the columns $E_1^{i,*}$ are the same as for prim maps for $i$ even and vanish for $i$ odd, which means that the whole first page of this spectral sequence is $0$.
\end{prf}

With the help of the spectral sequences described above, we are able to prove that in many cases the cobordism groups $\Prim_r(n,k)$ and $\Cob_r(n,k)$ are finite $2$-primary, i.e. they are elements of the Serre class $\CC_2$. Moreover, we are also able to compute the ranks of all cobordism groups $\Prim_r(P)$ and $\Cob_r(P)$ by considering the case $R:=\Q$; see remark \ref{rk}.

\begin{thm}\label{punorthm}
$~$
\begin{enumerate}
\item\label{u1} The group $\Prim_r(n,k)$ is finite $2$-primary, if $k$ is even and $r$ is either an odd number or $\infty$.
\item\label{u2} If both $k$ and $r$ are even, then we have
$$\Prim_r(n,k)\otimes\Q\cong H^{n-r(k+1)-k}(B\O(k);\Q),$$
hence its rank is the number of partitions $p_{\frac k2}\big(\frac{n-r(k+1)-k}4\big)$, in particular it is $0$ if $n-r(k+1)-k$ is not a multiple of $4$.
\item\label{u3} If $k=2l+1$ is odd and $r=\infty$, then we have
$$\Prim_\infty(n,k)\otimes\Q\cong H^{n-k-1}(B\O(k+1);\Q),$$
hence its rank is $p_{\frac{k+1}2}\big(\frac{n-k-1}4\big)$, in particular it is $0$ if $n-k-1$ is not a multiple of $4$.
\item\label{u4} If $k=2l+1$ is odd and $r<\infty$, then $\Prim_r(n,k)\otimes\Q$ is also the degree-$(n-k-1)$ part of a polynomial ring, namely
$$\Prim_r(n,k)\otimes\Q\cong\left[\Q[p_1,\ldots,p_{l+1}]\big/\big(p_{l+1}^{\lceil\frac{r}2\rceil}\big)\right]^{n-k-1}.$$
\end{enumerate}
\end{thm}

\begin{prf}
\emph{Proof of \ref{u1}.}\enspace\ignorespaces
The Kazarian spectral sequence from theorem \ref{punorsps} (together with part \ref{truncate} of remark \ref{crutial}) yields that in this case we have $H^*(\ol K_r;\tilde R)\cong0$ for any ring $R$ where $2$ is invertible. By the Thom isomorphism we have $H^{*}(\ol K_r;\tilde R)\cong H^{*+k}(T\ol\nu_r;R)$, hence $T\ol\nu_r$ has trivial (co)homology groups modulo the class $\CC_2$. Now Serre's mod-$\CC_2$ Hurewicz theorem from \cite{serre} implies that the stable homotopy groups $\pi_*^s(T\ol\nu_r)$ are also trivial modulo $\CC_2$, i.e. they are finite $2$-primary. Thus we have
$$\Prim_r(n,k)\cong\pi_{n+k}(\ol X_r)\cong\pi_{n+k}^s(T\ol\nu_r)\in\CC_2.$$

\medskip\noindent\emph{Proof of \ref{u2}.}\enspace\ignorespaces
Again we use part \ref{truncate} of remark \ref{crutial}; we obtain that in the twisted rational Kazarian spectral sequence $\ol E_2^{r,*}=\ol E_\infty^{r,*}$ is the only non-zero column in the second page and it is the ring $\Q\big[p_1,\ldots,p_{\frac k2}\big]=H^*(B\O(k);\Q)$ shifted by $(r+1)k$. Hence the twisted rational cohomology ring of the Kazarian space is $H^*(\ol K_r;\tilde\Q)\cong H^{*-r(k+1)-k}(B\O(k);\Q)$ and we get the statement of the theorem by remark \ref{rk}.

\medskip\noindent\emph{Proof of \ref{u3} and \ref{u4}.}\enspace\ignorespaces
In this case the twisted rational Kazarian spectral sequence yields
$$H^n(\ol K_\infty;\tilde\Q)\cong\bigoplus_{i=0}^\infty\ol E_\infty^{2i+1,n-2i-1}\cong\bigoplus_{i=0}^\infty\big[\Q[p_1,\ldots,p_l]\big]^{n-2i(k+1)-k-1}.$$
If we introduce a new variable $p_{l+1}$ of degree $2k+2$ and add it to the previous variable set $\{p_1,\ldots,p_l\}$, then we obtain that the direct sum above is the degree-$(n-k-1)$ part of the polynomial ring $\Q[p_1,\ldots,p_{l+1}]$. This ring is isomorphic to $H^*(B\O(k+1);\Q)$, thus $\Prim_\infty(n,k)\otimes\Q\cong H^n(\ol K_\infty;\tilde\Q)$ is of the form as we claimed. If we truncate the spectral sequence at the $r$-th column to obtain $\Prim_r(n,k)$, then the new variable $p_{l+1}$ has to be such that its $i$-th power is $0$ if $2i+1>r$, i.e. if $i\ge\lceil\frac r2\rceil$, hence we obtain the statement of the theorem again.
\end{prf}

Completely analogously one can prove the following.

\begin{thm}\label{unorthm}
$~$
\begin{enumerate}
\item The group $\Cob_r(n,k)$ is finite $2$-primary, if either $r=\infty$ or $k$ is odd or $k$ is even and $r$ is not divisible by $4$.
\item If $k$ is even and $r$ is divisible by $4$, then we have
$$\Cob_r(n,k)\otimes\Q\cong H^{n-r(k+1)-k}(B\O(k);\Q),$$
hence its rank is the number of partitions $p_{\frac k2}\big(\frac{n-r(k+1)-k}4\big)$, in particular it is $0$ if $n-k$ is not a multiple of $4$.
\end{enumerate}
\end{thm}

\subsection{Cooriented cobordism groups}

The present subsection is analogous to the previous subsection, the only change being that now we consider the cooriented case $L(k):=\SO(k)$ and with $R:=\Q$. We will use $e$ and $p_i$ respectively to denote the rational Euler class and Pontryagin classes of $\gamma_k^\SO$ and we put $A:=\Q[p_1,\ldots,p_l]$ (where the number $l$ will be defined later).

\begin{prop}\label{pcoorsps}
If $\ol E_*^{*,*}$ is the cohomological Kazarian spectral sequence for $k$-codimensional cooriented prim maps, then it has the following form:
\begin{enumerate}
\item If $k=2l$ is even, then for $i=0,1,\ldots$ the column $\ol E_1^{i,*}$ in the first page is the ring $(e^i\smallsmile A)\oplus(e^{i+1}\smallsmile A)$, the differential $d_1$ maps $e^{i+1}\smallsmile A\subset\ol E_1^{i,*}$ to $e^{i+1}\smallsmile A\subset\ol E_1^{i+1,*}$ isomorphically, hence the second page is $\ol E_2^{0,*}=\ol E_\infty^{0,*}\cong A$ and $\ol E_2^{i,*}=\ol E_\infty^{i,*}\cong0$ for $i>0$.
\item If $k=2l+1$ is odd, then all differentials vanish and for $i=0,1,\ldots$ the columns $\ol E_1^{i,*}=\ol E_\infty^{i,*}$ are the ring $A$ shifted by $ik$.
\end{enumerate}
\end{prop}

\begin{prf}
Analogously to the proof of proposition \ref{punorsps}, the first page of this spectral sequence is
$$\ol E_1^{i,j}=H^{i+j}(T\zeta_i^\SO;\Q)\cong H^j(T(i\gamma_k^\SO);\Q)\cong e^i\smallsmile H^{j-ik}(B\SO(k);\Q)$$
using that the universal bundle $\zeta_i^\SO$ coincides with the bundle $i\gamma_k^\SO\oplus\varepsilon^i$ over $B\SO(k)$ and the inclusion $B\SO(k)\into T\gamma_k^\SO$ induces an injection $H^*(T(i\gamma_k^\SO);\Q)\to H^*(B\SO(k);\Q)$ with image equal to the ideal generated by $e^i$.

For $k=2l+1$, the column $H^{*-ik}(B\SO(k);\Q)$ in the first page coincides with $A$ shifted by $ik$ and the sequence converges to $H^*(B\SO(k+1);\Q)=\Q[p_1,\ldots,p_{l+1}]\oplus(e\smallsmile \Q[p_1,\ldots,p_{l+1}])$. Now an easy computation shows that in this case the ranks of $\ol E_1^{*,*}$ and $H^*(B\SO(k+1);\Q)$ are the same, which implies that all differentials vanish. Hence the second statement is proved.

For $k=2l$, the column $H^{*-ik}(B\SO(k);\Q)$ is isomorphic to $A\oplus(e\smallsmile A)$ shifted by $ik$ and this shift is given by $e^i$, that is, this column is $(e^i\smallsmile A)\oplus(e^{i+1}\smallsmile A)$. This spectral sequence converges to $H^*(B\SO(k+1);\Q)=A$ and it can be identified with $A\subset\ol E_1^{0,*}$. Hence the Euler class $e$ from the lowest degree entry in the $0$-th column must be mapped to a non-zero multiple of $e$ in the next column, since this is the only chance for the multiples of $e$ to disappear. In the same way, the differential $d_1^{i,k}\colon\ol E_1^{i,k}\to\ol E_1^{i+1,k}$ maps $e^i$ to a non-zero multiple of $e^i$.

\medskip\begin{sclaim}
For any monomial $p_I$ of Pontryagin classes and $i=0,1,\ldots$, we have $d_1^{i,*}(e\smallsmile p_I)=d_1^{i,*}(e)\smallsmile p_I$.
\end{sclaim}

\begin{sprf}
The differential $d_1^{i,*}\colon\ol E_1^{i,*}\to\ol E_1^{i+1,*}$ is by definition the boundary map $\delta$ in the exact sequence of the triple $(\ol K_{i+1}^\SO,\ol K_{i}^\SO,\ol K_{i-1}^\SO)$. The Pontryagin classes $p_1,\ldots,p_l$ used in $\ol E_1^{i,*}$ can be identified with those of the universal bundle $\gamma_k^\SO$ restricted to the base space of $D(i\gamma_k^\SO\oplus\varepsilon^{i})\subset\ol K_\infty^\SO\cong B\SO(k+1)$ (recall that the Kazarian space $\ol K_\infty^\SO$ is constructed from blocks $D\zeta_r^\SO=D(r\gamma_k^\SO\oplus\varepsilon^r)$), thus the map $j^*$ induced by the inclusion $j\colon\ol K_{i}^\SO\into\ol K_{i+1}^\SO$ maps them into each other. Now we have
$$d_1^{i,*}(e\smallsmile p_I)=\delta(e\smallsmile p_I)=\delta(e\smallsmile j^*(p_I))=\delta(e)\smallsmile p_I+e\smallsmile\delta(j^*(p_I))=\delta(e)\smallsmile p_I,$$
as claimed.
\end{sprf}

This claim, combined with the observations above it, shows that the $e^{i+1}\smallsmile A$ part of $\ol E_1^{i,*}$ is mapped to that of $\ol E_1^{i+1,*}$ isomophically for all $i$.
\end{prf}

\begin{crly}\label{coorsps}
If $E_*^{*,*}$ is the cohomological Kazarian spectral sequence for $k$-codi- mensional cooriented Morin maps, then it has the following form:
\begin{enumerate}
\item If $k=2l$ is even, then for $i=0,1,\ldots$ the column $E_1^{4i,*}$ in the first page is the ring $(e^{4i}\smallsmile A)\oplus(e^{4i+1}\smallsmile A)$, the column $E_1^{4i+1,*}$ is the ring $e^{4i+1}\smallsmile A$, the column $E_1^{4i+2,*}$ vanishes and the column $E_1^{4i+3,*}$ is the ring $e^{4i+4}\smallsmile A$, the differential $d_1$ maps $e^{4i+1}\smallsmile A\subset\ol E_1^{4i,*}$ onto $e^{4i+1}\smallsmile A\subset\ol E_1^{4i+1,*}$ and $e^{4i+4}\smallsmile A\subset\ol E_1^{4i+3,*}$ onto $e^{4i+4}\smallsmile A\subset\ol E_1^{4i+4,*}$ isomorphically, hence the second page is $E_2^{0,*}=E_\infty^{0,*}=A$ and $E_2^{i,*}=E_\infty^{i,*}\cong0$ for $i>0$.
\item If $k=2l+1$ is odd, then all differentials vanish and for $i=0,1,\ldots$ the columns $E_1^{2i,*}=E_\infty^{2i,*}$ are the ring $A$ shifted by $2ik$ and the columns $E_1^{2i+1,*}=E_\infty^{2i+1,*}$ are $0$.
\end{enumerate}
\end{crly}

\begin{prf}
By lemma \ref{cover} (analogously to the proof of corollary \ref{unorsps}) we need to understand the $\Z_2$-action on the spectral sequence $\ol E_*^{*,*}$ corresponding to changing the orientation of the kernel line bundle to obtain the spectral sequence $E_*^{*,*}$.

The maximal compact subgroup $G_i^\SO$ is by \cite{rsz} and \cite{rrthesis} isomorphic to $\{(\varepsilon,B)\in\O(1)\times\O(k)\mid\varepsilon^i\cdot\det B>0\}$ with source representation
$$\lambda(\varepsilon,B)=\left\lceil\frac{i-1}2\right\rceil1\oplus\left\lfloor\frac{i+1}2\right\rfloor\varepsilon\oplus\left\lfloor\frac{i}2\right\rfloor B\oplus\left\lceil\frac{i}2\right\rceil\varepsilon B$$
and the restriction of this to $\SO(1)\times\SO(k)\cong\SO(k)\cong\ol G_i^\SO$ is the source representation for prim maps. The $\Z_2$-action is again given by the correspondence $\varepsilon\mapsto-\varepsilon$ on $\O(1)$, which now also changes the orientation of $B\in\O(k)$ iff $i$ is odd in order to keep $\varepsilon^i\cdot\det B$ positive.

Recall that for $k=2l$ we had $\ol E_1^{i,*}=u(\zeta_i^\SO)\smallsmile(A\oplus e\smallsmile A)$ shifted by $-i$. The $\Z_2$-action described above maps $e$ to $(-1)^ie$ and changes the orientation of $\lambda$ (and so the sign of $u(\zeta_i^\SO)$) exactly if $\lfloor\frac{i+1}2\rfloor+\lfloor\frac{i}2\rfloor+\lceil\frac{i}2\rceil=\lfloor\frac{i+1}2\rfloor+i$ is odd, i.e. for $i\equiv2,3~\mod4$. This means by remark \ref{coverr} that the columns $E_1^{i,*}$, which are the $\Z_2$-invariant parts of $\ol E_1^{i,*}$, are as we claimed. The differential $d_1$ is the restriction of the differential of $\ol E_1^{*,*}$, hence it is an isomorphism between the direct summands $e^{i+1}\smallsmile A\subset E_1^{i,*}$ and $e^{i+1}\smallsmile A\subset E_1^{i+1,*}$ where these exist.

For $k=2l+1$ and $i$ even (resp. odd), the above action changes the orientation of $\lambda$ exactly if $\lfloor\frac{i+1}2\rfloor+\lceil\frac{i}2\rceil$ (resp. $\lfloor\frac{i+1}2\rfloor+\lfloor\frac{i}2\rfloor$) is odd, i.e. the orientation changes iff $i$ is odd. Hence by remark \ref{coverr} the columns $E_1^{i,*}$ are the same as for prim maps for $i$ even and vanish for $i$ odd.
\end{prf}

Again the spectral sequences computed above can be used to obtain the ranks of the cobordism groups $\Prim_r^\SO(P)$ and $\Cob_r^\SO(P)$; exactly the same considerations as in the proof of theorem \ref{punorthm} yield the following two theorems.

\begin{thm}
$~$
\begin{enumerate}
\item If $r=\infty$, then we have
$$\Prim_\infty^\SO(n,k)\otimes\Q\cong H^n(B\SO(k+1);\Q).$$
\item If $k=2l$ is even and $r<\infty$, then $\Prim_r^\SO(n,k)\otimes\Q$ is also the degree-$n$ part of a polynomial ring, namely
$$\Prim_r^\SO(n,k)\otimes\Q\cong\big[\Q[p_1,\ldots,p_l]\oplus(e^{r+1}\smallsmile\Q[p_1,\ldots,p_l])\big]^n.$$
\item If $k=2l+1$ is odd and $r<\infty$, then we have
$$\Prim_r^\SO(n,k)\otimes\Q\cong\big[\Q[p_1,\ldots,p_l,e]/(e^{r+1})\big]^n.$$
\end{enumerate}
\end{thm}

\begin{thm}\label{coorthm}
$~$
\begin{enumerate}
\item If $k$ is even and either $r\equiv1,2~\mod4$ or $r=\infty$, then we have
$$\Cob_r^\SO(n,k)\otimes\Q\cong H^n(B\O(k);\Q),$$
hence its rank is the number of partitions $p_{\frac k2}\big(\frac n4\big)$, in particular it is $0$ if $n$ is not a multiple of $4$.
\item If $k=2l$ is even and $r\equiv0,3~\mod4$, then we have
$$\Cob_r^\SO(n,k)\otimes\Q\cong\big[\Q[p_1,\ldots,p_l]\oplus(e^{r+1}\smallsmile\Q[p_1,\ldots,p_l])\big]^n.$$
\item If $k$ is odd and $r=\infty$, then we have
$$\Cob_\infty^\SO(n,k)\otimes\Q\cong H^n(B\O(k+1);\Q),$$
hence its rank is $p_{\frac{k+1}2}\big(\frac n4\big)$, in particular it is $0$ if $n$ is not a multiple of $4$.
\item If $k=2l+1$ is odd and $r<\infty$, then we have
$$\Cob_r^\SO(n,k)\otimes\Q\cong\left[\Q[p_1,\ldots,p_{l+1}]\big/\big(p_{l+1}^{\lceil\frac r2\rceil}\big)\right]^n.$$
\end{enumerate}
\end{thm}

\subsection{General observations}

The computations in the previous two subsections allow us to obtain a few additional general properties of cobordism groups.

\begin{prop}
Let $k$ be an odd number and $r\in\N\cup\{\infty\}$. Then for any manifold $P^{n+k}$ we have
$$\Cob_r^\SO(n,P^{n+k})\otimes\Q\cong\displaystyle\bigoplus_{2i\le r}\Imm^{\tilde\xi_{2i}^\SO}(P)\otimes\Q.$$
\end{prop}

\begin{prf}
By corollary \ref{coorsps} we know that the spaces $T\xi_{2i+1}^\SO$ have trivial rational cohomologies, hence so do the spaces $T\tilde\xi_{2i+1}^\SO$ (using the Thom isomorphism for the virtual bundle $\nu_{2i+1}^\SO$), and so $\Gamma T\tilde\xi_{2i+1}^\SO$ is rationally trivial. Similarly, corollary \ref{coorsps} also implies $T\tilde\xi_{2i}^\SO\congq T\tilde\zeta_{2i}^\SO=S^{2i}T((2i+1)\gamma_k^\SO)$. Hence we can use theorem \ref{ratrivi1} recursively to obtain
$$X_r^\SO\congq\prod_{2i\le r}\Gamma T\tilde\xi_{2i}^\SO$$
and now the statement follows trivially.
\end{prf}

\begin{rmk}
We also proved here that $\Imm^{\tilde\xi_{2i}^\SO}(P)\otimes\Q$ is the rational cobordism group of immersions to $P$ with normal bundles induced from $(2i+1)\gamma_k^\SO\oplus\varepsilon^{2i}$, i.e. whose normal bundles split to the direct sum of $2i+1$ oriented $k$-bundles and a $2i$-dimensional trivial bundle.
\end{rmk}

\begin{prop}
Put $\tau:=\tau_r\cup\{\rm{III}_{2,2}\}$, where $r\ge2$ and $\rm{III}_{2,2}$ is the simplest $\Sigma^2$ singularity (see \cite{math}). Then for any manifold $P^{n+k}$ we have
$$\Cob_\tau^\SO(n,P^{n+k})\otimes\Q\cong\left(\Cob_r^\SO(P)\oplus\Imm^{\tilde\xi_{\rm{III}_{2,2}}^\SO}(P)\right)\otimes\Q.$$
\end{prop}

\begin{prf}
The global normal form of the singularity $\rm{III}_{2,2}$ can be found in \cite{rsz}. Now the case of odd $k$ follows from theorem \ref{ratrivi1} again similarly to the proof above. If $k$ is even, then theorem \ref{coorthm} implies that the conditions of theorem \ref{ratrivi2} are satisfied, hence the conclusion is the same again.
\end{prf}

Similar examples can be constructed to show details on the free parts of cobordism groups of other types of maps as well.

%-------------------------------------------------------------------------------------------------------------------------------

\section{Large codimensional cobordisms}\label{larco}

In the previous section we saw theorems on cobordism groups of Morin maps in general, but these were not complete computations, as we only determined the free parts of these groups (and in some cases we also saw that they are finite $2$-primary). Here we will completely compute cobordism groups of maps that only have the simplest type of singularity, i.e. maps with only regular ($\Sigma^0$) and fold ($\Sigma^{1,0}$) points.

Recall from subsection \ref{strata} that the submanifold $\Sigma^r(f)\subset M^n$ for a map $f\colon M^n\to P^{n+k}$ consists of the points where the corank of $df$ is $r$; and the submanifold $\ol{\Sigma^{1_r}(f)}\subset M^n$ is the set of singular points of $f|_{\ol{\Sigma^{1_{r-1}}(f)}}$. Now the jet transversality theorem (together with a bit of linear algebra) has the straighforward consequences $\codim\Sigma^r(f)=r(k+r)$ and $\codim\Sigma^{1_r}(f)=r(k+1)$.

Thus for $n<2k+2$, generic maps $M^n\to P^{n+k}$ only have fold singularities. If $n<2k+1$, then even a cobordism between two maps from $n$-manifolds to an $(n+k)$-manifold only has fold singularities, hence in this case we have
$$\Cob_1(n,P^{n+k})=\NN_n(P)$$
(and a similar equation holds for maps with given normal structures). Similarly for $n<2k+4$, generic maps $M^n\to P^{n+k}$ only have fold and cusp singularities. This means if $n=2k+1$ or $n=2k+2$, then the groups $\Cob_1^L(n,P^{n+k})$ are (a priori) not the same as the respective bordism groups of $P$, but the cobordisms are not too far away from being generic (since a generic map $W^{n+1}\to P\times[0,1]$ can only have fold and cusp singularities). These ($n=2k+1$ and $n=2k+2$) are the two cases we will consider in this section with normal structures $L=\O$ and $L=\SO$.

\subsection{The $(2k+1)$-dimensional case}

By the considerations above, we have natural forgetful homomorphisms
$$\varphi\colon\Cob_1(2k+1,P^{3k+1})\to\NN_{2k+1}(P)~~~~\text{and}~~~~\psi\colon\Cob_1^\SO(2k+1,P^{3k+1})\to\Omega_{2k+1}(P),$$
which send the cobordism class of a map $f\colon M\to P$ to the bordism class of $f$. To determine these, we have to understand the cusps of generic bordisms, which are maps from $(2k+2)$-manifolds to $(3k+2)$-manifolds. The set of cusp points of such a map is a $0$-dimensional manifold, i.e. a discrete set of points.

\begin{rmk}\label{cusproot}
According to Morin \cite{mor}, the root of cusp singularity is the equivalence class of the map 
$$\sigma_2\colon(\R^{2k+2},0)\to(\R^{3k+2},0);~(x,y,z)\mapsto(X,Y,Z)$$
where the coordinates are $x=(x_1,\ldots,x_{2k})\in\R^{2k}$, $y\in\R^1$, $z\in\R^1$ and $X=(X_1,\ldots,$ $X_{2k+1})\in\R^{2k+1}$, $Y=(Y_1,\ldots,Y_k)\in\R^k$, $Z\in\R^1$ and the map is given by the formula
\begin{alignat*}2
&X_1=x_1,\ldots,X_{2k}=x_{2k},X_{2k+1}=y,\\
&Y_1=zx_1+z^2x_2,Y_2=zx_3+z^2x_4,\ldots,Y_k=zx_{2k-1}+z^2x_{2k},\\
&Z=zy+z^3.
\end{alignat*}
In other words $\sigma_2$ has an isolated cusp point in the origin.
\end{rmk}

\begin{lemma}\label{removecusps}
Let $F\colon W^{2k+2}\to Q^{3k+2}$ be a map between manifolds (possibly with boundary), let $p$ and $q$ be two cusp points of $F$ and assume that there is an arc $\alpha\subset W$ connecting $p$ and $q$ such that $\Sigma(F)\cap\alpha=\{p,q\}$. Then there is a map $F'\colon W\to P$ that coincides with $F$ outside an arbitrarily small neighbourhood of $\alpha$ and has cusp set $\Sigma^{1,1}(F')=\Sigma^{1,1}(F)\setminus\{p,q\}$.
\end{lemma}

\begin{prf}
We can choose a coordinate neighbourhood of the arc $\alpha$ such that this neighbourhood corresponds to the disk
$$D:=\{(x,y,z)\in\R^{2k+2}\mid\lv x\rv\le\varepsilon,-\varepsilon\le y\le 1+\varepsilon,|z|\le1\}$$
(with the same coordinates as in remark \ref{cusproot} and some $\varepsilon>0$) and the map $F$ restricted to this neighbourhood is the composition of a map of the form
$$f\colon D\to\R^{3k+2};~(x,y,z)\mapsto(X,Y,Z)$$
given by the expression
\begin{alignat*}2
&X_1=x_1,\ldots,X_{2k}=x_{2k},X_{2k+1}=y,\\
&Y_1=zx_1+z^2x_2,Y_2=zx_3+z^2x_4,\ldots,Y_k=zx_{2k-1}+z^2x_{2k},\\
&Z=zy(1-y)+z^3,
\end{alignat*}
with an immersion $\R^{3k+2}\imto Q$. Clearly the $Z$-coordinate function of $f$, as a function $y\mapsto zy(1-y)+z^3$, has two critical points if $y<0$ or $y>1$ and these cancel at $y=0$ and $y=1$ respectively, hence the cusp points of $f$ are $(0,0,0)$ and $(0,1,0)$.

Now we define the map $f'\colon D\to\R^{3k+2}$ in the same way as $f$, except that the $Z$-coordinate function is
$$Z=z(g(y)+y(1-y))+z^3,$$
where $g\colon\R\to[-1,0]$ is a smooth function with $g|_{\R\setminus[-\varepsilon,1+\varepsilon]}\equiv0$, $g|_{[0,1]}\equiv-1$ and $g$ is decreasing on the negative numbers and increasing on the positives. Now the $Z$-coordinate function of $f'$ (again as a function of $y$) always has two critical points, hence $f'$ has no cusps. By a bump fuction we can make $f'$ coincide with $f$ on the boundary of $D$, and so if we compose it with the immersion $\R^{3k+2}\imto Q$ as above, then we get the desired map $F'$ that has the same cusp set as $F$, except that $p$ and $q$ are deleted.
\end{prf}

We will now first consider the simpler case of unoriented cobordisms, and after that get to the cooriented case.

\begin{lemma}\label{valamilyenlemma}
For any integers $c\ge0$ and $k>0$ there is a closed manifold $W^{2k+2}$ and a (stable) map $F\colon W^{2k+2}\to\R^{3k+2}$ with $c$ cusps.
\end{lemma}

\begin{prf}
If $F\colon W^{2k+2}\to\R^{3k+2}$ is a map with at least two cusps, then we can add $1$-handles to make $W$ connected and since the codimension of the singular set is $k+1>1$, we can use the previous lemma to remove pairs of cusps. Hence it is enough to show that there is a manifold $W^{2k+2}$ and a map $F\colon W\to\R^{3k+2}$ with an odd number of cusps.

Put $W:=(\RP^2)^{k+1}$, let $i\colon\RP^2\imto\R^3$ be an immersion and let $j\colon W\imto\R^{3k+3}$ be a (self-transverse) immersion regularly homotopic to $(i\times\ldots\times i)$. It is well-known that $j$ has an odd number of triple points (see \cite{2n3n}) and it was proved in \cite{prim} that the number of cusps of a generic hyperplane projection of such an immersion has the same parity as the number of triple points of the immersion. Thus if $F\colon W\to\R^{3k+2}$ is such a projection of $j$, then it has on odd number of cusps.
\end{prf}

\begin{thm}\label{cob=n}
For any manifold $P^{3k+1}$ we have
$$\Cob_1(2k+1,P)\cong\NN_{2k+1}(P).$$
\end{thm}

\begin{prf}
It is enough to show that the forgetful map $\varphi\colon\Cob_1(2k+1,P)\to\NN_{2k+1}(P)$ is an isomorphism. It is clearly surjective, hence it is sufficient to prove that it is also injective.

Let $f\colon M^{2k+1}\to P$ be a map that represents $0$ in $\NN_{2k+1}(P)$ and let $F\colon W^{2k+2}\to P\times\R_+$ be a null-bordism ($\partial W=M$ and $F|_M=f$). If $F$ has $c$ cusps, then by the previous lemma we can find a closed manifold $W'$ and a map $F'\colon W'\to P\times\R_+$ that also has $c$ cusps. Now by adding $1$-handles we connect the cusps of $W$ and $W'$ and then lemma \ref{removecusps} yields a manifold $W''$ with boundary $\partial W''=M$ and a map $F''\colon W''\to P\times\R_+$ with no cusps and for which $F''|_M=f$. Hence the map $f$ represents $0$ in $\Cob_1(2k+1,P)$ as well.
\end{prf}

In the following we will investigate the same cobordism groups, but with oriented normal structures. First we consider the case of maps to a Eucledian space (i.e. $P^{3k+2}:=\R^{3k+2}$).

\begin{rmk}
The global normal form of the cusp singularity (see \cite{rsz}) implies that for odd codimensional cooriented maps, the submanifold of cusp points has oriented normal bundle. Thus a cooriented map $F\colon W^{4l}\to Q^{6l-1}$ is such that $\Sigma^{1,1}(F)$ is a cooriented $0$-manifold, i.e. a discrete set of points equipped with signs; we call the sum of these signs the algebraic number of cusps of $F$ and denote it by $\#\Sigma^{1,1}(F)$. In this case (for $k=2l-1$) lemma \ref{removecusps} can be reformulated (with the same proof) to also include that the points $p$ and $q$ have opposite signs and then the map $F'$ is cooriented as well.
\end{rmk}

\begin{lemma}\label{nrcusps}
$~$
\begin{enumerate}
\item\label{la1} If $W^{4l}$ is a closed oriented manifold and $F\colon W\to\R^{6l-1}$ is a (stable) map, then $\#\Sigma^{1,1}(F)$ is the normal Pontryagin number $\ol p_l[W]$. Moreover, for any integer of the form $c=|\ol p_l[W]|+2r$ (for any $r\in\N$) there is a map $F\colon W\to\R^{6l-1}$ with $c$ cusps.
\item\label{la2} If $W^6$ is a closed oriented manifold and $F\colon W\to\R^8$ is a (stable) map, then $F$ has an even number of cusps and for any even number $c\ge0$ there is a map $F\colon W\to\R^8$ with $c$ cusps.
\item\label{la3} For any integers $c\ge0$ and $l\ge2$ there is a closed oriented manifold $W^{4l+2}$ and a (stable) map $F\colon W\to\R^{6l+2}$ with $c$ cusps.
\end{enumerate}
\end{lemma}

\begin{prf}
As in the proof of lemma \ref{valamilyenlemma}, we can use lemma \ref{removecusps} to delete pairs of cusps (with opposite signs in the case of \ref{la1}). Hence it is sufficient to prove for \ref{la1} that $\#\Sigma^{1,1}(F)=\ol p_l[W]$, for \ref{la2} that any map has an even number of cusps and for \ref{la3} that there is a map with an odd number of cusps.

\medskip\noindent\emph{Proof of \ref{la1}.}\enspace\ignorespaces
We have two homomorphisms
$$\#\Sigma^{1,1}\colon\Omega_{4l}\to\Z~~~~\text{and}~~~~\ol p_l[\cdot]\colon\Omega_{4l}\to\Z,$$
where $\ol p_l[\cdot]$ is the normal Pontryagin number assigned to each cobordism class and $\#\Sigma^{1,1}$ is defined for any cobordism class $[W]\in\Omega_{4l}$ as $\#\Sigma^{1,1}([W]):=\#\Sigma^{1,1}(F)$ for a generic map $F\colon W\to\R^{6l-1}$. This is well-defined, since $\Omega_{4l}$ can be identified with $\Omega_{4l}(\R^{6l-1})$ and if two maps, $F_0\colon W_0^{4l}\to\R^{6l-1}$ and $F_1\colon W_1^{4l}\to\R^{6l-1}$, represent the same bordism class, then the cusp set of a generic bordism between them is an oriented $1$-manifold with boundary the union of the cusp sets of $F_0$ and $F_1$ with opposite signs, hence $F_0$ and $F_1$ have the same algebraic number of cusps. Now it is enough to prove that the maps $\#\Sigma^{1,1}$ and $\ol p_l[\cdot]$ coincide.

If we have $F=\pr_{\R^{6l-1}}\circ i$ for an immersion $i\colon W^{4l}\imto\R^{6l}$, then by \cite{prim} $\#\Sigma^{1,1}(F)$ is the number of triple points of $i$. This number is by Herbert's formula \cite{herbert} $\ol p_l[W]$, hence we have $\#\Sigma^{1,1}|_G=\ol p_l[\cdot]|_G$, where $G$ is the subgroup of elements $[W^{4l}]\in\Omega_{4l}$ for which there is a prim map $W\to\R^{6l-1}$.

Now a theorem of Burlet \cite{burlet} claims that any cobordism class in $\Omega_{4l}$ has a non-zero multiple that contains a map which can be immersed to $\R^{6l}$, hence the index of $G$ is finite. This implies that $\#\Sigma^{1,1}$ and $\ol p_l[\cdot]$ coincide on the whole group $\Omega_{4l}$.

\medskip\noindent\emph{Proof of \ref{la2}.}\enspace\ignorespaces
For any map $F\colon W^6\to\R^8$, there is a map $G\colon V^7\to\R^8\times\R_+$ from a compact orientable manifold with boundary $\partial V=W$ for which $G|_W=F$, since $\Omega_6$ is trivial. Now $\Sigma^{1,1}(G)$ is a $1$-manifold with boundary $\Sigma^{1,1}(F)$, hence $\Sigma^{1,1}(F)$ has to contain an even number of points.

\medskip\noindent\emph{Proof of \ref{la3}.}\enspace\ignorespaces
Consider the Dold manifold $Y=(\CP^2\times S^1)/\Z_2$, where $\Z_2$ acts on $\CP^2$ by complex conjugation and on $S^1$ by multiplication by $-1$. The manifold $(Y\times(\RP^2)^{l-2})^2$ is a square, hence by \cite{wallcob} there is an oriented manifold cobordant to it. Let $W^{4l+2}$ be such a manifold and let $F\colon W\to\R^{6l+2}$ be a generic map. Our aim is to prove that the number of cusps of $F$ is odd.

We will use the normal Stiefel--Whitney classes $\ol w_i$ of $W$. We know from \cite{kazthom} (also using \cite[6.2. théorème]{bh}) that the Poincaré dual (with $\Z_2$ coefficients) of $\Sigma^{1,1}(F)$ is
$$\Tp_{\Sigma^{1,1}}=\ol w_{2l+1}^2+\ol w_{2l}\smallsmile\ol w_{2l+2},$$
hence the number of cusps has the same parity as $(\ol w_{2l+1}^2+\ol w_{2l}\smallsmile\ol w_{2l+2})[W]$. Now putting $V:=Y\times(\RP^2)^{l-2}$, we obtain
\begin{alignat*}2
\ol w_{2l+1}^2[V\times V]&=\sum_{i_1+j_1=2l+1}\sum_{i_2+j_2=2l+1}(\ol w_{i_1}\times\ol w_{j_1})(\ol w_{i_2}\times\ol w_{j_2})[V\times V]=\\
&=\sum_{i_1+j_1=2l+1}\sum_{i_2+j_2=2l+1}(\ol w_{i_1}\smallsmile\ol w_{i_2})[V](\ol w_{j_1}\smallsmile\ol w_{j_2})[V]=\\
&=2\sum_{i=0}^l(\ol w_{i}\smallsmile\ol w_{2l+1-i})[V](\ol w_{2l+1-i}\smallsmile\ol w_{i})[V]=0
\end{alignat*}
and
\begin{alignat*}2
(\ol w_{2l}\smallsmile\ol w_{2l+2})[V\times V]&=\sum_{i_1+j_1=2l+2}\sum_{i_2+j_2=2l}(\ol w_{i_1}\times\ol w_{j_1})(\ol w_{i_2}\times\ol w_{j_2})[V\times V]=\\
&=\sum_{i_1+j_1=2l+2}\sum_{i_2+j_2=2l}(\ol w_{i_1}\smallsmile\ol w_{i_2})[V](\ol w_{j_1}\smallsmile\ol w_{j_2})[V]=\\
&=\sum_{i=0}^{2l}(\ol w_{i+1}\smallsmile\ol w_{2l-i})[V](\ol w_{2l+1-i}\smallsmile\ol w_{i})[V]=\\
&=((\ol w_{l+1}\smallsmile\ol w_{l})[V])^2=(\ol w_{l+1}\smallsmile\ol w_{l})[V].
\end{alignat*}

Now it is sufficient to prove that $(\ol w_{l+1}\smallsmile\ol w_{l})[V]=1$. We use that $(\ol w_i\cs\ol w_j)[\RP^2]=1$ is equivalent to $i=j=1$, and so
\begin{alignat*}2
(\ol w_{l+1}\smallsmile\ol w_{l})[V]&=\\
=&\sum_{i_0+\ldots+i_{l-2}=l+1\atop j_0+\ldots+j_{l-2}=l}(\ol w_{i_0}\times\ldots\times\ol w_{i_{l-2}})(\ol w_{j_0}\times\ldots\times\ol w_{j_{l-2}})[Y\times(\RP^2)^{l-2}]=\\
=&\sum_{i_0+\ldots+i_{l-2}=l+1\atop j_0+\ldots+j_{l-2}=l}(\ol w_{i_0}\cs\ol w_{j_0})[Y](\ol w_{i_1}\cs\ol w_{j_1})[\RP^2]\ldots(\ol w_{i_{l-2}}\cs\ol w_{j_{l-2}})[\RP^2]=\\
=&(\ol w_3\cs\ol w_2)[Y].
\end{alignat*}
Hence we only have to prove that $(\ol w_3\cs\ol w_2)[Y]=1$. Note that $Y$ is odd-dimensional, hence its top Stiefel--Whitney class vanishes, i.e. $w_5(Y)=0$; and $Y$ is orientable, hence $w_1(Y)=0$. Now an easy calculation implies the equalities $w_2(Y)=\ol w_2(Y)$ and $w_3(Y)=\ol w_3(Y)$. The cobordism class of $Y$ generates $\Omega_5\cong\NN_5\cong\Z_2$ (see \cite{charclass}), so it is not null-cobordant, thus at least one of its Stiefel--Whitney numbers is not $0$. The only possibility for such a Stiefel--Whitney number is $(w_2\cs w_3)[Y]$, which means $(\ol w_2\cs\ol w_3)[Y]=(w_2\cs w_3)[Y]=1$ and the proof is finished.
\end{prf}

\begin{thm}\label{foldthm}
$~$
\begin{enumerate}
\item\label{a1} If $k=2l-1$ is odd, then we have
$$\Cob_1^\SO(4l-1,2l-1)\cong\Omega_{4l-1}\oplus\Z_{3^t},$$
where $t:=\min\{n\in\N\mid\alpha_3(2l+n)\le3n\}$ and $\alpha_3(n)$ denotes the sum of digits of $n$ in triadic system.
\item\label{a2} If $k=2$, then we have
$$\Cob_1^\SO(5,2)\cong\Omega_5\oplus\Z_2\cong\Z_2\oplus\Z_2.$$
\item\label{a3} If $k=2l\ge4$ is even, then we have
$$\Cob_1^\SO(4l+1,2l)\cong\Omega_{4l+1}.$$
\end{enumerate}
\end{thm}

\begin{prf}
\noindent\emph{Proof of \ref{a1}.}\enspace\ignorespaces
Consider the forgetful homomorphism
$$\psi\colon\Cob_1^\SO(4l-1,2l-1)\to\Omega_{4l-1}.$$

We want to determine the kernel of $\psi$ and in order to do this, we let $g$ be the greatest common divisor of the numbers $\#\Sigma^{1,1}(F)$ where $F$ ranges over all (stable) maps $W^{4l}\to\R^{6l-1}$ for all oriented closed manifolds $W$. Part \ref{la1} of the previous lemma together with a result of Stong \cite{stong} implies that $g=3^t$. We define the map
$$\gamma\colon\ker\psi\to\Z_g$$
for any $[f]\in\ker\psi$ in the following way: If $f\colon M^{4l-1}\to\R^{6l-2}$ is any representative, then $M$ is oriented null-cobordant, that is, there is a compact oriented manifold $W^{4l}$ with boundary $\partial W=M$. Let $F\colon W\to\R^{6l-2}\times\R_+$ be any map extending $f$ and put $\gamma([f]):=\#\Sigma^{1,1}(F)~\mod g$.

\medskip\begin{sclaim}
The map $\gamma$ is well-defined.
\end{sclaim}

\begin{sprf}
If $F'\colon W'\to\R^{6l-2}\times\R_+$ is another choice of extension of $f$, then we take the map $\tilde F\cup F'$, where $\tilde F$ is the composition of $F$ with the reflection in $\R^{6l-1}$ to the hyperplane $\R^{6l-2}$. This maps the closed oriented manifold $(-W)\cup W'$ to $\R^{6l-1}$ and its algebraic number of cusps is $\#\Sigma^{1,1}(F')-\#\Sigma^{1,1}(F)\equiv0~\mod g$, thus $\gamma([f])$ is independent of the choice of $F$.

Of course $\gamma([f])$ is also independent of the choice of representative of the cobordism class $[f]$, since if $f'$ is another representative, then $f$ and $f'$ can be connected by a cobordism without cusps.
\end{sprf}

\begin{sclaim}
The map $\gamma$ is an isomorphism.
\end{sclaim}

\begin{sprf}
I. \emph{$\gamma$ is surjective.}\medskip

Using a neighbourhood of an isolated cusp point (see remark \ref{cusproot} and figure \ref{kep1}) it is easy to construct a map $f\colon S^{4l-1}\to\R^{6l-2}$ that has an extension $F\colon W^{4l}\to\R^{6l-2}\times\R_+$ with only one cusp point. Hence $\gamma([f])=1$.

\medskip\noindent II. \emph{$\gamma$ is injective.}\medskip

If $\gamma([f])=0$ for a map $f\colon M^{4l-1}\to\R^{6l-2}$, then $f$ has an extension $F\colon W^{4l}\to\R^{6l-2}\times\R_+$ with an algebraic number of cusp points $c$, where $c$ is divisible by $g$. Now lemma \ref{nrcusps} yields a map $F'\colon W'\to\R^{6l-1}$ from a closed oriented manifold with an algebraic number of cusp points $-c$. By adding $1$-handles, we can connect $W$ and $W'$ and since the codimension of the singular set is $2l>1$, we can connect the cusp points with opposite signs by arcs and use lemma \ref{removecusps} to remove all cusps. Thus we have a manifold $W''$ and a map $F''\colon W''\to\R^{6l-2}\times\R_+$ that extends $f$ and has no cusp points, so $f$ represents $0$ in $\Cob_1^\SO(4l-1,2l-1)$.
\end{sprf}

Hence we have a short exact sequence
$$0\to\Z_g\to\Cob_1^\SO(4l-1,2l-1)\to\Omega_{4l-1}\to0$$
and the statement of the theorem follows from the facts that $g$ is a power of $3$ and $\Omega_{4l-1}$ is a direct sum of copies of $\Z_2$.

\medskip\noindent\emph{Proof of \ref{a2}.}\enspace\ignorespaces
Applying the same reasoning as above in the proof of \ref{a1} and with $g=2$, we obtain a short exact sequence
$$0\to\Z_2\to\Cob_1^\SO(5,2)\xra\psi\Omega_5\to0.$$
Now $\Omega_5\cong\Z_2$ leaves the two possibilities $\Cob_1^\SO(5,2)\cong\Z_4$ or $\Cob_1^\SO(5,2)\cong\Z_2\oplus\Z_2$. If this group was $\Z_4$, then any cobordism class $[f]\in\Cob_1^\SO(5,2)\setminus\ker\psi$ was an element of order $4$. Hence it is enough to show a manifold $Y^5$ that is not null-cobordant and a map $f\colon Y\to\R^7$ for which $[f]$ has order $2$.

We put $Y:=(\CP^2\times S^1)/\Z_2$ (the Dold manifold we used in the proof of lemma \ref{nrcusps}) and fix an orientation of it. By Cohen's result \cite{cohen}, there is an immersion $i\colon Y^5\imto\R^8$; we define $f:=\pr_{\R^7}\circ i\colon Y\to\R^7$. In the following we prove that $[f]$ has order $2$ in $\Cob_1^\SO(5,2)$.

Recall the definition of the inverse element in subsection \ref{grp}: For any cobordism class $[g]$ (represented by the map $g\colon M^5\to\R^7$), the inverse $-[g]$ is represented by a map $[\rho\circ g]$ where $\rho$ is the reflection to a hyperplane in $\R^7$. Now if the maps are cooriented, which is equivalent to the source manifolds being oriented, then this inverse also includes taking the opposite orientation of the source. If we want to indicate this in the notation, then we put $[g,M]:=[g]$ and we have $[g,M]+[\rho\circ g,-M]=0$.

\medskip\begin{sclaim}
For any fold map $g\colon Y\to\R^7$ we have $[g,Y]+[\rho\circ g,Y]=0$.
\end{sclaim}

\begin{sprf}
Observe that the manifold $Y$ admits an orientation reversing diffeomorphism $\alpha\colon Y\to Y$ induced by the complex conjugation on $S^1$ (as the unit circle in $\C$). Let $X$ be the mapping torus of $\alpha$, that is,
$$X:=Y\times[0,1]/(y,0)\sim(\alpha(y),1),$$
which is a closed non-orientable $6$-manifold (this is the $6$-dimensional Wall manifold used in \cite{wallcob}). The Stiefel--Whitney classes of $X$ were computed in \cite{wallcob}, in particular we have $(\ol w_3^2+\ol w_2\cs\ol w_4)[X]=0$. Again (as in the proof of theorem \ref{nrcusps}) we use that $\Tp_{\Sigma^{1,1}}=\ol w_3^2+\ol w_2\cs\ol w_4$, which implies that the number of cusps of any map $F\colon X\to\R^8$ is even.

We map $Y\times\big[0,\frac12\big]$ to the path of the rotation of the image of $g$ to the image of $\rho\circ g$ in $\R^7\times\R_+$ in the same way as in the proof of proposition \ref{inv}. Extend this to a map $G\colon X\to\R^8$ by mapping $Y\times\big[\frac12,1\big]$ to the negative half-space (such that its boundary components $Y$ are mapped by $g$ and $\rho\circ g$). Since $G$ has no cusps in $Y\times\big[0,\frac12\big]\subset X$, it must have an even number of cusps in $Y\times\big[\frac12,1\big]\subset X$.

Since $\Omega_5\cong\Z_2$, there is a compact orientable manifold $W^6$ with boundary $\partial W=Y\sqcup Y$. Let $F\colon W\to\R^7\times\R_+$ be a map such that the restriction to one boundary component $Y$ is $g$ and to the other is $\rho\circ g$. Joining the manifold $W$ with $Y\times\big[\frac12,1\big]$ and the map $F$ with $G|_{Y\times[\frac12,1]}$ on the common boundary yields a map of an orientable $6$-manifold to $\R^8$.

\begin{figure}[H]
\begin{center}
\centering\includegraphics[scale=0.1]{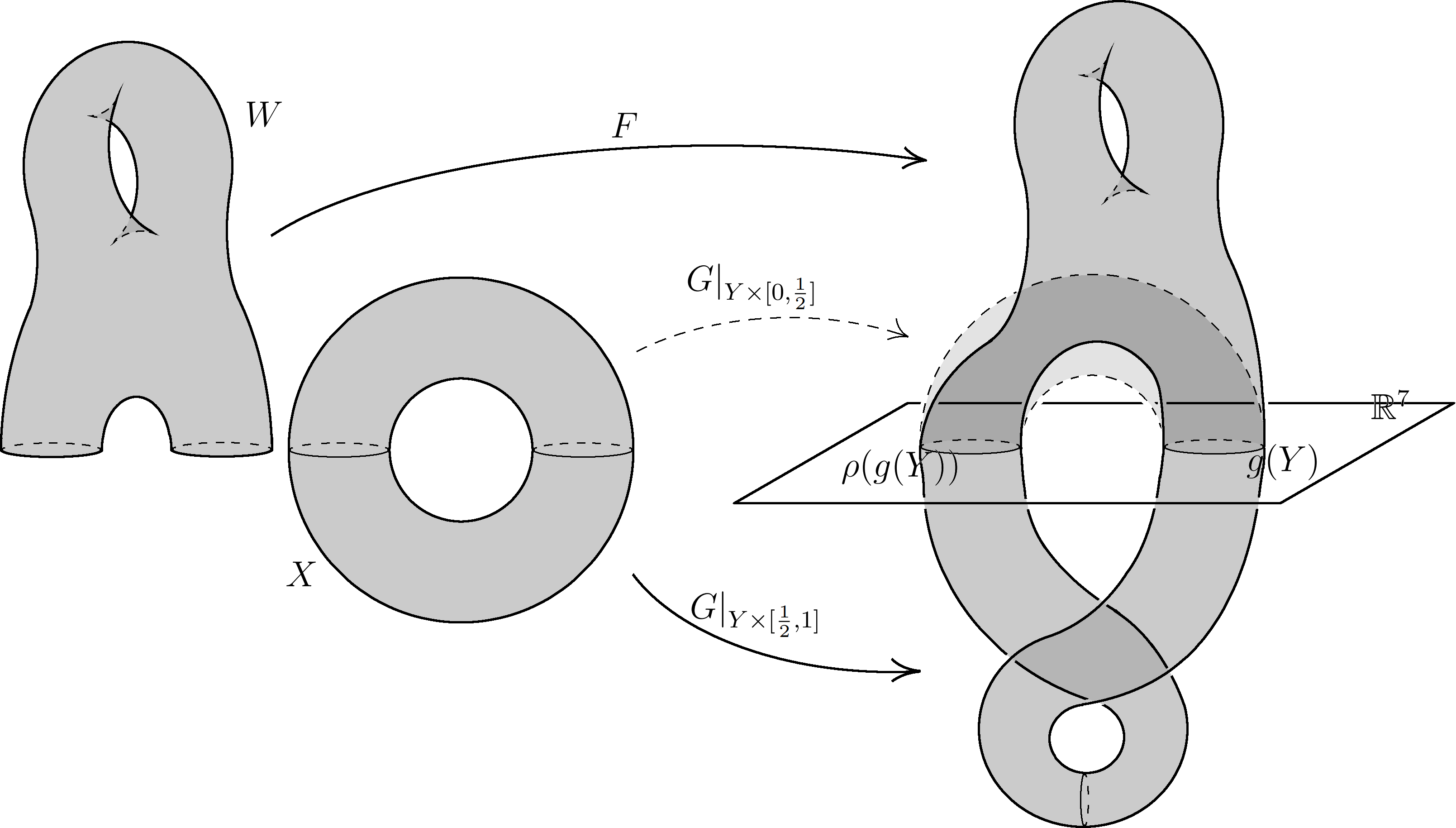}
\begin{changemargin}{2cm}{2cm} 
\caption{\hangindent=1.4cm\small We show the mapping torus $X$ mapped to $\R^8$ by $G$ and the manifold $W$ mapped to $\R^7\times\R_+$ by $F$. The images of $g$ and $\rho\circ g$ in $\R^7$ are represented by the two circles, they are rotated into each other by $G|_{Y\times[0,\frac12]}$ (indicated with dashed lines) and they are also connected by $G|_{Y\times[\frac12,1]}$ and $F$ in the two half-spaces.}\label{kep6}
\end{changemargin} 
\vspace{-1.3cm}
\end{center}
\end{figure}

This map has an even number of cusps by lemma \ref{nrcusps}, thus $F$ also has an even number of cusps. Now lemma \ref{removecusps} implies that these cusps can be eliminated, and then we have $[g,Y]+[\rho\circ g,Y]=0$.
\end{sprf}

\begin{sclaim}
For any prim map $f=\pr_{\R^7}\circ i$ we have $[f,Y]+[f,-Y]=0$.
\end{sclaim}

\begin{sprf}
We have
$$\Imm^\SO(5,3)\cong\pi^s_8(T\gamma^\SO_3)\cong\pi^s_9(ST\gamma^\SO_3),$$
where $\Imm^\SO(5,3)$ is the cobordism group of oriented immersions of $5$-manifolds to $\R^8$. The map on the cobordism group induced by changing the orientation of the source coincides with the map on the stable homotopy group induced by the involution $\iota\colon ST\gamma^\SO_3\to ST\gamma^\SO_3$ that is the reflection to $T\gamma^\SO_3$. The induced map $\iota_\#$ coincides with the multiplication by $-1$ in $\pi^s_9(ST\gamma^\SO_3)$, hence there is an oriented manifold $W^6$ with boundary $\partial W=Y\sqcup(-Y)$ and an immersion $I\colon W\to\R^8\times\R_+$ extending $i$ on both boundary components.

Define $\tilde W:=W\cup(Y\times[0,1])$ by attching the common boundary; hence $\tilde W$ is an orientable manifold. We can extend the normal bundle $\nu:=\nu_I$ to a bundle $\tilde\nu$ over $\tilde W$ in a canonical way (using the identity transition function where the attaching happened). Note that $T\tilde W\oplus\tilde\nu$ is the trivial bundle, so $\tilde\nu$ is a normal bundle for $\tilde W$.

Consider the constant vector field $\ua$ on $\R^9$ (the one we get from the standard basis vector on the last coordinate line) as a vector field on $I(W)$. Projecting $\ua$ to the normal bundle $\nu$ gives a section, which can be canonically continued to a section $s$ of $\tilde\nu$. Put $\Sigma:=s^{-1}(0)$, which is dual to $\ol w_3(\tilde W)$. Now $\ua$ also gives a vector field on $\Sigma$ by noting that the restriction of $\ua$ is always in $T\tilde W|_\Sigma$ and a point where $\ua$ is tangent to $\Sigma\cap W$ corresponds to a cusp point of $F:=\pr_{\R^8}\circ I$.

Generically $f=\pr_{\R^7}\circ i$ does not have any cusps, hence the section $\ua|_\Sigma$ is nowhere tangent to $\Sigma\cap(Y\times[0,1])$. Since the mod-$2$ number of tangency points is equal to the mod-$2$ self-intersection number $[\Sigma]\bullet[\Sigma]$ of $\Sigma\subset\tilde W$, this means that the number of cusps of $F$ has the same parity as $[\Sigma]\bullet[\Sigma]$. Now the Poincaré duality implies that $[\Sigma]\bullet[\Sigma]=\ol w_3^2(\tilde W)=0$, since $\tilde W$ is an orientable $6$-manifold and $\Omega_6\cong0$. We conclude that $F$ has an even number of cusps, which can be eliminated by lemma \ref{removecusps}, so it is a null-cobordism between $f$ on $Y$ and $f$ on $-Y$.
\end{sprf}

If we combine the above considerations, then we get
$$[f,Y]=-[\rho\circ f,-Y]=-[f,-Y]=-[f,Y],$$
thus $[f,Y]$ is indeed an element of order $2$, which is what we wanted to prove.

\medskip\noindent\emph{Proof of \ref{a3}.}\enspace\ignorespaces
This is a direct analogue of the proof of theorem \ref{cob=n} with the forgetful map $\psi$ in the place of $\varphi$.
\end{prf}

The following is the analogue of the previous theorem for more general target manifolds $P^{3k+2}$. The proof does not work for any target, hence there will be various restrictions on $P$.

\begin{thm}
$~$
\begin{enumerate}
\item\label{b1} If $k=2l-1$ is odd and $P^{6l-2}$ is an orientable manifold, then there is a short exact sequence
$$0\to\Z_{3^u}\to\Cob_1^\SO(4l-1,P)\to\Omega_{4l-1}(P)\to0,$$
where $u$ satisfies $1\le u\le t=\min\{n\in\N\mid\alpha_3(2l+n)\le3n\}$.
\item\label{b2} If $k=2$ and $P^7$ is an orientable manifold with $H_1(P;\Z_2)\cong H_2(P;\Z_2)\cong0$, then there is a short exact sequence
$$0\to\Z_2\to\Cob_1^\SO(5,P)\to\Omega_5(P)\to0.$$
\item\label{b3} If $k=2l\ge4$ is even and $P^{6l+1}$ is an orientable manifold, then we have
$$\Cob_1^\SO(4l+1,P)\cong\Omega_{4l+1}(P).$$
\end{enumerate}
\end{thm}

\begin{prf}
\noindent\emph{Proof of \ref{b1}.}\enspace\ignorespaces
Similarly to the proof of \ref{a1} in the previous theorem, we want to describe the kernel of the forgetful map
$$\psi\colon\Cob_1^\SO(4l-1,P)\to\Omega_{4l-1}(P).$$
If we define $g$ to be the greatest common divisor of the numbers $\la p_l(\nu_F),[W]\ra$, where $F$ (and $W$) ranges over all (stable) maps $F\colon W^{4l}\to P\times[0,1]$ for which $W$ is oriented. We have $g=3^u$ with some $1\le u\le t$ because of the same considerations as in the proof of \ref{a1} in the previous theorem. Again in exactly the same way as in the proof of \ref{a1} in the previous theorem, we have a well-defined isomorphism
$$\gamma\colon\ker\psi\to\Z_g$$
by assigning to any $[f]\in\ker\psi$ the algebraic number modulo $g$ of cusps of any map $F\colon W^{4l}\to P\times[0,1]$ extending $f$ (also using lemma \ref{nrcusps}). This yields the desired exact sequence.

\medskip\noindent\emph{Proof of \ref{b2}.}\enspace\ignorespaces
We first show the following.

\medskip\begin{sclaim}
For an orientable manifold $P^7$, the conditions $H_1(P;\Z_2)\cong H_2(P;\Z_2)\cong0$ and $\Hom(\Omega_6(P),\Z_2)\cong0$ are equivalent.
\end{sclaim}

\begin{sprf}
We know from \cite{conflo} that there is an isomorphism modulo odd torsion
$$\Omega_6(P)\congc{\rm{odd}}\bigoplus_{i=0}^6H_i(P;\Omega_{6-i})=H_1(P;\Z_2)\oplus H_2(P;\Z)\oplus H_6(P;\Z).$$
Hence we have to prove that the claimed homological condition is equivalent to all homomorphisms from the right-hand side to $\Z_2$ being trivial.

Write $H_1(P)=F\oplus T_{\rm{e}}\oplus T_{\rm{o}}$, where $F$ is the free part and $T_{\rm{e}}$ and $T_{\rm{o}}$ are respectively the even and odd torsion parts. The universal coefficient theorem implies that $H_1(P;\Z_2)=(F\oplus T_{\rm{e}})\otimes\Z_2$, so we have
$$\Hom(H_1(P;\Z_2),\Z_2)\cong0\iff H_1(P;\Z_2)\cong0\iff F\cong T_{\rm e}\cong0.$$
Again by the universal coefficient theorem and the Poincaré duality we get $H_6(P;\Z)\cong H^1(P;\Z)\cong F$, which means
$$\Hom(H_6(P;\Z),\Z_2)\cong0\iff F\cong0.$$
Finally, if $F\cong T_{\rm e}\cong0$, then yet again the universal coefficient theorem gives $H^2(P;\Z_2)\cong\Hom(H_2(P;\Z),\Z_2)$ and we obtain the equivalence
$$\Hom(H_2(P;\Z),\Z_2)\cong0\iff H_2(P;\Z_2)\cong H^2(P;\Z_2)\cong0.$$
This implies that the claim is indeed true.
\end{sprf}

Now we know that $\Hom(\Omega_6(P),\Z_2)=0$. The number of cusps modulo $2$ gives a homomorphism $\Omega_6(P)\to\Z_2$, since the cusps of generic bordisms are $1$-manifolds with boundary the cusps of the maps they connect. Thus the number of cusps of generic maps $W^6\to P\times[0,1]$ is always even (this is an analogue of part \ref{la2} in lemma \ref{nrcusps}).

Now in the same way as in the proof of \ref{a2} in theorem \ref{foldthm}, we get an isomorphism
$$\gamma\colon\ker(\psi\colon\Cob_1^\SO(5,P)\to\Omega_5(P))\to\Z_2$$
by counting the number of cusps modulo $2$ of any null-bordism of a given map and using lemma \ref{removecusps}. This yields the exact sequence as claimed.

\medskip\noindent\emph{Proof of \ref{b3}.}\enspace\ignorespaces
This is completely analogous to the proof of theorem \ref{cob=n} and to part \ref{a3} of theorem \ref{foldthm}.
\end{prf}

\subsection{The $(2k+2)$-dimensional case}

In this subsection we will use the key fibration corresponding to the singularity set $\tau_2=\{\Sigma^0,\Sigma^{1,0},\Sigma^{1,1,0}\}$ with a fixed codimension $k$, which is
$$\chi_2^L\colon X_2^L\xra{X_1^L}\Gamma T\tilde\xi_2^L$$
and using $L=\O$ and $L=\SO$. The long exact sequence in homotopy groups of this fibration is of the form
$$\ldots\to\pi_{n+1}^s(T\tilde\xi_2^L)\to\pi_n(X_1^L)\to\pi_n(X_2^L)\xra{(\chi_2^L)_\#}\pi_n^s(T\tilde\xi_2^L)\to\ldots$$
We are interested in the group $\pi_n(X_1^L)\cong\Cob_1^L(n-k,k)$ for $n=3k+2$. We start by proving lemmas on Thom spaces in general to understand the homotopy groups of $T\tilde\xi_2^L$. Recall from \cite{rsz} that $\tilde\xi_2^L$ is of dimension $3k+2$.

\begin{lemma}
If $\zeta$ is a vector bundle of dimension $n\ge1$ over a connected base space $B$, then $\pi_n(T\zeta)\cong\Z$ if $\zeta$ is orientable and $\pi_n(T\zeta)\cong\Z_2$ if $\zeta$ is non-orientable. Moreover, the map
$$\alpha\colon\pi_n(T\zeta)\to\begin{cases}
\Z;~[h]\mapsto\#|B\cap\im h|,&\text{ if }\zeta\text{ is orientable}\\
\Z_2;~[h]\mapsto|B\cap\im h|~\mod2,&\text{ if }\zeta\text{ is non-orientable}
\end{cases}$$
is an isomorphism (where we use such a representative of the homotopy class $[h]$ that is transverse to $B$; by $\#|B\cap\im h|$ we mean the number of intersection points taken with signs corresponding to an orientation of $\zeta$; this number is independent of the choice of the representative in both cases).
\end{lemma}

\begin{prf}
Since $T\zeta$ is $(n-1)$-connected, we have $\pi_n(T\zeta)\cong H_n(T\zeta)\cong H^n(T\zeta)$. This group is generated by the Thom class, which is a free generator if $\zeta$ is orientable and has order $2$ if it is non-orientable. Now the map $\alpha$ on any homotopy class $[h]\in\pi_n(T\zeta)$ is the evaluation of the Thom class on the image of $[h]$ under the Hurewicz homomorphism, hence it is an isomorphism.
\end{prf}

\begin{lemma}\label{pina}
If $\zeta$ is a vector bundle of dimension $n\ge3$ over a connected base space $B$, then the map
\begin{alignat*}2
\beta\colon\pi_{n+1}(T\zeta)&\to\begin{cases}
\Omega_1(B),&\text{ if }\zeta\text{ is orientable}\\
\{[f]\in\NN_1(B)\mid w_1(f^*\zeta)=0\},&\text{ if }\zeta\text{ is non-orientable}
\end{cases}\\
[h]&\mapsto[h|_{h^{-1}(B\cap\im h)}]
\end{alignat*}
is an epimorphism and its kernel is isomorphic to $\Z_2$ if $w_2(\zeta)$ vanishes and trivial if $w_2(\zeta)$ does not vanish.
\end{lemma}

\begin{prf}
We kill the homotopy group $\pi:=\pi_n(T\zeta)$ in the standard way, by pulling back the path fibration over $K(\pi,n)$ as shown in the square
$$\xymatrix{
K\ar[d]_{K(\pi,n-1)}\ar[r] & \ast\ar[d]^{K(\pi,n-1)} \\
T\zeta\ar[r]^(.4)u & K(\pi,n)
}$$
where the classifying map is the Thom class $u:=u(\zeta)$, which generates $\pi_n(T\zeta)\cong H^n(T\zeta)$. We will use the Serre spectral sequence to obtain the group $\pi_{n+1}(T\zeta)\cong\pi_{n+1}(K)\cong H_{n+1}(K)\cong H^{n+1}(K)$.

The only differentials in this spectral sequence that influence the group $H^{n+1}(K)$ are the transgressions $H^{n-1+i}(K(\pi,n-1))\to H^{n+i}(T\zeta)$ for $i=0,1,2$ (by dimensional reasons). For $i=0$ we know that the differential $H^{n-1}(K(\pi,n-1))\to H^n(T\zeta)$ is an isomorphism, since $H^{n-1}(K)\cong H^n(K)\cong0$. For $i=1$ and orientable $\zeta$, the group $H^n(K(\Z,n-1))$ is trivial; for non-orientable $\zeta$ it is $H^n(K(\Z_2,n-1))=\la\Sq^1\ra$ and the differential sends $\Sq^1$ to $\Sq^1(u)=u\cs w_1(\zeta)$, since transgressions commute with the Steenrod operations. For $i=2$ we have $H^{n+1}(K(\pi,n-1))=\la\Sq^2\ra$ and the image of $\Sq^2$ is $\Sq^2(u)=u\cs w_2(\zeta)$.

Combining these, we get that $H^{n+1}(K)$ is associated to the sum of $E_\infty^{n+1,0}=u\cs(H^1(B)/\la w_1(\zeta)\ra)$ and $E_\infty^{0,n+1}$, which is $0$ if $w_2(\zeta)$ is non-trivial and $\Z_2$ otherwise. We get the statement of the lemma by noting that $\Omega_1(B)\cong H_1(B)$ and $\{[f]\in\NN_1(B)\mid w_1(f^*\zeta)=0\}\cong\ker\la w_1(\zeta),\cdot\ra\subset H_1(B;\Z_2)$.
\end{prf}

\begin{crly}\label{pi2}
Let $G$ be a group with a fixed representation $\mu$ on $\R^n$ for an $n\ge3$ and put $\zeta:=EG\utimes\mu\R^n$. Then the map $\beta$ from the previous lemma is an isomorphism iff the image of $\mu_\#$ on the component of the neutral element of $G$ is the whole $\pi_1(\SO(n))\cong\Z_2$, that is, iff the image of $G$ contains a non-contractible loop.
\end{crly}

\begin{prf}
By the previous lemma it is enough to prove that $w_2(\zeta)$ does not vanish iff the image is $\pi_1(\SO(n))$. Now $w_2(\zeta)\ne0$ is equivalent to the existence of a map $h\colon S^2\to BG$ for which $w_2(h^*\zeta)=h^*(w_2(\zeta))\ne0$, and this latter condition means that the bundle $h^*\zeta$ is non-trivial.

The homotopy long exact sequence of the universal bundle $*\cong EG\xra GBG$ implies that we have an isomorphism $\pi_2(BG)\cong\pi_1(G)$ given by the boundary map $\partial$. For any homotopy class $[h]\in\pi_2(BG)$, the bundle $h^*\zeta$ over $S^2$ is in bijective correspondence with its gluing map, i.e. the difference between the trivialisations of the bundle over the two hemispheres, which is a map $S^1\to\SO(n)$ that maps to the image of $\mu_\#$. It is not hard to see that this map is $\mu_\#(\partial([h]))$. Now the existence of a map $h$ such that the bundle $h^*\zeta$ is non-trivial is equivalent to the existence of an $h$ such that $\mu_\#(\partial([h]))\in\Z_2$ is non-trivial and this is what we wanted to prove.
\end{prf}

Now we are ready to describe the segment of the homotopy long exact sequence of the key fibration containing $\Cob_1^L(2k+2,k)$, which can be seen as the short exact sequence
$$0\to\coker\chi_{3k+3}\to\Cob_1^L(2k+2,k)\to\ker\chi_{3k+2}\to0$$
using the notation $\chi_n:=(\chi_2^L)_\#\colon\pi_n(X_2^L)\to\pi^s_n(T\tilde\xi_2^L)$.

\begin{rmk}\label{chirmk}
The group $\ker\chi_{3k+2}$ has been computed in lemma \ref{valamilyenlemma} for $L=\O$ and in lemma \ref{nrcusps} for $L=\SO$: We have
$$\pi_{3k+2}(X_2^L)\cong\Cob_2^L(2k+2,k)\cong\begin{cases}
\Omega_{2k+2},&\text{ if }L=\SO\\
\NN_{2k+2},&\text{ if }L=\O
\end{cases}$$
and $\chi_{3k+2}$ maps the cobordism class of a generic map $f\colon M^{2k+2}\to\R^{3k+2}$ to the number $|\Sigma^{1,1}(f)|~\mod2$ if either $L=\O$ or $L=\SO$ and $k$ is even, and to the number $\#\Sigma^{1,1}(f)$ if $L=\SO$ and $k$ is odd. Now if $L=\O$, then lemma \ref{valamilyenlemma} implies that $\ker\chi_{3k+2}$ is an index-$2$ subgroup of $\NN_{2k+2}$; and if $L=\SO$, then lemma \ref{nrcusps} implies that $\ker\chi_{3k+2}$ is $\ker\ol p_l[\cdot]\subset\Omega_{4l}$ if $k=2l-1$ is odd, the whole $\Omega_6\cong0$ if $k=2$ and an index-$2$ subgroup of $\Omega_{2k+2}$ if $k\ge4$ is even. Hence it only remains to determine $\coker\chi_{3k+3}$.
\end{rmk}

\begin{lemma}\label{chilemma}
Consider the cases $L=\O$ and $L=\SO$ and the map $\beta$ (defined in lemma \ref{pina}) from $\pi_{3k+3}(T\tilde\xi_2^L)=\pi^s_{3k+3}(T\tilde\xi_2^L)$ to the respective bordism groups.
\begin{enumerate}
\item\label{beta1} $\coker\chi_{3k+3}\cong\coker(\beta\circ\chi_{3k+3})$.
\item\label{beta2} For any cobordism class $[f]\in\Cob_2^L(2k+3,k)\cong\pi_{3k+3}(X_2^L)$ represented by a map $f\colon M^{2k+3}\to\R^{3k+3}$, the image $\beta(\chi_{3k+3}([f]))$ depends only on the cobordism class of $M$ in $\Omega_{2k+3}$ or $\NN_{2k+3}$ respectively if $L=\SO$ or $L=\O$.
\end{enumerate}
\end{lemma}

\begin{prf}
\emph{Proof of \ref{beta1}.}\enspace\ignorespaces
Clearly it is enough to prove that $\beta$ is an isomorphism. Recall that the bundle $\tilde\xi_2^L$ is the universal vector bundle
$$EG_2^L\utimes{\tilde\lambda}\R^{3k+2},$$
where $\tilde\lambda$ is the target representation of $G_2^L$. We will check that the condition of corollary \ref{pi2} holds for $\tilde\xi_2^L$.

The group $G_2^L$ and its target representation are described in \cite{rsz} and \cite{rrthesis}; in the unoriented case we have $G_2\cong\O(1)\times\O(k)$ and in the cooriented case $G_2^\SO\cong\{(\varepsilon,A)\in\O(1)\times\O(k)\mid\varepsilon^2\det A>0\}\cong\O(1)\times\SO(k)$, the target representation is
$$\tilde\lambda(\varepsilon,A)=\varepsilon\oplus A\oplus 1\oplus A\oplus\varepsilon A$$
in both cases. The component of the unity is $\{1\}\times\SO(k)$ and if $\gamma$ is a non-contractible loop in $\SO(k)$, then clearly $\tilde\lambda(1,\cdot)\circ\gamma$ is a non-contractible loop in $\SO(3k+2)$. Thus corollary \ref{pi2} implies that $\beta$ is an isomorphism.

\medskip\noindent\emph{Proof of \ref{beta2}.}\enspace\ignorespaces
If $f\colon M^{2k+3}\to\R^{3k+3}$ represents an element $[f]\in\Cob_2^L(2k+3,k)$, then the element $\beta(\chi_{3k+3}([f]))$ in $\Omega_1(BG_2^L)$ or in $\NN_1(BG_2^L)$ (depending on $L$ and $k$) is represented by the inducing map $f(\Sigma^{1,1}(f))\to BG_2^L$ of the normal bundle of $f(\Sigma^{1,1}(f))$ (when we pull it back from $\tilde\xi_2^L$; see theorem \ref{univ}).

If we have an arbitrary cobordism of the manifold $M$ and generically map it to $\R^{3k+3}\times[0,1]$ (extending the map $f$ on the boundary part $M$ to $\R^{3k+3}\times\{0\}$ and mapping the rest of the boundary to $\R^{3k+3}\times\{1\}$), then this map will have isolated $\rm{III}_{2,2}$-points ($\rm{III}_{2,2}$ is the simplest $\Sigma^2$ singularity; see \cite{math}) apart from regular, fold and cusp points. The root of the singularity $\rm{III}_{2,2}$ (see \cite{rrthesis}) is such that the cusp points on its boundary form a circle with trivial normal bundle, hence a cobordism of the manifold $M$ yields the same $\beta(\chi_{3k+3}([f]))$. This is what we wanted to prove.
\end{prf}

\begin{thm}
There is a short exact sequence
$$0\to\Z_2\to\Cob_1(2k+2,k)\to G\to0,$$
where $G$ is an index-$2$ subgroup of $\NN_{2k+2}$.
\end{thm}

\begin{prf}
By remark \ref{chirmk} and lemma \ref{chilemma} we only have to prove that $\coker(\beta\circ\chi_{3k+3})$ is $\Z_2$. Since the elements of $\NN_{2k+3}$ are completely determined by their Stiefel--Whitney numbers, we want to express the map $\beta\circ\chi_{3k+3}$ in terms of these.

The image of $\beta$ can be identified with $\pi_1(BG_2)\cong\pi_1(B\O(1))\oplus\pi_1(B\O(k))\cong\Z_2\oplus\Z_2$. As we saw above, for a map $f\colon M^{2k+3}\to\R^{3k+3}$, the element $\beta(\chi_{3k+3}([f]))\in\pi_1(BG_2)$ is represented by the inducing map of the normal bundle of $f(\Sigma^{1,1}(f))$. Recall again that $G_2=\O(1)\times\O(k)$ was such that its source and target representations are
$$\lambda(\varepsilon,A)=\varepsilon\oplus 1\oplus A\oplus\varepsilon A~~~~\text{and}~~~~\tilde\lambda(\varepsilon,A)=\varepsilon\oplus A\oplus 1\oplus A\oplus\varepsilon A,$$
hence in $\pi_1(BG_2)\cong\Z_2\oplus\Z_2$ the projection to the first $\Z_2$ means the orientability of the kernel line bundle over each loop in $\Sigma^{1,1}(f)$ and the projection to the second $\Z_2$ is the orientability of the virtual normal bundle $\nu_f$ over each loop in $\Sigma^{1,1}(f)$.

Now we use \cite{kazthom} (together with \cite[6.2. théorème]{bh}) to obtain that the map $\beta\circ\chi_{3k+3}$ on a cobordism class $[M]\in\NN_{2k+3}$ is the evaluation of
$$(\ol w_{k+2}\cs\ol w_{k+1}+\ol w_{k+3}\cs\ol w_k,~\ol w_1\cs\ol w_{k+1}^2+\ol w_1\cs\ol w_k\cs\ol w_{k+2})$$
on $[M]$. The class $(\Sq^1+\ol w_1\cs)(\ol w_{k+1}^2+\ol w_k\cs\ol w_{k+2})$ is the second entry in this pair for $k$ odd and the sum of the two entries for $k$ even and it always evaluates to $0$ according to \cite{dold}. Hence the first entry determines $\im(\beta\circ\chi_{3k+3})$.

Similarly to the proof of part \ref{la3} of lemma \ref{nrcusps}, we use the Dold manifold $Y^5=(\CP^2\times S^1)/\Z_2$ multiplied by $(\RP^2)^{k-1}$. A computation analogous to the one for lemma \ref{nrcusps} yields that $\ol w_{k+2}\cs\ol w_{k+1}+\ol w_{k+3}\cs\ol w_k$ evaluates to $1$ on this manifold, hence we have $\im(\beta\circ\chi_{3k+3})\cong\Z_2$. Thus $\coker(\beta\circ\chi_{3k+3})$ is also $\Z_2$ and the statement of the theorem follows.
\end{prf}

\begin{thm}
$~$
\begin{enumerate}
\item\label{o1g} If $k=2l-1$ is odd, then $\Cob_1^\SO(4l,2l-1)$ is isomorphic to the kernel of the epimorphism $\ol p_l[\cdot]\colon\Omega_{4l}\twoheadrightarrow\Z$.
\item\label{o2g} If $k=2$, then we have $\Cob_1^\SO(6,2)\cong0$.
\item\label{o3g} If $k=2l\ge4$ is even, then $\Cob_1^\SO(4l+2,2l)$ is an index-$2$ subgroup of $\Omega_{4l+2}$.
\end{enumerate}
\end{thm}

\begin{prf}
By remark \ref{chirmk} and lemma \ref{chilemma} we only have to prove that $\coker(\beta\circ\chi_{3k+3})$ is trivial in all cases.

\medskip\noindent\emph{Proof of \ref{o1g}.}\enspace\ignorespaces
Now $\tilde\xi_2^\SO$ is orientable, hence $\beta$ maps to $\pi_1(BG_2^\SO)\cong\pi_1(B\O(1))\cong\Z_2$ and $\beta\circ\chi_{3k+3}$ can be identified with the evaluation of $\ol w_{k+2}\cs\ol w_{k+1}+\ol w_{k+3}\cs\ol w_k$. As in the proof above, we use that the manifold $Y^5\times(\RP^2)^{k-1}$ (which is cobordant to an orientable manifold by \cite{wallcob}) evaluates to $1$, and so $\beta\circ\chi_{3k+3}$ is surjective and its cokernel is trivial.

\medskip\noindent\emph{Proof of \ref{o2g} and \ref{o3g}.}\enspace\ignorespaces
Now $\tilde\xi_2^\SO$ changes orientation over all non-contractible loops in $BG_2^\SO$, hence the image of the isomorphism $\beta$ is trivial.
\end{prf}

%-------------------------------------------------------------------------------------------------------------------------------

\section{1-codimensional cobordisms}\label{1co}

Here we will consider cobordisms of cooriented Morin maps of codimension $1$ to Eucledian spaces, i.e. the groups $\Cob_r^\SO(n,1)$; the number $n$ will be fixed throughout this section. It will turn out that these groups are isomorphic to direct sums of stable homotopy groups of spheres modulo some torsion parts and in some cases a few of these torsion parts can also be described. First we will consider fold, then cusp and after that higher Morin maps.

\begin{rmk}
The analogue of this for unoriented Morin maps has been computed in theorem \ref{unorthm}, which includes as a special case that $\Cob_r(n,1)$ is finite $2$-primary for all $r\in\N\cup\{\infty\}$.
\end{rmk}

\begin{rmk}\label{morform}
We will need the global normal forms of $1$-codimensional Morin singularities. It follows from the computations for arbitrary codimensional singularities in \cite{rsz} and \cite{rrthesis} that we have $G_r^\SO\cong\Z_2$ for any positive integer $r$. The target representation $\tilde\lambda\colon\Z_2\to\O(2r+1)$ is such that its image is in $\SO(2r+1)$ iff $r$ is even, hence the universal bundle associated to it is
$$\tilde\xi_r^\SO=E\Z_2\utimes{\tilde\lambda}\R^{2r+1}=i\gamma_1\oplus\varepsilon^j$$
over $B\Z_2\cong\RP^\infty$ for some numbers for which $i+j=2r+1$ and $i\equiv r~\mod2$.

This implies that the Thom space $T\tilde\xi_r^\SO$ is $S^j(\RP^\infty/\RP^{i-1})$. It is not hard to see that for any odd prime $p$, the reduced cohomology $\tilde H^*(T\tilde\xi_r^\SO;\Z_p)$ isomorphic to $\tilde H^*(S^{2r+1};\Z_p)$ if $r$ is even (the isomorphism is induced by the inclusion $S^{2r+1}\subset T\tilde\xi_r^\SO$ as a ``fibre'') and it vanishes if $r$ is odd. Thus by Serre's mod-$\CC_2$ Whitehead theorem from \cite{serre}, the inclusion $\Gamma S^{2r+1}\subset\Gamma T\tilde\xi_r^\SO$ induces isomorphisms of the odd torsion parts of homotopy groups if $r$ is even and the homotopy groups of $\Gamma T\tilde\xi_r^\SO$ are finite $2$-primary if $r$ is odd.
\end{rmk}

\subsection{Fold maps}

The main tool in this subsection will be Koschorke's description of the Kahn--Priddy map from \cite{kos}, so first we recall this.

\begin{defi}
For any positive integer $m$, there are maps
\begin{itemize}
\item[(i)] $\RP^{m-1}\into\O(m)$ defined for a point $p\in\RP^{m-1}$, represented as a hyperplane in $\R^m$, by mapping $p$ to the reflection to this hyperplane,
\item[(ii)] $\O(m)\into\Omega^mS^m$ defined by mapping a point $A\in\O(m)$ to the element $s\in\Omega^mS^m$, represented as a map $S^m=\R^m\cup\{\infty\}\to\R^m\cup\{\infty\}=S^m$, for which $s(x)=A(x)$ if $x\in\R^m$ and $s(\infty)=\infty$.
\end{itemize}
The adjoint of the composition of these maps $\RP^{m-1}\into\O(m)\into\Omega^mS^m$ is a map $\lambda\colon S^m\RP^{m-1}\to S^m$. For any $n<m$ we have $\pi_{m+n}(S^m\RP^{m-1})\cong\pi^s_n(\RP^{m-1})\cong\pi^s_n(\RP^\infty)$ and $\pi_{m+n}(S^m)\cong\pi^s(n)$, hence $\lambda$ induces a homomorphism
$$\lambda_\#\colon\pi^s_n(\RP^\infty)\to\pi^s(n),$$
which is called the Kahn--Priddy map.
\end{defi}

\begin{rmk}
The stable homotopy groups $\pi^s_n(\RP^\infty)$ are $2$-primary (if $n\ge1$) and for small numbers $n$ they are completely computed; see \cite{liu}.
\end{rmk}

The following is due to Kahn and Priddy \cite{kahnpr} and we will not prove it here.

\begin{thm}\label{lambda}
The homomorphism $\lambda_\#\colon\pi^s_n(\RP^\infty)\to\pi^s(n)$ is onto the $2$-primary part of $\pi^s(n)$.
\end{thm}

Next we describe Koschorke's figure-$8$ construction that gives a very geometric understanding of the homomorphism $\lambda_\#$.

\begin{defi}
We define the figure-$8$ construction on an immersion $i\colon N^{n-1}\imto\R^n$ as follows: The composition of $i$ with the standard embedding $\R^n\into\R^{n+1}$ has normal bundle $\nu_i\oplus\varepsilon^1$; in each fibre of this bundle (considered as a plane in $\R^{n+1}$) we put a figure $8$ symmetrically to the $\nu_i$ factor (i.e. the image of $t\mapsto(\sin(2t),2\sin(t))$ with the first coordinate corresponding to $\nu_i$ and the second to $\varepsilon^1$). If we choose these figures $8$ smoothly, their union gives an immersion
$$8(i)\colon S(\nu_i\oplus\varepsilon^1)\imto\R^{n+1},$$
where the circle bundle $S(\nu_i\oplus\varepsilon^1)$ is an oriented $n$-manifold.
\end{defi}

\begin{rmk}
The unoriented cobordism group of immersions of $(n-1)$-manifolds to $\R^n$ is isomorphic to $\pi^s_n(\RP^\infty)$ and the oriented cobordism group of immersions of $n$-manifolds to $\R^{n+1}$ is isomorphic to $\pi^s(n)$. Clearly the figure-$8$ construction respects the cobordism relations (if $i_0\colon N_0^{n-1}\imto\R^n$ and $i_1\colon N_1^{n-1}\imto\R^n$ are unoriented cobordant, then $8(i_0)$ and $8(i_1)$ are oriented cobordant, since the figure-$8$ construction can be applied to cobordisms as well), hence we obtain a map
$$8_\#\colon\pi^s_n(\RP^\infty)\to\pi^s(n).$$
\end{rmk}

The following is a theorem of Koschorke \cite{kos} and again the proof will be omitted.

\begin{thm}\label{8}
The homomorphisms $\lambda_\#$ and $8_\#$ coincide.
\end{thm}

In the following we will see why the above theorems are improtant for the computation of cobordisms of fold maps. To do this, we use the key fibration for the singularity set $\{\Sigma^0,\Sigma^{1,0}\}$ that connects the classifying spaces $X_0^\SO$, $X_1^\SO$ and $\Gamma T\tilde\xi_1^\SO$ (for $1$-codimensional maps).

Now $X_0^\SO$ is the classifying space for $1$-codimensional cooriented immersions, which is clearly $\Gamma T\gamma_1^\SO=\Gamma S^1$; and $\tilde\xi_1^\SO$ is the bundle $\gamma_1\oplus\varepsilon^2$ over $\RP^\infty$, hence $\Gamma T\tilde\xi_1^\SO=\Gamma S^2\RP^\infty$. Thus the key fibration is
$$X_1^\SO\xra{\Gamma S^1}\Gamma S^2\RP^\infty,$$
and so its homotopy long exact sequence has the form
$$\ldots\xra\partial\pi^s(n)\to\Cob_1^\SO(n,1)\to\pi^s_{n-1}(\RP^\infty)\xra\partial\pi^s(n-1)\to\ldots$$
where we are using the isomorphisms $\pi_{n+1}(\Gamma S^1)\cong\pi^s(n)$ and $\pi_{n+1}(X_1^\SO)\cong\Cob_1^\SO(n,1)$ and $\pi_{n+1}(\Gamma S^2\RP^\infty)\cong\pi^s_{n+1}(S^2\RP^\infty)\cong\pi^s_{n-1}(\RP^\infty)$.

\begin{lemma}\label{d8}
The boundary map $\partial$ coincides with $8_\#$.
\end{lemma}

Before proving this lemma, we state (and prove) its most important corollary, which describes the cobordism groups $\Cob_1^\SO(n,1)$.

\begin{thm}\label{kpthm}
$\Cob_1^\SO(n,1)$ is a finite Abelian group with odd torsion part isomorphic to that of $\pi^s(n)$ and even torsion part isomorphic to the kernel of the Kahn--Priddy map $\lambda_\#\colon\pi^s_{n-1}(\RP^\infty)\twoheadrightarrow[\pi^s(n-1)]_2$. In other words, we have
$$\Cob_1^\SO(n,1)\cong\ker\lambda_\#\oplus\displaystyle\bigoplus_{p~\rm{odd}\atop\rm{prime}}[\pi^s(n)]_p.$$
\end{thm}

\begin{prf}
Trivial from theorems \ref{lambda} and \ref{8} and lemma \ref{d8}.
\end{prf}

\medskip\par\noindent\textbf{Proof of lemma \ref{d8}.\enspace\ignorespaces}
First observe that the homotopy group $\pi_{n+1}(\Gamma S^2\RP^\infty)\cong\pi_{n+1}(X_1^\SO,\Gamma S^1)$ is isomorphic to the cobordism group of cooriented fold maps
$$f\colon(M^n,\partial M^n)\to(D^{n+1},S^n)$$
for which $\partial f:=f|_{\partial M}$ is an immersion to $S^n$ (the cobordism relation and group structure can be defined analogously to section \ref{cobsec}). This $f$ will be a fixed representative throughout the proof. The boundary map
$$\partial\colon\pi_{n+1}(X_1^\SO,\Gamma S^1)\to\pi_n(\Gamma S^1)$$
maps this cobordism group to the cobordism group of cooriented immersions of $(n-1)$-manifolds to $S^n$ by assigning to the cobordism class of $[f]$ (as a cooriented fold map) the cobordism class of $[\partial f]$ (as a cooriented immersion).

Put $\Sigma:=\Sigma(f)=\Sigma^{1,0}(f)$, which is a $2$-codimensional submanifold of $M$. Note that $f|_\Sigma$ is an immersion and fix closed tubular nighbourhoods $T\subset M$ of $\Sigma$ and $\tilde T\subset D^{n+1}$ of $\tilde\Sigma:=f(\Sigma)$. Now theorem \ref{univ} (more precisely its analogue for multisingularities and cooriented maps) implies that there is a $D^3$-bundle $\hat T$ over $\Sigma$, an immersion $i\colon\hat T\imto D^{n+1}$ with image $\tilde T$ and a fibrewise map $\hat f\colon T\to\hat T$ such that $\hat f$ restricted to any fibre $D^2$ of $T$ is the Whitney umbrella
$$\sigma_1\colon(D^2,0)\to(D^3,0);~(x,y)\mapsto(x,xy,y^2)$$
(i.e. the normal form of $1$-codimensional fold singularity) and the following diagram is commutative:
$$\xymatrix@R=1.5pc{
T\ar[rr]^{f|_T}\ar[dr]^(.6){\hat f} && \tilde T\ar@{<-_)}[dd] \\
& \hat T\ar[dl]\ar[ur]^(.4)i &\\
\Sigma\ar@{_(->}[uu]\ar[rr]^{f|_\Sigma} && \tilde\Sigma
}$$

Put $N:=\ol{M\setminus T}$ and decompose its boundary as $\partial N=\partial_0N\sqcup\partial_1N$ with $\partial_0N=\partial M$ and $\partial_1N=\partial T$. Now the image of $f|_{\partial_2N}$ is the union of the images of the fibrewise Whitney umbrellas (composed with the immersion $i$) restricted to the boundary $S^1$ of each fibre. The image of such a restriction $\sigma_1|_{S^1}$ is a figure $8$ in a fibre $S^2$ of $\partial\hat T$.

\begin{figure}[H]
\begin{center}
\centering\includegraphics[scale=0.1]{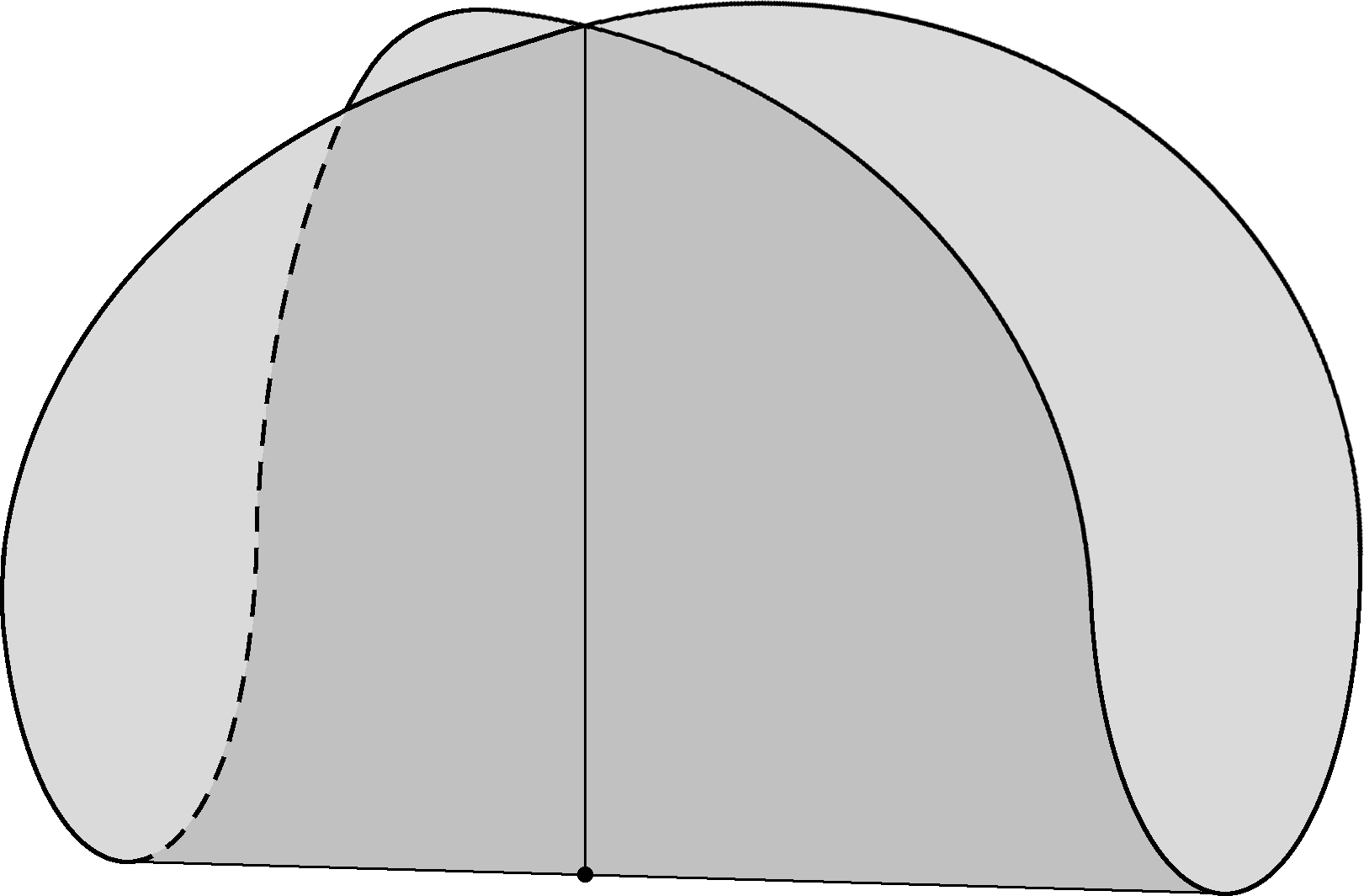}
\begin{changemargin}{2cm}{2cm} 
\caption{\hangindent=1.4cm\small Here we show the image of the Whitney umbrella; the indicated point is the origin and the curved figure $8$ is the image of the boundary.}\label{kep7}
\end{changemargin} 
\vspace{-1.5cm}
\end{center}
\end{figure}

Such a figure $8$ can canonically be endowed with a non-vanishing normal vector field in $S^2$ and with the inwards normal vector field in $D^3$. Thus the image of $f|_{\partial_2N}$ is a codimension-$2$ immersed submanifold in $D^{n+1}$ equipped with a normal framing.

We contract each figure $8$ in the corresponding $S^2$ to a small neighbourhood of its double point. We can identify the image of the singular set $\tilde\Sigma$ with the set of double points $\Delta$ of these figures $8$, hence $\Delta$ is also the image of an immersion with normal bundle induced from $\tilde\xi_1^\SO=\gamma_1\oplus2\varepsilon^1$. Over each point of $\Delta$, the first trivial line bundle $\varepsilon^1$ can be identified with the direction of the double point line of the Whitney umbrella; and the second $\varepsilon^1$ with the direction of the second coordinate in $t\mapsto(\sin(2t),2\sin(t))$ seen as the small figure $8$ around this point.

Of course we can view the immersion of $\Delta$ as a map to $\R^{n+1}$ (instead of $D^{n+1}$) and then the multi-compression theorem \ref{mct} (more precisely its version for immersions; see corollary \ref{ict}) can be applied to turn the two $\varepsilon^1$ directions in the normal bundle of $\Delta$ parallel to the last two coordinate lines. Projecting along these lines, we obtain an immersion $j$ to $\R^{n-1}$. If the figures $8$ were contracted to be small enough, then this $j$ extends to an immersion to $\R^n$ of a disk bundle over its source (induced from $D(\gamma_1\oplus\varepsilon^1)$) for which the image of each fibre contains the above described small figure $8$. Now the union of these figures $8$ is precisely the image of what we get by applying the figure-$8$ construction to the immersion $j$.

We obtained $j$ by a regular homotopy of the immersion of $\Delta$, which was regularly homotopic to that of $\tilde\Sigma$, i.e. the map $f|_\Sigma$. This implies that $j$ and $f|_\Sigma$ represent the same element of $\pi^s_{n-1}(\RP^\infty)\cong\pi^s_{n+1}(S^2\RP^\infty)$ noting that this suspension isomorphism is given between the cobordism groups of unoriented immersions to $\R^{n-1}$ and those to $\R^{n+1}$ with normal bundles induced from $\gamma_1\oplus\varepsilon^2$. Hence the immersion of the union of the figures $8$ on the boundary of $\tilde T$ (i.e. the map $f|_{\partial T}$) represents $8_\#([f|_\Sigma])$. This is the same as $8_\#([f])$, since $[f|_\Sigma]$ and $[f]$ correspond to each other in the isomorphism $\pi^s_{n+1}(S^2\RP^\infty)\cong\pi_{n+1}(X_1^\SO,\Gamma S^1)$.

Now the only thing left to prove is that $f|_{\partial T}\colon\partial T\imto D^{n+1}$, equipped with the normal framing of the figures $8$ mentioned above, represents the same element of $\pi^s(n-1)$ as the immersion $\partial f\colon\partial M\imto S^n$. Here $f|_{\partial T}$ is a $2$-codimensional framed immersion and $\partial f$ is just $1$-codimensional, so the first task is to turn $f|_{\partial T}$ $1$-codimensional as well.

Take a small neighbourhood of a point in $\tightoverset{_\circ~~~~}{D^{n+1}}$ disjoint from the image of $f$ and delete it from the disk. We obtain a space diffeomorphic to $S^n\times[0,1]$ and $f|_{\partial T}$ is a framed immersion to this. Now again the immersion version of the multi-compression theorem \ref{mct} yields a regular homotopy of $f|_{\partial T}$ that turns one of the framing vector fields parallel to $T[0,1]$. By combining this with the projection to $S^n\times\{1\}$, we get a $1$-codimensional immersion $\partial T\imto S^n$. This deformation can be identified with the suspension isomorphism for the stable homotopy groups of spheres, hence this immersion represents the same element of $\pi^s(n-1)$ as $f|_{\partial T}$.

But the above deformation can be extended to the manifold $N=\ol{M\setminus T}$ connecting $\partial T$ and $\partial M$, hence we get an immersion $N\imto S^n\times[0,1]$ which maps $\partial_0N$ to $S^n\times\{0\}$ by $\partial f$ and $\partial_1N$ to $S^n\times\{1\}$ by $f|_{\partial T}$ (up to regular homotopy). Hence we got a cobordism connecting $\partial f$ and $f|_{\partial T}$, which finishes the proof.
~$\square$\par\medskip

\begin{rmk}
Liulevicius \cite{liu} computed the stable homotopy groups of $\RP^\infty$ in dimensions at most $9$, which is shown below with the corresponding stable homotopy groups of spheres and the groups $\Cob_1^\SO(n,1)$ computed from these.
\begin{table}[H]\begin{center}\begin{tabular}{c||c|c|c|c|c|c|c|c|c|c}
$n$ & $1$ & $2$ & $3$ & $4$ & $5$ & $6$ & $7$ & $8$ & $9$ & $10$\\
\hline
$\pi^s_{n-1}(\RP^\infty)$ & $\Z$ & $\Z_2$ & $\Z_2$ & $\Z_8$ & $\Z_2$ & $0$ & $\Z_2$ & $\Z_{16}\oplus\Z_2$ & $(\Z_2)^3$ & $(\Z_2)^4$\\
\hline
$\pi^s(n)$ & $\Z_2$ & $\Z_2$ & $\Z_{24}$ & $0$ & $0$ & $\Z_2$ & $\Z_{240}$ & $(\Z_2)^2$ & $(\Z_2)^3$ & $\Z_6$\\
\hline
$\Cob_1^\SO(n,1)$ & $0$ & $0$ & $\Z_3$ & $0$ & $\Z_2$ & $0$ & $\Z_{15}$ & $\Z_2$ & $\Z_2$ & $\Z_6$
\end{tabular}\end{center}\end{table}\vspace{-.5cm}
\end{rmk}

\begin{rmk}
We also proved that the odd torsion parts of the groups $\Cob_1^\SO(n,1)$ can be represented by immersions. In particular for $n=3$ and $n=7$, these can be chosen to be immersions of $S^3$ and $S^7$ respectively, since in these dimensions the $J$ homomorphism is epi.
\end{rmk}

\subsection{Cusp maps}

Now the main tool will be the classifying space for the cobordisms of those coorinted codimension-$1$ cusp maps that are equipped with a trivialisation of the normal bundle of the immersion of the cusp stratum. This can be constructed completely analogously to the classifying spaces in chapter \ref{classp} and will be denoted by $\tilde X_2^\SO$.

By remark \ref{morform} we know that the inclusion $\Gamma S^5\subset\Gamma T\tilde\xi_2^\SO$ is a mod-$\CC_2$ homotopy equivalence. Consider the pullback of the key fibration by this inclusion
$$\xymatrix{
\tilde X_2^\SO\ar[r]\ar[d]_{X_1^\SO} & X_2^\SO\ar[d]^{X_1^\SO} \\
\Gamma S^5\ar@{^(->}[r] & \Gamma T\tilde\xi_2^\SO
}$$
The horizontal arrows here induce isomorphisms in the odd torsion parts of the homotopy groups; we will see that the left-hand fibration ``almost'' has a splitting.

\begin{defi}
Let $\pi\colon E\xra FB$ be a fibration and $l$ an integer.
\begin{enumerate}
\item We say that $\pi$ has an algebraic $l$-splitting, if for all $i\in\N$ there is a homomorphism $s_i\colon\pi_i(B)\to\pi_i(E)$ such that $\pi_\#\circ s_i$ is the multiplication by $l$ in the group $\pi_i(B)$.
\item We say that $\pi$ has a geometric $l$-splitting, if it has an algebraic one for which all homomorphisms $s_i$ are induced by a map $s\colon B\to E$.
\end{enumerate}
\end{defi}

\begin{lemma}\label{lspl}
If $l\in\Z$ is such that there is a map $f\colon M^4\to\R^5$ with algebraic number of cusps $\#\Sigma^{1,1}(f)=l$, then the fibration $\tilde X_2^\SO\xra{X_1^\SO}\Gamma S^5$ has an algebraic $l$-splitting.
\end{lemma}

\begin{prf}
Let $i\ge5$ be an integer and fix an element $[j]\in\pi_i(\Gamma S^5)\cong\pi^s(i-5)$ represented by an immersion $j\colon N^{i-5}\imto\R^i$ endowed with a normal framing. We may assume that $f$ maps to the disk $D^5$, and then the map
$${\id}_N\times f\colon N\times M\to N\times D^5$$
can be composed with an immersion of $N\times D^5$ to $\R^i$ onto a tubular neighbourhood of $j(N)$ using the normal framing. The composition of $\id_N\times f$ with this immersion is clearly a cusp map to $\R^i$ with framed cusp stratum, hence it represents an element in $\pi_i(\tilde X_2^\SO)$, and its restriction to the cusp stratum is $l\cdot[j]\in\pi^s(i-5)$. This defines the homomorphism $s_i$ and shows that it is an $l$-splitting.
\end{prf}

\begin{rmk}
The proof can be extended to also show that this $l$-splitting is actually geometric, but we will not need this.
\end{rmk}

Now the groups $\Cob_2^\SO(n,1)$ can be obtained modulo their $2$- and $3$-primary parts. Moreover, the $3$-primary parts can also be described up to a group extension; here we will use that the generator $\alpha_1\in\pi^s(3)\cong\Z_3\oplus\Z_8$ of the $\Z_3$ part defines a homomorphism of order $3$
$$(\alpha_1)_i\colon\pi^s(i)\to\pi^s(i+3);~[b]\mapsto[a\circ b]$$
(for all $i$), where $[b]$ is represented by $b\colon S^{m+i+3}\to S^{m+3}$ (for some large number $m$) and the map $a\colon S^{m+3}\to S^m$ represents $\alpha_1$.

\begin{thm}\label{cuspthm}
$\Cob_2^\SO(n,1)$ is such that there is a $\CC_{2,3}$-isomorphism
$$\Cob_2^\SO(n,1)\congc{2,3}\pi^s(n)\oplus\pi^s(n-4)$$
and a $\CC_2$-exact sequence
$$0\to\coker(\alpha_1)_{n-3}\to\Cob_2^\SO(n,1)\to\ker(\alpha_1)_{n-4}\to0.$$
\end{thm}

\begin{prf}
We will not prove the existence of the $\CC_2$-exact sequence here, as it relies on a spectral sequence for prim maps that will only be computed in the next section. For the proof we instead refer to remark \ref{c2ex}.

The $\CC_{2,3}$-isomorphism can be obtained from a sequence of isomorphisms
\begin{alignat*}2
\Cob_2^\SO(n,1)&\cong\pi_{n+1}(X_2^\SO)\toverset{^*}{\congc2}\pi_{n+1}(\tilde X_2^\SO)\toverset{^{**}}{\congc{2,3}}\\
&\toverset{^{**}}{\congc{2,3}}\pi_{n+1}(\Gamma S^5)\oplus\pi_{n+1}(X_1^\SO)\toverset{^{***}}{\congc2}\pi^s(n-4)\oplus\pi^s(n),
\end{alignat*}
where $^*$ follows from the observation in the beginning of this subsection, $^{***}$ is by theorem \ref{kpthm} and $^{**}$ will follow from lemma \ref{lspl} as soon as we prove that there is a map $f\colon M^4\to\R^5$ with algebraically $6$ cusps. This follows from a theorem of Eccles and Mitchell \cite{em} claiming that there is an immersion $i\colon M^4\imto\R^6$ with algebraically $2$ triple points, that is, $6$ such points in the source; now by \cite{prim} this equals the algebraic number of cusp points of $\pr_{\R^5}\circ i$.
\end{prf}

\begin{rmk}
Since $(\alpha_1)_i$ is a homomorphism of order $3$, the $\CC_2$-exact sequence above describes the $3$-primary part of the group $\Cob_2^\SO(n,1)$. As a consequence of Toda's work \cite{toda}, we know that for $i\le8$, the homomorphism $(\alpha_1)_i$ is non-trivial only for $i=0$. Hence for $3\ne n\le11$, we have a short exact sequence of $3$-primary parts
$$0\to[\pi^s(n)]_3\to\big[\Cob_2^\SO(n,1)\big]_3\to[\pi^s(n-4)]_3\to0.$$
\end{rmk}

\subsection{Higher Morin maps}

The computation in this subsection will be quite similar to the previous one. Namely we will use the $l$-splitting (for some $l$) of the bundle
$$\tilde X_{2r}^\SO\xra{X_{2r-1}^\SO}\Gamma S^{4r+1}$$
for all $r\ge1$, where $\tilde X_{2r}^\SO$ is the classifying space for the cobordisms of those coorinted codimension-$1$ $\Sigma^{1_{2r}}$-maps that are equipped with a normal framing of the immersion of the $\Sigma^{1_{2r}}$-stratum. Again we will obtain a $\CC_{\le2r+1}$-isomorphism between the cobordism groups and direct sums of certain stable homotopy groups of spheres. Before precisely stating and proving the main theorem, we need the following weak analogue of the Eccles--Mitchell theorem used above.

\begin{lemma}\label{lr}
Let $r$ be a natural number and $l(r)$ be the order of the cokernel of the stable Hurewicz homomorphism
$$\pi^s_{4r+2}(\CP^\infty)\to H_{4r+2}(\CP^\infty).$$
Then for any immersion of an oriented, closed $4r$-manifold to $\R^{4r+2}$, the algebraic number of $(2r+1)$-tuple points is divisible by $l(r)$, and there is an immersion for which this number is precisely $l(r)$.
\end{lemma}

\begin{prf}
By Herbert's formula \cite{herbert}, the algebraic number of $(2r+1)$-tuple points of an immersion $i\colon M^{4r}\imto\R^{4r+2}$ is the normal Pontryagin number $\ol p_1^r[M]$. Now $i$ represents an element $[i]$ in the cobordism group of oriented codimension-$2$ immersions to $\R^{4r+2}$, which is isomorphic to $\pi^s_{4r+2}(T\gamma_2^\SO)$. If the Pontryagin--Thom map of $i$ is $\alpha:=\alpha(i)\colon S^{4r+2}\nrightarrow T\gamma_2^\SO$, then the cobordism class $[i]$ corresponds to the homotopy class $[\alpha]$ by the above isomorphism.

Consider the composition of the maps
$$\pi^s_{4r+2}(T\gamma_2^\SO)\xra hH_{4r+2}(T\gamma_2^\SO)\xra\varphi H_{4r}(B\SO(2))\xra{\la p_1^r(\gamma_2^\SO),\cdot\ra}\Z,$$
where $h$ is the stable Hurewicz homomorphism and $\varphi$ is the homological Thom isomorphism. Since $\varphi$ and the evaluation map $x\mapsto\la p_1^r(\gamma_2^\SO),x\ra$ are isomorphisms, the cokernel of this composition is the same as $\coker h$, hence its image is $l(r)\Z$ (recall that the stable Hurewicz map is a rational isomorphism). Now it is enough to prove that $[\alpha]\in\pi^s_{4r+2}(T\gamma_2\SO)$ gets mapped to $\ol p_1^r[M]\in\Z$.

Note that $\alpha$ is a stable map, i.e. a map $S^{4r+2+m}\to S^mT\gamma_2^\SO$ for a sufficiently large number $m$. It is by definition given as a composition $\beta\circ q$, where $q\colon S^{4r+2+m}\to S^mT\nu_i$ is the quotient by the complement of a tubular neighbourhood and $\beta\colon S^mT\nu_i\to S^mT\gamma_2^\SO$ is a fibrewise map that restricts to $M$ as the inducing map $b\colon M\to B\SO(2)$ of the normal bundle $\nu_i$. Now using the notations $p_1:=p_1(\gamma_2^\SO)$ and $u:=u(\gamma_2^\SO)$, we have
\begin{alignat*}2
\la p_1^r,\varphi h([\alpha])\ra&=\la p_1^r,\varphi\alpha_*([S^{4r+2+m}])\ra=\la p_1^r,\varphi\beta_*q_*([S^{4r+2+m}])\ra=\la p_1^r,\varphi\beta_*([S^mT\nu_i])\ra=\\
&=\la p_1^r,\beta_*([S^mT\nu_i])\smallfrown S^mu\ra=\la p_1^r,b_*([M])\ra=\la b^*(p_1^r),[M]\ra=\\
&=\la p_1^r(\nu_i),[M]\ra=\ol p_1^r[M]
\end{alignat*}
and this is what we wanted to prove.
\end{prf}

\begin{thm}\label{hmor}
For all $r\ge1$ we have
$$\Cob_{2r+1}^\SO(n,1)\congc2\Cob_{2r}^\SO(n,1)\congc{\le2r+1}\displaystyle\bigoplus_{i=0}^r\pi^s(n-4i).$$
\end{thm}

\begin{prf}
By remark \ref{morform} we know that in the homotopy long exact sequence of the key fibration
$$X_{2r+1}^\SO\xra{X_{2r}^\SO}\Gamma\tilde\xi_{2r+1}^\SO$$
the homotopy groups of $\Gamma T\tilde\xi_{2r+1}^\SO$ vanish modulo $\CC_2$, hence the first $\CC_2$-isomorphism is proved.

The rest is completely analogous to the proof of theorem \ref{cuspthm}: We have a sequence of isomorphisms
\begin{alignat*}2
\Cob_{2r}^\SO(n,1)&\cong\pi_{n+1}(X_{2r}^\SO)\toverset{^*}{\congc2}\pi_{n+1}(\tilde X_{2r}^\SO)\toverset{^{**}}{\congc{\le2r+1}}\pi_{n+1}(\Gamma S^{4r+1})\oplus\pi_{n+1}(X_{2r-1}^\SO)\toverset{^{***}}{\congc{\le2r-1}}\\
&\toverset{^{***}}{\congc{\le2r-1}}\pi^s(n-4r)\oplus\displaystyle\bigoplus_{i=0}^{r-1}\pi^s(n-4i),
\end{alignat*}
where $^*$ follows again from remark \ref{morform} and $^{***}$ is obtained by an induction on $r$. To prove $^{**}$, we use lemma \ref{lr} and a direct analogue of lemma \ref{lspl} that yields an algebraic $l(r)$-splitting of the bundle $\tilde X_{2r}^\SO\xra{X_{2r-1}^\SO}\Gamma S^{4r+1}$. It remains to note that we can use Arlettaz's theorem \ref{arlthm} with the substitutions $V:=\CP^\infty$, $l:=2$, $m:=4r+2$ and get that $l(r)$ is a divisor of $\rho_1\ldots\rho_{4r-1}$, where $\rho_i$ is the exponent of $\pi^s(i)$. But Serre proved in \cite{serre} that $\rho_i$ is not divisible by the prime $p$ if $p>\frac i2+1$, which means that $l(r)$ does not have prime divisors greater than $2r+1$, hence we have the $\CC_{\le2r+1}$-isomorphism as claimed.
\end{prf}

%-------------------------------------------------------------------------------------------------------------------------------

\section{Prim cobordisms}\label{prico}

In the previous section we saw that the cobordism groups of $1$-codimensional Morin maps are closely related to the stable homotopy groups of spheres. In this section we establish even closer such connections; the main tool will be the so-called singularity spectral sequence for prim maps defined as follows.

\begin{defi}\label{sss}
Let $L$ be a stable linear group, $k$ a positive integer and $r\in\N\cup\{\infty\}$. The filtration $\ol X_0^L\subset\ol X_1^L\subset\ldots\subset\ol X_r^L$ of the classifying space for the cobordisms of $k$-codimensional prim maps defines a spectral sequence in homotopy groups with first page
$$\ol E^1_{i,j}=\pi_{i+j+k}(\ol X_i^L,\ol X_{i-1}^L)\cong\pi_{i+j+k}(\Gamma T\tilde\zeta_r^L)\cong\pi^s_{i+j+k}(T\tilde\zeta_r^L).$$
This spectral sequence converges to $\pi_{*+k}(\ol X_r^L)\cong\Prim_r^L(*,k)$, that is, $\underset{i+j=n}\bigoplus\ol E^\infty_{i,j}$ is associated to $\Prim_r^L(n,k)$.
\end{defi}

Similarly to the previous section, we will first consider the case of cooriented $1$-codimensional maps. We give a geometric description of the singularity spectral sequence in this case and compute it (modulo some torsion) for cooriented codimension-$1$ prim fold and cusp maps. After this we will see that in some cases (namely for $1$-codimensional cooriented prim maps and $3$-codimensional prim maps with spin normal structures) the singularity spectral sequence can be identified with spectral sequences arising from filtrations of projective spaces and starting with stable homotopy groups of spheres.

\subsection{Geometric description}\label{geode}

This subsection mainly contains remarks about how we can translate the singularity spectral sequence to cobordisms and transformations between them. In the following we consider prim maps of codimension $1$ and with oriented virtual normal bundles.

\begin{rmk}
$~$
\begin{enumerate}
\item For this sort of maps, the universal target bundle in the global normal form of $\Sigma^{1_r}$ singularity is easy to describe, namely it is $\tilde\zeta_r^\SO=(r+1)\gamma_1^\SO\oplus\varepsilon^r=\varepsilon^{2r+1}$, hence its Thom space is $T\tilde\zeta_r^\SO=S^{2r+1}$ and the key fibration has the form
$$\ol X_r^\SO\xra{\ol X_{r-1}^\SO}\Gamma S^{2r+1}.$$
This also implies that the starting page of the singularity spectral sequence is
$$\ol E^1_{i,j}=\pi_{i+j+1}(\ol X_i^\SO,\ol X_{i-1}^\SO)\cong\pi^s_{i+j+1}(S^{2i+1})=\pi^s(j-i).$$
\item It is not hard to see that a relative homotopy group $\pi_n(\ol X_r^\SO,\ol X_{r-1}^\SO)$ can be identified with a relative cobordism group, namely the cobordism group of such prim maps
$$f\colon(M^n,\partial M^n)\to(\R^n\times\R_+,\R^n)$$
for which $f$ is a cooriented $\tau_r$-map and $\partial f:=f|_{\partial M}$ has at most $\Sigma^{1_{r-1}}$ singularities. The cobordism of two such maps can be defined directly analogously to definition \ref{cobtau} (by requiring the cobordism to be a prim $(\tau_{r-1},\SO)$-cobordism between the boundaries and a prim $(\tau_r,\SO)$-cobordism altogether) and a complete analogue of Szűcs's Pontryagin--Thom type construction from chapter \ref{classp} gives that this cobordism group is indeed isomorphic to $\pi_n(\ol X_r^\SO,\ol X_{r-1}^\SO)$.
\end{enumerate}
\end{rmk}

Recall that for any fibration $\pi\colon E\xra FB$, the isomorphism $\pi_n(E,F)\cong\pi_n(B)$ is induced by the map $\pi$. Now the two observations in the remark above yield together a nice geometric description of the isomorphism $\pi_{n+1}(\ol X_r^\SO,\ol X_{r-1}^\SO)\cong\pi_{n+1}(\Gamma S^{2r+1})\cong\pi^s(n-2r)$. This is induced by the key fibration, so it assigns to the  cobordism class of a map
$$f\colon(M^n,\partial M^n)\to(\R^n\times\R_+,\R^n)$$
the cobordism class of its restriction to the top singularity stratum (which is disjoint from $\partial M$), that is, the map
$$f|_{\Sigma^{1_r}(f)}\colon\Sigma^{1_r}(f)\imto\R^{n+1}.$$

More precisely, this latter cobordism class is one in $\Imm^{\tilde\zeta_r^\SO}(\R^{n+1})$ (see remark \ref{keyrmk} or theorem \ref{terpthm}), which in the present case is the cobordism group of immersions of $(n-2r)$-dimensional manifolds with normal framings, i.e. the stable homotopy group of spheres $\pi^s(n-2r)$.

\begin{rmk}\label{relcob}
The above considerations also mean that the relative cobordism class of such a map $f$ is completely determined by the cobordism class of its restriction to $\Sigma^{1_r}(f)$.
\end{rmk}

Now that we know the geometric meaning of the groups in the first page of the singularity spectral sequence, we are ready to describe the differentials and what transformations they stand for in the cobordism groups. First consider the differential $d^1$, which has the form
$$d^1_{i,j}\colon\pi_{i+j+1}(\ol X_i^\SO,\ol X_{i-1}^\SO)\cong\ol E^1_{i,j}\to\ol E^1_{i-1,j}\cong\pi_{i+j}(\ol X_{i-1}^\SO,\ol X_{i-2}^\SO)\cong\pi^s(j-i+1).$$

\begin{rmk}\label{diffrmk}
Since $d^1_{i,*}$ is the boundary homomorphism $\partial$ in the homotopy long exact sequence of the triple $(\ol X_i^\SO,\ol X_{i-1}^\SO,\ol X_{i-2}^\SO)$ composed with the map induced by the key fibration $\ol\chi_{i-1}^\SO$, its meaning for cobordisms is that it assigns to the relative cobordism class of a map $f\colon(M^{i+j},\partial M^{i+j})\to(\R^{i+j}\times\R_+,\R^{i+j})$ the cobordism class of the restriction $f|_{\Sigma^{1_{i-1}}(\partial f)}$ (where we again put $\partial f:=f|_{\partial M}$). However, another description of the boundary map will be more useful, as we will soon see.
\end{rmk}

The following will be analogous to the proof of lemma \ref{d8}. Again represent an element of $\pi_{i+j+1}(\ol X_i^\SO,\ol X_{i-1}^\SO)$ by a map $f\colon(M^{i+j},\partial M^{i+j})\to(\R^{i+j}\times\R_+,\R^{i+j})$ as before and fix closed tubular nighbourhoods $T\subset M$ of $\Sigma^{1_r}(f)$ and $\tilde T\subset\R^{i+j}\times\R_+$ of $f(\Sigma^{1_r}(f))$ such that $f$ maps $(T,\partial T)$ to $(\tilde T,\partial\tilde T)$.

Put $N:=\Sigma^{1_{r-1}}(f|_{\ol{M\setminus T}})$ and decompose its boundary as $\partial N=\partial_0N\sqcup\partial_1N$ with $\partial_0N=\Sigma^{1_{r-1}}(\partial f)$ and $\partial_1N=\Sigma^{1_{r-1}}(f|_{\partial T})$. Since $f|_{\ol{M\setminus T}}$ has no $\Sigma^{1_r}$ singularities, we have that $f|_N$ is a framed immersion to $\ol{(\R^{i+j}\times\R_+)\setminus\tilde T}$. Note that the immersion $g:=f|_{\partial_1N}$ can be canonically equipped with a normal framing in $\R^{i+j+1}$ as well, since we can just add the image of the inwards normal vector field of $\partial T\subset T$ to the rest of the framing vectors. Thus we will think of $g$ as a framed immersion to $\R^{i+j+1}$, so it represents an element $[g]$ in $\pi^s(j-i+1)$.

\begin{prop}\label{d1f}
With the above notations we have $d^1_{i,j}([f])=[g]$.
\end{prop}

\begin{prf}
We may assume that $f$ maps to $\R^{i+j}\times[0,1]$, and so $g$ is a framed immersion of a $(j-i+1)$-manifold to the interior of $\R^{i+j}\times[0,1]$. Then the multi-compression theorem \ref{mct} (more precisely its analogue for immersions; see corollary \ref{ict}) yields a regular homotopy of $g$ that turns one of the framing vector fields parallel to $T[0,1]$ and by combining this with the projection to $\R^{i+j}\times\{1\}$, we get a framed immersion $\partial_1N\imto\R^{i+j}\times\{1\}$. This deformation can be identified with the suspension isomorphism for the stable homotopy groups of spheres, hence this immersion represents the same element of $\pi^s(j-i+1)$ as $g$.

But the above deformation can be extended to the whole manifold $N$, hence we get an immersion $N\imto\R^{i+j}\times[0,1]$ which maps $\partial_0N$ to $\R^{i+j}\times\{0\}$ by $f|_{\Sigma^{1_{i-1}}(\partial f)}$ and $\partial_1N$ to $\R^{i+j}\times\{1\}$ by $g$ (up to regular homotopy). Hence we got a cobordism connecting $f|_{\Sigma^{1_{i-1}}(\partial f)}$ and $g$, which is equivalent to our proposition by remark \ref{diffrmk}.
\end{prf}

All of the above describes what happens to the cobordism groups geometrically in the first page of the singularity spectral sequence.

In the rest of this subsection we restrict to the case of cusp maps and give a similar description of the differential $d^2$ as well. Note that now the spectral sequence only has three non-zero columns, so the only non-zero part of $d^2$ is the collection of the maps
$$d^2_{2,n-2}\colon\ol E^2_{2,n-2}\to\ol E^2_{0,n-1}.$$

We start by the analogue of remark \ref{diffrmk}.

\begin{rmk}\label{diffrmk2}
Since no differential maps to the second column, we have $\ol E^2_{2,n-2}=\ker d^1_{2,n-2}\subset\pi_{n+1}(\ol X_2^\SO,\ol X_1^\SO)$. The homomorphism $d^2_{2,n-2}$ is defined for any element $\alpha\in\ol E^2_{2,n-2}$ as follows: Since $d^1(\alpha)=0$, the element $\partial(\alpha)\in\pi_n(\ol X_1^\SO)$ is such that $j_\#(\partial(\alpha))$ vanishes in the homotopy long exact sequence
$$\ldots\to\pi_{n+1}(\ol X_1^\SO,\ol X_0^\SO)\xra{\partial'}\pi_n(\ol X_0^\SO)\xra{i_\#}\pi_n(\ol X_1^\SO)\xra{j_\#}\pi_n(\ol X_1^\SO,\ol X_0^\SO)\to\ldots$$
Hence $\partial(\alpha)$ is in the image of $i_\#$, i.e. we have $\partial(\alpha)=i_\#(\beta)$ where the coset
$$[\beta]\in\pi_n(\ol X_0^\SO)/\ker i_\#=\pi_n(\ol X_0^\SO)/\im\partial'=\ol E^2_{0,n-1}$$
is uniquely defined for $\alpha$. Now the definition is $d^2(\alpha):=[\beta]$.

Now turning to the geometric meaning of this definition, if we represent $\alpha\in\ol E^2_{2,n-2}\subset\pi_{n+1}(\ol X_2^\SO,\ol X_1^\SO)$ by a map $f\colon(M^n,\partial M^n)\to(\R^n\times\R_+,\R^n)$, then the vanishing of $d^1(\alpha)$ means that $f|_{\Sigma^{1,0}(\partial f)}\colon\Sigma^{1,0}(\partial f)\to\R^n$ is null-cobordant as a framed immersion. This is equivalent to saying that the classifying map $S^n\to\ol X_1^\SO$ composed with the key fibration $\ol\chi_1^\SO\colon\ol X_1^\SO\to\Gamma S^3$ is null-homotopic, which happens precisely if the classifying map itself can be deformed into a fibre $\ol X_0^\SO$ of $\ol\chi_1^\SO$. Since homotopy classes of maps to $\ol X_0^\SO$ correspond to cooriented cobordism classes of immersions, this deformation gives a prim $(\tau_1,\SO)$-cobordism between the prim $(\tau_1,\SO)$-map $\partial f$ and an immersion $g$ of an $(n-1)$-manifold to $\R^n$. Now the element $[g]\in\pi^s(n-1)$ is not uniquely defined for $\alpha$, but its coset in $\pi^s(n-1)/\im d^1$ is; this coset is defined to be $d^2(\alpha)$.
\end{rmk}

In order to understand the differential $d^2$ better, we will need the ring structure on the direct sum
$$\GG:=\bigoplus_{i=0}^\infty\pi^s(i),$$
so now we recall this.

\begin{defi}
The composition product operation $\circ$ on $\GG$ is defined for two elements $\alpha\in\pi^s(i)$ and $\beta\in\pi^s(j)$ in the following way: Represent $\alpha$ by a map $a\colon S^{m+j+i}\to S^{m+j}$ and $\beta$ by a map $b\colon S^{m+j}\to S^m$ and put $\alpha\circ\beta:=[a\circ b]$, that is, $\alpha\circ\beta\in\pi^s(i+j)$ is represented by the map $a\circ b\colon S^{m+i+j}\to S^m$.
\end{defi}

\begin{rmk}
This is a well-defined multiplication that makes $\GG$ a ring. Moreover, this multiplication is also skew-commutative, i.e. for homogeneous elements $\alpha,\beta\in\GG$ we have $\alpha\circ\beta=(-1)^{\dim\alpha\cdot\dim\beta}\beta\circ\alpha$ (see \cite{toda}).
\end{rmk}

\begin{rmk}\label{.com}
The composition product can also described in the language of cobordisms of immersed submanifolds with normal framings (which also represent the stable homotopy groups of spheres): If $\alpha\in\pi^s(i)$ and $\beta\in\pi^s(j)$ are represented respectively by $(A^i\imto\R^{m+i+j},U)$ and $(B^j\imto\R^{m+j},V)$, where $U=(u_1,\ldots,u_{m+j})$ and $V=(v_1,\ldots,v_m)$ are the normal framings of the corresponding immersions, then the framing $U$ identifies a tubular neighbourhood of $A$ with an immersion of $A\times\R^{m+j}$. By immersing $B$ to each fibre of this tubular neigbourhood in the natural way, we obtain a framed immersion of $A\times B$ into $\R^{m+i+j}$; this represents $\alpha\circ\beta$.
\end{rmk}

We will also need the normal forms of $1$-codimensional fold and cusp singularities, which are respectively
\begin{alignat*}2
&\sigma_1\colon(D^2,0)\to(D^3,0);~(x,y)\mapsto(x,xy,y^2),\\
&\sigma_2\colon(D^4,0)\to(D^5,0);~(x_1,x_2,x_3,y)\mapsto(x_1,x_2,x_3,x_1y+x_2y^2,x_3y+y^3)
\end{alignat*}
according to Morin \cite{mor} (where we choose the source disks to be the preimages of the target disks).

\begin{rmk}
The maps $\sigma_1$ and $\sigma_2$ are cooriented and prim, hence they represent elements of the groups $\pi^s(0)\cong\ol E^1_{i,i}$.
\end{rmk}

\begin{prop}\label{clajm}
For any $\alpha\in\pi^s(n)$ we have
\begin{enumerate}
\item\label{zsa} $d^1_{1,n+1}(\alpha)=d^1_{1,1}([\sigma_1])\circ\alpha$ and $d^1_{2,n+2}(\alpha)=d^1_{2,2}([\sigma_2])\circ\alpha=0$,
\item\label{he} $d^2_{2,n+2}(\alpha)$ is the coset of $d^2_{2,2}([\sigma_2])\circ\alpha$ whenever $\alpha\in\ker d^1_{2,n+2}$.
\end{enumerate}
\end{prop}

Note that in order for \ref{he} to make sense, we first need to show that $d^2_{2,2}([\sigma_2])$ is meaningful, i.e. that $d^1_{2,2}([\sigma_2])$ vanishes; this follows from \ref{zsa}. Moreover, we will also need the ambiguity of $d^2_{2,2}([\sigma_2])\circ\alpha$ (which is $(\im d^1_{1,3})\circ\pi^s(n)$) to be contained in the ambiguity of $d^2_{2,n+2}(\alpha)$ (which is $\im d^1_{1,n+3}$); this holds because \ref{zsa} implies
$$(\im d^1_{1,3})\circ\pi^s(n)=(\im d^1_{1,1})\circ\pi^s(2)\circ\pi^s(n)\subset(\im d^1_{1,1})\circ\pi^s(n+2)=\im d^1_{1,n+3}.$$

\medskip\begin{prf}
\emph{Proof of \ref{zsa}.}\enspace\ignorespaces
We will first prove the equality $d^1_{2,n+2}(\alpha)=d^1_{2,2}([\sigma_2])\circ\alpha$; the equality $d^1_{1,n+1}(\alpha)=d^1_{1,1}([\sigma_1])\circ\alpha$ can be obtained completely analogously, so we omit its proof. 

Let $f\colon(M^{n+4},\partial M^{n+4})\to(\R^{n+4}\times\R_+,\R^{n+4})$ be a representative of the element $\alpha\in\pi^s(n)\cong\pi_{n+5}(\ol X_2^\SO,\ol X_1^\SO)$. Recall proposition \ref{d1f} (and the construction that precedes it): $d^1_{2,n+2}(\alpha)$ can be represented by the framed immersion $f|_{\Sigma^{1,0}(f|_{\partial T})}$, where $T\subset M$ is a tubular neighbourhood of $\Sigma^{1,1,0}(f)$. Now the framed immersion $f|_{\Sigma^{1,1,0}(f)}$ represents $\alpha$ and the restriction of $f|_T$ to any fibre of $T$ can be identified with $\sigma_2$ (here we use theorem \ref{univ}), thus $f$ on the $\Sigma^{1,0}$-stratum on the boundary of this fibre represents $d^1_{2,2}([\sigma_2])$. Now taking $f|_{\Sigma^{1,0}(f)}$ on the boundary $\partial T$, we obtain a representative of $d^1_{2,2}([\sigma_2])\circ\alpha$ by remark \ref{.com}, hence the equality is proved.

Now the only thing left to prove that $d^1_{2,2}([\sigma_2])$ vanishes. The first page of the spectral sequence now looks as\\
\vspace{-.5cm}
\begin{center}
\begin{tikzpicture}
\matrix (m) [matrix of math nodes,nodes in empty cells,nodes={minimum width=5ex,minimum height=.5ex,text depth=0ex,inner sep=0pt,outer sep=0pt,anchor=base},column sep=5ex,row sep=1.7ex]{
  & & & & \\
  & & & & \\
\strut ~~~~~2 & \pi^s(2)\cong\Z_2 & \pi^s(1)\cong\Z_2 & \pi^s(0)\cong\Z & \strut\\
  & & & & \\
\strut ~~~~~1 & \pi^s(1)\cong\Z_2 & \pi^s(0)\cong\Z & 0 & \strut\\
  & & & & \\
\strut ~~~~~0 & \pi^s(0)\cong\Z & 0 & 0 & \strut\\
\strut  & 0 & 1 & 2 & \strut \\};
\begin{scope}[transparency group]
\draw[alt double,alt double distance=.7mm] ({$(m-1-1)!.4!(m-1-2)$} |- m-1-1.east) -- ({$(m-8-1)!.4!(m-8-2)$} |- m-8-1.south east);
\draw[alt double,alt double distance=.7mm] (m-8-1.north) -- (m-8-5.north west);
\end{scope}
\draw ({$(m-1-2)!.5!(m-1-3)$} |- m-1-2.east) -- ({$(m-8-2)!.5!(m-8-3)$} |- m-8-2.south east);
\draw ({$(m-1-3)!.5!(m-1-4)$} |- m-1-3.east) -- ({$(m-8-3)!.5!(m-8-4)$} |- m-8-3.south east);
\draw ({$(m-1-4)!.65!(m-1-5)$} |- m-1-4.east) -- ({$(m-8-4)!.65!(m-8-5)$} |- m-8-4.south east);
\draw (m-6-1.north) -- (m-6-5.north west);
\draw (m-4-1.north) -- (m-4-5.north west);
\draw (m-2-1.north) -- (m-2-5.north west);
\path[->,font=\scriptsize] (m-3-4) edge node [above] {$~d^1_{2,2}$} (m-3-3);
\path[->,font=\scriptsize] (m-3-3) edge node [above] {$~~d^1_{1,2}$} (m-3-2);
\path[->,font=\scriptsize] (m-5-3) edge node [above] {$~~d^1_{1,1}$} (m-5-2);
\end{tikzpicture}
\end{center}
In what follows, we refer to Toda \cite[chapter XIV]{toda} for informations on the stable homotopy groups of spheres, the composition product and the standard names of the generators. For the generator $\iota\in\pi^s(0)$, the element $d^1_{1,1}(\iota)$ is represented by $\partial\sigma_1$, which is the immersion of a circle to the sphere with one double point (see figure \ref{kep7}). This is the non-trivial element $\eta\in\pi^s(1)$, hence $d^1_{1,1}$ is epimorphic. Now we have
$$d^1_{1,2}(\eta)=d^1_{1,1}(\iota)\circ\eta=\eta\circ\eta\ne0\in\pi^s(2),$$
hence $d^1_{1,2}$ is an isomorphism. But since $d^1_{1,2}\circ d^1_{2,2}$ is trivial, this implies the vanishing of the differential $d^1_{2,2}$.

\medskip\noindent\emph{Proof of \ref{he}.}\enspace\ignorespaces
This will be quite similar to the above. First apply remark \ref{diffrmk2} to the map $\sigma_2\colon D^4\to D^5$: the boundary $\partial\sigma_2$ is cobordant to an immersion, i.e. there is a manifold $W^4$ with boundary $\partial W=S^3\sqcup N^3$ and a prim $(\tau_1,\SO)$-map $G\colon W\to S^4\times[0,1]$ which restricts to the boundary as $G|_{S^3}=\partial\sigma_2$ and with $G|_N$ being an immersion. We identify $S^4\times[0,1]$ with $D^5_2\setminus\tightoverset{_\circ~}{D^5}$ (where $D^5_2$ is the $5$-disk with radius $2$) and put
$$\tilde\sigma_2:=\sigma_2\cup G\colon D^4\usqcup{S^3}W\to D^5_2.$$
Then clearly $d^2_{2,2}([\sigma_2])$ is represented by $\partial\tilde\sigma_2=\tilde\sigma_2|_N$.

Now if we again represent $\alpha\in\pi^s(n)$ by the map $f$ as above, then again the restriction of $f|_T$ to each fibre of the tubular neighbourhood $T$ can be identified with $\sigma_2$, hence the whole $f|_T\colon T\to\tilde T$ can be identified with $f|_{\Sigma^{1,1,0}(f)}\times\sigma_2\colon\Sigma^{1,1,0}(f)\times D^4\to f(\Sigma^{1,1,0}(f))\times D^5$. Put
$$\tilde f:=f|_{\Sigma^{1,1,0}(f)}\times\tilde\sigma_2\colon\Sigma^{1,1,0}(f)\times\left(D^4\usqcup{S^3}W\right)\to f(\Sigma^{1,1,0}(f))\times D^5_2$$
and denote by $T_2$ and $\tilde T_2$ the source and target spaces of $\tilde f$.

Here $\tilde T_2\subset\R^{n+4}\times\R_+$ is a tubular neighbourhood of $f(\Sigma^{1,1,0}(f))$ that contains $\tilde T$ and $\partial\tilde f$ is a $1$-codimensional cooriented immersion to $\partial\tilde T_2$, hence it can be endowed with a non-vanishing normal vector field in $\tilde T_2$ and the inwards normal vector field of $\tilde T_2$. This means that $\partial\tilde f$ is a framed immersion that clearly represents $d^2_{2,2}([\sigma_2])\circ\alpha$.

We now use the immersion version of the multi-compression theorem \ref{mct} (see corollary \ref{ict}) on $\partial\tilde f$ to turn the normal vector field of $\tilde T_2$ restricted to $\partial\tilde f$ parallel to the direction of $T\R_+$ in $\R^{n+4}\times\R_+$. After this, the projection of the image of this immersion to $\R^{n+4}$ yields an immersion to $\R^{n+4}$, which is regular homotopic to $\partial\tilde f$ and in particular it also represents $d^2_{2,2}([\sigma_2])\circ\alpha$.

Now proposition \ref{ct+} ensures that this regular homotopy can be chosen such that its ``time'' derivative vectors are nowhere tangent to the current image of the source manifold $\partial T_2$. Hence we can glue $T_2$ to $\partial T_2\times[0,1]$ along one common boundary part and map the obtained manifold to $\R^{n+4}\times\R_+$ by the union of $\tilde f$ and the above regular homotopy of $\partial\tilde f$. This way we get a map with the same cusp set as $f$, which means it also represents $\alpha$ (by remark \ref{relcob}), and such that its boundary immersion, which is by definition a representative of $d^2_{2,n+2}(\alpha)$, represents $d^2_{2,2}([\sigma_2])\circ\alpha$. This is what we wanted to prove.
\end{prf}

\subsection{Fold and cusp prim maps}

The previous subsection gave a complete geometric description of the singularity spectral sequence for codimension-$1$ cooriented prim fold and cusp maps. Here we give (almost complete) computations on these cobordism groups, first for the fold case, then for the cusp case. We will use that for all $i$ the generator $\eta\in\pi^s(1)\cong\Z_2$ and the generator $\alpha_1\in\pi^s(3)\cong\Z_3\oplus\Z_8$ of the $\Z_3$ part define homomorphisms
\begin{alignat*}2
\eta_i\colon\pi^s(i)\to\pi^s(i+1);&~\beta\mapsto\eta\circ\beta,\\
(\alpha_1)_i\colon\pi^s(i)\to\pi^s(i+3);&~\beta\mapsto\alpha_1\circ\beta.
\end{alignat*} 

\begin{thm}\label{primfold}
$\Prim_1^\SO(n,1)$ is such that there is a $\CC_2$-isomorphism
$$\Prim_1^\SO(n,1)\congc2\pi^s(n)\oplus\pi^s(n-2)$$
and an exact sequence
$$0\to\coker\eta_{n-1}\to\Prim_1^\SO(n,1)\to\ker\eta_{n-2}\to0.$$
\end{thm}

\begin{prf}
The $\CC_2$-isomorphism can be proved quite similarly to the $\CC_{2,3}$-isomorphism in theorem \ref{cuspthm} and it will also follow from theorem \ref{ploc} later, so we only sketch its proof. The singularity spectral sequence for fold maps degenerates modulo $\CC_2$, since by proposition \ref{clajm} the differential $d^1_{1,n+1}$ is just the multiplication by the order-$2$ element $\eta=d^1_{1,1}([\sigma_1])$. The fact that the cobordism group $\Prim_1^\SO(n,1)$ is a direct sum modulo $\CC_2$ is due to a $2$-splitting of the key fibration
$$\ol\chi_1^\SO\colon\ol X_1^\SO\xra{\ol X_0^\SO}\Gamma S^3.$$

The exact sequence is just part of the homotopy long exact sequence of the key fibration, using that the boundary homomorphism in this long exact sequence is $d^1_{1,n+1}$. By proposition \ref{clajm}, this homomorphism is the multiplication by $d^1_{1,1}([\sigma_1])=\eta$, i.e. the map $\eta_n$.
\end{prf}

To obtain an analogous theorem for cusp maps as well, we need a bit more information on the differential $d^2$ of the singularity spectral sequence. Since $d^1_{1,3}$ maps the generator $\eta\circ\eta\in\pi^s(2)\cong\Z_2$ to $d^1_{1,1}([\sigma_1])\circ\eta\circ\eta=\eta\circ\eta\circ\eta\ne0\in\pi^s(3)$ (see proposition \ref{clajm} and \cite[chapter XIV]{toda}), we have $\ol E^2_{0,3}\cong\Z_{24}/\Z_2=\Z_{12}$. The second page now has the form
\vspace{-.3cm}
\begin{center}
\begin{tikzpicture}
\matrix (m) [matrix of math nodes,nodes in empty cells,nodes={minimum width=5ex,minimum height=.5ex,text depth=0ex,inner sep=0pt,outer sep=0pt,anchor=base},column sep=5ex,row sep=1.7ex]{
  & & & & \\
  & & & & \\
\strut ~~~~~3 & \Z_{12} & 0 & \Z_2 & \strut\\
  & & & & \\
\strut ~~~~~2 & 0 & 0 & \Z & \strut\\
  & & & & \\
\strut ~~~~~1 & 0 & \Z & 0 & \strut\\
  & & & & \\
\strut ~~~~~0 & \Z & 0 & 0 & \strut\\
\strut  & 0 & 1 & 2 & \strut \\};
\begin{scope}[transparency group]
\draw[alt double,alt double distance=.7mm] ({$(m-1-1)!.45!(m-1-2)$} |- m-1-1.east) -- ({$(m-10-1)!.45!(m-10-2)$} |- m-10-1.south east);
\draw[alt double,alt double distance=.7mm] (m-10-1.north) -- (m-10-5.north west);
\end{scope}
\draw ({$(m-1-2)!.5!(m-1-3)$} |- m-1-2.east) -- ({$(m-10-2)!.5!(m-10-3)$} |- m-10-2.south east);
\draw ({$(m-1-3)!.5!(m-1-4)$} |- m-1-3.east) -- ({$(m-10-3)!.5!(m-10-4)$} |- m-10-3.south east);
\draw ({$(m-1-4)!.5!(m-1-5)$} |- m-1-4.east) -- ({$(m-10-4)!.5!(m-10-5)$} |- m-10-4.south east);
\draw (m-8-1.north) -- (m-8-5.north west);
\draw (m-6-1.north) -- (m-6-5.north west);
\draw (m-4-1.north) -- (m-4-5.north west);
\draw (m-2-1.north) -- (m-2-5.north west);
\path[->,font=\scriptsize] (m-5-4) edge node [above,pos=.27] {$d^2_{2,2}$} (m-3-2);
\end{tikzpicture}
\end{center}

\begin{lemma}\label{teritettloboleny}
The differential $d^2_{2,2}\colon\Z\to\Z_{12}$ maps the generator $\iota\in\Z$ to an element of order $6$.
\end{lemma}

\begin{prf}
Since $d^1_{2,2}\colon\ol E^1_{2,2}=\pi_5(\ol X_2^\SO,\ol X_1^\SO)\to\pi_4(\ol X_1^\SO,\ol X_0^\SO)=\ol E^1_{1,2}$ was trivial, we have $\ol E^1_{2,2}=\ol E^2_{2,2}$. Consider the following commutative diagram with exact row and column and with $\partial$ being the boundary map in the homotopy long exact sequence of the key fibration $\ol\chi_1^\SO$:
$$\xymatrix{
&& \pi_4(\ol X_0^\SO)\big/\im\partial \\
\pi_5(\ol X_2^\SO)\ar[r]^(.45){j_\#} & \pi_5(\ol X_2^\SO,\ol X_1^\SO)\ar[ur]^{d^2}\ar[r]\ar[dr]_{d^1} & \pi_4(\ol X_1^\SO)\ar[d]\ar@{<-<}[u] \\
&& \pi_4(\ol X_1^\SO,\ol X_0^\SO)
}$$

A diagram chasing shows that the generator $\iota\in\pi_5(\ol X_2^\SO,\ol X_1^\SO)\cong\Z$ (represented by the map $\sigma_2$) is such that $d^2_{2,2}(\iota)$ has the same order as the order of $\coker j_\#$. The latter is the minimal positive algebraic number of cusps that a cooriented prim cusp map $f\colon M^4\to\R^5$ can have, since $j_\#$ assigns to such a map $f$ the number $\#\Sigma^{1,1,0}(f)$. This minimal number is known to be $6$ from \cite{prim}.
\end{prf}

Combining this lemma with proposition \ref{clajm} yields the following.

\begin{crly}
The differential $d^2_{2,n+2}$ acts on the $3$-primary component as the homomorphism $(\alpha_1)_n$ up to sign.
\end{crly}

\begin{thm}\label{primcusp}
$\Prim_2^\SO(n,1)$ is such that there is a $\CC_{2,3}$-isomorphism
$$\Prim_2^\SO(n,1)\congc{2,3}\pi^s(n)\oplus\pi^s(n-2)\oplus\pi^s(n-4)$$
and an exact sequence of $3$-primary parts
$$0\to[\coker(\alpha_1)_{n-3}]_3\oplus[\pi^s(n-2)]_3\to\big[\Prim_2^\SO(n,1)\big]_3\to[\ker(\alpha_1)_{n-4}]_3\to0.$$
\end{thm}

\begin{prf}
Again as in theorem \ref{primfold}, the $\CC_{2,3}$-isomorphism can be proved analogously to the $\CC_{2,3}$-isomorphism in theorem \ref{cuspthm} and theorem \ref{ploc} will also imply it, so we only sketch its proof. The singularity spectral sequence for cusp maps degenerates modulo $\CC_{2,3}$, since by proposition \ref{clajm} and lemma \ref{teritettloboleny} the differential $d^2_{2,n+2}$ is the multiplication by an order-$6$ element. The fact that the cobordism group $\Prim_2^\SO(n,1)$ is a direct sum modulo $\CC_{2,3}$ is due to a $6$-splitting of the key fibration
$$\ol\chi_2^\SO\colon\ol X_2^\SO\xra{\ol X_1^\SO}\Gamma S^5.$$

Getting to the exact sequence between the $3$-primary parts, we first note that the singularity spectral sequence stabilises at the third page and the $3$-primary part of the last differential $d^2$ can be identified up to sign with the multiplication by $\alpha_1$. Thus we have
$$\big[\ol E^\infty_{0,j}\big]_3=[\coker(\alpha_1)_{j-3}]_3,~~~~\big[\ol E^\infty_{1,j}\big]_3=[\pi^s(j-1)]_3~~~~\text{and}~~~~\big[\ol E^\infty_{2,j}\big]_3=[\ker(\alpha_1)_{j-2}]_3.$$
By general properties of spectral sequences, if we define the groups
\begin{alignat*}2
&F_{2,n}:=\Prim_2^\SO(n,1)=\pi_{n+1}(\ol X_2^\SO),\\
&F_{1,n}:=\im(\Prim_1^\SO(n,1)\to\Prim_2^\SO(n,1))=\im(\pi_{n+1}(\ol X_1^\SO)\to\pi_{n+1}(\ol X_2^\SO)),\\
&F_{0,n}:=\im(\pi^s(n)\to\Prim_2^\SO(n,1))=\im(\pi_{n+1}(\ol X_0^\SO)\to\pi_{n+1}(\ol X_2^\SO))
\end{alignat*}
(with the unmarked arrows being the forgetful maps), then we have $F_{2,n}/F_{1,n}=\ol E^\infty_{2,n-2}$, $F_{1,n}/F_{0,n}=\ol E^\infty_{1,n-1}$ and $F_{0,n}=\ol E^\infty_{0,n}$.

Our plan is now to show that the short exact sequence
$$0\to[F_{0,n}]_3\to[F_{1,n}]_3\to[F_{1,n}/F_{0,n}]_3\to0$$
splits. If this is true, then $[F_{1,n}]_3$ has the form $[F_{0,n}]_3\oplus[F_{1,n}/F_{0,n}]_3$ and the analogous exact sequence with the indices increased can be written as
$$0\to[F_{0,n}]_3\oplus[F_{1,n}/F_{0,n}]_3\to[F_{2,n}]_3\to[F_{2,n}/F_{3,n}]_3\to0.$$
Hence by substituting the definitions of the groups $F_{i,j}$ and the groups $\ol E^\infty_{i,j}$, we obtain the desired exact sequence for $\big[\Prim_2^\SO(n,1)\big]_3$. So it only remains to show that the sequence above indeed splits.

In the following commutative diagram we show $F_{0,n}$ and $F_{1,n}$ with the two left-hand squares being clear from the definitions; the right-hand column and the top and bottom rows are segments of homotopy long exact sequences and the map $s$ is the $2$-splitting of the key fibration $\ol\chi_1^\SO$ used in theorem \ref{primfold} (defined in \cite{nszt}).
$$\xymatrix{
&&& \pi_{n+2}(\ol X_2^\SO,\ol X_1^\SO)\ar[d]^{\partial=d^1_{2,n-1}} \\
\pi_{n+1}(\ol X_0^\SO)\ar[r]\ar@{->>}[d] & \pi_{n+1}(\ol X_1^\SO)\ar[rr]\ar@{->>}[d] && \pi_{n+1}(\ol X_1^\SO,\ol X_0^\SO)\ar[dd]^{i_\#}\ar@/^-1pc/[ll]_s \\
F_{0,n}\ar[r] & F_{1,n}\ar[r] & F_{1,n}/F_{0,n}\ar@{-->}[dr]^a\ar@{-->}[ur]^b\ar@{-->}[ul]_c &\\
\pi_{n+1}(\ol X_2^\SO)\ar@{<-<}[u]\ar@{=}[r] & \pi_{n+1}(\ol X_2^\SO)\ar@{<-<}[u]\ar[rr]^{j_\#} && \pi_{n+1}(\ol X_2^\SO,\ol X_0^\SO)
}$$

The kernel of the composition $F_{1,n}\rightarrowtail\pi_{n+1}(\ol X_2^\SO)\xra{j_\#}\pi_{n+1}(\ol X_2^\SO,\ol X_0^\SO)$ is clearly $\ker j_\#\cap F_{1,n}$. But $\ker j_\#$ is the image of $\pi_{n+1}(\ol X_0^\SO)$, which is $F_{0,n}$, thus the map $j_\#$ uniquely defines a map $a\colon F_{1,n}/F_{0,n}\to\pi_{n+1}(\ol X_2^\SO,\ol X_0^\SO)$. The boundary homomorphism $\partial$ in the left-hand column can be identified with $d^1_{2,n-1}$, which is trivial by proposition \ref{clajm}, hence $i_\#$ is injective. But $\im i_\#$ contains $\im a$ because of the commutativity of the left-hand large square and the definition of $a$. It follows that $a$ can be lifted to a map $b\colon F_{1,n}/F_{0,n}\to\pi_{n+1}(\ol X_1^\SO,\ol X_0^\SO)$.

Now defining $c:=s\circ b$ and composing it with the map $\pi_{n+1}(\ol X_1^\SO)\twoheadrightarrow F_{1,n}$ yields a $2$-splitting of the short exact sequence
$$0\to F_{0,n}\to F_{1,n}\to F_{1,n}/F_{0,n}\to0.$$
When considering the $3$-primary parts, this $2$-splitting becomes just a splitting and this proves the statement.
\end{prf}

\begin{rmk}\label{c2ex}
The singularity spectral sequence $E^*_{*,*}$ can be defined for non-prim Morin maps by a trivial modification of definition \ref{sss} and in the starting page we again have the stable homotopy groups of the Thom spaces of the universal target bundles.

When considering this sequence for cooriented cusp maps, remark \ref{morform} implies that the groups $E^1_{1,j}$ are trivial modulo $\CC_2$ and the groups $E^1_{0,j}$ and $E^1_{2,j}$ coincide with $\ol E^1_{0,j}$ and $\ol E^1_{2,j}$ respectively modulo $\CC_2$ again. We obtained that the natural forgetting map $\ol E^1_{i,j}\to E^1_{i,j}$ (see lemma \ref{cover} and remark \ref{coverr}) induces a $\CC_2$-isomorphism for $i=0,2$ and $E^1_{1,j}\in\CC_2$.

Since the differential $d^1$ is trivial modulo $\CC_2$ in both spectral sequences, the forgetting map $\ol E^2_{i,j}\to E^2_{i,j}$ is also a $\CC_2$-isomorphism for $i=0,2$. Hence the differential $d^2$ restricted to the $3$-primary part can be identified in the two spectral sequences, and so we have $\big[E^\infty_{i,j}\big]_3\cong\big[\ol E^\infty_{i,j}\big]_3$ for $i=0,2$ and $\big[E^\infty_{1,j}\big]_3$ vanishes. Now (the proof of) the previous theorem implies the existence of a $\CC_2$-exact sequence
$$0\to\coker(\alpha_1)_{n-3}\to\Cob_2^\SO(n,1)\to\ker(\alpha_1)_{n-4}\to0,$$
as claimed in theorem \ref{cuspthm}.
\end{rmk}

\subsection{Relation to projective spaces}\label{rps}

In this subsection we give (among others) an alternative description of the singularity spectral sequence considered so far. We will consider two types of prim maps with fixed codimensions and normal structures, namely we will investigate the groups $\Prim_r^L(n,k)$ in the following cases:
\begin{enumerate}
\renewcommand{\labelenumi}{\theenumi}
\renewcommand{\theenumi}{\rm{(\roman{enumi})}}
\item\label{cplx} $L:=\U$, $k:=1$ and in this case we put $\F:=\C$,
\item\label{quat} $L:=\Sp$, $k:=3$ and in this case we put $\F:=\HH$.
\end{enumerate}
\renewcommand{\labelenumi}{\theenumi}
\renewcommand{\theenumi}{\rm{(\arabic{enumi})}}

\begin{rmk}
Here $\Prim_r^\U(n,1)$ is the cobordism group of those prim maps $f\colon M^n\to\R^{n+1}$ that have immersion lifts $i_f\colon M^n\imto\R^{n+2}$ with normal bundles induced from $\gamma_2^\U$ (where the number $2$ now stands for the real dimension of the universal bundle); see remark \ref{primrmk}. Since $\gamma_2^\U$ can be identified with $\gamma_2^\SO$, we now have
$$\Prim_r^\U(n,1)=\Prim_r^\SO(n,1).$$
Similarly, since $\Sp(1)\cong\Spin(3)$, we also have
$$\Prim_r^\Sp(n,3)=\Prim_r^{\Spin}(n,3).$$
\end{rmk}

Throughout this subsection we fix the notations $L$, $k$ and $\F$ either as in \ref{cplx} or as in \ref{quat}. We will need the following simplification of the vector bundle $\zeta_S^r$ (see definition \ref{zetasdef}) in the present two cases, which connects it with $\gamma_{k+1}^L$ that is either the complex or the quaternionic tautological line bundle.

\begin{lemma}
The bundle $\zeta_S^r$ is homotopy equivalent to $\gamma_{k+1}^L|_{\FP^r}$ in the sense that there is a homotopy equivalence of base spaces $h\colon\FP^r\to S((r+1)\gamma_{k+1}^L)$ such that $h^*\zeta_S^r\cong\gamma_{k+1}^L|_{\FP^r}$.
\end{lemma}

\begin{prf}
Consider the space $S^\infty\times S(\F^{r+1})\times\F^1$, where $S(\F^{r+1})$ is the set of unit length elements in $\F^{r+1}$. The group $S^k=S(\F^1)$ naturally acts on all three factors of this space and factoring by the diagonal action yields an $\F^1$-bundle $E$ over $B:=S^\infty\utimes{S^k}S(\F^{r+1})$.

The base space $B$ is an $S(\F^{r+1})$-bundle over $S^\infty/S^k=\FP^\infty=BL(k+1)$ and it coincides with the sphere bundle $S((r+1)\gamma_{k+1}^L)$. But we can also view $B$ as an $S^\infty$-bundle over $S(\F^{r+1})/S^k=\FP^r$, hence we have
$$S((r+1)\gamma_{k+1}^L)\cong B\cong\FP^r.$$
Moreover, this homotopy equivalence takes the pullback of the tautological line bundle over $\FP^\infty$ to the tutological line bundle over $\FP^r$, since it extends to the whole orbit space $E$ which is the total space of these bundles. This is what we wanted to prove.
\end{prf}

\begin{thm}
$\ol X_r^L\cong\Omega\Gamma\FP^{r+1}.$
\end{thm}

\begin{prf}
Trivial from theorem \ref{zetasthm} and the above lemma.
\end{prf}

This theorem has several applications, we will just outline a few of them. The first one is that now the filtration $\ol X_0^L\subset\ol X_1^L\subset\ldots\subset\ol X_r^L$ (for any $r$) can be identified with the filtration $\Omega\Gamma\FP^1\subset\Omega\Gamma\FP^2\subset\ldots\subset\Omega\Gamma\FP^{r+1}$.

\begin{crly}
The singularity spectral sequence corresponding to the homotopy groups of $\ol X_r^L$ coincides with the spectral sequence corresponding to the stable homotopy groups of $\FP^{r+1}$ after an index shift.
\end{crly}

In particular this spectral sequence has first page
\begin{alignat*}2
\ol E^1_{i,j}&=\pi_{i+j+k}(\ol X_i^L,\ol X_{i-1}^L)\cong\pi_{i+j+k}(\Omega\Gamma\FP^{i+1},\Omega\Gamma\FP^i)\cong\pi^s_{i+j+k+1}(\FP^{i+1},\FP^i)\cong\\
&\cong\pi^s_{i+j+k+1}(S^{(k+1)(i+1)})=\pi^s(j-ki).
\end{alignat*}

\begin{ex}
In the case of \ref{cplx} the first page $\ol E^1_{i,j}$ is such that the groups $\ol E^1_{i,i+n}$ are all $\pi^s(n)$ and the differentials for small indices are shown in the following table (see \cite{stabhom} and \cite{nehez}):
\vspace{-.3cm}
\begin{center}
\begin{tikzpicture}
\matrix (m) [matrix of math nodes,nodes in empty cells,nodes={minimum width=5ex,minimum height=.5ex,text depth=0ex,inner sep=0pt,outer sep=0pt,anchor=base},column sep=5ex,row sep=1.7ex]{
  & & & & & & & \\
  & & & & & & & \\
\strut ~~~~~10 & \Z_2\la\eta\circ\mu\ra & (\Z_2)^3 & (\Z_2)^2 & \Z_{240} & \Z_2 & \strut\\
  & & & & & & & \\
\strut ~~~~~9 & \Z_2\la\nu^3\ra\oplus\Z_2\la\mu\ra\oplus\Z_2\la\eta\circ\varepsilon\ra & (\Z_2)^2 & \Z_{240} & \Z_2 & 0 & \strut\\
  & & & & & & & \\
\strut ~~~~~8 & \Z_2\la\ol\nu\ra\oplus\Z_2\la\varepsilon\ra & \Z_{240} & \Z_2 & 0 & 0 & \strut\\
  & & & & & & & \\
\strut ~~~~~7 & \Z_{240}\la\sigma\ra & \Z_2 & 0 & 0 & \Z_{24} & \strut\\
  & & & & & & & \\
\strut ~~~~~6 & \Z_2\la\nu^2\ra & 0 & 0 & \Z_{24} & \Z_2 & \strut\\
  & & & & & & & \\
\strut ~~~~~5 & 0 & 0 & \Z_{24} & \Z_2 & \Z_2 & \strut\\
  & & & & & & & \\
\strut ~~~~~4 & 0 & \Z_{24} & \Z_2 & \Z_2 & \Z & \strut\\
  & & & & & & & \\
\strut ~~~~~3 & \Z_{24}\la\nu\ra & \Z_2 & \Z_2 & \Z & 0 & \strut\\
  & & & & & & & \\
\strut ~~~~~2 & \Z_2\la\eta^2\ra & \Z_2 & \Z & 0 & 0 & \strut\\
  & & & & & & & \\
\strut ~~~~~1 & \Z_2\la\eta\ra & \Z & 0 & 0 & 0 & \strut\\
  & & & & & & & \\
\strut ~~~~~0 & \Z\la\iota\ra & 0 & 0 & 0 & 0 & \strut\\
\strut  & 0 & 1 & 2 & 3 & 4 & \strut \\};
\begin{scope}[transparency group]
\draw[alt double,alt double distance=.7mm] ({$(m-1-1)!.25!(m-1-2)$} |- m-1-1.east) -- ({$(m-24-1)!.25!(m-24-2)$} |- m-24-1.south east);
\draw[alt double,alt double distance=.7mm] (m-24-1.north) -- (m-24-7.north west);
\end{scope}
\draw ({$(m-1-2)!.75!(m-1-3)$} |- m-1-2.east) -- ({$(m-24-2)!.75!(m-24-3)$} |- m-24-2.south east);
\draw ({$(m-1-3)!.5!(m-1-4)$} |- m-1-3.east) -- ({$(m-24-3)!.5!(m-24-4)$} |- m-24-3.south east);
\draw ({$(m-1-4)!.5!(m-1-5)$} |- m-1-4.east) -- ({$(m-24-4)!.5!(m-24-5)$} |- m-24-4.south east);
\draw ({$(m-1-5)!.5!(m-1-6)$} |- m-1-5.east) -- ({$(m-24-5)!.5!(m-24-6)$} |- m-24-5.south east);
\draw ({$(m-1-6)!.5!(m-1-7)$} |- m-1-6.east) -- ({$(m-24-6)!.5!(m-24-7)$} |- m-24-6.south east);
\draw (m-22-1.north) -- (m-22-7.north west);
\draw (m-20-1.north) -- (m-20-7.north west);
\draw (m-18-1.north) -- (m-18-7.north west);
\draw (m-16-1.north) -- (m-16-7.north west);
\draw (m-14-1.north) -- (m-14-7.north west);
\draw (m-12-1.north) -- (m-12-7.north west);
\draw (m-10-1.north) -- (m-10-7.north west);
\draw (m-8-1.north) -- (m-8-7.north west);
\draw (m-6-1.north) -- (m-6-7.north west);
\draw (m-4-1.north) -- (m-4-7.north west);
\draw (m-2-1.north) -- (m-2-7.north west);
\path[->,font=\scriptsize] (m-3-3) edge node [above,pos=.4] {$\eta\circ\mu$} (m-3-2);
\path[->,font=\scriptsize] (m-3-4) edge node [above,pos=.5] {$0$} (m-3-3);
\path[->,font=\scriptsize] (m-3-5) edge node [above,pos=.5] {$\ol\nu+\varepsilon$} (m-3-4);
\path[->,font=\scriptsize] (m-3-6) edge node [above,pos=.5] {$0$} (m-3-5);
\path[->,font=\scriptsize] (m-5-3) edge node [above,pos=.4] {$?$} (m-5-2);
\path[->,font=\scriptsize] (m-5-4) edge node [above,pos=.5] {$0$} (m-5-3);
\path[->,font=\scriptsize] (m-5-5) edge node [above,pos=.5] {$0$} (m-5-4);
\path[->,font=\scriptsize] (m-7-3) edge node [above,pos=.5] {$\ol\nu+\varepsilon$} (m-7-2);
\path[->,font=\scriptsize] (m-7-4) edge node [above,pos=.5] {$0$} (m-7-3);
\path[->,font=\scriptsize] (m-9-3) edge node [above,pos=.5] {$0$} (m-9-2);
\path[->,font=\scriptsize] (m-11-6) edge node [above,pos=.5] {$0$} (m-11-5);
\path[>->,font=\scriptsize] (m-13-5) edge node [above,pos=.5] {$~$} (m-13-4);
\path[->,font=\scriptsize] (m-13-6) edge node [above,pos=.5] {$0$} (m-13-5);
\path[->,font=\scriptsize] (m-15-4) edge node [above,pos=.5] {$0$} (m-15-3);
\path[->,font=\scriptsize] (m-15-5) edge node [above,pos=.5] {$\cong$} (m-15-4);
\path[->,font=\scriptsize] (m-15-6) edge node [above,pos=.5] {$0$} (m-15-5);
\path[>->,font=\scriptsize] (m-17-3) edge node [above,pos=.5] {$~$} (m-17-2);
\path[->,font=\scriptsize] (m-17-4) edge node [above,pos=.5] {$0$} (m-17-3);
\path[->>,font=\scriptsize] (m-17-5) edge node [above,pos=.5] {$~$} (m-17-4);
\path[->,font=\scriptsize] (m-19-3) edge node [above,pos=.5] {$\cong$} (m-19-2);
\path[->,font=\scriptsize] (m-19-4) edge node [above,pos=.5] {$0$} (m-19-3);
\path[->>,font=\scriptsize] (m-21-3) edge node [above,pos=.5] {$~$} (m-21-2);
\end{tikzpicture}
\end{center}
\end{ex}

The following is analogous to theorem \ref{hmor} and is a common generalisation of the $\CC_2$-isomorphism part of theorem \ref{primfold} and the $\CC_{2,3}$-isomorphism part of theorem \ref{primcusp}.

\begin{thm}\label{ploc}
For all $r\ge1$ we have
$$\Prim_r^L(n,k)\congc{\le\frac{k+1}2r+1}\displaystyle\bigoplus_{i=0}^r\pi^s(n+1-(k+1)i-k).$$
\end{thm}

In the two cases \ref{cplx} and \ref{quat} this theorem states respectively
\begin{gather*}
\Prim_r^\SO(n,1)=\Prim_r^\U(n,1)\congc{\le r+1}\displaystyle\bigoplus_{i=0}^r\pi^s(n-2i)\\
\text{and}~~~~\Prim_r^{\Spin}(n,3)=\Prim_r^\Sp(n,3)\congc{\le 2r+1}\displaystyle\bigoplus_{i=0}^r\pi^s(n-4i-2).
\end{gather*}

\medskip\begin{prf}
Fix any prime $p>\frac{k+1}2r+1$. We prove by induction on $r$ that the stable homotopy type of the $p$-localised projective space $\big[\FP^{r+1}\big]_p$ coincides with the $p$-localisation of $S^{k+1}\vee S^{2(k+1)}\vee\ldots\vee S^{(r+1)(k+1)}$. This is clearly true for $r=0$, hence we only have to prove the induction step.

Using the induction hypothesis, the stable homotopy class of the attaching map of the top dimensional cell of $\FP^{r+1}$ after $p$-localisation is given by a collection of $p$-localised stable maps from $S^{(k+1)r+k}$ to $S^{i(k+1)}$ for $i=1,\ldots,r$. Such a map represents an element in the $p$-primary part of $\pi^s((k+1)(r-i)+k)$, but by a theorem of Serre from \cite{serre} this $p$-primary part is trivial, since $p>\frac{(k+1)(r-i)+k}2+1$. Hence after the $p$-localisation all stable attaching maps are null-homotopic and this proves the stated form of $\big[\FP^{r+1}\big]_p$.

Thus using that $\Omega\Gamma S=\Gamma$ and $\Gamma(A\vee B)\cong\Gamma A\times\Gamma B$, we have
\begin{alignat*}2
\big[\ol X_r^L\big]_p&\cong\Omega\Gamma\big[S^{k+1}\vee\ldots\vee S^{(r+1)(k+1)}\big]_p\cong\\
&\cong\Gamma[S^{k}\vee\ldots\vee S^{r(k+1)+k}\big]_p\cong\prod_{i=1}^r\big[\Gamma S^{(k+1)i+k}\big]_p,
\end{alignat*}
and now applying the functor $\pi_{n+1}$ yields the statement.
\end{prf}

\begin{rmk}
All of the above also works with the substitutions $L:=\O$, $k:=0$ and $\F:=\R$, we just omitted this case because we considered only positive codimensional maps throughout this thesis.
\end{rmk}

\chapter{Final remarks and open questions}

In the present thesis we described results on computing the cobordism groups of singular maps. What we did not include is the consequences and applications of these, so we now just mention some of them.

The first application of these groups (the initial reason why Szűcs introduced them) is computing cobordism groups of immersions and embeddings in dimensions where the classical theory did not succeed; this resulted in \cite{analog}, \cite{immemb} and \cite{immor}. A few orientability properties could also be described as consequences of the cobordism theory of singular maps; see \cite{rsz}. Yet another application is the answer to (a version of) a question posed by Arnold on eliminating the most complicated singularity of a map by some deformation, which was answered by the key bundle; see \cite{elimcob} and remark \ref{keyrmk}. A very vague, yet important problem is to find other applications of this theory as well.

Turning again to the cobordism groups themselves, observe (which I also mentioned in the introduction) that many basic notions in this thesis were introduced more generally than they were in the original papers (such as the constructions of classifying spaces, the Kazarian conjecture, the key fibration, etc.). These were all such that the original proofs of theorems could be applied in the general setting too, however, there were a few occasions where I had trouble implementing these proofs; let us recall these now. The analogue of theorem \ref{kazc} is probably not true for unoriented cobordisms; theorem \ref{charcob} is most certainly true in some version in the unoriented case as well; and I expect that the orientability conditions in theorems \ref{ratrivi1} and \ref{ratrivi2} can also be relaxed.

Even so, there probably are a few concrete computations from chapter \ref{compcobgr} that can be modified to also work for maps which have different (for example spin) normal structures than the ones considered there. So a project can be to find such computations for maps with various additional structures.

In the following we describe questions in a few more related topics.

\section{Multiplicative structures}

In \cite{multmor} and \cite{szszt} two different ways were introduced to define product operations on direct sums of the free parts of certain cobordism groups of singular maps (tensored by $\Q$), which make these direct sums rings. An interesting question to consider is whether we can generalise (one of) them to obtain multiplicative structures on direct sums of other types of cobordism groups as well. If not, can we find another way to define such product operations? Can we define a ring structure that also includes the torsion parts of these groups?

If the answer to (one of) the above questions is affirmative, then it is natural to ask how we can describe such a ring. For example, how can we express the homology class represented by the singularities of a map in the product of two cobordism classes using the singularities of representatives of the two initial cobordism classes?

\section{Double spectral sequences}

Suppose we have a stable linear group $L$ and a set of $k$-codimensional singularities $\tau=\{\eta_0,\eta_1,\ldots\}$ with a complete order of type $\omega$ or less extending the natural partial order. Put $\xi_i^L:=\xi_{\eta_i}^L$, $\tilde\xi_i^L:=\tilde\xi_{\eta_i}^L$, $c_i:=\dim\xi_i^L$ and $G_i^L:=G_{\eta_i}^L$.

It was noted in \cite{hosszu} that there are two double spectral sequences converging to the cobordism groups $\Cob_\tau^L(n,k)$ as shown on the diagram below. The double spectral sequences here mean systems of groups $E^1_{i,j,l}$ which form the starting page of a spectral sequence with variables $i,j$ and fixed $l$ converging to groups that form the starting page of another spectral sequence with variables $i+j,l$ (this is shown on the diagram as the left-hand vertical and the bottom horizontal arrows); similarly the groups $E^1_{i,j,l}$ are also the starting page of a spectral sequence with variables $i,l$ and fixed $j$ converging to the starting page of a spectral sequence with variables $i+l,j$ (this is the top horizontal and the right-hand vertical arrows).
$$\xymatrix{
H_i\big(BG_l^L;(\widetilde{\pi^s(j)})_{\tilde\xi_l^L}\big)\ar[r]^(.47)\cong\ar[d]_\cong & H_{i+c_l}\big(T\xi_l^L;(\widetilde{\pi^s(j)})_{\nu_\tau^L}\big)\ar@{=>}[r]^(.48){\rm{K}} & H_{i+c_l+l}\big(K_\tau^L;(\widetilde{\pi^s(j)})_{\nu_\tau^L}\big)\ar[d]^\cong \\
H_{i+c_l+k}(T\tilde\xi_l^L;\pi^s(j))\ar@{=>}[d]_{\rm{AH}} && H_{i+c_l+l+k}(V_\tau^L;\pi^s(j))\ar@{=>}[d]^{\rm{AH}} \\
\pi^s_{i+c_l+k+j}(T\tilde\xi_l^L)\ar@{=>}[rr]^(.43){\rm{Sz}} && *!<1cm,0cm>{~\pi_{*+k}(X_\tau^L)\cong\pi^s_{*+k}(V_\tau^L)}
}$$

The isomorphisms indicated in the diagram are all Thom isomorphisms and the double arrows marked by $\rm{K}$, $\rm{AH}$ and $\rm{Sz}$ respectively mean the Kazarian spectral sequence, the Atiyah--Hirzebruch spectral sequence and Szűcs's singularity spectral sequence (the analogue of definition \ref{sss} for non-prim maps). Hence we obtained the following.

\begin{prop}
There are two double spectral sequences
$$E^1_{i,j,l}:=H_i\big(BG_l^L;(\widetilde{\pi^s(j)})_{\tilde\xi_i^L}\big)\implies\pi_{*+k}(X_\tau^L)\cong\Cob_\tau^L(*,k).$$
\end{prop}

Although this gives a way in principle to compute the groups $\Cob_\tau^L(n,k)$ completely, the practical computation seems to be rather difficult; to my knowledge there was no concrete result on cobordism groups derived from these double spectral sequences so far. So it may be interesting to study these and (for example) trying to obtain estimates on the torsion of the cobordism groups of cooriented Morin maps (recall that the ranks of these groups were computed in theorem \ref{coorthm}).

\section{Mosher's spectral sequence}

In subsection \ref{rps} we saw that in some cases (namely for $1$-codimensional cooriented prim maps and $3$-codimensional prim maps with spin normal structures) the singularity spectral sequence can be identified with spectral sequences in stable homotopy groups arising from filtrations of projective spaces. In one case this was the filtration of $\CP^\infty$ by the subspaces $\CP^n$, which was studied by Mosher \cite{stabhom} (his paper is rather compressed and misses many proofs and details, but these are clarified and completed by Szűcs and Terpai in \cite{nehez}).

In the following we recall (without proofs) a few results on this spectral sequence, which will be denoted by $E^*_{*,*}$. We will also view its $p$-localised version denoted by $^p\! E^*_{*,*}$ (where $p$ is a prime), which starts with the $p$-components of $E^1_{*,*}$ and the first differential $^p\! d^1$ is the $p$-component of $d^1$.

Firstly there is a periodicity theorem on $E^*_{*,*}$ and its analogue on $^p\! E^*_{*,*}$ using the Atiyah--Todd numbers $M_n$ (see \cite{at}) that have the form
$$M_n=\prod_{p~\rm{prime}}p^{\nu_p(M_n)}~~~~\text{with}~~~~\nu_p(M_n)=\max\big\{m+\nu_p(m)\mid 1\le m\le\textstyle\big\lfloor\frac{n-1}{p-1}\big\rfloor\big\}.$$
We will also put $[M_n]_p:=p^{\nu_p(M_n)}$ for the $p$-component of $M_n$. The first few Atiyah--Todd numbers are listed below.
\begin{table}[H]\begin{center}\begin{tabular}{c||c|c|c|c|c|c}
$n$ & $1$ & $2$ & $3$ & $4$ & $5$ & $6$\\
\hline
$M_n$ & $1$ & $2$ & $2^3\cdot 3$ & $2^3\cdot 3$ & $2^6\cdot3^2\cdot5$ & $2^6\cdot3^2\cdot5$
\end{tabular}\end{center}\end{table}\vspace{-.5cm}

\begin{thm}[\cite{stabhom},\cite{nehez}]
$~$
\begin{enumerate}
\item For any $r<n$ we have $E^r_{i,j}\cong E^r_{i+M_n,j+M_n}$, moreover, this isomorphism commutes with the differential $d^r$.
\item For any $r<n$ and any prime $p$ we have $^p\! E^r_{i,j}\cong {^p\! E^r_{i+[M_n]_p,j+[M_n]_p}}$, moreover, this isomorphism commutes with the differential $^p\! d^r$.
\end{enumerate}
\end{thm}

Next we use the $J$ homomorphism and its ananlogue in the $p$-localised case (for a prime $p$) denoted by $J_p$ and defined in \cite{nehez}. We remark that $J_p$ had two different definitions in \cite{nehez}, an algebraic and a geometric one; these coincide in the present case, but generally it is an open question whether the two definitions are equivalent.

\begin{thm}[\cite{stabhom},\cite{nehez}]
$~$
\begin{enumerate}
\item If $\iota\in E^1_{i,i}\cong\pi^s(0)\cong\Z$ is a generator for which $d^r(\iota)$ vanishes for all $r<n$, then $d^n(\iota)\in E^n_{i-n,i+n-1}$ belongs to the image of $\im J\subset\pi^s(2n-1)\cong E^1_{i-n,i+n-1}$ in $E^n_{i-n,i+n-1}$.
\item If $\iota\in {^p\! E^1_{i,i}}\cong\pi^s(0)\otimes\Z_{(p)}\cong\Z_{(p)}$ is a generator for which $^p\! d^r(\iota)$ vanishes for all $r<n$, then $^p\! d^n(\iota)\in {^p\! E^n_{i-n,i+n-1}}$ belongs to the image of $\im J_p\subset\pi^s(2n-1)\otimes\Z_{(p)}\cong {^p\! E^1_{i-n,i+n-1}}$ in $^p\! E^n_{i-n,i+n-1}$.
\end{enumerate}
\end{thm}

By the geometric interpretation of the spectral sequence $E^*_{*,*}$ in subsection \ref{geode}, the first claim above translates to the following: The vanishing of the first $n-1$ differentials on the class represented by the normal form $\sigma_{i-1}\colon(D^{2i-2},S^{2i-3})\to(D^{2i-1},S^{2i-2})$ of the $\Sigma^{1_{i-1}}$ singularity implies that the boundary $\partial\sigma_{i-1}$ can be chosen (up to cobordism) to have at most $\Sigma^{1_{i-n-1}}$ singularities and such that its restriction to the $\Sigma^{1_{i-n-1}}$-stratum is a framed immersion of $S^{2n-1}$.

Lastly there is also a way to determine the first non-zero differential on a generator $\iota\in E^1_{i,i}$. This also uses the following.

\begin{lemma}[\cite{stabhom},\cite{nehez}]
For a generator $\iota\in E^1_{i,i}$, the differential $d^n$ is defined on $\iota$ iff $i+1$ is divisible by $M_n$ and if this is the case, then $d^n(\iota)$ vanishes iff $i+1$ is divisible by $M_{n+1}$ as well.
\end{lemma}

Now the tool will be the so-called $e$-invariants of Adams, more precisely the homomorphism $e_C\colon\pi^s(2n-1)\to\Q/\Z$ (see \cite{adams}). One important property of this is that any two representatives of an element in $[\im J]\subset E^n_{*,*}$ (the image of $\im J\subset E^1_{*,*}$) have the same $e_C$-value, hence $e_C$ is well-defined on $[\im J]$, moreover, it is also injective on $[\im J]$ (see \cite{adams} and \cite{stabhom}). We saw in the previous theorem that the vanishing of $d^1,\ldots,d^{n-1}$ on $\iota$ means that we have $d^n(\iota)\in[\im J]\subset E^n_{i-n,i+n-1}$, thus the value $e_C(d^n(\iota))$ uniquely determines $d^n(\iota)$.

\begin{thm}[\cite{stabhom},\cite{nehez}]
If $i+1\equiv kM_n~\mod M_{n+1}$ and $\iota\in E^1_{i,i}$ is a generator, then we have $e_C(d^n(\iota))=ku_n$, where $u_n$ is the coefficient of $z^n$ in the Taylor expansion
$$\left(\frac{\log(1+z)}z\right)^{M_n}=\sum_{j=0}^\infty u_nz^n.$$
\end{thm}

The theorems above are strong results on the spectral sequence arising from the natural filtration of $\CP^\infty$, which coincides with the singularity spectral sequence for cooriented codimension-$1$ prim maps. Now what can we say about the other case considered in subsection \ref{rps}? Can we obtain similar theorems on the singularity spectral sequence for codimension-$3$ prim maps with spin normal structures, i.e. the spectral sequence arising from the natural filtration of $\HP^\infty$? If so, what geometric meaning do they yield for singular maps? What other applications do these spectral sequences have for singularities besides those described in \cite{nehez}, \cite{ctrl} and subsection \ref{rps}?

\section{Bordism groups}

There is a different notion of cobordism relation between singular maps, which is also interesting to consider. This is called (left-right) bordism and is defined in the following way.

\begin{defi}
Let $\ul\tau$ be a set of (multi)singularities of a fixed codimension $k$. We call two $\ul\tau$-maps $f_0\colon M_0^n\to P_0^{n+k}$ and $f_1\colon M_1^n\to P_1^{n+k}$ (with closed source and target manifolds) left-right $\ul\tau$-bordant (or simply $\ul\tau$-bordant) if there is
\begin{itemize}
\item[(i)] a compact manifold with boundary $W^{n+1}$ such that $\partial W=M_0\sqcup M_1$,
\item[(ii)] a compact manifold with boundary $Z^{n+k+1}$ such that $\partial Z=P_0\sqcup P_1$,
\item[(iii)] a $\ul\tau$-map $F\colon W^{n+1}\to Z^{n+k+1}$ such that for $i=0,1$ we have $F^{-1}(P_i)=M_i$ and $F|_{M_i}=f_i$.
\end{itemize}
The set of (left-right) $\ul\tau$-bordism classes of $\ul\tau$-maps of $n$-manifolds to $(n+k)$-manifolds is denoted by $\Bord_{\ul\tau}(n,k)$.
\end{defi}

The disjoint union $f_0\sqcup f_1\colon M_0\sqcup M_1\to P_0\sqcup P_1$ now obviously defines an Abelian group operation on $\Bord_{\ul\tau}(n,k)$ and it is not hard to see (by a modified version of the Pontryagin--Thom construction) that the isomorphism
$$\Bord_{\ul\tau}(n,k)\cong\NN_{n+k}(X_{\ul\tau})$$
holds. Moreover, additional normal structures can be defined for these bordism groups similarly to those on cobordism groups and analogous isomorphisms hold for bordism groups equipped with extra structures.

Some theorems and computations on these bordism groups were obtained in \cite{analog}, \cite{eszt}, \cite{hosszu}, \cite{szszt} and \cite{bordfold}. Apart from the obvious problem to generalise the existing results on them, we can also ask which of the theorems and conjectures on cobordism groups have analogues for bordism groups.

\section{Some more questions}

The following are a few more interesting problems related to the topic of this thesis:
\begin{enumerate}
\item What connections can we find between the Thom polynomials describing homology classes represented by singularity strata and the key fibration connecting the classifying spaces?
\item We know as a result of Grant and Szűcs \cite{grr} that given a manifold $P^{n+k}$, not all homology classes in $H_n(P;\Z_2)$ can be represented by immersions $M^n\imto P^{n+k}$, moreover, no finite set of multisingularities $\ul\tau$ is enough to represent all homology classes by $\ul\tau$-maps. Is the analogue of this true or false for infinite multisingularity sets or sets of singularities with no restrictions on multiplicities? What can we say about representing homology classes with integer coefficients?
\item If we have an immersion equipped with a non-vanishing normal vector field, then projecting it to hyperplanes yields prim maps. On the other hand, the compression theorem implies that after a deformation such a projection will be an immersion (see corollary \ref{ict}). How and why do all obstructions represented by the singularity strata vanish in this case?
\end{enumerate}

%-------------------------------------------------------------------------------------------------------------------------------
%-------------------------------------------------------------------------------------------------------------------------------

\chapter*{Appendix}
\addcontentsline{toc}{chapter}{Appendix}

\newcounter{kaki}
\renewcommand\thesection{\Alph{kaki}}
\renewcommand\thefigure{\thesection.\arabic{figure}}

%-------------------------------------------------------------------------------------------------------------------------------

\refstepcounter{kaki}  
\section{The compression theorem}

We recall the compression theorem of Rourke and Sanderson from \cite{rs}, then see some addenda and extensions due to Rourke, Sanderson, Szűcs and Terpai. Throughout this section $P$ will always denote a manifold with a fixed Riemannian metric, the distance of two points is always denoted by $d(\cdot,\cdot)$, the inner product by $\la\cdot,\cdot\ra$ and the length of a vector by $\lv\cdot\rv$. The basic compression theorem is the following.

\begin{prop}\label{ct}
Let $n,k$ be natural numbers such that $k\ge1$ and let $i\colon M^n\hookrightarrow P^{n+k}\times\R^1$ be the embedding of a compact manifold (possibly with boundary) equipped with a nowhere zero normal vector field $u$ (i.e. a section of the bundle $T(P\times\R^1)|_{i(M)}\setminus Ti(M)$). Then there is an (ambient) isotopy $\varphi_t~(t\in[0,1])$ of $i$ such that $\varphi_0=i$ and $d\varphi_1(u)$ is the vertically upwards (i.e. everywhere on the positive ray of $\R^1$ in $T(P\times\R^1)$) unit vector field, which will be denoted by $\ua$.
%\item[\rm{(2)}] If there is a compact subset $C\subset M$ for which $u|_C$ is already vertical, then the diffeotopy $\varphi_t$ can be chosen such that it fixes $i(C)$ and $u|_C$.
\end{prop}

For the sake of simplicity, throughout this section we will call the summands $TP$ and $\R^1$ in $T(P\times\R^1)$ horizontal and vertical respectively; the positive and negative directions on $\R^1$ will be called upwards and downwards respectively; and the parameter $t$ in one-parameter families of diffeomorphisms will be thought of as the time.

\begin{rmk}
The name ``compression theorem'' refers to the fact that an embedding $M\hookrightarrow P\times\R^1$ equipped with a vertical normal vector field can be ``compressed'' to an immersion $M\looparrowright P$ by the vertical projection.
\end{rmk}

\begin{rmk}
By the isotopy extension theorem, an isotopy of the embedding $i$ is equivalent to a diffeotopy of $P\times\R^1$, hence we will not differentiate between these two notions.
\end{rmk}

Before getting to the proof of proposition \ref{ct} we show the following preliminary observations.

\begin{lemma}\label{ctlemma}
We can assume that the vector field $u$ is
\begin{enumerate}
\item\label{ct1} everywhere of unit length and orthogonal to the tangent space of $i(M)$,
\item\label{ct2} nowhere vertically downwards (i.e. nowhere equal to $\da:=-\ua$).
\end{enumerate}
\end{lemma}

\begin{prf}
\emph{Proof of \ref{ct1}.}\enspace\ignorespaces We can apply the Gram--Schmidt process pointwise in $T(P\times\R^1)|_{i(M)}$ to turn $u$ into a vector field of the required form. This is a smooth deformation, hence it results in a smooth vector field.

\medskip\noindent\emph{Proof of \ref{ct2}.}\enspace\ignorespaces Now because of the assumption \ref{ct1}, both $u$ and $\da$ are sections of the unit sphere bundle of $T(P\times\R^1)|_{i(M)}$. This sphere bundle is a $(2n+k)$-dimensional manifold and the images of $u$ and $\da$ are both $n$-dimensional submanifolds, hence if $u$ is transverse to $\da$ (which we can assume), then they are disjoint (because $k\ge1$).
\end{prf}

\par\noindent\textbf{Proof of proposition \ref{ct}.\enspace\ignorespaces}
Denote by $\alpha(p)$ the angle of $u(p)$ and $\ua(p)$ for all $p\in i(M)$; this defines a smooth map $\alpha\colon i(M)\to[0,\pi)$. Take an $\alpha_0$ for which $0<\alpha_0<\pi-\max\alpha$ and let $v$ be the vector field that we get for all $p\in i(M)$ by rotating $u(p)$ in the plane of $u(p)$ and $\ua(p)$ upwards by angle $\frac\pi2-\alpha_0$ if $\alpha(p)\ge\frac\pi2-\alpha_0$, otherwise only rotating until $\ua(p)$. This $v$ is not necessarily smooth but we can assume that it is, since we can approximate it with a smooth vector field. Now $v$ is such a normal vector field that $\la v(p),\ua(p)\ra>0$ for all $p\in i(M)$.

Let $V$ be a closed tubular neighbourhood of $i(M)$ of radius $r$ and define an extension $\tilde v$ of the vector field $v$ to $P\times\R^1$ in the following way: Let $\beta\colon i(M)\to[0,\frac\pi2)$ be the angle of $v$ and $\ua$ pointwise and let $f\colon[0,r]\to[0,1]$ be a smooth function with $f(0)=0, f(r)=1$ and $\underset{0+0}\lim f'=\underset{r-0}\lim f'=0$. Any point $q\in V$ is in a fibre $V_p$ for some $p\in i(M)$; if the distance of $q$ and this $p$ is $t\in[0,r]$, then we define $\tilde v(q)$ as follows: Parallel translate $v(p)$ to $q$ along the minimal geodesic $[p,q]$, then rotate it upwards by angle $f(t)\beta(p)$. In the complement of $V$ we define $\tilde v$ to be everywhere $\ua$.

If the flow of $\tilde v$ is $\{\Phi_t\in\Diff(P\times\R^1)\mid t\in\R_+\}$, then every integral curve leaves the compact $V$ in finite time, so there is an $s>0$ for which $\Phi_t(p)\notin V$ for all $t\ge s$ and $p\in i(M)$. Therefore if we define $\varphi_t$ as $\Phi_{ts}$ (for $t\in[0,1]$), then the vector field on $\varphi_1(M)$ is $d\varphi_1(v)=\ua$ (now we also used that the preliminary rotation of $u$ to $v$ can be achieved by deforming a neighbourhood of $i(M)$ leaving $i(M)$ fixed).
~$\square$\par\medskip

\subsection{Local compression}

In this subsection we prove the key extension of proposition \ref{ct} (again by Rourke and Sanderson) which is called local compression theorem. In order to prove it we will need a few remarks and lemmas first; to make notations simpler we will use the convention $(p,x)+y:=(p,x+y)$ for all $(p,x)\in P\times\R^1$ and $y\in\R^1$ and identify the tangent spaces $T_{(p,x)}(P\times\R^1)$ for all $x\in\R^1$. We are in the same setting as before and we again assume the conditions in lemma \ref{ctlemma}.

\begin{rmk}\label{lctrmk}
We can change the vector field $\tilde v$ in the proof of proposition \ref{ct} to the time-dependent vector field
$$\tilde v_t\colon P\times\R^1\to T(P\times\R^1);~p\mapsto\tilde v(p+t)+\da~~~~(t\in\R_+).$$
Let the flow of this vector field be $\{\tilde\Phi_t\in\Diff(P\times\R^1)\mid t\in\R_+\}$, i.e. for any $p\in P\times\R^1$ the curve $\gamma\colon t\mapsto\tilde\Phi_t(p)$ is the unique solution of the differential equation $\gamma(0)=p,\gamma'(t)=\tilde v_t(\gamma(t))$.

Then we have the identity $\tilde\Phi_t(p)=\Phi_t(p)-t$ for all $p\in P\times\R^1$; in other words, $\tilde\Phi_t$ is obtained by flowing along $\tilde v$ for time $t$ and then flowing back down the unit downwards flow again for time $t$. This way, if the vector $u(p)$ was initially vertically upwards at a point $p\in i(M)$, then the flow $\tilde\Phi_t$ fixes $p$ and $u(p)$ for sufficiently small numbers $t$ (for larger numbers $t$ it may happen that $\tilde v(p+t)$ is not vertically upwards, and so $\tilde\Phi_t$ eventually moves the point $p$).
\end{rmk}

\begin{defi}
The set of points in $i(M)$ where the orthogonal complement of the tangent space is horizontal is called the horizontal set, i.e. the horizontal set is
$$H:=H(i):=\{(p,x)\in i(M)\mid T_{(p,x)}(\{p\}\times\R^1)\subset T_{(p,x)}i(M)\}$$
\end{defi}

\begin{lemma}
We can assume that $H$ is a submanifold of $i(M)$.
\end{lemma}

\begin{prf}
Consider the $1$-jet bundle $J^1(M,P\times\R^1)\to M\times(P\times\R^1)\to M$. We can choose in each fibre of $J^1(M,P\times\R^1)\to M\times(P\times\R^1)$ (which can be identified with the space of linear maps $\R^n\to\R^{n+k}\times\R^1$) the space of those linear maps $A$ for which $\im A$ contains the last coordinate line $\R^1$. The union of these in all fibres is a ($(k+1)$-codimensional) submanifold of $J^1(M,P\times\R^1)$ and $H(i)$ is the $J^1i$-preimage of it.

By the jet transversality theorem, we can perturb $i$ slightly in $C^\infty(M,P\times\R^1)$ to get a new map $j\colon M\to P\times\R^1$ for which $H(j)$ is a submanifold. By the openness of the space of embeddings in $C^\infty(M,P\times\R^1)$, this $j$ is an embedding isotopic to $i$.
\end{prf}

\begin{defi}
For all $p\in i(M)\setminus H$ there is a unique unit vector in $N_pi(M)$ (the orthogonal complement of $T_pi(M)$) with maximal vertical coordinate; it will be denoted by $x(p)$. The set of points in $i(M)\setminus H$ where the vector $u$ is downmost is called the downset, i.e. the downset is
$$D:=D(i,u):=\{p\in i(M)\setminus H\mid u(p)=-x(p)\}$$
\end{defi}

\begin{lemma}
We can assume that $\overline D$ is a submanifold of $i(M)$ with boundary $H$.
\end{lemma}

\begin{prf}
Both $u|_{i(M)\setminus H}$ and $-x$ are sections of the unit sphere bundle of $N(i(M)\setminus H)$. We may assume that $u|_{i(M)\setminus H}$ is transverse to $-x$, and then the preimage of their intersection in $i(M)\setminus H$ (i.e. $D$) is a ($k$-codimensional) submanifold. 

If $T\subset i(M)$ is a tubular neighbourhood of $H$, then we can canonically identify the normal spaces $N_pi(M)$ over each fibre of $T$. We may also homotope the vector field $u$ around $H$ such that it becomes constant on each fibre of $T$.

Moreover, we can identify the fibres of $T(P\times\R^1)|_{i(M)}$ with $\R^{n+k+1}$ over each fibre of $T$, and then we may assume that the bundle of the ($(k+1)$-dimensional) subspaces $N_pi(M)$ intersects the submanifold of the horizontal $(k+1)$-planes in the Grassmannian $G_{k+1}(\R^{n+k+1})$ transversally over each fibre of $T$. Now if $T$ was chosen sufficiently small, then for all $p\in H$ any non-horizontal $(k+1)$-plane in a neighbourhood of $N_pi(M)$ takes place exactly once among the spaces $N_qi(M)$ where $q$ is in the fibre of $T$ over $p$. This implies that the vector field $-x|_{S_p}$, where $S_p\approx S^k$ is the fibre over $p$ of the sphere bundle of $T$ with any sufficiently small radius, composed with the projection to the fibre of the unit sphere bundle of $Ni(M)$ (which is also an $S^k$) is a diffeomorphism.

Hence $-x$ and $u$ intersect exactly once in $S_p$ for all $p\in H$ and so the closure of $D$ will be a manifold with boundary $H$.
\end{prf}

\begin{defi}
Let $y$ be the gradient field of the projection $\pr_{\R^1}|_{i(M)}$ to the vertical coordinate line, i.e. $y(p)$ is the projected image of $\ua(p)$ in $T_pi(M)$ for all $p\in i(M)$.
\end{defi}

\begin{lemma}\label{lctlemma}
For any tubular neighbourhood $T\subset i(M)$ of $H$ and any number $\delta>0$ we can assume that $\overline{D\setminus T}$ has a tubular neighbourhood $U\subset i(M)$ such that each component of the intersection of an integral curve of $y$ with $U$ has length less than $\delta$.
\end{lemma}

\begin{figure}[H]
\begin{center}
\centering\includegraphics[scale=0.1]{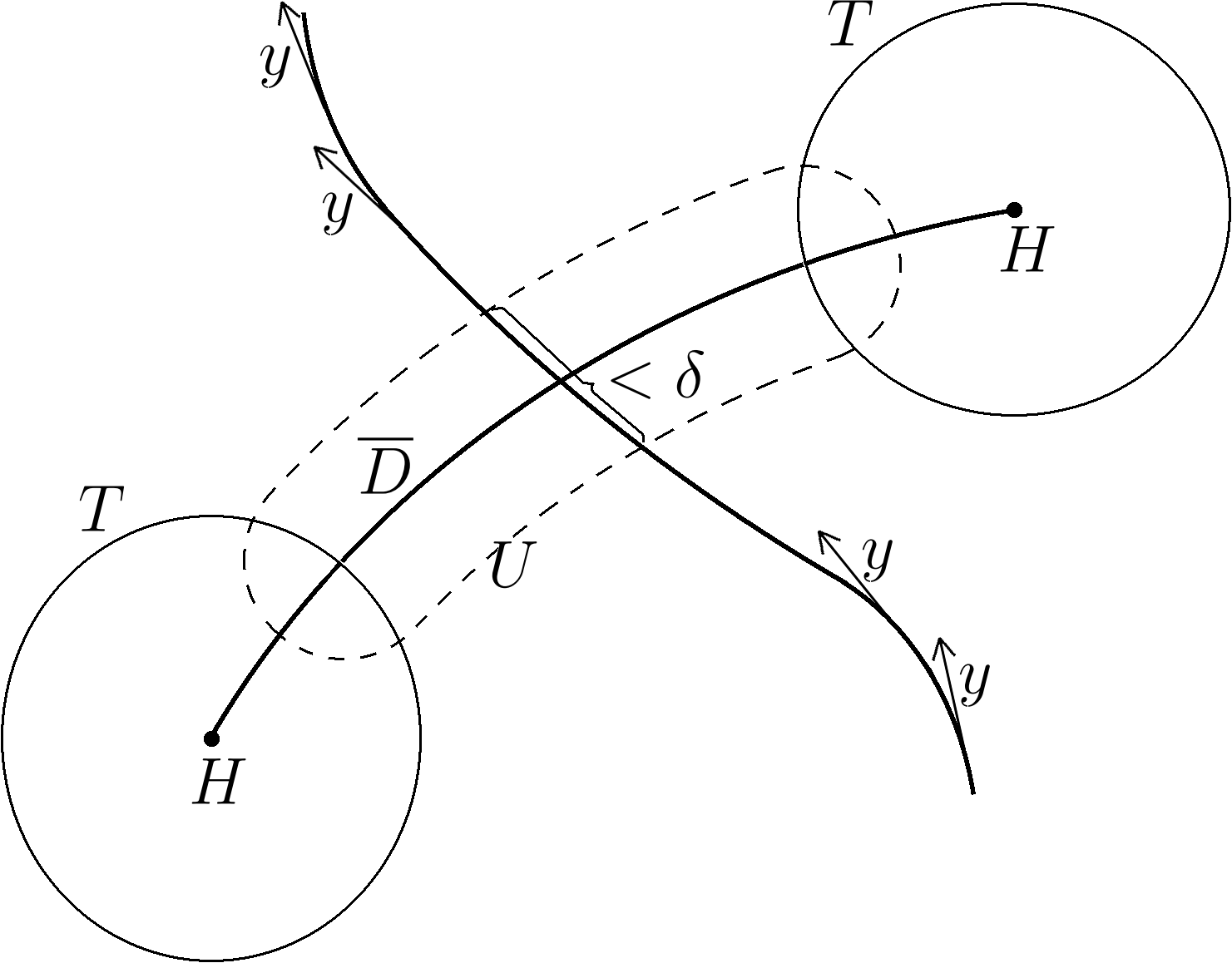}\label{kepa1}
\begin{changemargin}{2cm}{2cm} 
\caption{\hangindent=1.4cm\small We represent $H$ as two points, the arc connecting them is $\overline D$ and the curve with the indicated tangent vectors is an integral curve of $y$.}
\end{changemargin} 
\vspace{-1.5cm}
\end{center}
\end{figure}

\begin{prf}
Observe that a small homotopy of the vector field $u$ fixed in $H$ has the effect of moving $\overline D$ inside a small tubular neighbourhood while keeping the boundary $H$ fixed. Moreover, if $D'$ is the result of a small isotopy of $\overline D$ fixed in $H$, then it is not hard to construct a vector field $u'$ homotopic to $u$ such that $\overline{D(i,u')}=D'$. Hence it is enough to obtain the statement of the lemma after a small isotopy of $\overline D$ fixed in $H$.

Choose a small subset $B\subset\overline{D\setminus T}$ diffeomorphic to $\tightoverset{_\circ~~~~}{D^{n-k}}$ and a transverse slice $S$ of the integral curves of $y$ passing through $B$. Now the projection $B\to S$ along the integral curves is a $0$-codimensional map, hence in a generic position the preimage of each point is discrete. Now we can perturb the embedding of $B$ into $i(M)$ such that this projection is generic, which means the integral curves of $y$ intersect $B$ in a discrete set.

We cover $\overline{D\setminus T}$ by finitely many such subsets $B$ and obtain a small isotopy of $\overline D$ fixed in $H$ that makes all integral curves of $y$ intersect $\overline{D\setminus T}$ in a discrete set. Now the compactness of $\overline{D\setminus T}$ implies that we can choose a small neighbourhood $U$ of $\overline{D\setminus T}$ with the desired property.
\end{prf}

Now we can get to the local compression theorem. We will follow the proof of Rourke and Sanderson \cite{rs}, which is rather hard to understand in their paper, as much of it is only given as a visual intuition. We attempt to make it clearer by being more precise.

\begin{thm}\label{lct}
For any $\varepsilon>0$, the isotopy $\varphi_t~(t\in[0,1])$ in proposition \ref{ct} can be chosen such that for all $p\in P\times\R^1$ and $t\in[0,1]$ we have $d(p,\varphi_t(p))<\varepsilon$.
\end{thm}

\begin{prf}
Let $\delta>0$ be a small number to be specified later, $T\subset i(M)$ a tubular neighbourhood of $H$ and $U\subset i(M)$ a closed tubular neighbourhood of $\overline{D\setminus T}$ according to lemma \ref{lctlemma} (such that $U\cap H=\varnothing$). We will turn the vector field $u$ vertically upwards in three steps (step \ref{st2} will be the only considerably hard part).

\newcounter{step}
\renewcommand\thestep{\Roman{step}}
\medskip\noindent\refstepcounter{step}\thestep\label{st1}. \emph{Moving upwards on $i(M)\setminus(T\cup U)$.}\medskip

For all points $p\in i(M)\setminus(T\cup U)$ we have $u(p)\ne-x(p)$, hence we can canonically rotate $u(p)$ to $x(p)$ in the plane of $u(p)$ and $x(p)$. This rotation clearly extends to a small neighbourhood of $i(M)\setminus(T\cup U)$ in $i(M)$, hence it also extends to $i(M)$ such that on $\overline D$ it remains the initial $u|_{\overline D}$.

This way we can assume that $u$ was initially such that $u|_{\overline{i(M)\setminus(T\cup U)}}=x|_{\overline{i(M)\setminus(T\cup U)}}$. Now for all $p\in\overline{i(M)\setminus(T\cup U)}$ the angle of $u(p)$ and $\ua(p)$ is $\alpha(p)<\frac\pi2$, so by the compactness of $\overline{i(M)\setminus(T\cup U)}$ there is an $\alpha_1>0$ such that $\alpha(p)<\frac\pi2-\alpha_1$ is also true. Now we can use the upwards rotation of $u$ into $v$ that we defined in the proof of proposition \ref{ct} with $0<\alpha_0<\alpha_1$, hence we get that the new vector field $v$ is such that $v|_{\overline{i(M)\setminus(T\cup U)}}=\ua|_{\overline{i(M)\setminus(T\cup U)}}$.

\medskip\noindent\refstepcounter{step}\thestep\label{st2}. \emph{Moving upwards on $U$.}\medskip

Let $V$ be the closed tubular neighbourhood of $i(M)$ of radius $r$ (where $r>0$ is to be specified later) as in the proof of proposition \ref{ct}. For all $(p,x)\in U$ we put
$$A_{(p,x)}:=N_{(p,x)}i(M)\cap T_{(p,x)}(P\times\{x\}),$$
that is, $A_{(p,x)}$ is the horizontal subspace in the normal space $N_{(p,x)}i(M)$ (the orthogonal complement of $T_{(p,x)}i(M)$), which is a hyperplane in $N_{(p,x)}i(M)$. Let $N'_{(p,x)}i(M)$ be the subspace in $T_{(p,x)}(P\times\R^1)$ generated by $A_{(p,x)}$ and $\ua(p,x)$.

Now for any point $p\in U$ we can canonically rotate the fibre $V_p\subset N_pi(M)$ of $V$ around the axis $A_p$ to a closed $(k+1)$-dimensional disk in $N'_pi(M)$; let us denote this disk by $V'_p$ (here we identified the disks $V_p$ and $V'_p$ with their exponential images). This rotation clearly extends to a small neighbourhood of $U$ in $i(M)\setminus H$, hence it also extends to $i(M)$ such that at the points $p$ in the complement of a small neighbourhood of $U$ the disk $V'_p$ remains the initial $V_p$. If the radius $r$ is small enough, then the disks $V'_p$ are disjoint and the union
$$V':=\bigcup_{p\in i(M)}V'_p$$
is again a tubular neighbourhood of $i(M)$. In other words, we get the neighbourhood $V'$ by rotating $V$ fibrewise around the horizontal hyperplane to be vertical over $U$, then extending this rotation to the whole $i(M)$.

We now redefine the extension $\tilde v$ of $v$ as follows: Recall the proof of proposition \ref{ct} where we defined $\tilde v$ on a fibre $V_p$ by rotating (the parallel translate of) the vector $v(p)$ upwards as we go along the minimal geodesic from $p$ towards the boundary of $V_p$. Now we proceed in exactly the same way except that we use the disks $V'_p$ instead of $V_p$. Hence the vector field $\tilde v$ is defined on $V'$ and is vertically upwards on $\partial V'$; in the complement of $V'$ we define it to be everywhere $\ua$.

As in the proof of proposition \ref{ct}, let $\{\Phi_t\in\Diff(P\times\R^1)\mid t\in\R_+\}$ be the flow of $\tilde v$; and as in remark \ref{lctrmk}, let $\{\tilde\Phi_t\in\Diff(P\times\R^1)\mid t\in\R_+\}$ be the flow of the time-dependent vector field
$$\tilde v_t\colon P\times\R^1\to T(P\times\R^1);~p\mapsto\tilde v(p+t)+\da~~~~(t\in\R_+).$$

We need one final modification of the flow $\tilde\Phi_t$ in order to terminate the effect of the flow once the vector field near $U$ is vertically upwards. As we mentioned in remark \ref{lctrmk}, it may happen that the flow $\tilde\Phi_t$ leaves a point fixed for some time but begins to move it later; the purpose of this final modification is to not let this happen to points of $U\cup(i(M)\setminus T)$.

Observe that the vertical projection $\pr_P|_{i(M)\setminus H}$ to $P$ is an immersion and $U\cup(i(M)\setminus T)$ is compact, hence we can fix a small number $d$ (to be specified later) such that
$$0<d<\frac13\min\{d(p,q)\mid p\in U\cup(i(M)\setminus T),q\in i(M),p\ne q,{\pr}_P(p)={\pr}_P(q)\}.$$
Let $g\colon[0,\infty)\to[0,1]$ be a smooth decreasing function with $g|_{[0,d]}\equiv1$ and $g|_{[2d,\infty)}\equiv0$. We define the final version of the time-dependent vector field as
$$\hat v_t\colon P\times\R^1\to T(P\times\R^1);~p\mapsto g(t)\tilde v_t(p)~~~~(t\in\R_+)$$
and denote its flow by $\{\hat\Phi_t\in\Diff(P\times\R^1)\mid t\in\R_+\}$. Our aim in the rest of step \ref{st2} will be to prove that the flow $\hat\Phi_t$ deforms $v$ into $\ua$ on $U$ while moving all points of $i(M)$ by less than $\frac\varepsilon2$ if we define $\delta,r$ and $d$ appropriately.

Let $\Gamma$ be a component of the intersection of an integral curve of $y$ (the gradient field of $\pr_{\R^1}|_{i(M)}$) with $U$ (and so by lemma \ref{lctlemma}, shorter than $\delta$). Denote by $\pi$ the projection $V'_p\to\{p\}$ of the tubular neighbourhood $V'$ of $i(M)$ and define
$$V'_\Gamma:=\bigcup_{x\in\R^1}(\pi^{-1}(\Gamma)+x),$$
that is, the union of all vertical translates of the restriction of $V'$ to $\Gamma$.

\medskip\begin{sclaim}
The flow $\hat\Phi_t$ can take any point of $\Gamma$ only to points of $V'_\Gamma$.
\end{sclaim}

\begin{sprf}
Fix any $p\in\Gamma$. We rotated $u(p)$ to $v(p)$ in the plane of $u(p)$ and $\ua(p)$, then $v(p)$ to the vectors $\tilde v(q)~(q\in V'_p)$ in the parallel translates of the same plane. The vector $u(p)$ was in the normal space $N_pi(M)$ which is generated by the horizontal subspace $A_p$ and the upmost vector $x(p)$, hence the vectors $\tilde v(q)~(q\in V'_p)$ are in the parallel translate of the subspace generated by $A_p,x(p)$ and $\ua(p)$.

\begin{figure}[H]
\begin{center}
\centering\includegraphics[scale=0.1]{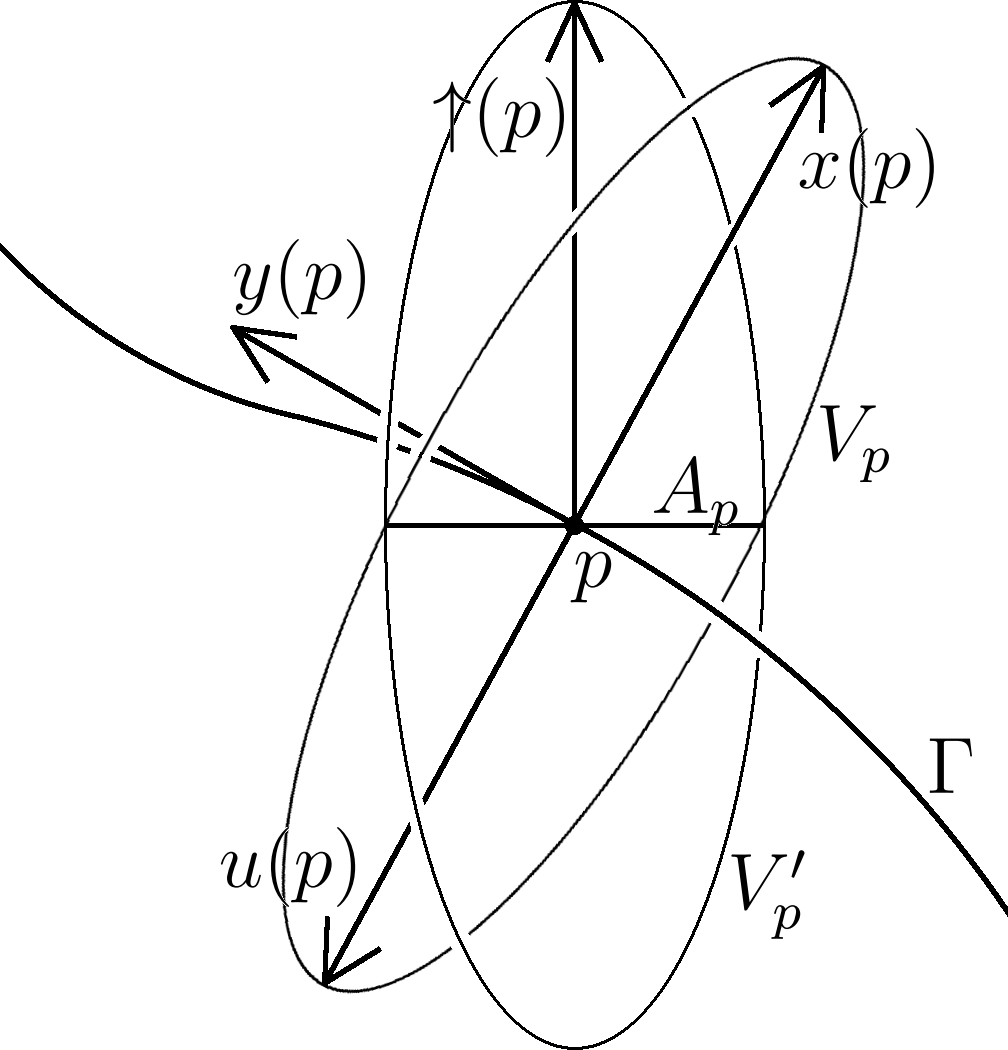}\label{kepa2}
\begin{changemargin}{2cm}{2cm} 
\caption{\hangindent=1.4cm\small Here we show the point $p$ on the curve $\Gamma$ with the vectors at $p$ and the normal disks $V_p$ and $V'_p$.}
\end{changemargin} 
\vspace{-1.3cm}
\end{center}
\end{figure}

Now $p\notin H$ implies that $y(p)\ne\ua(p)$. Moreover, $y(p)\ne0$ because $y(p)=0$ would mean that $T_pi(M)$ is horizontal, thus the upmost unit vector orthogonal to it is $x(p)=\ua(p)$, and so $-x(p)=\da(p)$; now if $p$ is a point of the downset $D$, then this cannot happen because of lemma \ref{ctlemma}, hence it also cannot happen in a small neighbourhood of $D$ and we can choose $U$ to be a neighbourhood small enough.

Therefore $y(p)$ and $\ua(p)$ generate a plane and it is not hard to see that $x(p)$ is in this plane as well. Indeed, $x(p)$ and $\frac{y(p)}{\lv y(p)\rv}$ are the upmost unit vectors in $N_pi(M)$ and $T_pi(M)$ respectively, hence the plane they generate has to contain the upmost unit vector of $N_pi(M)\oplus T_pi(M)=T_p(P\times\R^1)$ which is $\ua(p)$.

This way the vectors $\tilde v(q)~(q\in V'_p)$ are in the parallel translate of the subspace generated by $A_p,y(p)$ and $\ua(p)$, where $A_p$ and $\ua(p)$ generate the linear space $N'_pi(M)$ by definition. Now $V'_p$ is a small closed ball in $N'_pi(M)$ and if $q$ is in $\partial V'_p$ or $p$ is an endpoint of $\Gamma$, then $\tilde v(q)=\ua(q)$, thus all vectors $\tilde v(q)~(q\in V'_p)$ are tangent to $V'_\Gamma$.

If we add the downwards unit vector field $\da$ to these vectors, then we still get vectors tangent to $V'_\Gamma$ and the same is true if we add $\da$ to the vectors $\tilde v(q+t)$ where $q+s\in V'_p$ and $0\le s\le t\le3d$ and multiply them by small numbers. We get the vectors of the time-dependent vector field $\hat v_t=g(t)\tilde v_t$ this way, hence the integral curves of $\hat v_t$ that start in $\Gamma$ remain in $V'_\Gamma$.
\end{sprf}

In the following $q$ will denote a point and $t$ will denote a time such that $0\le t\le3d$ and $q+s\in V'_p$ (for some $0\le s\le t$ and $p\in U$). Based on the above proof we can decompose the vector $\tilde v_t(q)$ in the form
$$\tilde v_t(q)=a_t(q)+b_t(q)+c_t(q),$$
where $a_t(q)$ is the vertical component, $b_t(q)$ is the horizontal component in the (parallel translate of the) plane of $\ua(p)$ and $y(p)$ and $c_t(q)$ is the component in the (parallel translate of the) horizontal subspace $A_p$ orthogonal to $T_pi(M)$. These three vectors are orthogonal to each other.

Observe that we get the vector $\tilde v_t(q)$ by adding the dowards unit vector $\da$ to some other unit vector. One can easily check that this implies that exactly one of the following conditions hold:
\begin{enumerate}
\renewcommand{\labelenumi}{\theenumi}
\renewcommand{\theenumi}{\rm{(\roman{enumi})}}
\item\label{un1} The length of the vertical component is $\lv a_t(q)\rv\le1-\frac{\sqrt2}{2}$.
\item\label{un2} The length of the horizontal component is $\lv b_t(q)+c_t(q)\rv>\frac{\sqrt2}{2}$.
\end{enumerate}

Let $\gamma\colon\R_+\to P\times\R^1$ be an integral curve of $\tilde v_t$ starting at a point $p\in U$ and put $t_0:=\min\{t>0\mid\gamma'(t)=0\}$. Denote by $I_1$ (resp. $I_2$) the union of those subintervals $I\subset[0,t_0]$ for which at each time $t\in I$ the condition \ref{un1} (resp. \ref{un2}) holds with $q=\gamma(t)$. Let $t_i$ denote the length (i.e. Lebesgue measure) of $I_i$ (for $i=1,2$).

\medskip\begin{sclaim}
For $i=1,2$ we have $t_i<\sqrt2(\delta+r)$.
\end{sclaim}

\begin{sprf}
If $p\in D$, then $u(p)$ is in the plane generated by $\ua(p)$ and $y(p)$, which implies $c_0(p)=0$. Therefore, given any number $\zeta>0$, we can choose the neighbourhood $U$ so small that we have $\lv c_0(p)\rv<\zeta$ for all $p\in U$. Then $\lv c_t(q)\rv<\zeta$ is also true for any $q$ and $t$ (where $q+s\in V'_p,0\le s\le t\le3d$ and $p\in U$), since this component can only decrease when we rotate the vector upwards. 

Let $\Gamma$ be the component of the intersection of the integral curve of $y$ with $U$ for which $p\in\Gamma$ (and so $\gamma|_{[0,t_0]}$ only maps to points of $V'_\Gamma$). Recall that $\delta$ is an upper estimate on the length of $\Gamma$ and $r$ is the radius of the disks $V'_p$.

By choosing $\delta$ and $r$ small enough, we can assure that $\pi^{-1}(\Gamma)$ (where $\pi$ is the projection of $V'$) is in one Eucledian neighbourhood, hence so is $V'_\Gamma$. Now we can identify all tangent spaces in the points of $V'_\Gamma$, hence it makes sense to talk about the change of direction of the vector $c_t~(t\in[0,t_0])$ along the curve $\gamma$. By setting $\delta$ sufficiently small, we can make this change of direction arbitrarily small along $\gamma$; moreover, because of the compactness of $U$, this can be done universally for all integral curves of $\tilde v_t$ starting in $U$ and for $t\le 3d$.

The above considerations imply that the integral curves of $\tilde v_t=a_t+b_t+c_t$ are arbitrarily close to the integral curves of $a_t+b_t$ (for $t\le t_0$), thus it is enough to prove the claim so that we forget about the component $c_t$ and assume $\tilde v_t=a_t+b_t$.

\begin{figure}[H]
\begin{center}
\centering\includegraphics[scale=0.1]{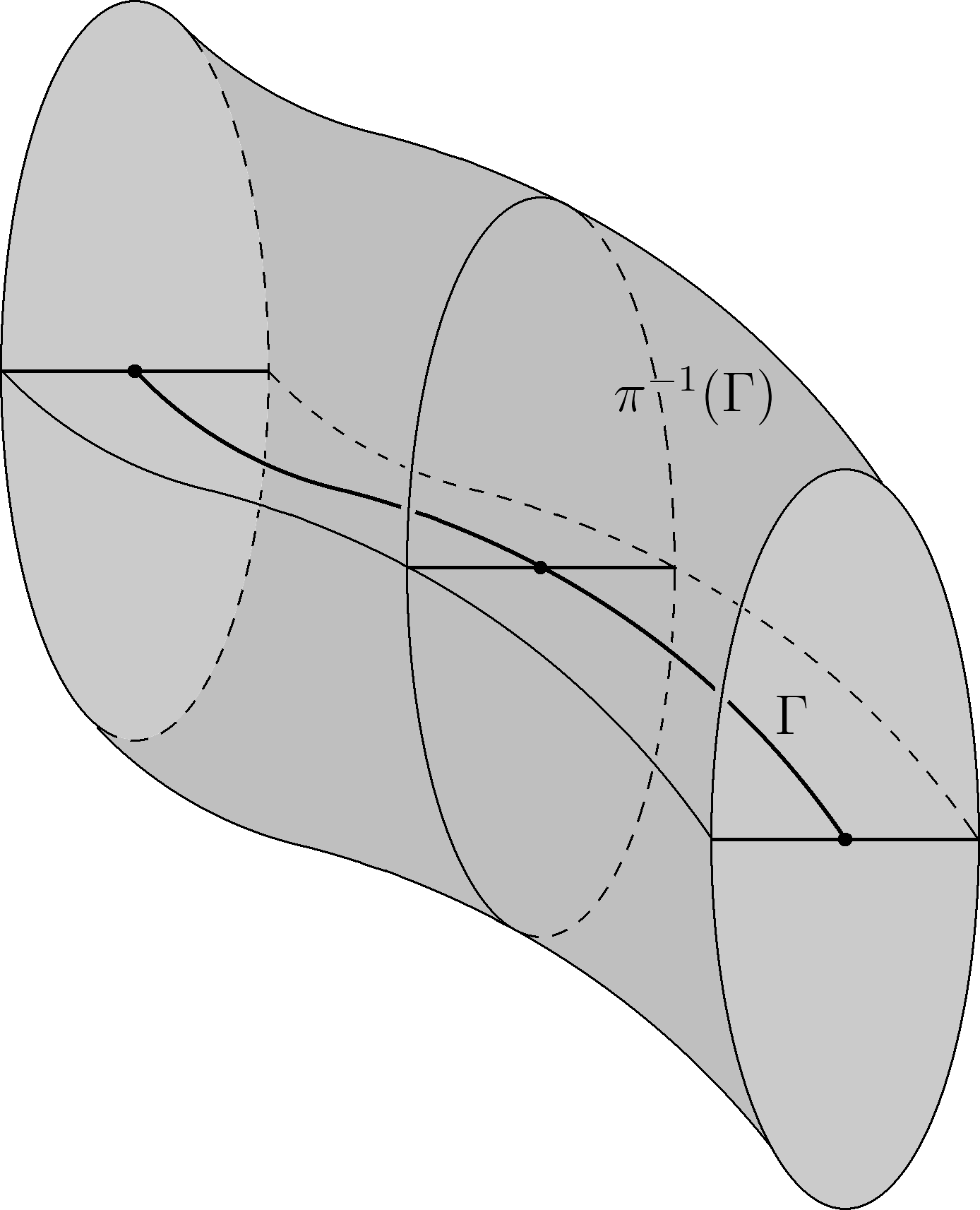}\label{kepa3}
\begin{changemargin}{2cm}{2cm} 
\caption{\hangindent=1.4cm\small Here we show $\pi^{-1}(\Gamma)$; we get the set $V'_\Gamma$ by taking the union of all vertical translates of this.}
\end{changemargin} 
\vspace{-1.5cm}
\end{center}
\end{figure}

We get the flow of $\tilde v_t$ by combining the flow of $\tilde v$ (the non-time-dependent vector field) with the unit downwards flow, in other words, flowing along $\tilde v_t$ is the same as flowing along $\tilde v$, except that in the meantime everything is flowing downwards with unit speed. Hence $\tilde v_t(\gamma(t))=0$ is equivalent to $\tilde v(\gamma(t)+t)=\ua$. The length of the vertical component $a_t(\gamma(t))$ being at most $1-\frac{\sqrt2}2$ is equivalent to the vertical component $a(\gamma(t)+t)$ of $\tilde v(\gamma(t)+t)$ being at least $\frac{\sqrt2}2$; the length of the horizontal component $b_t(\gamma(t))$ being more than $\frac{\sqrt2}2$ is equivalent to the length of the horizontal component $b(\gamma(t)+t)$ of $\tilde v(\gamma(t)+t)$ being more than $\frac{\sqrt2}2$.

\medskip\noindent\emph{Proof for $t_1$.}\enspace\ignorespaces The difference of the vertical coordinate of $p$ and the vertical coordinate of the upmost point of $\pi^{-1}(\Gamma)$ is clearly smaller than $\delta+r$. The flow of $\tilde v$ moves points vertically only upwards and it moves $p$ to the boundary of $\pi^{-1}(\Gamma)$ where $\tilde v$ is vertically upwards.

Now the time the integral curve of $\tilde v$ starting at $p$ spends at points where the vertical component $a$ of the velocity vector is at least $\frac{\sqrt2}{2}$ before reaching the boundary of $\pi^{-1}(\Gamma)$ is clearly an upper bound for $t_1$. But this time cannot be more than the time it takes to get to the vertical level of the upmost point of $\pi^{-1}(\Gamma)$ on a curve for which the velocity vector always has vertical component at least $\frac{\sqrt2}2$. This is trivially less than $\sqrt2(\delta+r)$, hence $\sqrt2(\delta+r)$ is also an upper bound for $t_1$.

\medskip\noindent\emph{Proof for $t_2$.}\enspace\ignorespaces The curve $\gamma$ cannot reach the boundary of $V'_\Gamma$ (where the vector field $\tilde v_t$ is $0$) later than the integral curve of $\tilde v$ starting at $p$ reaches the horizontal level of the boundary of $\pi^{-1}(\Gamma)$.

We assumed that the horizontal component of $\tilde v$ is in the plane generated by (the parallel translates of) vectors of $\ua$ and $y$, hence the horizontal movement along integral curves of $\tilde v$ happens in the directions of (parallel translates of) vectors of $y$. The vectors of $y$ are tangent to $\Gamma$ by definition, therefore the length of the hoizontal component of such an integral curve is smaller than $\delta+r$.

Now the time the integral curve of $\tilde v$ starting at $p$ spends at points where the horizontal component $b$ of the velocity vector is more than $\frac{\sqrt2}{2}$ before reaching the horizontal level of the boundary of $\pi^{-1}(\Gamma)$ is clearly an upper bound for $t_2$. But this time cannot be more than the time it would take on a curve which always has velocity vectors with vertical components more than $\frac{\sqrt2}2$. This is again less than $\sqrt2(\delta+r)$, hence $\sqrt2(\delta+r)$ is an upper bound for $t_2$ as well.
\end{sprf}

The claim we have just proved implies the inequality $t_0=t_1+t_2<2\sqrt2(\delta+r)$, hence for any integral curve $\gamma$ of $\tilde v_t$ starting in $U$ there is $0<t<2\sqrt2(\delta+r)$ for which $\gamma'(t)=0$. We can choose $\delta$ and $r$ such that $2\sqrt2(\delta+r)\le d$, which implies the same statement where $\gamma$ is an integral curve of $\hat v_t=g(t)\tilde v_t$ instead of $\tilde v_t$. But if $\gamma$ is an integral curve of $\hat v_t$, then $\gamma'(t)=0$ implies that the curve $\gamma$ always remains stationary after $t$, because the bump function $g$ was chosen such that the vectors of $\hat v_t$ that are far above have no effect on $\gamma$.

What we proved by this is that all integral curves of $\hat v_t$ that start in $U$ are constant after $2\sqrt2(\delta+r)$ time. Now setting $d$ to be $2\sqrt2(\delta+r)$ we can even make the whole flow $\hat\Phi_t$ terminate after $4\sqrt2(\delta+r)$ time. It is easy to see that for all $p\in P\times\R^1$ and $t\in\R_+$ we have $\lv\hat v_t(p)\rv\le\sqrt2$, thus $d(p,\hat\Phi_t(p))\le8(\delta+r)$ which is smaller than $\frac\varepsilon2$ if $\delta$ and $r$ are sufficiently small. 

\medskip\noindent\refstepcounter{step}\thestep\label{st3}. \emph{Moving upwards on $T$.}\medskip

In the course of step \ref{st1} and step \ref{st2} we deformed the embedding of $M$ to a new one which we will denote by $j(M):=\hat\Phi_{2d}(M)$. The normal vector field $u$ was also deformed to a normal vector field $w$ on $j(M)$ which is vertically upwards in the complement of $T':=\hat\Phi_{2d}(T)$.

Observe that the original vector field $u$ was horizontal on $H$, hence when we rotated it into $v$, the angle of $v(p)$ and $\ua(p)$ was $\beta(p)=\alpha_0$ for all $p\in H$. Hence by choosing the tubular neighbourhood $T$ small enough, we can assure the inequality $\beta(p)\le2\alpha_0$ for all $p\in T$. The deformations of steps \ref{st1} and \ref{st2} could only decrease this angle, thus the angle of $w(p)$ and $\ua(p)$ is again at most $2\alpha_0$ for all $p\in T'$.

Now we extend $w$ to a vector field $\tilde w$ on $P\times\R^1$ using a closed tubular neighbourhood $W$ of $j(M)$ in exactly the same way as we did in the proof of proposition \ref{ct}. We denote the flow of $\tilde w$ by $\{\Psi_t\in\Diff(P\times\R^1)\mid t\in\R_+\}$ and, as in remark \ref{lctrmk}, we denote by $\{\tilde\Psi_t\in\Diff(P\times\R^1)\mid t\in\R_+\}$ the flow of the time-dependent vector field
$$\tilde w_t\colon P\times\R^1\to T(P\times\R^1);~p\mapsto\tilde w(p+t)+\da~~~~(t\in\R_+).$$

By setting the angle $\alpha_0$ small enough, we can make the horizontal components of all vectors of $\tilde w$ arbitrarily small. The compactness of $W$ implies that there is an $s>0$ for which $\Psi_t(p)\notin W$ for all $t\ge s$ and $p\in j(M)$, this way turning the vector field vertical on the image of $T'$ as well, while moving all points arbitrarily close to the unit upwards flow. Hence combining it with the unit downwards flow to get $\tilde\Psi_t$, the distance $d(p,\tilde\Psi_t(p))$ can be made arbitrarily small for all $p\in P\times\R^1$ and $t\in\R_+$, in particular we can set it smaller than $\frac\varepsilon2$.

\medskip Note that each step did not affect the achievements of the previous one, thus we can define the desired isotopy $\varphi_t~(t\in[0,1])$ to apply first $\hat\Phi_t~(t\in[0,2d])$ and then $\tilde\Psi_t~(t\in[0,s])$.
\end{prf}

\subsection{Addenda}

In this subsection we will prove some of the numerous extensions and corollaries of the local compression theorem \ref{lct}. The first one is called the multi-compression theorem and is as follows.

\begin{thm}\label{mct}
Let $n,k,m$ be natural numbers such that $k\ge1$ and let $i\colon M^n\hookrightarrow P^{n+k}\times\R^m$ be the embedding of a compact manifold (possibly with boundary) equipped with $m$ pointwise independent normal vector fields $u_1,\ldots,u_m$. Then there is an (ambient) isotopy $\varphi_t~(t\in[0,1])$ of $i$ such that $\varphi_0=i$ and $d\varphi_1(u_1),\ldots,d\varphi_1(u_m)$ are vertical (i.e. $d\varphi_1(u_j)$ is the positive unit vector on the $j$-th coordinate line of $\R^m$ in $T(P\times\R^m)$). Moreover, for any $\varepsilon>0$ this $\varphi_t$ can be chosen such that for all $p\in P\times\R^m$ and $t\in[0,1]$ we have $d(p,\varphi_t(p))<\varepsilon$.
\end{thm}

\begin{prf}
We proceed by induction on the number $m$. The $m=1$ case is proposition \ref{ct} and theorem \ref{lct}. Now suppose that we know the theorem for $m-1$ and want to prove for $m$.

We may assume by induction and lemma \ref{ctlemma} that the vector fields $u_1,\ldots,u_{m-1}$ are already vertical and $u_m$ is a unit vector field pointwise orthogonal to $Ti(M)$ and $u_1,\ldots,u_{m-1}$. Now if we put $Q:=P\times\R^1$ (where this $\R^1$ is the last real coordinate line), then the projection $\pr_Q|_{i(M)}\colon i(M)\to Q$ is an immersion because we projected along (the first) $m-1$ normal vector fields. A tubular neighbourhood of the image of this immersion can be lifted to an $(n+k+1)$-dimensional strip $S\subset P\times\R^m$ containing $i(M)$; the fibres of $S$ are everywhere parallel to $Q$ and $u_m$ is tangent to $S$.

Let $U_j\subset V_j\subset W_j\subset S$ (where $j$ is in some finite index set) be such that $U_j\approx\tightoverset{_\circ~~~~~~~}{D^{n+k+1}}$, $V_j$ is the $\varepsilon$-neighbourhood of $U_j$ and $W_j$ is the $\varepsilon$-neighbourhood of $V_j$ in $S$ (for some $\varepsilon>0$) and the projection $\pr_Q|_{W_j}$ is an embedding. Now we can use the local compression theorem \ref{lct} inside $W_j\approx\R^{n+k+1}$, composed with a bump function that is $0$ outside $V_j$, to obtain an isotopy of $i$ that keeps the image of $M$ inside the strip $S$ and only moves the points inside $W_j$.

By applying this process for all indices $j$, we obtain an isotopy of $i$ that moves $i(M)$ inside the strip $S$ and turns the vector field $u_m$ vertical. We can keep the vector fields $u_1,\ldots,u_{m-1}$ (which are orthogonal to the fibres of $S$) vertical throughout this process, hence the theorem is proved.
\end{prf}

\begin{crly}\label{ict}
In theorem \ref{mct} we can replace the word ``embedding'' by ``immersion'' and the word ``isotopy'' by ``regular homotopy''.
\end{crly}

\begin{prf}
If $i\colon M\looparrowright P\times\R^m$ is an immersion, then fix any embedding $j\colon M\hookrightarrow\R^l$ (where $l$ is a large number) and obtain an embedding $i\times j\colon M\hookrightarrow P\times\R^{m+l}$. We can lift the normal vector fields $u_1,\ldots,u_m$ to this embedding and define the normal vector fields $u_{m+1},\ldots,u_{m+l}$ to be everywhere vertical. We can use the multi-compression theorem \ref{mct} to this to obtain an isotopy of $i\times j$, which composed with the projection to $P\times\R^m$ is a regular homotopy of $i$.
\end{prf}

\begin{rmk}
It follows from the proofs that the relative versions of the statements above (and below) are also true, that is, if the vector fields are already vertical on a neighbourhood of a compact subset $C\subset i(M)$, then the isotopy (or regular homotopy) $\varphi_t~(t\in[0,1])$ is fixed on a neighbourhood of $C$.
\end{rmk}

The following is the stratified version of the multi-compression theorem, which was proved in \cite{hosszu}.

\begin{thm}\label{sct}
Let $n,k,m$ be natural numbers such that $k\ge1$ and let $f\colon M^n\to P^{n+k}\times\R^m$ be the map of a compact manifold (possibly with boundary). Suppose that $M$ is equipped with a finite stratification by submanifolds $S_i~(i=1,\ldots,r)$, where $M=S_1\cup\ldots\cup S_r$ and $S_{i-1}\subset\overline{S_i}$ for all $i$, moreover, the restriction $f|_{S_i}$ is an embedding for all $i$ and $f(S_i)\cap f(S_j)=\varnothing$ for all $i,j$. Let $u_1,\ldots,u_m$ be pointwise independent vector fields along $f(M)$ transverse to the stratification, that is, sections of the bundle $T(P\times\R^m)|_{f(M)}$ such that for all $i$ the vector fields $u_1|_{f(S_i)},\ldots,u_m|_{f(S_i)}$ are nowhere tangent to $f(S_i)$. Then there is a diffeotopy $\varphi_t~(t\in[0,1])$ of $P\times\R^m$ such that $\varphi_0=\id_{P\times\R^m}$ and $d\varphi_1(u_1),\ldots,d\varphi_1(u_m)$ are vertical. Moreover, for any $\varepsilon>0$ this $\varphi_t$ can be chosen such that for all $p\in P\times\R^m$ and $t\in[0,1]$ we have $d(p,\varphi_t(p))<\varepsilon$.
\end{thm}

\begin{prf}
We prove by induction on the number $r$ of the strata. The $r=1$ case is theorem \ref{mct}. Now suppose that we know the theorem for $r-1$ and want to prove for $r$.

Now $\overline{S_{r-1}}$ is stratified by $r-1$ submanifolds, hence by the induction hypothesis we can deform the map $f|_{\overline{S_{r-1}}}$ to turn the vector fields $u_1|_{f(\overline{S_{r-1}})},\ldots,u_m|_{f(\overline{S_{r-1}})}$ vertical by a diffeotopy of $P\times\R^m$. This way we obtained a new map $g\colon M\to P\times\R^m$ with the same properties as $f$ and new vector fields $v_1,\ldots,v_m$ along $g(M)$ that are transverse to the stratification and vertical on $g(\overline{S_{r-1}})$.

Of course, we can homotope the vector fields $v_1,\ldots,v_m$ to be vertical on a small neighbourhood $U\subset g(M)$ of $g(S_{r-1})$ as well. Now we can choose an even smaller neighbourhood $V\subset g(M)$ of $g(S_{r-1})$ such that $\overline V\subset U$ and apply the relative version of the multi-compression theorem \ref{mct} to $g(M)\setminus V$ to turn the vector fields $v_1|_{g(M)\setminus V},\ldots,v_m|_{g(M)\setminus V}$ vertical by only moving the manifold in the complement of $V$ and leaving a neighbourhood of $\partial(g(M)\setminus V)$ fixed. This way the vector fields become vertical on the whole manifold.
\end{prf}

The last statement in this section is a simple addendum to the multi-compression theorem, which was proved in \cite{nszt}.

\begin{prop}\label{ct+}
If in theorem \ref{mct} we do not require that the points stay in their $\varepsilon$-neighbourhoods, then the isotopy $\varphi_t~(t\in[0,1])$ can be chosen such that for any point $p\in M$ and any time $t_0\in[0,1]$ the tangent vector of the curve $t\mapsto\varphi_t(p)$ at the point $\varphi_{t_0}(p)$ is not tangent to $\varphi_{t_0}(M)$.
\end{prop}

\begin{prf}
We will slightly modify the steps of the induction that proves the multi-compression theorem \ref{mct}. In the starting step we just use the vector field $\tilde v$ construsted in the proof of proposition \ref{ct} (instead of the one in the local compression theorem \ref{lct}); now the derivative in $t_0$ of each curve $t\mapsto\varphi_t(p)$ is just the image of the normal vector $v(p)$ under the diffeomorphism $\varphi_{t_0}$, hence it is independent of $T_{\varphi_{t_0}(p)}\varphi_{t_0}(M)$.

%It is clear from the proofs above that we only have to check that the flow $\{\tilde\Phi_t\in\Diff(P\times\R^1)\mid t\in\R_+\}$ of the time-dependent vector field
%$$\tilde v_t\colon P\times\R^1\to T(P\times\R^1);~p\mapsto\tilde v(p+t)+\da~~~~(t\in\R_+)$$
%has the desired property for small numbers $t$ (see remark \ref{lctrmk} and the proof of theorem \ref{lct}). If $\gamma\colon t\mapsto\tilde\Phi_t(p)$ is the integral curve for which $\gamma(0)=p$ and $\gamma'(t)=\tilde v_t(\gamma(t))$, then we have to prove that the vector $\tilde v_t(\gamma(t))$ is not in $T_{\gamma(t)}\tilde\Phi_t(M)$.

%We get the point $\gamma(t)$ by flowing along the non-time-dependent vector field $\tilde v$ for time $t$, then flowing downwards with unit speed for time $t$. The same goes for the whole manifold, hence if we translate the tangent space $T_{\gamma(t)}\tilde\Phi_t(M)$ upwards by $t$, then we get the tangent space $T_{\Phi_t(p)}\Phi_t(M)$ (which is obtained by applying the flow of $\tilde v$ to $T_pi(M)$). Therefore we only have to prove that $\tilde v(\Phi_t(p))+\da(\Phi_t(p))$ is not in $T_{\Phi_t(p)}\Phi_t(M)$.

%Now saying $\tilde v(\Phi_t(p))+\da(\Phi_t(p))\in T_{\Phi_t(p)}\Phi_t(M)$ is equivalent to saying that the components of $\tilde v(\Phi_t(p))$ and $\da(\Phi_t(p))$ in $N_{\Phi_t(p)}\Phi_t(M)$ (the orthogonal complement of $T_{\Phi_t(p)}\Phi_t(M)$) are each other's negatives. If we look at the plane generated by $\tilde v(\Phi_t(p))$ and $\da(\Phi_t(p))$ (this is really a plane because we assumed $\tilde v_t(\gamma(t))\ne0$), then 

Note that since the vector field we used in the first step is not time-dependent, the image of $M$ keeps flowing with unit speed upwards in the direction of the first real coordinate of $P\times\R^m$, which we know is now independent of the tangent spaces of the image of $M$. Thus by the compactness of $M$, the minimal length of those vectors which added to this upwards unit vector get into the tangent space of the image of $M$ is bounded from below. In the induction step we are bound to use the time-dependent vector fields defined in the proof of theorem \ref{lct}, but by reparameterising time for these, we can make them arbitrarily small. Hence if we add them to the upwards unit vectors remained from the first step, we get that the derivatives of the curves $t\mapsto\varphi_t(p)$ are again not tangent to the image of $M$.
\end{prf}

\begin{rmk}
Much of the versions and addenda of the compression theorem (such as corollary \ref{ict} or proposition \ref{ct+}) can be obtained from Hirsch--Smale immersion theory. But for example the local versions (to my knowledge) do not follow from such general theories.
\end{rmk}

%-------------------------------------------------------------------------------------------------------------------------------

%\section{Tools from stable homotopy theory}

%-------------------------------------------------------------------------------------------------------------------------------

\refstepcounter{kaki}
\section{Virtual complexes}\label{virtc}

Here we introduce the notion of virtual complexes based on \cite{hosszu}. In short, a virtual complex is a ``CW-complex'' where the gluing maps of the cells are only defined in stable sense.

First we recall the definition of a standard gluing procedure that results in a CW-complex, and then we give a similar definition of a ``stable gluing procedure'' that results in a virtual complex. This way the analogies and differences between usual and virtual CW-complexes will be easy to see.

\begin{defi}
Let $\{(A_i,B_i)\mid i=0,1,\ldots\}$ be a sequence of CW-pairs and suppose that (each component of) the pair $(A_i,B_i)$ is $c_i$-connected where $c_i\to\infty$ for $i\to\infty$.
\begin{itemize}
\item[(0)] Put $X_0:=A_0$ and fix a map $\rho_1\colon B_1\to X_0$.
\item[(1)] Put $X_1:=X_0\usqcup{\rho_1}A_1$ and fix a map $\rho_2\colon B_2\to X_1$.
\item[($i$)] Put $X_i:=X_{i-1}\usqcup{\rho_i}A_i$ and fix a map $\rho_{i+1}\colon B_{i+1}\to X_i$.
\end{itemize}
Using this recursion we define the space $X:=\liminfty{i}X_i$.
\end{defi}

\begin{ex}
If we set in the above definition $A_i$ to be a disjoint union $\underset{j}{\bigsqcup}D_j^i$ (where $D_j^i\approx D^i$ for all $j$) and $B_i$ to be its boundary $\underset j\bigsqcup S_j^{i-1}$, then we get the standard definition of CW-complexes.
\end{ex}

\begin{defi}\label{vc}
Let $\{(A_i,B_i)\mid i=0,1,\ldots\}$ be a sequence of CW-pairs and suppose that (each component of) the pair $(A_i,B_i)$ is $c_i$-connected where $c_i\to\infty$ for $i\to\infty$.
\begin{itemize}
\item[(0)] Put $X_0:=A_0$ and fix a stable map $\rho_1\colon B_1\nrightarrow X_0$, that is, for some $n_1\in\N$ a map $S^{n_1}B_1\to S^{n_1}X_0$ denoted by $S^{n_1}\rho_1$ (although $\rho_1$ may not exist).
\item[(1)] Put $S^{n_1}X_1:=S^{n_1}X_0\usqcup{S^{n_1}\rho_1}S^{n_1}A_1$ (although $X_1$ may not exist) and fix a stable map $\rho_2\colon B_2\nrightarrow X_1$, that is, for some $n_2\ge n_1$ a map $S^{n_2}B_2\to S^{n_2}X_1$ denoted by $S^{n_2}\rho_2$ (although $\rho_2$ may not exist).
\item[($i$)] Put $S^{n_i}X_i:=S^{n_i}X_{i-1}\usqcup{S^{n_i}\rho_i}S^{n_i}A_i$ (although $X_i$ may not exist) and fix a stable map $\rho_{i+1}\colon B_{i+1}\nrightarrow X_i$, that is, for some $n_{i+1}\ge n_i$ a map $S^{n_{i+1}}B_{i+1}\to S^{n_{i+1}}X_i$ denoted by $S^{n_{i+1}}\rho_{i+1}$ (although $\rho_{i+1}$ may not exist).
\end{itemize}
Using this recursion we define the symbol $X:=\liminfty{i}X_i$.
\end{defi}

\begin{prop}
Although the limit space $X$ does not necessarily exist, for any $r$ there is an $n(r)$ such that for all $n\ge n(r)$ the $(n+r)$-homotopy type of $S^nX$ is well-defined.
\end{prop}

\begin{prf}
If $i(r)$ is such a number that $c_i>r$ for all $i\ge i(r)$, then we can set $n(r):=n_{\max(i(r),c_{i(r)})+1}$, hence the spaces $S^nX_j$ exist for all $n\ge n(r)$ and $j=0,\ldots,i(r)$; we can define this $S^nX_{i(r)}$ to be the $(n+r)$-type of $S^nX$. This is well-defined, as for any $i\ge i(r)$ and $m\ge n$, if the $(m+r)$-type of $S^mX_j$ is well-defined for $j=0,\ldots,i$, then the $(m+r)$-type of $S^mX_i$ is equal to that of $S^{m-n}S^nX_{i(r)}$.
\end{prf}

\begin{crly}
\begin{itemize}
\item[]
\item[\rm{(1)}] For any virtual complex $X$ the space $\Gamma X=\Omega^\infty S^\infty X$ is well-defined.
\item[\rm{(2)}] The suspension functor has an inverse for virtual complexes, i.e. for any virtual complex $X$ there is a virtual complex $S^{-1}X$.
\item[\rm{(3)}] For any virtual complex $X$ and topological space $Y$ the stable homotopy classes $\{Y,X\}$ are well-defined.
\end{itemize}
\end{crly}

A direct application of the notion of virtual complexes is the following.

\begin{crly}\label{gammatnu}
If $\nu$ is a virtual vector bundle over a space $A=\liminfty iA_i$, which is the union of compact subspaces $A_1\subset A_2\subset\ldots$, then $T\nu$ exists as a virtual complex and so the space $\Gamma T\nu$ exists.
\end{crly}

\begin{prf}
If $\nu$ is a virtual bundle of dimension $k$, then it is an equivalence class of formal differences $\alpha-\beta$ where $\alpha,\beta$ are vector bundles over $A$ and $\dim\alpha-\dim\beta=k$. If $A$ is compact, then $\nu$ has a representative of the form $\alpha-\varepsilon^m$ for some $m\in\N$. The Thom space of a virtual bundle of this form can be defined as $S^{-m}T\alpha$ which makes sense as a virtual complex, hence $\Gamma T\nu:=\Gamma S^{-m}T\alpha$ is a well-defined space.

Now even if $A$ is not compact, $\nu$ can be represented by a sequence of formal differences $\alpha_i-\varepsilon^{m_i}$ where $\alpha_i$ is a vector bundle over $A_i$ of dimension $m_i+k$ and $\alpha_i|_{A_{i-1}}=\alpha_{i-1}\oplus\varepsilon^{m_i-m_{i-1}}$. This way $\Gamma T\nu|_{A_i}$ exists for all $i$ and contains $\Gamma T\nu|_{A_{i-1}}$, thus we can define $\Gamma T\nu$ as $\liminfty i\Gamma T\nu|_{A_i}$.
\end{prf}

\begin{rmk}\label{spec}
We should point out the connection between the notion of virtual complexes and the notion of CW-spectra. In order to obtain a spectrum from a virtual complex, one has to choose finite dimensional approximations of the blocks $A_i$ in definition \ref{vc} attached stably to the previously obtained spaces $X_{i-1}$. Two different sequences of approximations give different spectra, but these are naturally equivalent, in particular they define the same extraordinary cohomology theory. So the advantage of the notion of the virtual complex is that one does not have to choose data (the sequence of approximations) that finally turn out to be unimportant. On the other hand a virtual complex has a filtration by the $X_i$ in the construction, so one can say, that a virtual complex defines an equivalence class of spectra without pointing out a particular spectrum of this equivalence class.
\end{rmk}

%-------------------------------------------------------------------------------------------------------------------------------
%-------------------------------------------------------------------------------------------------------------------------------

\end{document}